\documentclass[11pt]{article}
\usepackage[utf8]{inputenc}
\usepackage{geometry}
\usepackage{etoolbox}
\usepackage{lmodern}
\usepackage[affil-it]{authblk} 
\usepackage{wrapfig}
\usepackage{graphicx}
\usepackage{listings}
\usepackage[dvipsnames]{xcolor}
\usepackage{tcolorbox}
\usepackage{subfiles}
\usepackage{esvect}

\usepackage{systeme}
\usepackage{enumerate}
\usepackage{amssymb, amsmath, amsthm, bm}
\usepackage{mathtools}
\usepackage{commath}
\usepackage[export]{adjustbox}
\usepackage{caption}
\usepackage{subcaption}
\usepackage[english]{babel}
\usepackage{dsfont}
\usepackage{hyperref}
\usepackage{cleveref}
\usepackage[numbers]{natbib}
\usepackage{thmtools}
\usepackage{thm-restate}

\theoremstyle{plain}
\newtheorem{thm}{Theorem}[section]
\newtheorem{lemma}[thm]{Lemma}
\newtheorem{cor}[thm]{Corollary}
\newtheorem{prop}[thm]{Proposition}

\newtheorem{rmk}[thm]{Remark}
\newtheorem{eg}[thm]{Example}
\newtheorem{claim}[thm]{Claim}
\theoremstyle{definition}
\newtheorem{defn}[thm]{Definition}

\newtheorem{obs}[thm]{Observation}

\DeclareMathOperator{\trace}{tr}
\DeclareMathOperator{\nortrace}{\bar{\trace}}

\DeclareMathOperator{\val}{val}

\DeclareMathOperator{\spann}{span}

\DeclareMathOperator{\perm}{perm}

\DeclareMathOperator{\var}{Var}

\def \a{\alpha}
\def \b{\beta}
\def \c{\gamma}

\def \A{\mathcal{A}}
\def \B{\mathcal{B}}
\def \C{\mathcal{C}}
\def \D{\mathcal{D}}
\def \E{\mathcal{E}}

\def \I{\mathcal{I}}

\def \L{\mathcal{L}}
\def \M{\mathcal{M}}
\def \N{\mathcal{N}}

\def \P{\mathcal{P}}
\def \Q{\mathcal{Q}}
\def \R{\mathcal{R}}

\def \X{\mathcal{X}}

\def \Z{\mathcal{Z}}

\def \Cbb{\mathbb{C}}
\def \Rbb{\mathbb{R}}

\def \Zbb{\mathbb{Z}}
\def \Nbb{\mathbb{N}}
\def \Ebb{\mathbb{E}}

\def \mh{\hat{m}}

\def \opr{\circ_{R}}
\def \o{\Omega}
\def \oset{\mathbf{\o_{0,1}}}
\def \oa{\o_{E\(\a\)}} 
\def \oc{\o_{E\(\c\)}} 

\def \oams{\o_{\mzshape, s}}
\def \oz{\o_{\zshape, 1}}
\def \ozm{\o_{Z(m)}} 
 
\def \ozmm{\o_{Z\(m'\)}}
\def \ozms{\o_{Z(m), \o^{(1)},\dots,\o^{(s)}}}

\def \zshape{\a_Z}
\def \mzshape{\a_{Z(m)}}

\def \mzshapesdistr{M_{Z(m),s}}
\def \mzshapesdistrg{M^{(G)}_{Z(m),s}}

\def \({\left(}
\def \){\right)}
\def \[{\left[}
\def \]{\right]}
\def \<{\left\langle}
\def \>{\right\rangle}
\def \lcurb{\left\{}
\def \rcurb{\right\}}
\def \={\longleftrightarrow}

\def \vaa{\vec{\a}}
\def \vbb{\vec{\b}}
\def \vcc{\vec{\c}}
\def \voo{\vv{\o}}
\def \vd{\vec{d}}
\def \ve{\vec{e}}
\def \vi{\vec{i}}
\def \vj{\vec{j}}

\def \vvv{\vec{v}}
\def \vw{\vec{w}}

\def \np{\N\P}
\def \npt{\np^{(T)}}

\def \nplx{\np\L_{\X}}

\def \pla{\P\L_{A}}
\def \npla{\np\L_{A}}

\def \nplax{\np\L_{A,\X}}

\def \npm{\np_m}
\def \npmla{\npm\L_{A}}

\def \npmlax{\npm\L_{A,\X}}

\def \nptm{\npt_m}
\def \nptmla{\nptm\L_{A}}
\def \nptmlx{\nptm\L_{\X}}
\def \nptmlax{\nptm\L_{A,\X}}

\def \ams{A_m^{(s)}}
\def \bms{B_m^{(s)}}

\newcommand{\comm}[1]{}

\comm{
\begin{figure}[hbt!]
    \centering
    \begin{subfigure}[t]{.45\textwidth}
        \includegraphics[width=1\linewidth]{}
        \caption{}
        \label{fig:}
    \end{subfigure}\hspace{1cm}
    \begin{subfigure}[t]{.44\textwidth}
        \includegraphics[width=1\linewidth]{}
        \caption{}
        \label{fig:}
    \end{subfigure}
    \caption{}
    \label{fig:}
\end{figure}
}

\title{On Mixing Distributions Via Random Orthogonal Matrices and the Spectrum of the Singular Values of Multi-Z Shaped Graph Matrices}
\author{Wenjun Cai and Aaron Potechin}
\date{\today}

\begin{document}

\maketitle

\setlength{\parskip}{1.3mm}
\setlength{\baselineskip}{1.3em}

\begin{abstract}
In this paper, we determine the limiting distribution of the singular values of Z-shaped and multi-Z-shaped graph matrices as $n \to \infty$ for arbitrary input distributions whose odd moments are $0$. This extends the results of our previous paper \cite{CP20} which analyzed the limiting distribution of the singular values of Z-shaped and multi-Z-shaped graph matrices but only for the case when the input has $\pm{1}$ entries.

For our results, we use the following operation $\opr$ which mixes two distributions $\o$ and $\o'$ via a random orthogonal matrix. Given $\o$ and $\o'$, we take $\o \opr \o'$ to be the limit as $n \to \infty$ of the distribution of the singular values of $DRD'$ where $D$ and $D'$ are $n \times n$ diagonal matrices whose diagonal entries have distributions $\o$ and $\o'$ respectively and $R$ is a random $n \times n$ orthogonal matrix.

Using this operation, we show that the limiting distribution of the singular values of Z-shaped graph matrices can be described as follows. Let $\o_{Z}$ be the limiting distribution of the singular values for the Z-shaped graph matrix with $\pm{1}$ input entries (with the appropriate scaling). Given an arbitrary input distribution $\o$ with variance $1$ whose odd moments are $0$, the limiting distribution of the singular values for the Z-shaped graph matrix with input entries $\o$ is $\o_{Z} \opr \o$. We show that this result generalizes very nicely for multi-Z-shaped graph matrices.

Our results are closely connected to free probability theory and the combinatorics of non-crossing partitions. To analyze the operation $\opr$, we give a new direct formula for the moments of the product of two freely independent random variables. We also prove a number of results on non-crossing partitions which may be of independent interest.

\end{abstract}

Acknowledgement: This research was supported by NSF grant CCF-2008920. We thank Mark Sellke for pointing out the connection between our results and free probability theory.
\newpage
\tableofcontents

\section{Introduction}\label{section:introduction}

\setlength{\parskip}{1.5mm}
\setlength{\baselineskip}{1.3em}

Graph matrices are a type of matrix which plays a key role in analyzing the sum of squares hierarchy on average-case problems \cite{doi:10.1137/17M1138236,9317995,9719766,DBLP:journals/corr/abs-2011-09416,DBLP:journals/corr/abs-2011-04253} and is also useful for analyzing other methods involving higher moments \cite{AMP20, https://doi.org/10.48550/arxiv.2208.00122, https://doi.org/10.48550/arxiv.2208.09493, https://doi.org/10.48550/arxiv.2209.02655} 
(for background on graph matrices, see \cite{AMP20}).
However, the behavior of graph matrices is only partially understood. At present, we have rough norm bounds on graph matrices \cite{AMP20, 9719766, https://doi.org/10.48550/arxiv.2209.02655} which are sufficient for many applications, but we can hope to analyze graph matrices much more precisely.

In our previous paper \cite{CP20}, we took a first step towards analyzing graph matrices more precisely by determining the limiting distribution of the spectrum of the singular values of the Z-shaped graph matrix $\frac{1}{n}M_{\alpha_Z}$ as $n \to \infty$. This result can be seen as an analogue of Wigner's Semicircle Law \cite{wigner1958distribution, wigner1993characteristic}. However, unlike Wigner's Semicircle Law, which holds as long as the entries of the symmetric random matrix have mean $0$ and variance $1$, the results of our previous paper only hold if the input distribution for the spoke\footnote{the input distributions for the other two edges can be arbitrary as long as they have mean $0$ and variance $1$} is $\{-1,1\}$. This raised a natural question. Can we determine the spectrum of the singular values of $\frac{1}{n}M_{\alpha_Z}$ for other input distributions?

In this paper, we resolve this question using an operation $\opr$ (see Definition \ref{defn:circ-R}) which mixes two distributions $\o$ and $\o'$ via a random orthogonal matrix. We show that if we take $\Omega_{Z}$ to be the limiting distribution of the singular values of $\frac{1}{n}M_{\alpha_Z}$ where the input distribution is $\{-1,1\}$ then for any distribution $\Omega$ with variance $1$ whose odd moments are $0$, the limiting distribution of the singular values of $\frac{1}{n}M_{\alpha_Z}$ where the input distribution is $\Omega$ is $\Omega_{Z} \opr \Omega$. We further show that this result generalizes very nicely for multi-Z shaped graph matrices.

This operation $\opr$ is closely connected to free probability theory. In particular, for any distributions $\o$ and $\o'$, if $x$ and $y$ are two freely independent random variables with the same moments as ${\o}^2$ and ${\o'}^2$ then the kth moment of $x * y$ is equal to the $2k^{th}$ moment of $\o \opr \o'$. To determine the even moments of $\o \opr \o'$, we prove a new direct formula (see \Cref{thm:main-free-prob}) for the moments of the product of two freely independent random variables.

To prove our results, we express the moments of our distributions in terms of various types of non-crossing partitions and then analyze the combinatorics of these non-crossing partitions. Non-crossing partitions, which are also important for analyzing free probability theory (see \cite{nica2006lectures}), were first studied systematically by Kreweras \cite{KREWERAS1972333} and Poupard \cite{POUPARD1972279}. In particular, Kreweras showed that the number of non-crossing partitions of $[n]$ with $\alpha_i$ parts of size $i$ is $\frac{n!}{(n - a)!}\prod\limits_{j=1}^{n}{\frac{1}{{\alpha_j}!}}$ where $a = \sum\limits_{j=1}^{n}{\alpha_j}$ is the total number of parts. For a simple and direct proof of this result, see Liaw, Yeh, Hwang, and Chang \cite{liaw1998simple}. In our analysis, we prove various generalizations of this result where there are further restrictions on the non-crossing partitions. 

\subsection{Free Probability Theory and the Moments of the Product of Two Freely Independent Random Variables}\label{subsection:intro-free-prob}

\setlength{\parskip}{1.5mm}
\setlength{\baselineskip}{1.3em}
In order to state our results, we need some definitions for graph matrices and free probability theory. We start by briefly introducing free probability theory and giving our formula, \Cref{thm:main-free-prob}, for the moments of the product of two freely independent random variables. We give further results on free probability theory which are needed for our analysis in Section \ref{subsection:prelim-free-prob}.

\begin{defn}\label{defn:noncom-prob-space}
    We say that $\(\A,\varphi\)$ is a \emph{non-commutative probability space} if $\A$ is a unital algebra over $\Cbb$ and $\varphi:\A\to \Cbb$ is a unital linear functional.
\end{defn}
In this paper, we focus on the following type of non-commutative probability space. For further examples of non-commutative probability spaces, see \cite{nica2006lectures}.
\begin{eg}
$\(M_n\(\Cbb\), \nortrace\)$ is a non-commutative probability space where $M_n(\Cbb)$ is the algebra of $n \times n$ complex matrices and $\nortrace$ is the normalized trace $\nortrace(A) = \frac{1}{n}\(\sum_{j=1}^{n}{A_{jj}}\)$.
\end{eg}
        
        

\begin{defn}[Free Independence]\label{defn:free-independence}
    Let $(\A,\varphi)$ be a non-commutative probability space. Let $\{\A_i:i\in \I\}$ be a collection of unital subalgebras of $\A$. We say that $\left\{\A_i: i\in\I \right\}$ are \emph{freely independent} if for any $k\in \Zbb^{+}$, $\varphi(a_1\dots a_k) = 0 $ whenever 
    \begin{enumerate}
        \item $a_j\in \A_{i_j}$ for all $j\in [k]$
        \item $\varphi(a_j) = 0$ for all $j\in [k]$
        \item $i_j\neq i_{j+1}$ for any $j\in [k-1]$
    \end{enumerate}
    
    Let $X_i \subset \A$ for $i\in \I$. We say that $\{X_i\}_{i\in\I}$ are \emph{freely independent} if $\{\A_i\}_{i\in \I}$ are freely independent where $\A_i := \text{alg}\(1, X_i\)$ is the unital algebra generated by the set $X_i$. 
    
    Let $a_i\in \A$ for $i\in\I$. We say that $\{a_i\}_{i\in\I}$ are \emph{freely independent} if $\{\A_i\}_{i\in \I}$ are freely independent where $\A_i := \text{alg}(1,a_i)$ is the unital algebra generated by $a_i$.
\end{defn}
Given freely independent random variables $a$ and $b$, the moments of the product of $a$ and $b$ (i.e., $\{\varphi((ab)^k): k \in \mathbb{N}\}$) are completely determined by the moments of $a$ and $b$. These moments can be computed directly using the definition of free independence but this quickly gets complicated.
\begin{eg}\label{eg:phi-a-b}
If $\(\A,\varphi\)$ is a non-commutative probability space and $a,b\in \A$ are freely independent then: 
\begin{enumerate}
    \item $\varphi(ab) = \varphi(a)\varphi(b)$ since $0 = \varphi\(\(a-\varphi(a)\) \(b-\varphi(b)\)\) = \varphi\(ab - \varphi(a)b - a\varphi(b) + \varphi(a)\varphi(b)\) = \varphi(ab) - \varphi(a)\varphi(b)$. 
    More generally, $\varphi(a^mb^n) = \varphi(a^m)\varphi(b^n)$.
    
    \item $\varphi(aba) = \varphi(a^2)\varphi(b)$ since
    \begin{align*}
       0 &= \varphi\(\(a-\varphi(a)\) \(b-\varphi(b)\) \(a-\varphi(a)\)\) \\
         &= \varphi\(a\(b-\varphi(b)\) \(a-\varphi(a)\)\) - \varphi(a)\varphi\( \(b-\varphi(b)\) \(a-\varphi(a)\)\)\\
         &= \varphi(aba) - \varphi(a)\varphi(ba) - \varphi(b)\varphi(a^2) + \varphi(a)\varphi(b)\varphi(a)-0\\
         &= \varphi(aba) -  \varphi(a^2)\varphi(b)
    \end{align*}
    
    \item $\varphi(abab) = -\varphi(a)^2\varphi(b)^2 + \varphi(a^2)\varphi(b)^2 + \varphi(a)^2\varphi(b^2)$ since
    \begin{align*}
        0 &= \varphi\(\(a-\varphi(a)\) \(b-\varphi(b)\) \(a-\varphi(a)\) \(b-\varphi(b)\)\) \\
          &= \varphi\(a\(b-\varphi(b)\) \(a-\varphi(a)\) \(b-\varphi(b)\)\) - \varphi(a)\varphi\(\(b-\varphi(b)\) \(a-\varphi(a)\) \(b-\varphi(b)\)\) \\
          &=\varphi\(abab\) - \varphi(b)\varphi\(aba\) -  \varphi(a)\varphi\(ab^2\) + \varphi(a)\varphi(b)\varphi\(ab\) - \varphi(b)\varphi\({a^2}b\) + \varphi(b)^{2}\varphi\({a^2}\)\\
          &+\varphi(a)\varphi(b)\varphi\(ab\) - \varphi(a)^2\varphi(b)^2 - 0\\
          &=\varphi\(abab\) - \varphi(b)^{2}\varphi\(a^2\) -  \varphi(a)^2\varphi\(b^2\) + \varphi(a)^2\varphi(b)^2 - \varphi(b)^2\varphi\(a^2\) + \varphi(b)^{2}\varphi\(a^2\)\\
          &+\varphi(a)^2\varphi(b)^2- \varphi(a)^2\varphi(b)^2\\
          &= \varphi(abab) + \varphi(a)^2\varphi(b)^2 - \varphi(a)^2\varphi(b^2) - \varphi(a^2)\varphi(b)^2
    \end{align*}
\end{enumerate}
\end{eg}
We now give our direct formula for the moments of the product of two freely independent random variables.
\begin{defn}
    Let $\displaystyle P_k = \left\{\vaa=\(\a_1,\dots,\a_k\)\subseteq [k]^k: \sum_{i=1}^k i\cdot\a_i = k\right\}$.
\end{defn}

\begin{defn}
    Let $\(\A,\varphi\)$ be a non-commutative probability space. Given $\vaa\in P_k$ and $a\in\A$, we denote $\vec{\varphi_a}^{\vaa}$ to be 
    \begin{equation}
        \vec{\varphi_a}^{\vaa}
        = \prod_{i=1}^k\, \varphi(a^{i})^{\a_i} 
        = \varphi(a)^{\a_1}\varphi(a^2)^{\a_2}\dots \varphi(a^k)^{\a_k}
    \end{equation}
\end{defn}

\begin{defn}\label{defn:C-vaa-vbb}
    Given $\vaa,\vbb\in P_k$, we denote $C\(\vaa,\vbb\)$ to be
    \begin{equation}
        C\(\vaa,\vbb\) = (-1)^{A+B-k-1}\cdot k\cdot \binom{A+B-2}{k-1} \cdot \dfrac{(A-1)!}{\a_1!\dots\a_k!}\cdot \dfrac{(B-1)!}{\b_1!\dots\b_k!}\;.
    \end{equation}
    where $A = \a_1+ \dots + \a_k$ and $B = \b_1+ \dots + \b_k$.
    
    Note that when $A+B\leq k$, $C\(\vaa,\vbb\) = 0$.
\end{defn}

\begin{restatable}{thm}{mainfreeprob}
\label{thm:main-free-prob}
    If $\(\A,\varphi\)$ is a non-commutative probability space and $a,b\in\A$ are freely independent then
    \begin{equation}
        \varphi\((ab)^k\) = \sum_{\vaa, \vbb\in P_k} C\(\vaa,\vbb\)\cdot \vec{\varphi_a}^{\vaa} \cdot \vec{\varphi_b}^{\vbb}
    \end{equation}
\end{restatable}

\begin{eg}
    We list below the first few cases when $k=1,2,3,4$.
    \begin{enumerate}
        \item $\varphi(ab) = \varphi(a)\varphi(b)$.
        \item $\varphi(abab) = -\varphi(a)^2\varphi(b)^2 + \varphi(a)^2\varphi(b^2) + \varphi(a^2)\varphi(b)^2$.
        \item $\varphi\(ababab\) = 2\varphi(a)^3\varphi(b)^3 - 3\varphi(a)^3\varphi(b)\varphi(b^2) - 3\varphi(a)\varphi(a^2)\varphi(b)^3 + 3\varphi(a)\varphi(a^2)\varphi(b)\varphi(b^2) + \varphi(a)^3\varphi(b^3) + \varphi(a^3)\varphi(b)^3$.
        \item $\varphi\(abababab\) = -5\varphi(a)^4\varphi(b)^4 + 10\varphi(a)^4\varphi(b)^2\varphi(b^2) + 10 \varphi(a)^2\varphi(a^2)\varphi(b)^4 - 4\varphi(a)^4\varphi(b)\varphi(b^3) - 4\varphi(a)\varphi(a^3)\varphi(b)^4 - 2\varphi(a)^4\varphi(b^2)^2 - 2\varphi(a^2)^2\varphi(b)^4 + \varphi(a)^4\varphi(b^4) + \varphi(a^4)\varphi(b)^4 - 16\varphi(a)^2\varphi(a^2)\varphi(b)^2\varphi(b^2) + 4\varphi(a)^2\varphi(a^2)\varphi(b)\varphi(b^3) + 4\varphi(a)\varphi(a^3)\varphi(b)^2\varphi(b^2) + 2\varphi(a)^2\varphi(a^2)\varphi(b^2)^2 + 2\varphi(a^2)^2\varphi(b)^2\varphi(b^2)$.
    \end{enumerate}
    
     Note that we computed the first two results directly in \Cref{eg:phi-a-b}.
\end{eg}

\subsection{The \texorpdfstring{$\opr$}{ºR} Operation and its Connection to Free Probability Theory}\label{subsection:oproperation}

\setlength{\parskip}{1.5mm}
\setlength{\baselineskip}{1.3em}
We now formally define the operation $\opr$ and describe its connection with the moments of the product of two freely independent random variables.
\begin{defn}
    We say $D$ is an \emph{$\o$-distribution diagonal matrix} if $D$ is diagonal and each diagonal entry of $D$ is drawn independently from the distribution $\o$. 
\end{defn}
\begin{defn}
    Given $n \in \mathbb{N}$, we choose a random orthogonal matrix $R$ as follows:
    \begin{enumerate}
        \item Choose the first row $\vec{r}_1$ to be a random vector in $S^{n-1} = \left\{\vec{x} \in \mathbb{R}^n: \norm{\vec{x}} = 1 \right\}$. One way to do this is to choose the entries of $\vec{r}_1$ to be random Gaussian variables and then rescale $\vec{r}_1$ so that $\norm{\vec{r}_1} = 1$.
        \item For each $j \in \{2,\ldots,n\}$, choose the jth row $\vec{r}_j$ to be a random vector in the $(n-j)$-dimensional sphere $\left\{\vec{x} \in \mathbb{R}^n: \norm{\vec{x}} = 1, \vec{x} \cdot \vec{r}_{i} = 0 \text{ for all } i \in [j-1] \right\}$. One way to do this is to choose the entries of $\vec{r}_j$ to be random Gaussian variables, remove the components of $\vec{r}_j$ which are parallel to $\vec{r}_1,\ldots,\vec{r}_{j-1}$, and then rescale $\vec{r}_j$ so that $\norm{\vec{r}_j} = 1$.
    \end{enumerate}
\end{defn}

\begin{defn}\label{defn:circ-R}
    Given two distributions $\o$ and $\o'$, we define $\o \opr \o'$ to be the limiting distribution as $n \to \infty$ of the singular values of the random matrix $M=DRD'$ where $D$ and $D'$ are $n\times n$ $\o$ and $\o'$-distribution diagonal matrices, respectively, and $R$ is a random $n\times n$ orthogonal matrix. 
\end{defn}

\begin{restatable}{prop}{main}
\label{prop:main-0}
    The operation $\opr$ is commutative and associative.
\end{restatable}

\begin{proof}
$ $
    \begin{enumerate}
        \item $\o\opr\o' = \o'\opr\o$: Let $D$ and $D'$ be $\o$ and $\o'$-distribution diagonal matrices and $R$ be a random orthognal matrix. $\(DRD'\)^T = D'R^TD$ where $R^T$ is also a random orthogonal matrix. Since $M$ and $M^T$ have the same singular values, $\{DRD':R \text{ is a random orthogonal matrix}\}$ and $\{D'RD:R \text{ is a random orthogonal matrix}\}$ have the same singular value distributions.
        
        \item $\(\o\opr\o'\)\opr\o'' = \o\opr\(\o'\opr\o''\)$: Let $D$, $D'$ and $D''$ be $\o$, $\o'$ and $\o''$-distribution diagonal matrices. Let $D_{\o\opr\o'}$ and $D_{\o'\opr\o''}$ be $\o\opr\o'$ and $\o'\opr\o''$-distribution diagonal matrices. Since $D_{\o\opr\o'}$ and $DRD'$ have the same singular values distributions in the limit as $n \to \infty$, 
        $\left\{ R_1D_{\o\opr\o'}R_2: R_1, R_2\text{ orthogonal} \right\} = \left\{\hat{R_1}\(DRD'\)\hat{R_2}: \hat{R_1}, \hat{R_2}\text{ orthogonal}\right\}$.
        
        We have that in the limit as $n \to \infty$,
        \begin{align*}
             &\left\{R_1\(D_{\o\opr\o'}R_2D''\)R_3: R_1, R_2, R_3\text{ orthogonal} \right\} \\
            =& \left\{\(R_1D_{\o\opr\o'}R_2\)D''R_3: R_1, R_2, R_3\text{ orthogonal} \right\}\\
            =& \left\{\(\hat{R_1}\(DRD'\)\hat{R_2}\)D''R_3: R, \hat{R_1}, \hat{R_2}, R_3\text{ orthogonal} \right\}\\
            =& \left\{R_1DR_2D'R_3D''R_4: R_1, R_2, R_3, R_4\text{ orthogonal} \right\}\\
            =& \left\{R_1D\(R_2\(D'R_3D''\)R_4\): R_1, R_2, R_3, R_4\text{ orthogonal} \right\}\\
            =& \left\{R_1D\(\hat{R_2}D_{\o'\opr\o''}\hat{R_4}\): R_1, \hat{R_2}, \hat{R_4}\text{ orthogonal}\right\}\\
            =& \left\{R_1\(DR_2D_{\o'\opr\o''}\)R_4: R_1, R_2, R_4\text{ orthogonal} \right\}
        \end{align*}
        
        For any matrix $M$, $M$ and $R_1MR_2$ have the same singular values if $R_1$ and $R_2$ are orthogonal. Thus, in the limit as $n \to \infty$, $\left\{D_{\o\opr\o''}RD'': R\text{ orthogonal}\right\}$ and $\left\{DRD_{\o'\opr\o''}: R\text{ orthogonal}\right\}$ have the same singular value distributions.
    \end{enumerate}
\end{proof}
To see why the $\opr$ operation is connected with free probability theory, we make the following observation.
\begin{prop}
Let $R_n$ be an $n\times n$ random orthogonal matrix, $D_n$ and $D_n'$ be $n\times n$ random diagonal matrices where the diagonal elements are drawn independently from distributions $\Omega$ and $\Omega'$, respectively. Let $M_n=D_n R_n D_n'$, $A_n = {D_n}^2$, $B_n = R_n{D_n'}^2R_n^T$, and $\varphi = \nortrace\otimes \Ebb$. Then
\begin{equation}
 \dfrac{\,1\,}{n}\cdot \Ebb\[\trace\(\(M_nM_n^T\)^k\)\] = \varphi\(\({A_n}B_n\)^{k}\) 
\end{equation}
\end{prop}

\begin{proof}
Observe that
    \begin{align*}
         \dfrac{\,1\,}{n}\cdot \Ebb\[\trace\(\(MM^T\)^k\)\] 
        =& \Ebb\[\dfrac{\,1\,}{n}\trace\(DR{D'}^2R^T D^2 R{D'}^2R^T \dots D^2 R{D'}^2R^T D\)\] \\
        =& \Ebb\[\dfrac{\,1\,}{n}\trace\(D^2 \(R{D'}^2R^T\) D^2 \(R{D'}^2R^T\) \dots D^2\(R{D'}^2R^T\)\)\]\\
        =& \Ebb\[\dfrac{\,1\,}{n}\trace\( A_nB_n A_nB_n\dots A_n B_n\)\]\\
        =& \varphi\(\(A_nB_n\)^{k}\)
    \end{align*}
where the second equality uses the fact that $\trace(AB) = \trace(BA)$.
\end{proof}

As we discuss in \Cref{subsection:prelim-free-prob}, an important result in free probability theory is that in the limit as $n \to \infty$, the variables $A_n$ and $B_n$ are freely independent. Together with \Cref{thm:main-free-prob}, this gives the following corollary.

\begin{defn}
    Given a distribution $\o$ and $\vaa\in P_k$, we denote
    \begin{enumerate}
        \item $\o_{2m}=\Ebb_{x\sim \o}\[\; x^{2m}\]$, and
        \item $\voo^{\vaa} = \o_2^{\a_1}\dots\o_{2k}^{\a_k}$.
    \end{enumerate}
\end{defn}


\begin{restatable}{cor}{mainonecirc}\label{thm:main-1-circ}
Let $\o$ and $\o'$ be two distributions. Then for all $k\in\Nbb$,
    \begin{equation}
        \(\o\opr\o'\)_{2k} = \sum_{\vec{\a}, \vec{\b}\in P_k} C\(\vaa,\vbb\)\cdot \voo^{\,\vaa}\cdot \vv{\o'}^{\,\vbb}
    \end{equation}
    where if we let $a=a_1+\dots+a_k$ and $b=b_1+\dots+b_k$, then
    \begin{equation}
        C\(\vaa,\vbb\) = (-1)^{k+a+b-1}\cdot k\cdot \binom{a+b-2}{k-1} \cdot \dfrac{(a-1)!}{\a_1!\dots\a_k!}\cdot \dfrac{(b-1)!}{\b_1!\dots\b_k!}\;.
    \end{equation}
\end{restatable}

\subsection{Graph Matrices}\label{subsection:intro-graph-matrices}

\setlength{\parskip}{1.5mm}
\setlength{\baselineskip}{1.3em}

We now define graph matrices and state our main results on the limiting distribution of the singular values of multi-Z-shaped graph matrices. These definitions are equivalent to the definitions of graph matrices in \cite{AMP20} except that instead of having a single random input matrix, we instead have a random input matrix $M^{e}$ for each edge $e$ of our shape.
\begin{defn}[Shapes]
We define a \emph{shape} $\alpha$ to be a graph with vertices $V(\alpha)$, edges $E(\alpha)$, and distinguished tuples of vertices $U_{\alpha} = \(u_1,\ldots,u_{|U_{\alpha}|}\)$ and $V_{\alpha} = \(v_1,\ldots,v_{|V_{\alpha}|}\)$.
\end{defn}
In this paper, we associate each shape $\a$ with a set of distributions $\oa=\{\o_e:e\in E(\a)\}$, one for each edge of $\a$, where each $\o_e$ has mean $0$ and variance $1$.
\begin{defn}[Bipartite Shapes]
We say that a shape $\alpha$ is \emph{bipartite} if $U_{\alpha} \cap V_{\alpha} = \emptyset$, $V(\alpha) = U_{\alpha} \cup V_{\alpha}$, and all edges in $E(\alpha)$ are between $U_{\alpha}$ and $V_{\alpha}$.
\end{defn}

\begin{defn}[Random Input Matrix $M^e$]
    Let $\a$ be a shape associated with distributions $\oa$. For each $e\in E(\a)$, we define the \emph{random input matrix $M^{e}$} to be a symmetric random matrix where each entry $M^{e}(i,j)$ is drawn independently from $\o_e$.
\end{defn}

\begin{defn}[Fourier Characters $\chi_{\sigma\(E(\a)\)}$]
    Let $\a$ be a shape associated with distributions $\oa$. Given an injective map $\sigma: V(\alpha) \to [n]$, we define
    \begin{equation}
        \chi_{\sigma\(E(\a)\)} = \prod_{e=(u,v) \in E(\a)} M^{e}\(\sigma(u),\sigma(v)\).
    \end{equation}
\end{defn}

\begin{defn}[Graph Matrices with input distributions $\oa$]\label{def:graphmatrices}
    Given a shape $\a$ associated with distributions $\oa$, we define $M_{\alpha,\oa}$, the \emph{graph matrix with input distributions $\oa$}, to be the $\frac{n!}{\(n-\abs{U_{\a}}\)!}\times\frac{n!}{\(n-\abs{V_{\a}}\)!}$ matrix with rows indexed by tuples of $\abs{U_{\alpha}}$ distinct vertices, columns indexed by tuples of $\abs{V_{\alpha}}$ distinct vertices, and entries
\begin{equation}
    M_{\alpha,\oa}(A,B) = \sum_{\substack{\sigma: V(\alpha) \to [n]: \sigma \text{ is injective}, \\ \sigma(U_{\alpha}) = A, \sigma(V_{\alpha}) = B }} \chi_{\sigma\(E(\a)\)}.
\end{equation}
\end{defn}
\begin{rmk}
In this paper, we only consider bipartite shapes so for all $A$ and $B$, we will either have that $A \cap B = \emptyset$ and there is exactly one injective map $\sigma$ such that $\sigma(U_{\alpha}) = A$ and $\sigma(V_{\alpha}) = B$ or $A \cap B \neq \emptyset$ and $M_{\alpha,\oa}(A,B) = 0$ as there are no such injective maps $\sigma$.
\end{rmk}
\begin{defn}[Graph Matrices with a single input distribution $\o$]
    If $\o_e=\o$ for all $e\in E(\a)$ then we denote the corresponding graph matrix $M_{\a,\oa}$ as $M_{\a,\o}$. Note that even though all $\o_e$ are the same, there is an independent input matrix $M^{e}$ with the same underlying distribution $\o$ for each edge $e$.
\end{defn}

\begin{defn}[$\pm 1$ distribution]
    We define the \emph{$\pm1$ distribution} $\o_{\pm1}$ to be the distribution which is $1$ with probability $1/2$ and $-1$ with probability $1/2$.
\end{defn}

\begin{defn}
    Given a shape $\a$, if $\oa$ consists of only $\o_{\pm1}$, then we denote $M_{\a,\oa}$ as just $M_{\a}$. In other words, we take graph matrices to have input distributions $\o_{\pm 1}$ by default.
\end{defn}



\begin{defn}\label{defn:m-zshape}
    Let $\mzshape$ be the bipartite shape with vertices $V\(\mzshape\)=\{u_1,\dots,u_m,v_1,\dots, v_m\}$ and edges $E\(\mzshape\)=\left\{\{u_i,v_i\}: i\in[m]\right\}\cup \left\{\{u_{i+1},v_{i}\}: i\in[m-1]\right\}$ with distinguished tuples of vertices $U_{\mzshape}=(u_1,\dots,u_m)$ and $V_{\mzshape}=(v_1,\dots, v_m)$. We call the slanted edges $\left\{\{u_{i+1},v_{i}\}: i\in[m-1]\right\}$ \emph{spokes} and call $\{u_{i+1},v_{i}\}$ the $i^{th}$ spoke. See \Cref{fig:m-zshape} for an illustration.
    
    We refer to $\mzshape$ as the \emph{$m-$layer Z-shape} or the \emph{Z(m)-shape}. When $m=2$, the shape is exactly a Z shape.
    
    
    \begin{figure}[hbt!]
        \centering
        \includegraphics[scale=0.35]{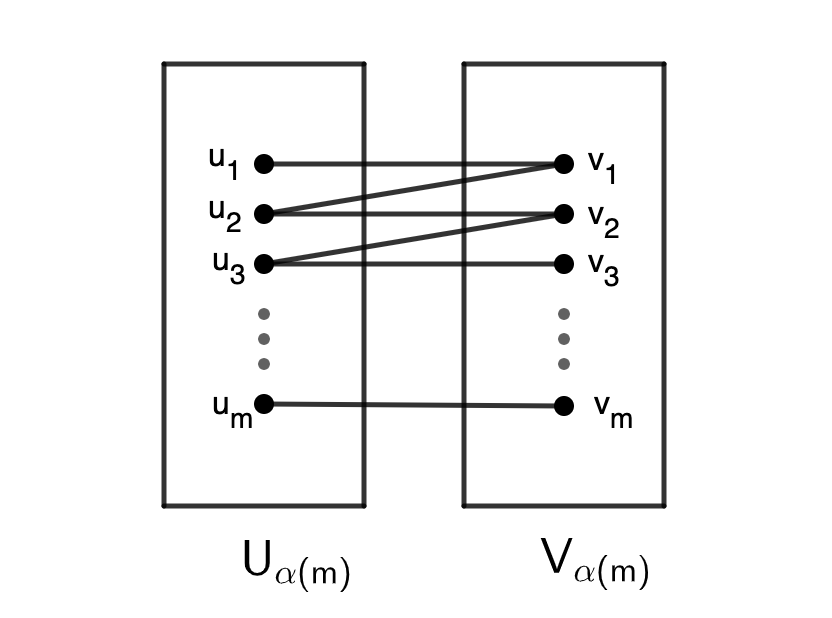}
        \caption{The $m-$layer Z-shape $\mzshape$}
        \label{fig:m-zshape}
    \end{figure}
\end{defn} 

\begin{defn}\label{defn:mzshape-sdistr-graph-matrix}
    Given $\o^{(1)},\dots,\o^{(s)}$ and $m \geq s+1$, we define
    \begin{enumerate}
        \item $\oams$ to be the set of distributions associated with $\mzshape$ where the $i^{th}$ spoke has distribution $\o^{(i)}$ for $i=1,\dots,s$ and all other edges have distribution $\o_{\pm 1}$, and
        \item $M_{\mzshape,s}$ to be the graph matrix with shape $\mzshape$ and input distributions $\oams$, and
        \item $\displaystyle \mzshapesdistrg = \dfrac{1}{n^{m/2}}\cdot M_{\mzshape,s}$ to be the normalized graph matrix.
    \end{enumerate}
\end{defn}

\begin{defn}\label{defn:ozm-ozms}
    Given $\o^{(1)},\dots,\o^{(s)}$ and $m \geq s+1$, we define
    \begin{enumerate}
        \item $\ozm$ to be the limiting distribution of the singular values of $\displaystyle \dfrac{1}{n^{m/2}}\cdot M_{\mzshape}$, or equivalently, $M^{(G)}_{Z(m),0}$.
        \item $\ozms$ to be the limiting distribution of the singular values of $\mzshapesdistrg$.
    \end{enumerate}
\end{defn}

\comm{
\begin{defn}\label{defn:mzshape-sdistr-grid-matrix}
    Given $\o^{(1)},\dots,\o^{(s)}$, consider the random matrix $\mzshapesdistr$ obtained inductively on $s$ as the following:
\begin{enumerate}
    \item $M_{Z(m),0} = D$ where $\o=\ozm$.
    \item $M_{Z(m),1} = D'RD$ where $\o'=$ the limiting distribution of singular values of $M_{Z(m),0}$ and $\o=\o^{(1)}$. 
    \item $M_{Z(m),i} = D'RD$ where $\o'=\o_{Z(m),\o^{(1)},\dots,\o^{(i-1)}}:=$ the limiting distribution of singular values of $M_{Z(m),i-1}$ and $\o=\o^{(i)}$ for $i=1,\dots,s$.
\end{enumerate}
\end{defn}
}

With these definitions, we can now state our main results which show that the spectrum of the singular values of multi-Z-shaped graph matrices can be easily described using the $\opr$ operation.
\begin{defn}
    We define $\mathbf{\o}$ to be the set of all random distributions and $\oset$ to be the set of distributions $\o$ with odd moments $0$ and variance $1$. 
\end{defn}
\begin{defn}
A distribution $\Omega$ satisfies Carleman's condition if $ \sum\limits_{k=1}^{\infty}{\beta_{2k}^{-\frac{1}{2k}}} = \infty$ where $\beta_{2k} = \Ebb_{\Omega}\[x^{2k}\]$ is the 2k-th moment of $\Omega$.
\end{defn}
\begin{restatable}{thm}{maintwocirc}
\label{thm:main-2-circ}
For all $\o^{(1)},\dots,\o^{(s)}\in \oset$ which satisfy Carleman's condition,
    \begin{equation}
        \o_{Z(m)}\opr \o^{(1)}\opr\dots \opr \o^{(s)} = \ozms
    \end{equation}
\end{restatable}

\begin{restatable}{thm}{mainthreecirc}
\label{thm:main-3-circ}
For any $m,m'\in\Nbb$,
\begin{equation}
    \ozm\opr \ozmm = \o_{Z\(m+m'\)}
\end{equation}
\end{restatable}
\begin{rmk}
Carleman's condition is needed as it ensures that the resulting distribution can be recovered from its moments.
\end{rmk}

\section{Preliminaries}\label{section:preliminary}

\setlength{\parskip}{1.5mm}
\setlength{\baselineskip}{1.3em}
Before proving our results, we need some results about non-crossing partitions and some results from free probability theory.
\subsection{Non-crossing Partitions} \label{subsection:prelim-np}

\setlength{\parskip}{1.5mm}
\setlength{\baselineskip}{1.3em}
We first discuss non-crossing partitions as they play a central role in our analysis. Many of the following definitions and theorems are from the lecture notes of Nica and Speicher \cite{nica2006lectures}.

\begin{defn}[Noncrossing Partitions]
$ $

    A \emph{partition $\pi$ of $[n]$} is $\pi=\lcurb P_1,\dots, P_m \rcurb$ where $\displaystyle\bigsqcup_{i=1}^m \; P_i = [n]$ for each $i\in[m]$.
    
    We say a partition $\pi = \lcurb P_1,\dots, P_m \rcurb$ is \emph{non-crossing} if there do not exist $a<b<c<d$ such that $a,c\in P_i$ and $b,d\in P_j$ for some $i\neq j\in[m]$.

    For $n\in\Nbb$, we denote $P(n)$ to be the set of partitions of $[n]$ and $NC(n)$ to be the noncrossing ones.
\end{defn}

\begin{rmk}
    Given a non-crossing partition $\pi$ of $[n]$, we can visualize the partition $\pi = \{P_1,\dots,P_m\}$ as placing polygons of size $\abs{P_i}$ within a cycle of length $n$ such that no two polygons intersect each other. When $\abs{P_i}=1$, the corresponding polygon is simply a point.
\end{rmk}

\begin{eg}
    Consider the partition $\pi=\left\{\{2\},\{4\},\{6\}, \{1,3\}, \{5,7,8\}\right\}$ of $[8]$. Then $\pi$ can be represented as placing the line $\{4,6\}$ and the triangle $\{5,7,8\}$ on the cycle $\C_8$ as shown in \Cref{fig:np-polygons}.
    
    \begin{figure}[hbt!]
        \centering
        \begin{subfigure}[t]{.4\textwidth}
            \includegraphics[width=1\linewidth]{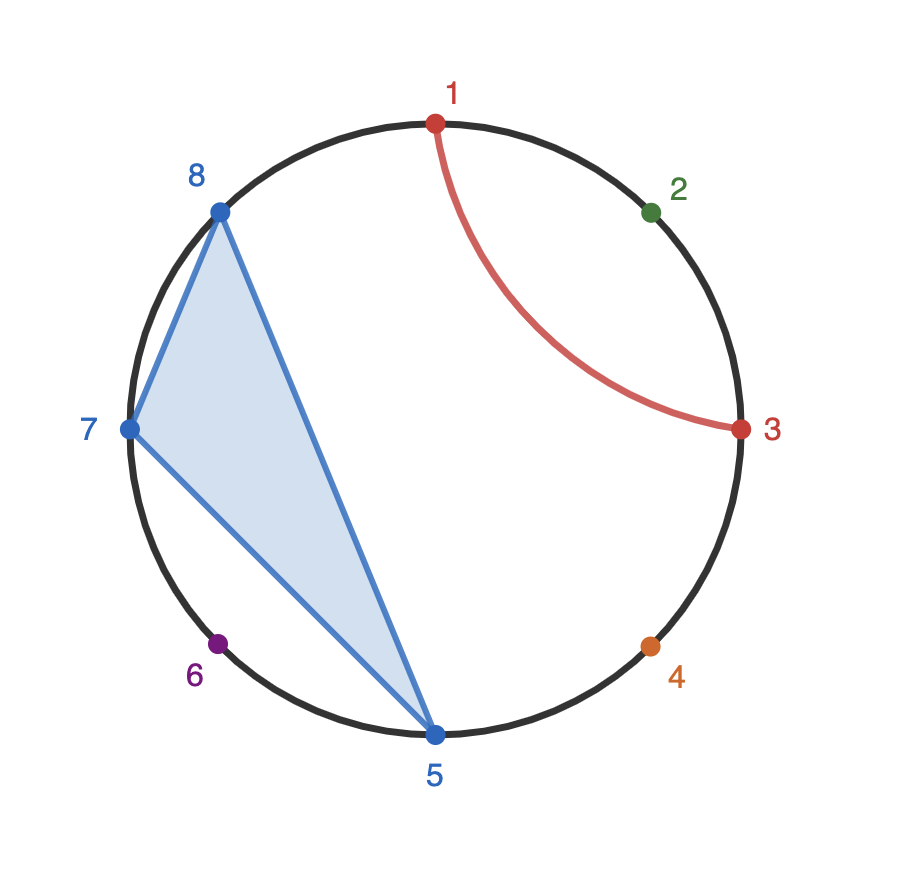}
            \caption{$\{1,3\}$ represents a line, and $\{5,7,8\}$ represents a triangle and all others are points.}
            \label{fig:np-polygons}
        \end{subfigure}\hspace{0.7cm}
        \begin{subfigure}[t]{.43\textwidth}
            \includegraphics[width=1\linewidth]{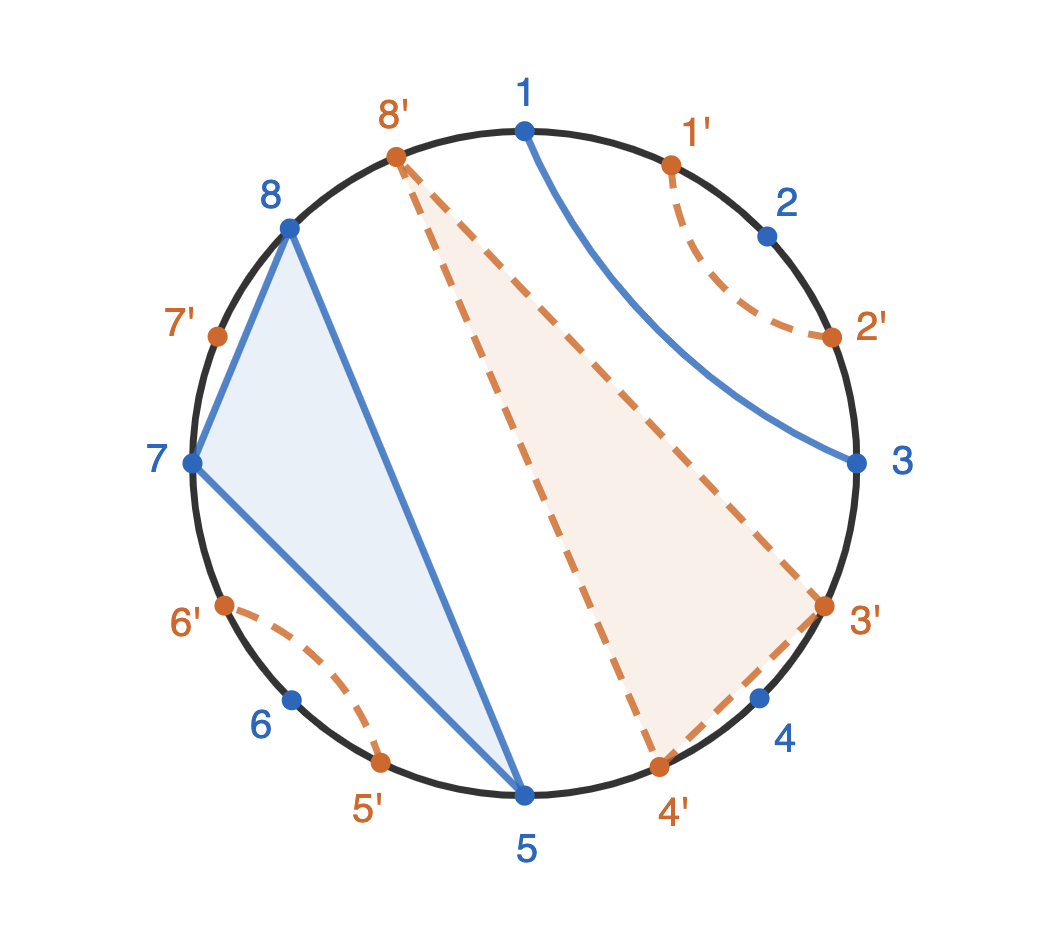}
            \caption{The blue solid parts are $\pi$, and the orange dashed parts are $K(\pi)$.}
            \label{fig:np-complement}
        \end{subfigure}
        \caption{Illustration of $\pi$ and $K(\pi)$. Here $\pi=\left\{\{2\},\{4\},\{6\}, \{1,3\}, \{5,7,8\}\right\}$.}
        \label{fig:np}
    \end{figure}
\end{eg}

\begin{defn}
    We say a partition $\pi = \lcurb P_1,\dots, P_m \rcurb$ is \emph{of parts with sizes $(n_1,\dots,n_m)$} if $\abs{P_i}=n_i$ for all $i\in[m]$.
    
    Let $\vaa = (\a_1,\dots,\a_n) \in P_n$. We denote $\np(\vaa)$ to be the set of noncrossing partitions of parts with sizes $\(\underbrace{1,\dots,1}_{\a_1\text{ times}}, \underbrace{2,\dots,2}_{\a_2 \text{ times }},\dots, \underbrace{n,\dots,n}_{\a_n \text{ times}}\)$. In other words, $\pi$ has $\a_i$ number of parts with size $i$ for each $i\in[n]$.
\end{defn}


\begin{defn}
    Let $\pi=\{V_1,\dots, V_r\}, \sigma=\{W_1,\dots, W_s\} \in NC(n)$. We define $\pi\leq \sigma$ if for each $i\in[r]$, $V_i\subset W_j$ for some $j\in[s]$. 
    
    For $\pi,\sigma\in NC(n)$, we denote $[\pi,\sigma] = \{\tau\in NC(n): \pi\leq \tau\leq \sigma\}$
\end{defn}

\begin{prop}
    $\(NC(n), \leq\)$ is a poset (partially ordered set). In particular, $1_n = \{\{1,2,\dots,n\}\}$ is the maximal element and $0_n = \{\{1\},\dots,\{n\}\}$ is the minimal element.
\end{prop}

\begin{defn}[Kreweras Complement Map, Definition 9.21 from \cite{nica2006lectures}]
    We define the \emph{Kreweras complement map $K: NC(n) \to NC(n)$} as the following: 
    \begin{enumerate}
        \item take $\pi\in NC(n)$.
        \item expand the vertex set $[n]$ into $\{1,\bar{1},2,\bar{2},\dots, n, \bar{n}\}$.
        \item $K(\pi) = \sigma$ where $\sigma$ is the biggest element in $NC(\bar{1},\dots,\bar{n}) \cong NC(n)$ where $\pi\cup \,\sigma \in NC(1,\bar{1},\dots, n,\bar{n})$.
    \end{enumerate}
\end{defn}

\begin{eg}
    Let $\pi=\left\{\{2\},\{4\},\{6\}, \{1,3\}, \{5,7,8\}\right\}$ of $[8]$ as in the previous example. Then $K(\pi) = \{\{1,2\}, \{3,4,8\}, \{5,6\}, \{7\}\}$. 
    
    See \Cref{fig:np-complement} for an illustration.
\end{eg}

\begin{prop}
    Let $\pi\leq \sigma\in NC(n)$. Then $[\pi, \sigma] \cong [K(\sigma), K(\pi)]$.
\end{prop}

\begin{thm}[Canonical Factorization, Theorem 9.29 from \cite{nica2006lectures}]\label{thm:canonical-factorization}
    Let $\pi\leq \sigma\in NC(n)$. Then there exists a canonical choice of $(k_1,\dots, k_n)\in \(\Zbb_{\geq 0}\)^{n}$ such that 
    \begin{equation}
        [\pi, \sigma] \cong NC(1)^{k_1}\times \dots\times NC(n)^{k_n}
    \end{equation}
\end{thm}

Following the proof for \Cref{thm:canonical-factorization}, the canonical factorization for a special interval $[0,\pi]$ is as follows. 
\begin{cor}
    Let $\pi\in NC(n)$. Assume $\pi\in\np(\vaa)$ for some $\vaa\in P_n$. Then the canonical factorization gives
    \begin{equation}
        [0,\pi] \cong NC(1)^{\a_1}\times \dots\times NC(n)^{\a_n} = \prod_{V\in\pi} NC\(\,\abs{V}\)
    \end{equation}
    and 
    \begin{equation}
        [\pi, 1_n] \cong \[0_n, K(\pi)\] \cong \prod_{V\in K(\pi)} NC\(\,\abs{V}\)
    \end{equation}
\end{cor}

\begin{eg}
    Let $\pi = \left\{\{2\},\{4\},\{6\}, \{1,3\}, \{5,7,8\}\right\}$ be as in previous example. Then $[0,\pi]\cong NC(1)^3\times NC(2)\times NC(3)$.
\end{eg}

As noted in the introduction, Kreweras \cite{KREWERAS1972333} proved the following result on the number of non-crossing partitions.

\begin{thm}\label{thm:num-of-np}
    Let $\vaa\in P_k$. Then 
    \begin{equation}
        \abs{\np\(\vaa\)} = \binom{k}{a-1}\cdot \dfrac{(a-1)!}{\a_1!\dots\a_k!}\,.
    \end{equation}
\end{thm}

We now give an extended definition of $\np\(\vaa\)$.
\begin{defn}
    Given $\vaa\in P_k$, we define $\np\(m\vaa\)$ to be the set of non-crossing partitions of $[mk]$ corresponding to augmented $\vaa$: there are $\a_i$ parts of size $mi$ in a partition $P\in\np\(m\vaa\)$.
\end{defn}

A direct corollary is the following.
\begin{cor}\label{cor:num-of-m*np}
    Let $\vaa\in P_k$. Then 
    \begin{equation}
        \abs{\np\(m\vaa\)} = \binom{mk}{a-1} \cdot\dfrac{(a-1)!}{\a_1!\dots\a_k!}\,.
    \end{equation}
\end{cor}
\begin{proof}
    This is an application of \Cref{thm:num-of-np}. Since $\np\(m\vaa\)$ is the set of non-crossing partitions of $[mk]$ with $a$ parts, and there are $\a_i$ number of parts with size $mi$, $\displaystyle\abs{\np\(m\vaa\)}= \binom{mk}{a-1}\cdot \dfrac{(a-1)!}{\a_1!\dots\a_k!}$, as needed.
\end{proof}

Using \Cref{cor:num-of-m*np}, the generalized Catalan number can be interpreted as the number of non-crossing partitions in the following way.
\begin{defn}
    Let $\displaystyle C\(k,m\)=\dfrac{1}{mk+1}\binom{(m+1)k}{k}$ denote the level-$m$ $k^{th}$ Catalan number. In particular, when $m=1$, $C(k,1)=C_k$ is the $k^{th}$ Catalan number.
\end{defn}

\begin{cor}\label{cor:m-catalan-nnp}
Let $\vaa = \(k,0,\dots,0\)\in P_k$. Then
\begin{equation}
    C(k,m) = \abs{\np\((m+1)\vaa\)}
\end{equation}

i.e. $C(k,m)$ is the number of non-crossing partitions of $[(m+1)k]$ where all the parts are of size $(m+1)$.
\end{cor}
\begin{proof}
    By \Cref{cor:num-of-m*np},
    \begin{align*}
        \abs{\np\((m+1)\vaa\)} 
        &= \binom{(m+1)k}{k-1}\cdot \dfrac{(k-1)!}{k!}
         = \dfrac{\((m+1)k\)!}{(k-1)!(mk+1)!}\cdot\dfrac{\,1\,}{k}
         = \dfrac{\((m+1)k\)!}{k!(mk)!}\cdot \dfrac{1}{mk+1}\\
        &= \binom{(m+1)k}{k}\cdot \dfrac{1}{mk+1} = C(k,m).
    \end{align*}
    as needed.
\end{proof}

\subsection{Results from Free Probability Theory} \label{subsection:prelim-free-prob}

\setlength{\parskip}{1.5mm}
\setlength{\baselineskip}{1.3em}


In this subsection we give several results from free probability theory which we need for our analysis. The definitions and theorems are from the lecture notes of Nica and Speicher \cite{nica2006lectures}.


\begin{defn}
    Let $P$ be a finite poset and denote $P^{(2)} = \{(\pi,\sigma)\in P\times P: \pi\leq \sigma\}$. Let $F,G: P^{(2)}\to \Cbb$. We define the \emph{convolution of $F$ and $G$, $F*G: P^{(2)}\to \Cbb$} to be
    \begin{equation}
        (F*G)(\pi, \sigma) = \sum_{\tau\in P:\, \pi\leq \tau\leq \sigma} F(\pi, \tau)G(\tau, \sigma)
    \end{equation}
\end{defn}

\begin{defn}
    Let $P$ be a finite poset. Define $\delta: P^{(2)}\to \Cbb$ to be 
    \begin{equation}
        \delta(\pi, \sigma) = 
        \begin{cases}
            \; 1 & \text{ if } \pi = \sigma\\
            \; 0 & \text{ if } \pi < \sigma
        \end{cases}
    \end{equation}
\end{defn}

\begin{prop}
    $\delta$ is a unit element for the convolution operation. i.e. $\delta * F = F*\delta = F$. 
\end{prop}

\begin{defn}
    Let $P$ be a finite poset. We define 
    \begin{enumerate}
        \item the \emph{Zeta function $\zeta: P^{(2)}\to \Cbb$} to be $\zeta(\pi, \sigma) = 1$ for any $\(\pi, \sigma\) \in P^{(2)}$
        \item the \emph{M\"{o}bius function} $\mu = \zeta^{-1}$ under the convolution operation. i.e. $\mu * \zeta = \zeta * \mu = \delta$.
    \end{enumerate}
\end{defn}

\begin{prop}
$ $
    \begin{enumerate}
        \item Let $P,Q$ be finite posets and $\Phi: P\to Q$ be an isomorphism. Then $\mu_P(\pi,\sigma) = \mu_Q\(\Phi(\pi), \Phi(\sigma)\)$.
        \item Let $P_1,\dots, P_n$ be finite posets and let $P= P_1\times\dots \times P_n$. Then 
        \begin{equation}
            \mu_P\( (\pi_1,\dots,\pi_n), (\sigma_1, \dots, \sigma_n)\) = \mu_{P_1}(\pi_1,\sigma_1)\dots \mu_{P_n}(\pi_n, \sigma_n).
        \end{equation}
    \end{enumerate}
\end{prop}

\begin{prop}\label{prop:mu-NC-n}
    Let $n\in\Nbb$ and let $\mu_n$ be the M\"{o}bius function for $NC(n)$. Then $\mu_n(0_n, 1_n) = (-1)^{n-1}C_{n-1}$ where $\displaystyle C_k = \dfrac{1}{k+1}\binom{2k}{k}$ is the $k^{th}$ Catalan number.
\end{prop}

\begin{cor}
    Let $\pi\leq \sigma\in NC(n)$. Then $\displaystyle \mu_n(\pi,\sigma) = \prod_{i=1}^{n} \((-1)^{i-1}C_{i-1}\)^{k_i}$ where $[\pi, \sigma] \cong NC(1)^{k_1}\times \dots\times NC(n)^{k_n}$ is the canonical factorization.
\end{cor}


\begin{defn}
   Let $\(\A, \varphi\)$ be a non-commutative probability space. For $n\in \Nbb$ and $V\subset [n]$, we define $\varphi(V)[a_1,\dots, a_n]$ to be
   \begin{equation}
       \varphi(V)[a_1,\dots, a_n] = \varphi(a_{i_1}\dots a_{i_m}) \text{ where } V= \{i_1,\dots, i_m\}\subset [n]
   \end{equation}
   
   Let $\pi = \{V_1,\dots, V_r\} \in NC(n)$, we further define $\varphi_{\pi}[a_1,\dots, a_n]$ to be
    \begin{equation}
        \varphi_\pi[a_1,\dots, a_n] = \prod_{i=1}^r \varphi(V_i)[a_1,\dots, a_n]
    \end{equation}
    
    When $\pi = 1_n$, we denote $\varphi_{1_n}[a_1,\dots, a_n]$ as $\varphi_n(a_1\dots a_n)$ which is the same as $\varphi(a_1\dots a_n)$.
\end{defn}

\begin{defn}[Free Cumulants]
    Let $\(\A, \varphi\)$ be a non-commutative probability space. For $n\in\Nbb$, and $\pi\in NC(n)$, the \emph{free cumulant} $\kappa_{\pi}$ is a multilinear functional $\kappa_\pi : \A^n\to \Cbb$ where 
    \begin{equation}
        \kappa_\pi[a_1,\dots,a_n] = \sum_{\sigma\in NC(n):\, \sigma\leq \pi} \varphi_{\sigma} [a_1,\dots, a_n]\cdot \mu(\sigma, \pi)
    \end{equation}
    
    In particular, when $\pi = 1_n$, we denote $\kappa_{1_n}[a_1,\dots, a_n]$ as $\kappa_n(a_1,\dots, a_n)$. We have that 
    \begin{equation}
        \kappa_n(a_1,\dots, a_n) = \sum_{\sigma\in NC(n)} \varphi_{\sigma} [a_1,\dots, a_n]\cdot \mu(\sigma, 1_n)
    \end{equation}
    
    Let $V=\{i_1,\dots, i_m\}\subset [n]$. We define $\kappa(V)[a_1,\dots, a_n] = \kappa_m(a_{i_1},\dots, a_{i_m})$.
\end{defn}

\begin{prop}[Proposition 11.4 of \cite{nica2006lectures}]
\label{prop:Mobius-Inversion}
Let $\(\A, \varphi\)$ be a non-commutative probability space. Then
    \begin{enumerate}
        \item $\displaystyle \kappa_{\pi}[a_1,\dots, a_n] = \prod_{V\in\pi} \kappa(V)[a_1,\dots, a_n]$
        \item (M\"{o}bius Inversion) $\displaystyle \varphi(a_1\dots a_n) = \sum_{\pi\in NC(n)} \kappa_\pi[a_1,\dots, a_n]$
    \end{enumerate}
\end{prop}

\begin{thm}[Theorem 14.4 of \cite{nica2006lectures}]
\label{thm:cumulants}
    Let $\(\A,\varphi\)$ be a non-commutative probability space. Let $a_1,\dots, a_n, b_1,\dots, b_n \in \A$ where $\{a_1,\dots,a_n\}$ and $\{b_1,\dots,b_n\}$ are freely independent. Then 
    \begin{equation}
        \varphi(a_1b_1\dots a_nb_n) = \sum_{\pi\in NC(n)} \kappa_{\pi}[a_1,\dots, a_n]\cdot \varphi_{K(\pi)} [b_1,\dots, b_n]
    \end{equation}
    
    In particular, we can apply this to $a_i=a$ and $b_i=b$ for all $i\in [n]$, then 
    \begin{equation}
        \varphi\((ab)^n\) = \varphi(ab\dots ab) = \sum_{\pi\in NC(n)} \kappa_{\pi}[a,\dots, a]\cdot \varphi_{K(\pi)} [b,\dots, b]
    \end{equation}
\end{thm}

\begin{defn}
    Let $\mu$ and $\nu$ be compactly supported probability measure on $\Rbb^{+}$. The \emph{multiplicative free convolution} $\mu \boxtimes \nu$ is the distribution of $\sqrt{x}y\sqrt{x}$ where $x,y$ are positive elements in some $C^*$-probability space, $x$ and $y$ are free and $x\sim \mu$, $y\sim \nu$.
\end{defn}

\begin{thm}[Theorem 23.14 of \cite{nica2006lectures}]\label{thm:free-independence-matrices}
    Let $(A_n)_{n\in\Nbb}$ and $(B_n)_{n\in\Nbb}$ be sequences of $n\times n$ matrices such that $A_n$ converges in distribution (with respect to tr) for $n\to\infty$ and $B_n$ converges in distribution (with respect to tr) for  $n\to\infty$. Furthermore, let $(U_n)_{n\in\Nbb}$ be a sequence of Haar unitary $n\times n$ random matrices. Then, $U_n A_n U_n^{*}$ and $B_n$ are asymptotically free for $n\to\infty$.
\end{thm}




\section{Moments for \texorpdfstring{$ab$}{ab} when \texorpdfstring{$a,b$}{a, b} are Freely Independent}
\label{section:moments}

\setlength{\parskip}{1.5mm}
\setlength{\baselineskip}{1.3em}

In this section, we prove \Cref{thm:main-free-prob}. For convenience, we restate \Cref{thm:main-free-prob} here.


\setlength{\parskip}{1.5mm}
\setlength{\baselineskip}{1.3em}

Recall that
\begin{equation}\nonumber
    C\(\vaa,\vbb\) = (-1)^{A+B-k-1}\cdot k\cdot \binom{A+B-2}{k-1} \cdot \dfrac{(A-1)!}{\a_1!\dots\a_k!}\cdot \dfrac{(B-1)!}{\b_1!\dots\b_k!}\;
\end{equation}
where $A = \a_1+ \dots + \a_k$ and $B = \b_1+ \dots + \b_k$.

\mainfreeprob*

\begin{defn}
    Let $\pi=\lcurb V_1,\dots, V_m \rcurb \in NC(k)$. Let $\C_k$ be the cycle graph with vertices $\{1,\dots,k\}$. Define $\C_k/\pi$ to be the graph obtained by identifying vertices together under $\pi$. 
\end{defn}

\begin{figure}[hbt!]
    \centering
    \includegraphics[scale=0.35]{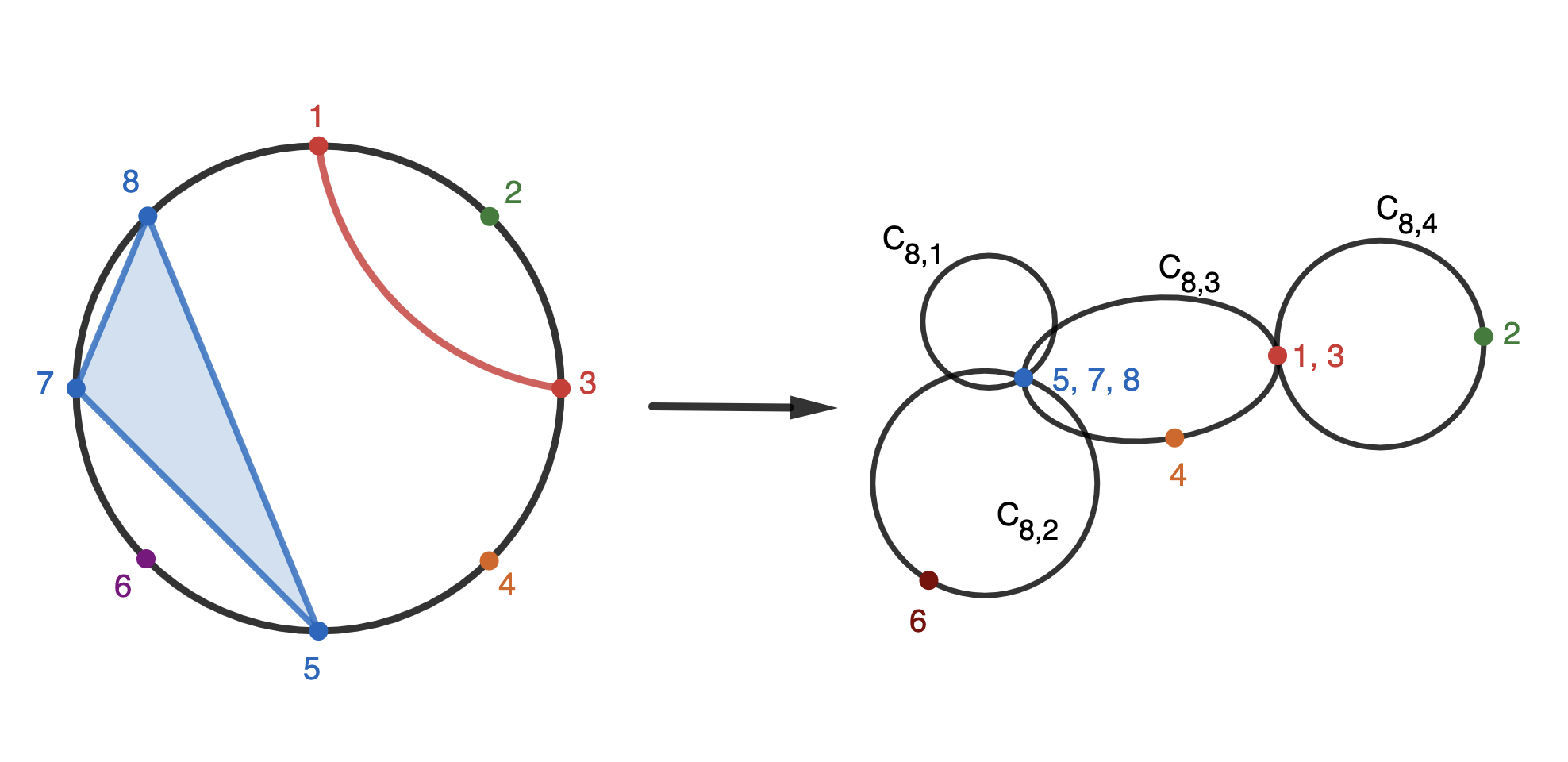}
    \caption{Illustration of \Cref{prop:noncross-partition-cycles} and \Cref{defn:np-cycle-size-Sp}: $\C_8/\pi = \C_{8,1}\cup \C_{8,2}\cup \C_{8,3}\cup \C_{8,4}$ where $\pi = \left\{\{2\},\{4\},\{6\},\{1,3\}, \{5,7,8\}\right\}$. Here $S_{\pi}=\{1,2,3,2\}$.}
    \label{fig:np-poly-contract}
\end{figure}

\begin{prop}\label{prop:noncross-partition-cycles}
    Let $\pi=\lcurb V_1,\dots, V_m \rcurb \in NC(k)$. Then $\C_k/\pi$ is a union of cycles $\C_{k,1},\dots,\C_{k,p}$ where
    \begin{enumerate}
        \item $\abs{V\(\C_{k,i}\)\cap V\(\C_{k,i+1}\)} = 1$ for each  $i\in[p-1]$.
        \item $\displaystyle p = k-m+1$.
        \item $\displaystyle\sum_{i=1}^{p}\, \abs{\C_{k,i}} = k$.
    \end{enumerate} 
\end{prop}

\begin{defn}\label{defn:np-cycle-size-Sp}
    Let $\pi\in NC(k)$. Assume $\C_k/\pi = \C_{k,1}\cup\dots\cup \C_{k,p}$. We define $S_{\pi}$ to be the unordered sequence $\lcurb i_1,\dots,i_p\rcurb$ where $i_j$ is the size of the cycle $\C_{k,j}$ for each $j\in[p]$.
\end{defn}

\begin{prop}\
Let $k\in\Nbb$ and $\pi\in NC(k)$. Then
\begin{equation}
    [\pi, 1_k] \cong \[0_k, K(\pi)\] \cong \prod_{V\in K(\pi)} NC\(\,\abs{V}\) = \prod_{x\in S_{\pi}} NC(x)
\end{equation}
\end{prop}

\begin{prop}
Let $k\in\Nbb$ and $\pi=\{V_1,\dots, V_m\}\in NC(k)$. Let $\mu$ be the M\"{o}bius function of $NC(k)$. Then
\begin{equation}
    \mu\(\pi, 1_k\) = \prod_{x\in S_{\pi}}(-1)^{x-1}\cdot C_{x-1} = (-1)^{m-1} \cdot \prod_{x\in S_{\pi}} C_{x-1}
\end{equation}
\end{prop}

\begin{proof}
By the canonical factorization of $[\pi,1_k]$ and \Cref{prop:mu-NC-n},
\begin{align*}
    \mu\(\pi, 1_k\) = \mu\(\prod_{x\in S_{\pi}} NC(x)\) = \prod_{x\in S_{\pi}} \mu\(NC(x)\) = \prod_{x\in S_{\pi}} (-1)^{x-1}\cdot C_{x-1}
\end{align*}
The last equality is because $\displaystyle \sum_{x\in S_\pi} x-1 = k-p = m-1$.
\end{proof}

\begin{figure}[hbt!]
    \centering
    \includegraphics[scale=0.32]{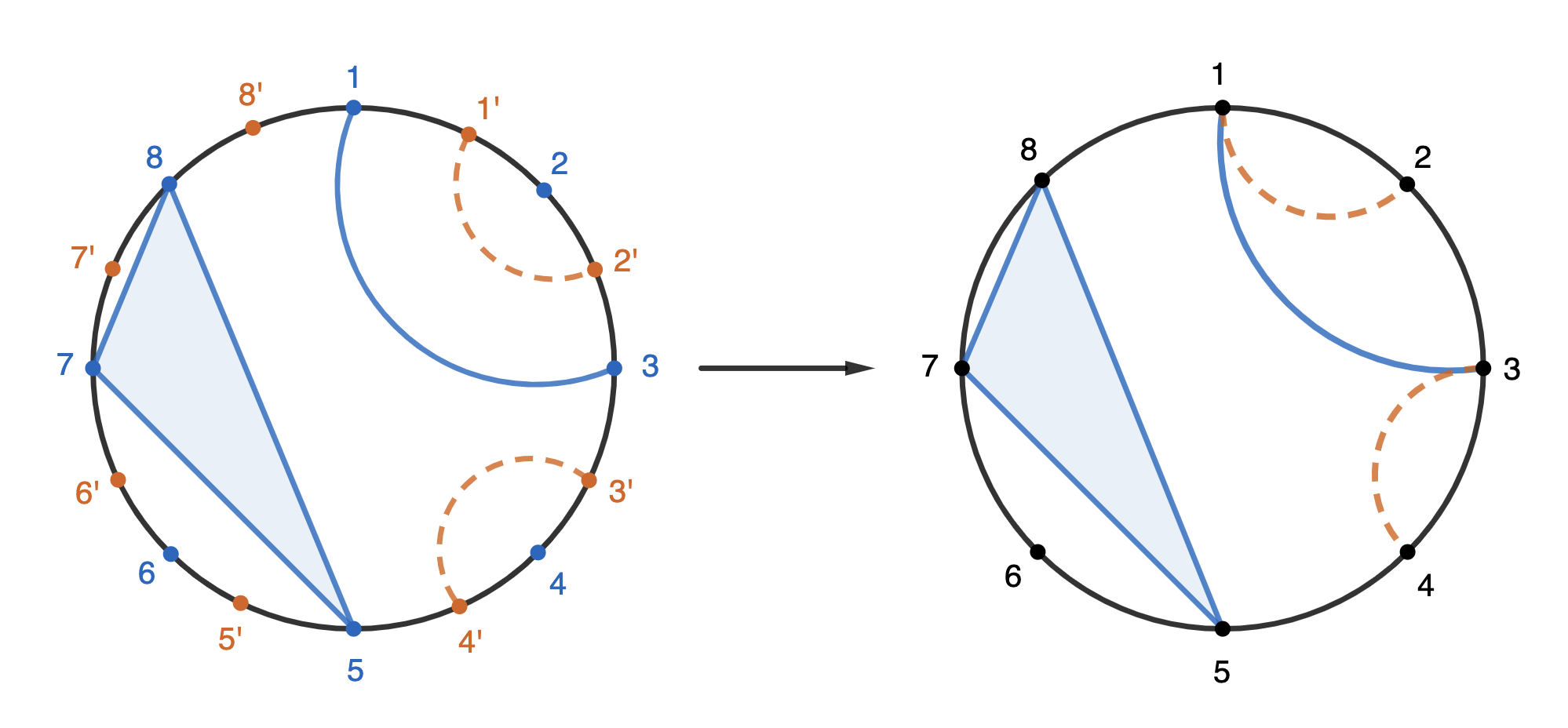}
    \caption{Illustration of \Cref{defn:np-oplus} and \Cref{defn:size-oplus}: Blue solid represents $\pi$ and orange dashed represents $\sigma$. Here $\pi\oplus \sigma\in NC(k)$ since $\pi\cup\bar{\sigma}$ is noncrossing in the left figure. Moreover, $S_{\pi\oplus \sigma} = \{1,1,1,1,2,2\}$.}
    \label{fig:np-oplus}
\end{figure}

\begin{defn}\label{defn:np-oplus}
    Given $\pi,\sigma\in NC(k)$, we say that $\pi\oplus \sigma \in NC(k)$ if $\pi\cup \bar{\sigma}\in NC(1,\bar{1}, \dots, k,\bar{k})$ where $\bar{\sigma} := \left\{ \left\{ \bar{v_1},\dots \bar{v_m} \right\}: W=\{v_1,\dots, v_m\}\in \sigma \right\}$.
    
    See \Cref{fig:np-oplus} for illustration.
\end{defn}

\begin{defn}\label{defn:size-oplus}
    Let $\pi,\sigma\in NC(k)$ such that $\pi\oplus \sigma \in NC(k)$. Let $\C_k/\pi\oplus \sigma$ be the graph obtained by identifying vertices together under $\pi\cup\sigma$. Then $\C_k/\pi\oplus \sigma = \C_{k,1}\cup \dots\cup \C_{k,p}$ is a union of cycles. We define $S_{\pi\oplus \sigma} = \{i_1,\dots, i_p\}$ where $i_j$ is the size of $\C_{k,j}$. See \Cref{fig:np-oplus} for illustration.
\end{defn}

\begin{prop}\label{prop:moments-ab-step1}
Let $\(\A,\varphi\)$ be a non-commutative probability space. Let $\{a_1,\dots, a_n\}, \{b_1,\dots, b_n\}\in\A$ be freely independent. Then
\begin{equation}
    \varphi\(a_1b_1\dots a_nb_n\) = \sum_{\substack{\pi,\sigma\in NC(n): \\ \pi\oplus \sigma \in NC(n)}} (-1)^{\,\abs{\pi} + \abs{\sigma} - k-1} \cdot \prod_{x\in S_{\pi\oplus \sigma}}  C_{x-1} \cdot \prod_{V\in\pi} \varphi(V)[a_1,\dots, a_n] \cdot \prod_{W\in\sigma} \varphi(W)[b_1,\dots, b_n]
\end{equation}

In particular, for $a,b\in\A$ freely independent,
\begin{equation}
    \varphi\((ab)^n\) = \sum_{\substack{\pi,\sigma\in NC(n): \\ \pi\oplus \sigma \in NC(n)}} (-1)^{\,\abs{\pi} + \abs{\sigma}-k-1} \cdot \prod_{x\in S_{\pi\oplus \sigma}}  C_{x-1} \cdot \prod_{V\in\pi} \varphi\(a^{\,\abs{V}}\) \cdot \prod_{W\in\sigma} \varphi\(b^{\,\abs{W}}\)
\end{equation}

\end{prop}

\begin{proof}
By \Cref{thm:cumulants},
    \begin{align*}
        &\varphi\(a_1b_1\dots a_nb_n\) 
        = \sum_{\pi\in NC(n)} \kappa_{\pi}[a_1,\dots, a_n]\cdot \varphi_{K(\pi)} [b_1,\dots, b_n] \\
        =&  \sum_{\pi\in NC(n)} \kappa_{K(\pi)}[a_1,\dots, a_n]\cdot \varphi_{\pi} [b_1,\dots, b_n] \\
        =& \sum_{\pi\in NC(n)} \prod_{V\in K(\pi)} \kappa(V)[a_1,\dots, a_n]\cdot \prod_{W\in\pi} \varphi(W) [b_1,\dots, b_n] \\
        =& \sum_{\pi\in NC(n)} \prod_{V\in K(\pi)} \kappa_{\,\abs{V}}\(a_{i_1}, \dots, a_{i_{\,\abs{V}}} \)\cdot \prod_{W\in\pi} \varphi(W) [b_1,\dots, b_n] \text{ where } V=\{i_1,\dots, i_{\,\abs{V}}\} \\
        =& \sum_{\pi\in NC(n)} \(\prod_{V\in K(\pi)} \sum_{\sigma\in NC(\,\abs{V})} \varphi_{\sigma}\[a_{i_1}, \dots, a_{i_{\,\abs{V}}} \] \cdot \mu\(\sigma, 1_{\,\abs{V}}\) \)\cdot \prod_{W\in\pi} \varphi(W) [b_1,\dots, b_n] \\
        =& \sum_{\pi\in NC(n)} \(\prod_{V_j\in K(\pi)} \sum_{\sigma_j\in NC(\,\abs{V_j})} \varphi_{\sigma_j}\[a_{i_{j1}}, \dots, a_{i_{j\,\abs{V_j}}} \] \cdot \prod_{x\in S_{\sigma_j}} (-1)^{x-1} C_{x-1}\) \cdot \prod_{W\in\pi} \varphi(W) [b_1,\dots, b_n] \\
        =& \sum_{\substack{\pi,\sigma\in NC(n):\\ \pi\oplus \sigma\in NC(n)}} \(\prod_{x\in S_{\pi\oplus \sigma}} (-1)^{x-1} C_{x-1}\) \cdot \prod_{V\in \sigma} \varphi(V) [a_1,\dots, a_n] \cdot \prod_{W\in\pi} \varphi(W) [b_1,\dots, b_n] \\
        =& \sum_{\substack{\pi,\sigma\in NC(n):\\ \pi\oplus \sigma\in NC(n)}} (-1)^{\,\abs{\pi} + \abs{\sigma}-k-1} \cdot \prod_{x\in S_{\pi\oplus \sigma}} C_{x-1} \cdot \prod_{V\in \sigma} \varphi(V) [a_1,\dots, a_n] \cdot \prod_{W\in\pi} \varphi(W) [b_1,\dots, b_n] \\
    \end{align*}
\end{proof}

\begin{defn}
    For $\vaa, \vbb\in P_k$, we will denote $\np\(\vaa, \vbb\)$ to be 
    \begin{equation}
        \np\(\vaa, \vbb\) = \left\{\(\pi,\sigma\): \pi\in\np(\vaa), \sigma\in \np(\vbb), \pi\oplus \sigma \in NC(k)\right\}
    \end{equation}
\end{defn}

\begin{cor}\label{cor:moments-ab-step1}
    Let $\(\A,\varphi\)$ be a non-commutative probability space. Let $a, b\in\A$ be freely independent. Then
\begin{equation}
    \varphi\((ab)^k\) = \sum_{\vaa, \vbb\in P_k}\, (-1)^{\,A+B-k-1} \cdot\( \sum_{(\pi,\sigma)\in \np\(\vaa,\vbb\)}  \prod_{x\in S_{\pi\oplus \sigma}}  C_{x-1}\) \cdot \vec{\varphi_a}^{\vaa} \cdot \vec{\varphi_b}^{\vbb}
\end{equation}
where $A=\a_1+\dots + \a_k$ and $B=\b_1+\dots + \b_k$ for $\vaa,\vbb\in P_k$.
\end{cor}

\begin{proof}
    This is because for $\pi\in \np(\vaa)$ and $\sigma\in \np(\vbb)$, $\abs{\pi} = A$ and $\abs{\sigma} = B$, and $\displaystyle \prod_{V\in\pi} \varphi\(a^{\,\abs{V}}\) = \prod_{V\in\pi} \varphi(a^i)^{\a_i} = \vec{\varphi_a}^{\vaa}$, $\displaystyle \prod_{W\in\sigma} \varphi\(b^{\,\abs{W}}\) = \prod_{W\in\sigma} \varphi(b^i)^{\b_i} = \vec{\varphi_b}^{\vbb}$.
\end{proof}

\begin{restatable}{thm}{momentMiddleStep} \label{thm:moments-ab-step2}
    Let $\vaa, \vbb\in P_k$ Then
    \begin{equation}
         \sum_{(\pi,\sigma)\in \np\(\vaa,\vbb\)} \prod_{x\in S_{\pi\oplus\sigma}} C_{x-1} = k\cdot \binom{A+B-2}{k-1}\cdot \dfrac{(A-1)!}{\a_1!\dots\a_k!} \cdot \dfrac{(B-1)!}{\b_1!\dots\b_k!}
    \end{equation}
\end{restatable}

Note that in the special case when $\vaa=\vbb=(k,0,\dots,0)\in P_k$, $\displaystyle C(\vaa,\vbb) = (-1)^{k-1}\cdot \dfrac{1}{k} \binom{2k-2}{k-1} = (-1)^{k-1}\cdot C_{k-1}$, and $LHS = (-1)^{\abs{0_k} + \abs{0_k}-k-1}\cdot C_{k-1} = (-1)^{k-1}C_{k-1}$.

\Cref{thm:main-free-prob} follows directly from \Cref{cor:moments-ab-step1} and \Cref{thm:moments-ab-step2}.

\begin{proof}
\begin{align*}
    \varphi\((ab)^k\) 
    &= \sum_{\vaa, \vbb\in P_k}\, (-1)^{\,A+B-k-1} \cdot\( \sum_{(\pi,\sigma)\in \np\(\vaa,\vbb\)}  \prod_{x\in S_{\pi\oplus \sigma}}  C_{x-1}\) \cdot \vec{\varphi_a}^{\vaa} \cdot \vec{\varphi_b}^{\vbb} \\
    &= \sum_{\vaa, \vbb\in P_k}\, (-1)^{\,A+B-k-1} \cdot k\cdot \binom{A+B-2}{k-1}\cdot \dfrac{(A-1)!}{\a_1!\dots\a_k!} \cdot \dfrac{(B-1)!}{\b_1!\dots\b_k!} \cdot \vec{\varphi_a}^{\vaa} \cdot \vec{\varphi_b}^{\vbb} \\
    &= \sum_{\vaa, \vbb\in P_k} C\(\vaa,\vbb\)\cdot \vec{\varphi_a}^{\vaa} \cdot \vec{\varphi_b}^{\vbb}
\end{align*}
\end{proof}
 
So to prove \Cref{thm:main-free-prob}, it suffices to prove \Cref{thm:moments-ab-step2}. We will first prove \Cref{thm:moments-ab-step2} when $\vbb= (k,0,\dots,0)\in P_k$ and then prove the theorem in the general case.


\subsection{ Case when \texorpdfstring{$\vbb = (k,0,\dots,0)$}{ß = (k,0,...,0)}} \label{subsection:moments-base-2}

\setlength{\parskip}{1.3mm}
\setlength{\baselineskip}{1.3em}

When $\vbb=(k,0,\dots,0)$, we want to show the following.
\begin{restatable}{thm}{FreeProbCoefStepTwo}
\label{thm:free-prob-coefficient-step-2}
    Let $\vaa = \(\a_1,\dots,\a_k\)\in P_k$ and $a= \a_1+\dots+\a_k$. Then
\begin{equation}
    \begin{aligned}
        & \sum_{\pi\in\np(\vaa)}\, \prod_{x\in S_{\pi}} C_{x-1} = \binom{a+k-2}{k-1}\cdot \dfrac{(a-1)!}{\a_1!\dots \a_k!}
    \end{aligned}
\end{equation}
\end{restatable}

\begin{defn}
    Given a set of elements $\I=\lcurb i_1,\dots,i_p\rcurb$ where $i_j$'s are not necessarily distinct, we denote $\perm\(i_1,\dots,i_p\)$, or $\perm\(\I\)$, to be the number of permutations of the elements in $\I$. 
\end{defn}

\begin{eg}
    $\perm\(\{1,2,3,4\}\)=4!=24$. $\perm\(\{1,1,2,2\}\)=\dfrac{4!}{2!2!}=6$.
\end{eg}

\begin{defn}
    Let $\pi$ be a non-crossing partition of $[k]$. Assume $\C_k/\pi = \C_{k,1}\cup\dots\cup \C_{k,p}$. Let $\X=\{x_1,\dots,x_p\}$ be such that $x_1+\dots+x_p=k$. 
    We say a permutation $\sigma:[p]\to[p]$ is a \emph{label of $\C_k/\pi$ corresponding to $\X$} if $\abs{\C_{k,j}}=x_{\sigma[j]}$ for all $j\in[p]$.
\end{defn}

\begin{defn}\label{defn:np-vaa-x}
    Let $\vaa\in P_k$. Let $p=k-a+1$ and $\X=\lcurb x_1,\dots,x_p\rcurb$ be such that $x_1+\dots+x_p=k$. We denote
    \begin{enumerate}
        \item $\np\(\vaa,\X\)$, or $\np\(\vaa,x_1,\dots,x_p\)$ to be 
        \begin{equation}
            \np\(\vaa,\X\) = \lcurb \pi\in\np\(\vaa\):S_{\pi}=\{x_1,\dots,x_p\} \rcurb,
        \end{equation}
        reads ``the set of $\N$oncrossing $\P$artitions corresponding to $\vaa$ that partitions $\C_k$ into parts of sizes $\X$ ".
        
        \item $\npla\(\vaa,\X\)$, or $\npla\(\vaa,x_1,\dots,x_p\)$ to be 
        \begin{equation}
            \npla\(\vaa,\X\) = \lcurb
            \begin{aligned}
            \(\pi,\sigma\): &\pi=\{P_1,\dots,P_a\}\in\np\(\vaa,\X\), \\
            &\sigma:[a]\to[a] \text{ where } \abs{P_i}=\abs{P_{\sigma(i)}} \text{ for all } i\in[a] 
            \end{aligned}
            \rcurb,
        \end{equation}
        reads ``the set of $\N$oncrossing $\P$artitions corresponding to $\L$abeled polygons $\vaa$, which partitions $\C_k$ into parts of sizes $\X$ ".
        
        \item $\nplx\(\vaa,\X\)$, or $\nplx\(\vaa,x_1,\dots,x_p\)$ to be 
        \begin{equation}
            \nplx\(\vaa,\X\)=\lcurb \(\pi,\sigma\):\pi\in\np\(\vaa,\X\), \sigma \text{ is a label of }\C_k/\pi \text{ corresponding to } \X \rcurb,
        \end{equation}
        reads ``the set of $\N$oncrossing $\P$artitions corresponding to $\vaa$ that partitions $\C_k$ into parts $\L$abeled with $\X$ ".
        
        \item $\nplax\(\vaa,\X\)$, or $\nplax\(\vaa,x_1,\dots,x_p\)$ to be 
        \begin{equation}
            \nplax\(\vaa,\X\)=\lcurb \(\pi,\sigma_1,\sigma_2\): \(\pi,\sigma_1\)\in\npla\(\vaa,\X\), \sigma_2 \text{ is a label of }\C_k/\pi \text{ corresponding to } \X \rcurb,
        \end{equation}
        reads ``the set of $\N$oncrossing $\P$artitions corresponding to $\L$abeled polygons $\vaa$, which partitions $\C_k$ into parts $\L$abeled with $\X$ ".
    \end{enumerate}
\end{defn}

\begin{eg}
\begin{enumerate}
$ $
    \item Let $\vaa\in P_k$ correspond to $\o_2^k$ and $\X=\{k\}$. Then $\abs{\np\(\vaa,\X\)} = 1$, $\abs{\npla\(\vaa,\X\)} = k!$, $\abs{\nplx\(\vaa,\X\)} = 1$ and $\abs{\nplax\(\vaa,\X\)} = k!$.

    \item Let $\vaa\in P_k$ correspond to $\o_{2k}$ and $\X=\{1,1,\dots,1\}$. Then $\abs{\np\(\vaa,\X\)} = 1$, $\abs{\npla\(\vaa,\X\)} = 1$, $\abs{\nplx\(\vaa,\X\)} = k!$ and $\abs{\nplax\(\vaa,\X\)} = k!$.
    
    \item Let $\vaa\in P_5$ correspond to $\o_2\o_4^2$ and $\X = \{1,1,3\}$. Then $\abs{\np\(\vaa,\X\)} = 5$, $\abs{\npla\(\vaa,\X\)} = 10$ (can exchange the position of the two lines), $\abs{\nplx\(\vaa,\X\)} = 10$ (can exchange the labels $x_1$ and $x_2$ which are both of size $1$) and $\abs{\nplax\(\vaa,\X\)} = 20$. See \Cref{fig:np-two-lines}, \Cref{fig:np-two-lines-la} and \Cref{fig:np-two-lines-lx} for illustration.
    
    \begin{figure}[hbt!]
        \centering
        \includegraphics[scale = 0.35]{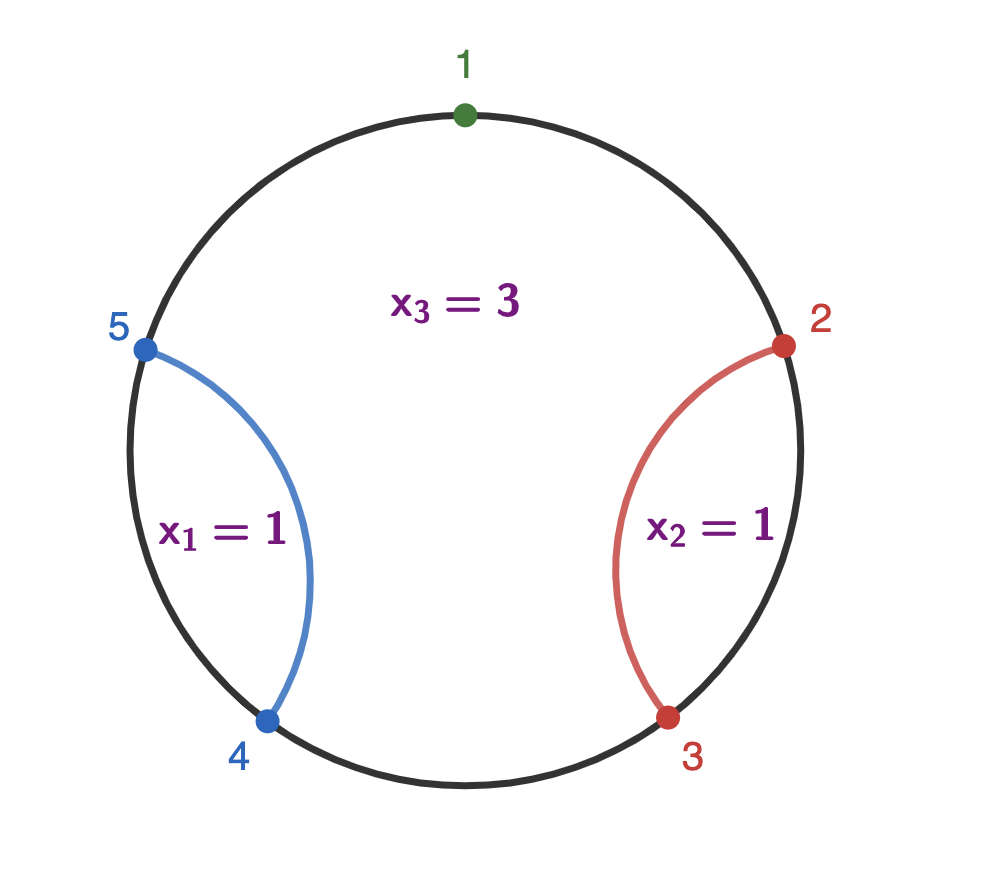}
        \caption{Example of $\protect \voo^{\vaa} = \o_2\o_4^2 \protect$ and $\X=\{1,1,3\}$ where $\pi=\left\{\{1\},\{2,3\}, \{4,5\}\right\}\in \np\(\vaa,\X\)$. The other four partitions in $\np\(\vaa,\X\)$ are $\left\{\{2\},\{3,4\}, \{5,1\}\right\}$, $\left\{\{3\},\{4,5\}, \{1,2\}\right\}$, $\left\{\{4\},\{5,1\}, \{2,3\}\right\}$ and $\left\{\{5\},\{1,2\}, \{3,4\}\right\}$.}
        \label{fig:np-two-lines}
    \end{figure}
    
    \begin{figure}[hbt!]
        \centering
        \includegraphics[scale=0.35]{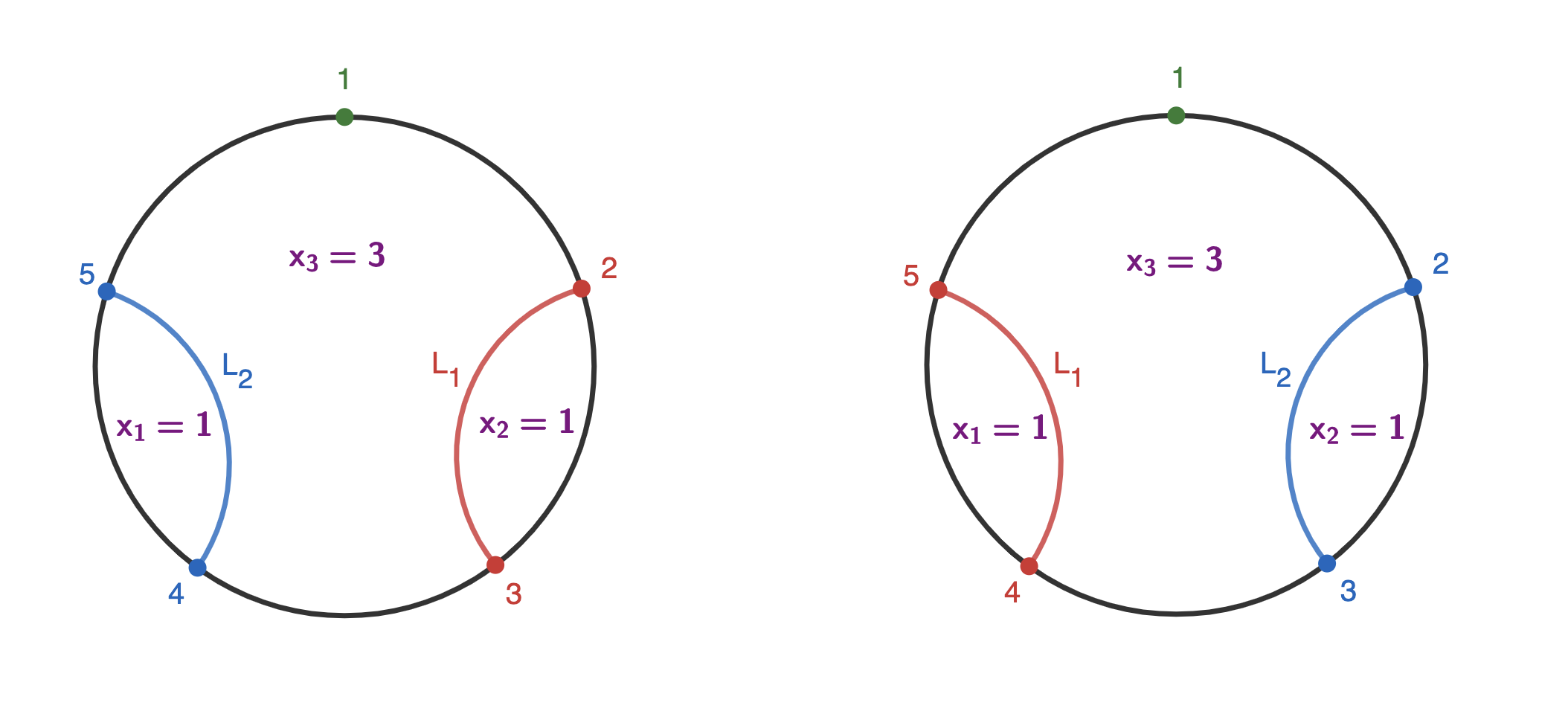}
        \caption{For $\npla\(\vaa,\X\)$, we can swap the position of the lines $L_1$ and $L_2$, doubling the counting for $\np\(\vaa,\X\)$.}
        \label{fig:np-two-lines-la}
    \end{figure}
    
    \begin{figure}[hbt!]
        \centering
        \includegraphics[scale=0.35]{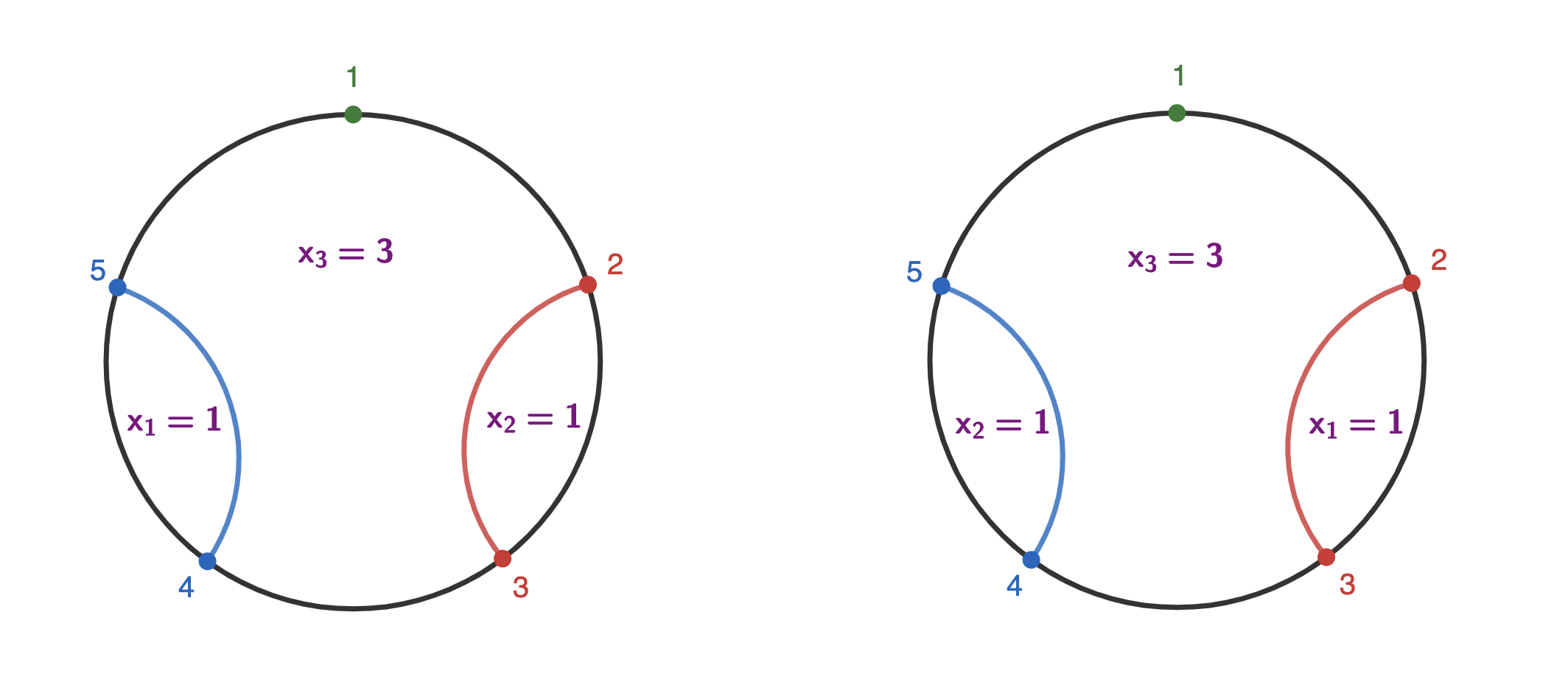}
        \caption{For $\nplx\(\vaa,\X\)$, we can swap the position of the $x_1$ and $x_2$, doubling the counting for $\np\(\vaa,\X\)$.}
        \label{fig:np-two-lines-lx}
    \end{figure}
\end{enumerate}
\end{eg}

\begin{obs}
    Let $\vaa=(\a_1,\dots,\a_k)\in P_k$, $a=\a_1+\dots+\a_k$ and $p=k-a+1$. Let $\X=\{x_1,\dots,x_p\}$ be such that $x_1+\dots+x_p=k$. Then
    \begin{equation}
        \abs{\np\(\vaa,\X\)} = \dfrac{\perm\(\X\)}{p!}\cdot\abs{\nplx\(\vaa,\X\)}\,
    \end{equation}
    and 
    \begin{equation}
        \abs{\np\(\vaa,\X\)} = \dfrac{1}{\a_1!\dots\a_k!}\cdot\abs{\npla\(\vaa,\X\)}\,
    \end{equation}
\end{obs}

\begin{thm}\label{thm:num-of-partitions-1}
    Let $\vaa = \(\a_1,\dots,\a_k\)\in P_k$ and $a=\a_1+\dots+\a_k$, $p=k-a+1$. Let $x_1,\dots,x_p$ be such that $x_1+\dots+x_p=k$. Then 
    \begin{equation}
        \abs{\nplax\(\vaa,\X\)} =  k\cdot (p-1)!\cdot (a-1)! \,.
    \end{equation}
\end{thm}

A direct corollary is the following.
\begin{cor}
    Let $\vaa = \(\a_1,\dots,\a_k\)\in P_k$ and $a=\a_1+\dots+\a_k$, $p=k-a+1$. Let $x_1,\dots,x_p$ be such that $x_1+\dots+x_p=k$. Then 
    \begin{equation}
        \abs{\nplx\(\vaa,\X\)} = k\cdot (p-1)!\cdot \dfrac{(a-1)!}{\a_1!\dots\a_k!}\,.
    \end{equation}
\end{cor}

To prove \Cref{thm:num-of-partitions-1}, we will start with the case of a single polygon.
\begin{lemma}\label{lem:num-of-partitions-polygon}
    Let $\vaa$ corresponds to $\o_{2}^{\,k-t}\o_{2t}$. Let $\X=\{x_1,\dots,x_t\}$ be such that $x_1+\dots+x_t=k$. Then 
    \begin{equation}
        \abs{\nplax\(\vaa,\X\)} = k\cdot(t-1)!\cdot(k-t)!\,.
    \end{equation}
\end{lemma}
\begin{proof}
\begin{figure}[htb!]
    \centering
    \includegraphics[scale=0.3]{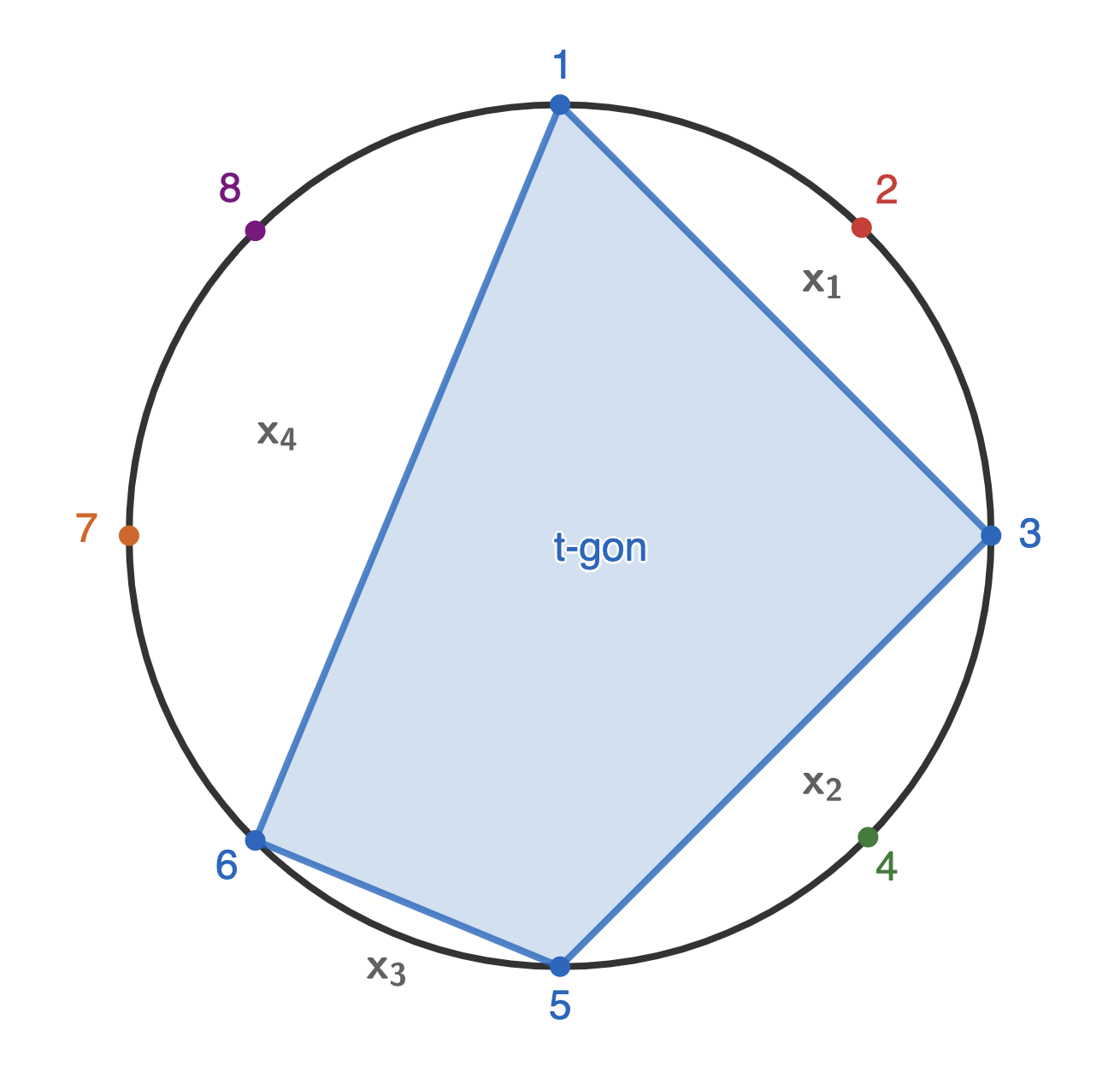}
    \caption{Illustration of \Cref{lem:num-of-partitions-polygon}. In this case we choose $1$ to be the starting point of the $t$-gon, and $x_1,x_2,x_3,x_4$ as the order of $x_i$'s. On the other hand, choosing $3$ to be the starting point of the $t$-gon and $x_2,x_3,x_4,x_1$ as order of the $x_i$'s results in the exact same configuration.}
    \label{fig:np-one-polygon}
\end{figure}
    $\o_{2t}$ corresponds to a $t$-gon on the cycle $\C_k$ and divides $\C_k$ into $t$ parts of size $x_1,\dots,x_t$. There are $k$ ways to choose where to place the first vertex of the $t$-gon on $\C_k$, $t!$ ways to arrange the order of the $x_i$ (corresponding to ways to put the remaining $t-1$ vertices on $\C_k$). For each placement of the $t$-gon on $\C_k$, there are $t$ choices for the first vertex to result in this placement. Since all the dots corresponding to $\o^{k-t}$ are considered to be distinct, there are $(k-t)!$ ways rearrange the order of these dots. Thus $\abs{\nplax\(\vaa,\X\)} = k\cdot t!/t\cdot (k-t)! = k\cdot (t-1)!\cdot (k-t)!$ as wanted.
\end{proof}

\begin{lemma}\label{lem:num-of-partitions-reduction}
    Let $\vaa$ be corresponding to $\o_2^{\,\a_1}\,\o_{2t_1}\dots\o_{2t_m}$ and $\vbb$ be corresponding to \\ $\o_2^{\,\a_1+1}\,\o_{2t_1}\dots\o_{2t_{m-2}}\o_{2(t_{m-1}+t_m-1)}$. Notice that $a = a = \a_1+m$. Let $p = k-a+1= k-\a_1-m+1$ and $\X=\{x_1,\dots,x_p\}$ be such that $x_1+\dots+x_p=k$. Then 
    \begin{equation}
        \abs{\nplax\(\vaa, \X\)} = \abs{\nplax\(\vbb, \X\)}\,.
    \end{equation}
\end{lemma}

\begin{proof}
    We will show that there is a bijection between $\nplax\(\vaa, \X\)$ and $\nplax\(\vbb, \X\)$.
    \begin{enumerate}
        \item $\nplax\(\vaa, \X\) \to \nplax\(\vbb, \X\)$: Let $\pi=\{P_1,\dots,P_a\}$ be a noncrossing partition of the cycle $\C_k$ corresponding to distinct $\vaa$ and $\X$. Assume $P_{a-1}$ and $P_{a}$ are the $t_{m-1}$-gon and $t_m$-gon, respectively. Assume  $v_1,\dots,v_{t_m},w_1,\dots,w_{t_{m-1}}$ are ordered clockwise on $\C_k$ where $P_{a} = \{v_1,\dots,v_{t_m}\}$ and $P_{a-1} = \{w_1,\dots,w_{t_{m-1}}\}$.
        
        We now merge $P_{a-1}$ and $P_a$ to get $P_b$ of size $(t_{m-1}+t_{m}-1)$ and shift the other polygons in the following way:
        \begin{enumerate}[i.]
            \item Let $P_b=\big\{v_1,\dots,v_{t_m} \big\}\bigcup \big\{w_i-\(w_1-v_{t_m}\) \mod k: i\in[t_{m-1}] \big\}$. i.e. We shift the $P_{a-1}$ polygon $\(w_1-v_{t_m} \mod k\)$ units in the counter-clockwise direction so that the original vertex $w_1$ touches $v_{t_m}$.
            
            \item For any $i\in[a-2]$, let $P_i' = \big\{ p+\delta_1(p)\cdot \(w_{t_m}-w_{1}\) - \delta_2(p)\cdot \(w_1-v_{t_m}\) \mod k: p\in P_i\big\}$ where
            \begin{equation*}
                \delta_1(p) = 
                \begin{cases}
                     \, 1 & \text{if $p$ is strictly in between $v_{t_m}$ and $w_1$ in the clockwise direction} \\
                     \, 0 & \text{ otherwise}
                \end{cases}
            \end{equation*}
            and 
            \begin{equation*}
                \delta_2(p) = 
                \begin{cases}
                     \, 1 & \text{if $p$ is strictly in between $w_1$ and $w_{t_{m-1}}$ in the clockwise direction} \\
                     \, 0 & \text{ otherwise}
                \end{cases}
            \end{equation*}
            
            i.e. For each of the remaining polygons $P_i$, we shift its vertices that are in between $v_{t_m}$ and $w_1$ \textbf{clockwise} $\(w_{t_m}-w_{1}\)$ units and those in between $w_1$ and $w_{t_{m-1}}$  $\(w_1-v_{t_m} \mod k\)$ units in the \textbf{counter-clockwise} direction.
            
            \item Let $s=w_{t_{m-1}}$. 
        \end{enumerate}
        
        
        The resulting $\pi' = \big\{P_1',\dots,P_{a-2}',\{s\}, P_b\big\}$ is in $\nplax\(\vbb,\X\)$. See \Cref{fig:np-two-poly-combine} for illustration.
        
        \begin{figure}[hbt!]
            \centering
            \includegraphics[scale=0.35]{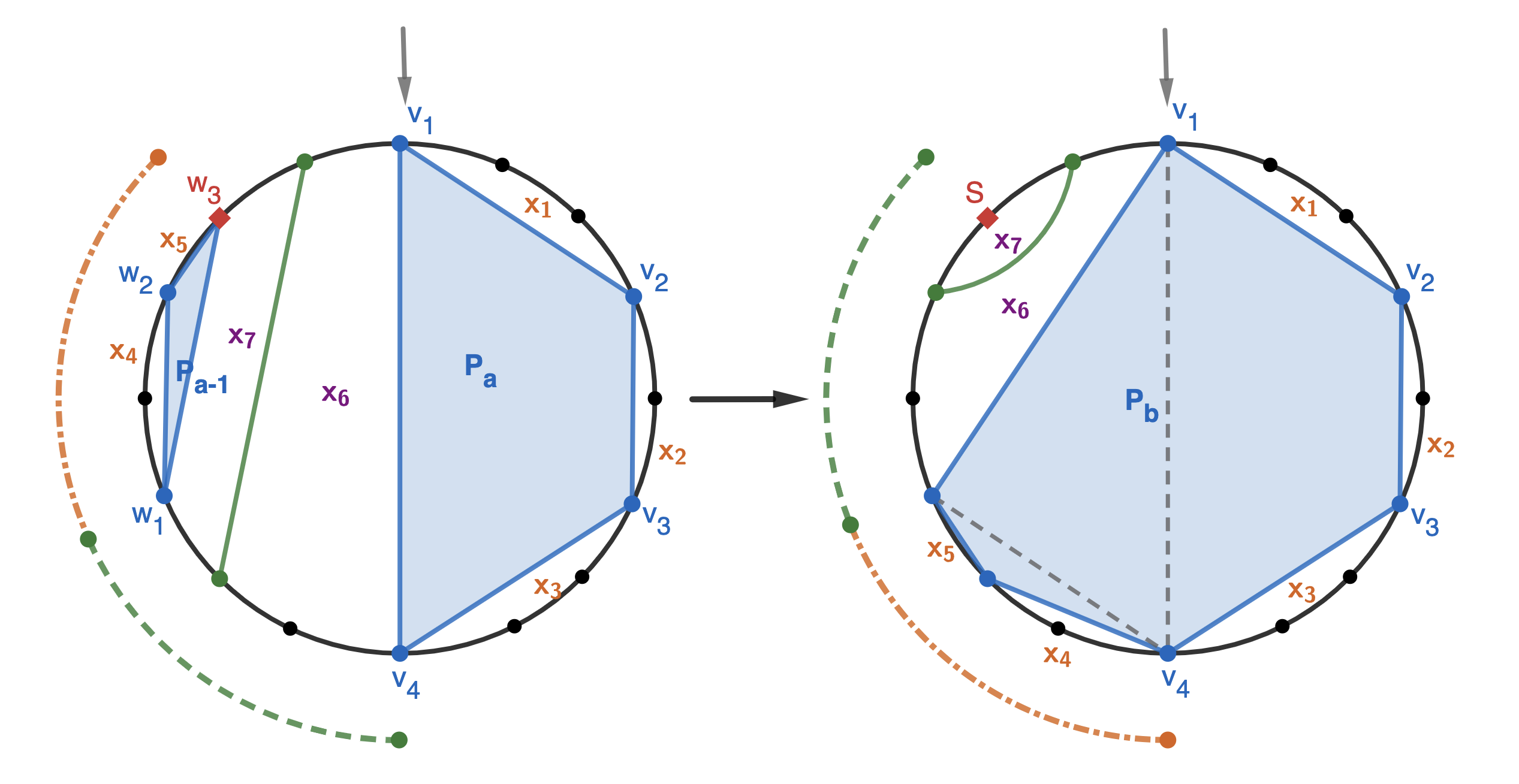}
            \caption{Illustration of $\nplax\(\vaa, \X\) \to \nplax\(\vbb, \X\)$. }
            \label{fig:np-two-poly-combine}
        \end{figure}
        
        \item $\nplax\(\vbb, \X\) \to \nplax\(\vaa, \X\)$: Let $\pi=\{P_1,\dots,P_b\}$ be a noncrossing partition of the cycle $\C_k$ corresponding to distinct $\vbb$ and $\X$. W.O.L.G assume $P_1=\{s\}$ is a point, call $s$ the special point, and $P_b=\{p_1,\dots,p_{t_{m-1}+t_{m}-1}\}$ is the $(t_{m-1}+t_m-1)$-polygon. Now we will split $P_b$ into $P_{b,t_m}$ and $P_{b,t_{m-1}}$, a $t_m$-gon and a $t_{m-1}$-gon respectively, in the following way:
        \begin{enumerate}[i.]
            \item Let $P_{b,t_m}=\{p_1,\dots,p_{t_m}\}\subseteq P_b$ be the first $t_m$ points in the clockwise direction starting from $s$.
            \item Let $P_{b,t_{m-1}} = \{q_1,\dots,q_{t_{m-1}}\}$ where $q_i = p_{t_m+i-1}+\(s-p_{t_{m-1}+t_{m}-1}\) \mod k$. i.e. We shift the $t_{m-1}$-gon $\{p_{t_m},\dots, p_{t_{m-1}+t_{m}-1}\}$ \, $\(s-p_{t_{m-1}+t_{m}-1} \mod k\)$ units so that the last vertex of the $t_{m-1}$-gon replace the original $s$.
            
            \item For any $i\in\{2,3,\dots,b-1\}$, let $P_i' = \big\{ p+ \delta_1(p)\cdot \(s-p_{t_m+t_{m-1}-1}\) - \delta_2(p)\cdot \(p_{t_m+t_{m-1}-1}-p_{t_m}\) \mod k: p\in P_i\big\}$ where 
            \begin{equation*}
                \delta_1(p) = 
                \begin{cases}
                     \, 1 & \text{if $p$ is strictly in between $p_{t_m}$ and $p_{t_m+t_{m-1}-1}$ in the clockwise direction} \\
                     \, 0 & \text{ otherwise}
                \end{cases}
            \end{equation*}
            and 
            \begin{equation*}
                \delta_2(p) = 
                \begin{cases}
                     \, 1 & \text{if $p$ is strictly in between $p_{t_m+t_{m-1}-1}$ and $s$ in the clockwise direction} \\
                     \, 0 & \text{ otherwise}
                \end{cases}
            \end{equation*}
            
            i.e. For each of the remaining polygons $P_i$, we shift its vertices that are in between $p_{t_m}$ and $p_{t_m+t_{m-1}-1}$ \textbf{clockwise} $\(s-p_{t_m+t_{m-1}-1}\)$ units and those in between $p_{t_m+t_{m-1}-1}$ and $s$  $\(p_{t_m+t_{m-1}-1}-p_{t_m} \mod k\)$ units in the \textbf{counter-clockwise} direction.
            
        \end{enumerate}
        
        
        The resulting partition $\pi'=\{P_2',\dots,P_{b-1}', P_{b,t_{m-1}}, P_{b,t_m}\}$ is in $\nplax\(\vaa,\X\)$. See \Cref{fig:np-two-poly-split-2} for illustration.
        
        
        \begin{figure}[hbt!]
            \centering
            \includegraphics[scale=0.35]{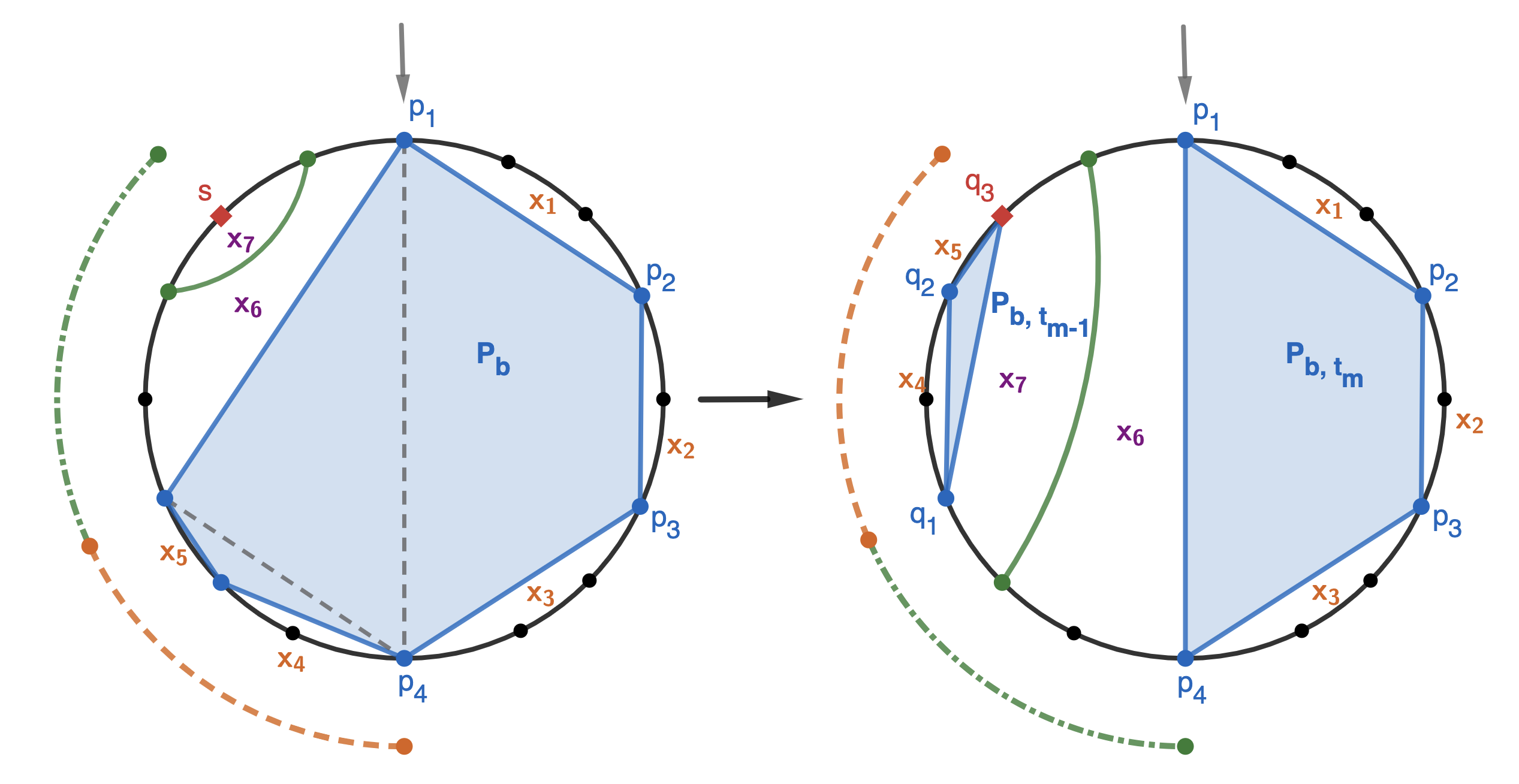}
            \caption{Illustration of $\nplax\(\vbb, \X\) \to \nplax\(\vaa, \X\)$. }
            \label{fig:np-two-poly-split-2}
        \end{figure}
        
        \item It is not hard to see that the above two operations are reverses of each other.
    \end{enumerate}
\end{proof}

\begin{cor}\label{cor:num-of-partitions-multipoly}
    Let $\vaa\in P_k$. Let $a=\a_1+m$ and $p=k-a+1$. Let $\X=\lcurb x_1,\dots,x_p\rcurb$ be such that $x_1+\dots+x_p=k$. Then 
    \begin{equation}
        \abs{\nplax\(\vaa,\X\)} = k\cdot (p-1)!\cdot (a-1)!\,.
    \end{equation}
\end{cor}

\begin{proof}
    We can write $\voo^{\vaa}$ as $\o_2^{\,\a_1}\,\o_{2t_1}\dots\o_{2t_m}$ where $t_i$'s are not necessarily distinct. eg. $\o_2^{2}\o_4^3\o^{6}$ can be written as $\o_2^2\o_4\o_4\o_4\o_6$.
    
    We can prove this by induction on $m$.
    \begin{enumerate}
        \item Base case $m=0$: $\nplax\(\vaa,\X\) = k!$ and $k\cdot \(p-1\)!\cdot (a-1)! = k\cdot 0!\cdot (k-1)! = k!$.
        
        \item Base case $m=1$: $a=k-t+1$. By \Cref{lem:num-of-partitions-polygon}, $\abs{\nplax\(\vaa,\X\)} = k\cdot (t-1)!\cdot(k-t)!$ and $RHS = k\cdot (t-1)!\cdot (a-1)! = k\cdot (t-1)!\cdot (k-t)!$.
        
        \item Inductive case $(m-1) \implies m$: By \Cref{lem:num-of-partitions-reduction} and the inductive hypothesis on the term $\voo^{\vbb} = \o_2^{\,\a_1+1}\o_{2t_1}\dots\o_{2t_{m-2}}\o_{2(t_{m-1}+t_m-1)}$, which has $b=\a_1+1+m-1=\a_1+m = a$ and $p'=k-b+1 = k-a+1 = p$, we get
        \begin{align*}
            \abs{\nplax\(\vaa,\X\)}
            &= \abs{\nplax\(\o_2^{\,\a_1+1}\o_{2t_1}\dots\o_{2t_{m-2}}\o_{2(t_{m-1}+t_m-1)}, \X\)}\\
            &= k\cdot (p'-1)!\cdot (b-1)!
            = k\cdot (p-1)!\cdot (a-1)!\,.
    \end{align*}
    \end{enumerate}
\end{proof}


To prove \Cref{thm:free-prob-coefficient-step-2}, we still need the following Catalan number identity.

\begin{lemma}\label{lem:sum-of-prod-of-catalan}
Let $\displaystyle C_{n} = \dfrac{1}{n+1}\binom{2n}{n}$ denote the $n^{th}$ Catalan number. Then
    \begin{equation}
        \sum_{i_j\geq 0:\, i_1+\dots+i_k=n} C_{i_1}C_{i_2}\dots C_{i_k} = k\cdot\dfrac{(2n+k-1)!}{n!(n+k)!}
    \end{equation}
\end{lemma}

\begin{defn}
    We denote $\displaystyle A_n^{(k)}=\sum_{i_j\geq 0: i_1+\dots+i_k=n} C_{i_1}C_{i_2}\dots C_{i_k}$.
\end{defn}

\begin{prop}\label{prop:A-n-k-base-case}
    $A_n^{(1)}=C_n$ and $A_n^{(2)}=C_{n+1}$.
\end{prop}

\begin{proof}
    $A_n^{(1)}=C_n$ by the definition of $A_n^{(k)}$. $\displaystyle A_n^{(2)} = \sum_{i_1+i_2=n} C_{i_1}C_{i_2} = C_{n+1}$.
\end{proof}

\begin{proof}[Proof of \Cref{lem:sum-of-prod-of-catalan}]
    We first prove that $\displaystyle A_n^{(k)}=A_{n+1}^{(k-1)}-A_{n+1}^{(k-2)}$, and then prove the result by induction on $k$.
    
    \begin{claim}
        $A_n^{(k)}=A_{n+1}^{(k-1)}-A_{n+1}^{(k-2)}$.
    \end{claim}
    \begin{proof}
        \begin{align*}
            A_n^{(k)}
            &=\sum_{i_j\geq 0: i_1+\dots+i_k=n} C_{i_1}C_{i_2}\dots C_{i_k}
            =\sum_{\substack{i_j\geq 0, r\geq 0:\\ i_1+\dots+i_{k-2}+r=n}} C_{i_1}C_{i_2}\dots C_{i_{k-2}}\cdot\(\sum_{i_{k-1}+i_{k}=r} C_{i_{k-1}}C_{i_k}\)\\
            &=\sum_{\substack{i_j\geq 0,r\geq 0:\\ i_1+\dots+i_{k-2}+r=n}} C_{i_1}C_{i_2}\dots C_{i_{k-2}}\cdot C_{r+1}
            = \sum_{\substack{i_1,\dots,i_{k-2}\geq 0, i_{k-1}'\geq 1:\\i_1+\dots+i_{k-2}+i_{k-1}'=n+1}} C_{i_1}C_{i_2}\dots C_{i_{k-2}}\cdot C_{i_{k-1}'}\\
            &=\(\sum_{\substack{i_j\geq 0:\\ i_1+\dots+i_{k-2}+i_{k-1}'=n+1}} C_{i_1}C_{i_2}\dots C_{i_{k-2}}\cdot C_{i_{k-1}'}\) - \(\sum_{\substack{i_j\geq 0:\\i_1+\dots+i_{k-2}=n+1}} C_{i_1}C_{i_2}\dots C_{i_{k-2}}\cdot C_0\) \\
            &=A_{n+1}^{(k-1)}-A_{n+1}^{(k-2)}.
        \end{align*}
    \end{proof}
    
    Now we prove the lemma by induction on $k$.
    \begin{enumerate}
        \item Base case: $A_n^{(0)}=0=0\cdot\dfrac{(2n-1)!}{n!n!}$. $A_n^{(1)}=C_n=\dfrac{1}{n+1}\cdot\dfrac{(2n)!}{n!n!}$ and $RHS=\dfrac{(2n)!}{n!(n+1)!}$.
        \item Assume the result holds for $k-1$, then 
        \begin{align*}
            A_n^{(k)}
            &=(k-1)\cdot\dfrac{(2n+k)!}{(n+1)!(n+k)!}-(k-2)\cdot\dfrac{(2n+k-1)!}{(n+1)!(n+k-1)!}\\
            &\\
            &=\dfrac{(2n+k-1)!}{n!(n+k)!}\cdot\(\dfrac{(k-1)(2n+k)}{n+1}-\dfrac{(k-2)(n+k)}{n+1}\)\\
            &\\
            &=\dfrac{(2n+k-1)!}{n!(n+k)!}\cdot\(\dfrac{nk+k}{n+1} \)
            = k\cdot\dfrac{(2n+k-1)!}{n!(n+k)!}.
        \end{align*}
    \end{enumerate}
\end{proof}

\begin{cor}\label{cor:sum-of-prod-of-catalan}
    We can rewrite \Cref{lem:sum-of-prod-of-catalan} as
    \begin{equation}
        \sum_{i_j\geq 0\,:\,i_1+\dots+i_k=n-k} C_{i_1}\dots C_{i_k} = k\cdot\dfrac{(2n-k-1)!}{(n-k)!n!}
    \end{equation}
    or 
    \begin{equation}
        \sum_{i_j\geq 1\,:\,i_1+\dots+i_k=n} C_{i_1-1}\dots C_{i_k-1} = k\cdot\dfrac{(2n-k-1)!}{(n-k)!n!}.
    \end{equation}
\end{cor}


\comm{
\begin{lemma}
Assume $n\geq k$. Then
    \begin{equation}
        \sum_{i_j\geq 1\,:\,i_1+\dots+i_k=n} C_{i_1}\dots C_{i_k} = \dfrac{\,k\,}{n}\cdot \dfrac{(2n)!}{(n-k)!\,(n+k)!}
    \end{equation}
\end{lemma}
}


Now we are ready to prove \Cref{thm:free-prob-coefficient-step-2}.

\FreeProbCoefStepTwo*

\begin{proof}
    Recall that $p = k - a + 1$. By \Cref{thm:num-of-partitions-1} and \Cref{cor:sum-of-prod-of-catalan},
    \begin{align*}
        \sum_{\pi\in\np\(\vaa\)}\, \prod_{i_j\in S_{\pi}} C_{i_j-1}
        &= \sum_{i_j\geq 1:\, i_1+\dots+i_p=k}\, \dfrac{\abs{\np\(\vaa,i_1,\dots,i_p\)}}{\perm\(i_1,\dots,i_p\)} \cdot C_{i_1-1}\cdot\dots\cdot C_{i_p-1} \\
        &= \sum_{i_j\geq 1:\, i_1+\dots+i_p=k}\, \dfrac{\abs{\nplx\(\vaa,i_1,\dots,i_p\)}}{p\,!}\cdot C_{i_1-1}\cdot\dots\cdot C_{i_p-1} \\
        &= \sum_{i_j\geq 1:\, i_1+\dots+i_p=k}\, \(\dfrac{(a-1)!}{\a_1!\dots\a_k!} \cdot k\cdot (p-1)!/ p\,!\) \cdot C_{i_1-1}\cdot\dots\cdot C_{i_p-1} \\
        &= \(\dfrac{(a-1)!}{\a_1!\dots\a_k!} \cdot k\cdot (p-1)!/ p\,!\) \cdot \(p\cdot \dfrac{(2k-p-1)!}{(k-p)!\, k!}\)\\
        &= \dfrac{(a-1)!}{\a_1!\dots\a_k!} \cdot \dfrac{(k+a-2)!}{(a-1)!(k-1)!}\\
        &= \dfrac{(a-1)!}{\a_1!\dots\a_k!} \cdot \binom{k+a-2}{k-1}\,.
    \end{align*}
\end{proof}

\subsection{The General Case}\label{subsection:moments-general}

\setlength{\parskip}{1.5mm}
\setlength{\baselineskip}{1.3em}

Recall that
\begin{equation}
    \np\(\vaa, \vbb\) := \left\{\(\pi,\sigma\): \pi\in\np(\vaa), \sigma\in \np(\vbb), \pi\oplus \sigma \in NC(k)\right\}
\end{equation}

We want to prove the following.

\begin{equation}
    \sum_{(\pi,\sigma)\in \np\(\vaa,\vbb\)} \prod_{x\in S_{\pi\oplus\sigma}} C_{x-1} = k\cdot \binom{A+B-2}{k-1}\cdot \dfrac{(A-1)!}{\a_1!\dots\a_k!} \cdot \dfrac{(B-1)!}{\b_1!\dots\b_k!}
\end{equation}

\begin{defn}
    Let $\vaa,\vbb\in P_k$ and $p=k-(a+b)+1$. Let $\X=\lcurb x_1,\dots,x_p\rcurb$ be such that $x_1+\dots+x_p=k$. We denote
    \begin{enumerate}
        \item $\np\(\vaa,\vbb,\X\)$ to be 
        \begin{equation}
            \np\(\vaa,\vbb, \X\) = \lcurb \(\pi, \sigma\) \in \np\(\vaa,\vbb\): S_{\pi\oplus \sigma}=\{x_1,\dots,x_p\} \rcurb,
        \end{equation}
        reads ``the set of $\N$oncrossing $\P$artitions corresponding to $\vaa$ and $\vbb$ that partition $\C_k$ into parts of sizes $\X$ ".
        
        \item $\npla\(\vaa,\vbb,\X\)$ to be 
        \begin{equation}
            \npla\(\vaa, \vbb, \X\) = \lcurb
            \begin{aligned}
            \(\pi, \sigma, \tau_1, \tau_2\): & \(\pi, \sigma\) \in\np\(\vaa,\vbb, \X\), \\
            &\tau_1: [a]\to[a] \text{ where } \abs{P_i}=\abs{P_{\tau(i)}} \text{ for all } i\in[a]\\ 
            &\tau_2: [b]\to[b] \text{ where } \abs{Q_i}=\abs{Q_{\tau(i)}} \text{ for all } i\in[b]
            \end{aligned}
            \rcurb.
        \end{equation}
        
        \item $\nplx\(\vaa,\vbb,\X\)$ to be 
        \begin{align}
            \lcurb (\pi,\sigma,\rho): \(\pi, \sigma\)\in \np\(\vaa,\vbb,\X\),
            \rho \text{ is a label of }\C_k/\pi\oplus\sigma \text{ corresponding to } \X \rcurb.
        \end{align}
        reads ``the set of $\N$oncrossing $\P$artitions corresponding to $\vaa$ and $\vbb$ that partition $\C_k$ into parts $\L$abeled with $\X$ ".
        
        \item $\nplax\(\vaa,\vbb, \X\)$ to be 
        \begin{equation}
            \nplax\(\vaa,\vbb,\X\)=\lcurb
            \begin{aligned}
                \(\pi,\sigma, \tau_1, \tau_2, \rho\): &\(\pi,\sigma, \tau_1, \tau_2\) \in\npla\(\vaa,\vbb, \X\), \\
                &\rho \text{ is a label of } \C_k/\pi\oplus\sigma \text{ corresponding to } \X 
            \end{aligned}
            \rcurb.
        \end{equation}
    \end{enumerate}
\end{defn}

\begin{lemma}\label{lem:num-of-partitions-polygon-2}
    Let $\vaa$ and $\vbb$ be corresponding to $\o_2^{k-t_a}\o_{2t_a}$ and $\o_2^{k-t_b}\o_{2t_b}$, respectively. Let $a=k-t_a+1$, $b=k-t_b+1$ and $p=t_a+t_b-1 = 2k-a-b+1$. Let $\X=\{x_1,\dots,x_p\}$ be such that $x_1+\dots+x_p=k$. Then
    \begin{equation}
        \abs{\nplax\(\vaa,\vbb,\X\)} = k^2(p-1)!\cdot (a-1)!\cdot (b-1)!\,.
    \end{equation}
\end{lemma}

\begin{proof}
    Let $\vcc$ corresponds to $\o_2^{k-t_a-t_b+1}\o_{2(t_a+t_b-1)}$. Note that $p_c = t_a+t_b-1 = p$. We will prove that there is a bijection between $\nplx\(\vaa,\vbb,\X\)$ and $\nplx\(\vcc, \X\)\times [k]$.
    \begin{enumerate}
        \item $\nplx\(\vaa,\vbb,\X\)\to \nplx\(\vcc, \X\)\times [k]$: let $\P\cup\Q$ be a noncrossing partition of $[k]$ corresponding to $\vaa$ and $\vbb$ and labels $\X$. In this case $\P$ contains a polygon $P_a$ of size $t_a$ and $\Q$ contains a polygon $P_b$ of size $t_b$. Assume $P_b=\{p_1,\dots,p_{t_b}\}$ and $P_a=\{q_1,\dots,q_{t_a}\}$ such that $p_1,\dots,p_{t_b},q_1,\dots,q_{t_a}$ are ordered clockwise on $\C_k$.
        
        
        We now merge $P_a$ and $P_b$ to get $P_c$ of size $(t_a+t_b-1)$ in the following way: let $P_c=\big\{p_1,\dots,p_{t_b} \big\}\bigcup \big\{q_i-\(q_1-p_{t_b}\) \mod k: i\in[t_a] \big\}$. i.e. We shift the $P_a$ polygon $\(q_1-p_{t_b} \mod k\)$ units in the counter-clockwise direction so that the original vertex $q_1$ touches $p_{t_b}$. We further let $s= q_{t_a}$. The resulting $P_c$ corresponds to a partition in $\nplx\(\vcc,\X\)$ and $s\in[k]$. See \Cref{fig:np-two-poly-ab-combine} for an illustration.
        
        \begin{figure}[hbt!]
            \centering
            \includegraphics[scale=0.3]{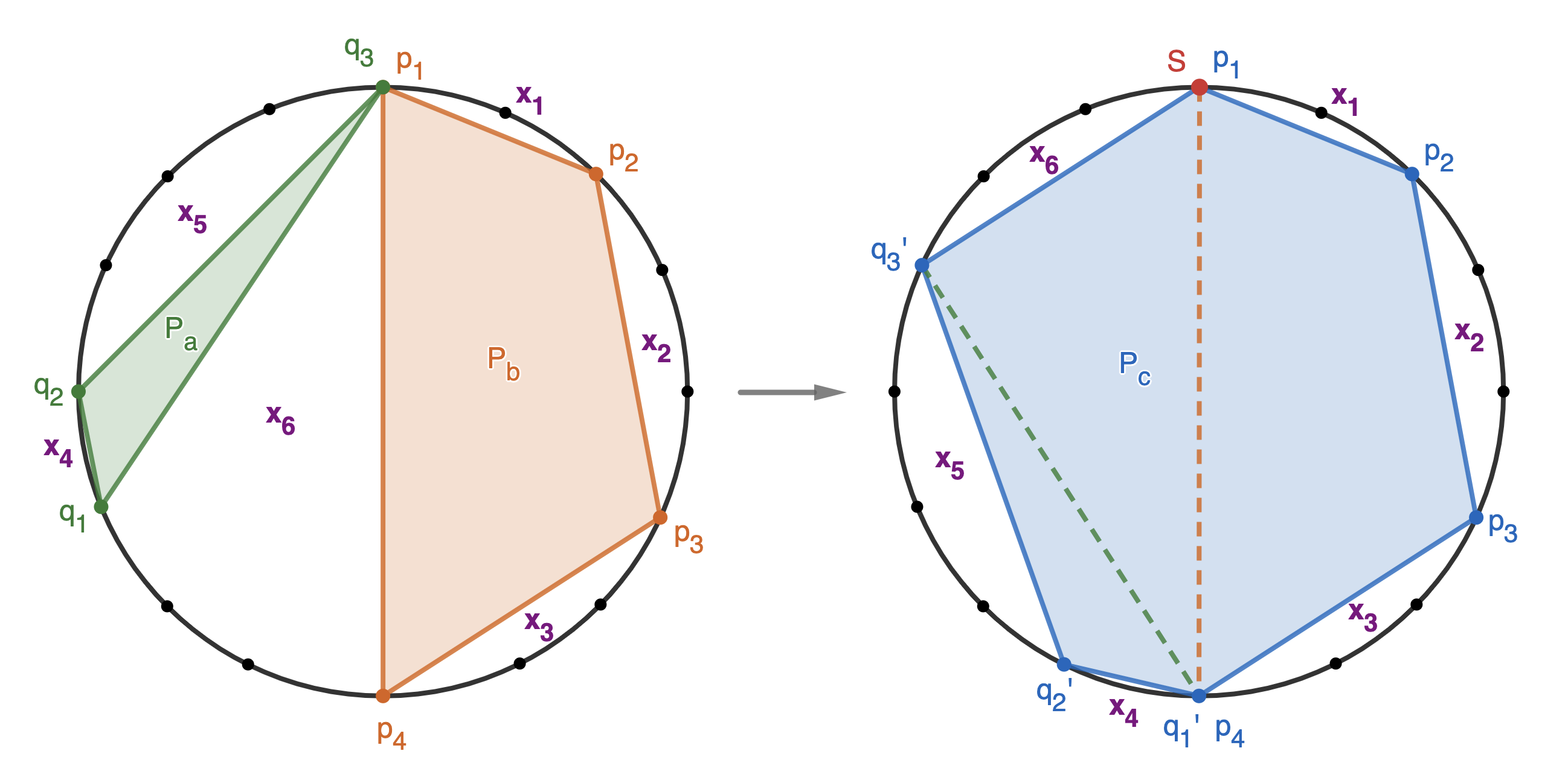}
            \caption{Illustration of $\nplx\(\vaa,\vbb,\X\)\to \nplx\(\vcc, \X\)\times [k]$: here $p_1=q_3=s$, $t_a=3$, $t_b=4$ and $t_c=6$.}
            \label{fig:np-two-poly-ab-combine}
        \end{figure}
        
        \item $\nplx\(\vcc, \X\)\times [k]\to \nplx\(\vaa,\vbb,\X\)$: Let $\P$ be a noncrossing partition of $[k]$ corresponding to $\vcc$ and labels $\X$, and $s$ be a point in $[k]$. In this case $\P$ contains a polygon $P_c$ of size $t_c = t_a+t_b-1$. Assume $P_c= \{p_1,\dots,p_{t_c}\}$ such that $p_i$'s are ordered in the clockwise direction starting from $s$. Note that $s$ might be a point in $P_c$ and then in that case $p_1=s$. Now we will split $P_c$ into $P_a$ and $P_b$ in the following way:
        \begin{enumerate}[i.]
            \item Let $P_b = \{p_1,\dots,p_{t_b}\}$.
            \item Let $P_a = \{q_1,\dots,q_{t_a}\}$ where $q_i = p_{t_b+i-1}+\(s-p_{t_c}\) \mod k$ for each $i\in[t_a]$. i.e. We shift the $\{p_{t_b},\dots, p_{t_c}\}$ polygon $\(s-p_{t_c} \mod k\)$ units in the clockwise direction so that the last vertex of the $t_a$-gon touches $s$.
        \end{enumerate}
        
        Note that in the case of $s=p_1$, $P_a$ and $P_b$ touches at $s$, which is allowed for $\np\(\vaa,\vbb\)$. The resulting $P_a$ and $P_b$ correspond to a partition in $\nplx\(\vaa,\vbb,\X\)$. See \Cref{fig:np-two-poly-ab-split} for illustration.
        
        \begin{figure}[hbt!]
            \centering
            \includegraphics[scale=0.3]{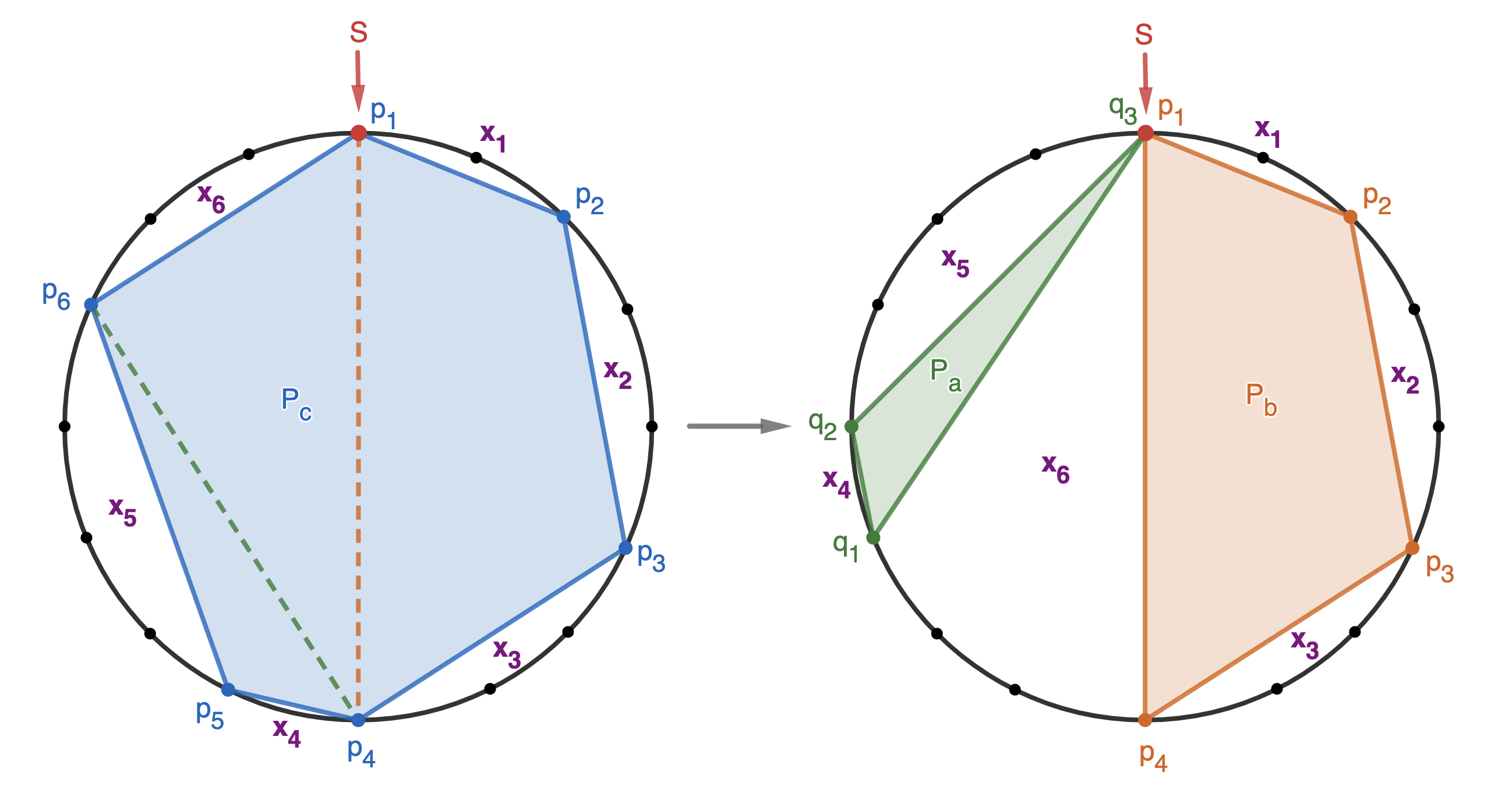}
            \caption{Illustration of $\nplx\(\vcc, \X\)\times [k]\to \nplx\(\vaa,\vbb,\X\)$: here $s=p_1$, $t_a=3$, $t_b=4$ and $t_c=6$.}
            \label{fig:np-two-poly-ab-split}
        \end{figure}
        
        \item It is not hard to see that the above two operations are inverses of each other.
    \end{enumerate}
    
    Thus $\abs{\nplx\(\vaa,\vbb,\X\)} = \big| \nplx\(\vcc, \X\)\times [k] \big| = \big| \nplx\(\vcc,\X\) \big| \cdot k = k\cdot(p-1)!\cdot k = k^2(p-1)!$. To consider $\nplax\(\vaa,\vbb,\X\)$, we take into account the fact that all the $(a-1)$ points in $\vaa$ and $(b-1)$ points in $\vbb$ can be permuted, thus $\nplax\(\vaa,\vbb,\X\) = k^2(p-1)!\cdot(a-1)!\cdot(b-1)!$.
\end{proof}

\begin{thm}\label{thm:num-of-partitions-2}
    Let $\vaa = \(\a_1,\dots,\a_k\), \vbb=\(\b_1,\dots,\b_k\)\in P_k$ and $a=\a_1+\dots+\a_k$, $b=\b_1+\dots,\b_k$, $p=2k-a-b+1$. Let $\X=\{x_1,\dots,x_p\}$ be such that $x_1+\dots+x_p=k$. Then 
    \begin{equation}
        \abs{\nplax\(\vaa,\vbb,\X\)} = k^2(p-1)!\cdot (a-1)!\cdot (b-1)!.
    \end{equation}
\end{thm}

\begin{proof}
     We will applying the inductive combining argument for $\nplax\(\vaa,\X\)$ from the last section to $\nplax\(\vaa,\vbb,\X\)$, until we reach the case when $\o$ has only one polygon. Each combining step of a $t_1$-gon and $t_2$-gon gives a $t_1+t_2-1$-gon and an extra point. So in the end we will obtain a polygon of size $2\a_2+3\a_3+\dots+k\a_k-(\a_2+\dots+\a_k-1) = (k-\a_1)-(a-\a_1-1) = k-a+1$ and $\a_1+(\a_2+\dots+\a_k-1) = a-1$ points. Let $\vaa'$ corresponds to $\o_2^{a-1}\o_{2(k-a+1)}$. Thus $\nplax\(\vaa,\vbb,\X\) = \nplax\(\vv{\a'},\vbb,\X\)$. 
     
     Similarly we can apply the same argument to $\vbb$ and result in $\o_2^{b-1}\o_{2(k-b+1)}$. Let $\vbb'$ correspond to $\o_2^{b-1}\o_{2(k-b+1)}$. Now $\nplax\(\vaa,\vbb,\X\) = \nplax\(\vv{\a'},\vbb,\X\) = \nplax\(\vv{\a'},\vv{\b'},\X\)$.
     
     Note that $a'=a-1+1 = a$, $b'=b-1+1 = b$ and $p'=(k-a+1)+(k-b+1)-1 = 2k-a-b+1 = p$. By \Cref{lem:num-of-partitions-polygon-2}, $\nplax\(\vv{\a'},\vv{\b'},\X\) = k^2(p'-1)!(a'-1)!(b'-1)! = k^2(p-1)!(a-1)!(b-1)!$.
\end{proof}

\begin{cor}
    Let $\vaa = \(\a_1,\dots,\a_k\), \vbb=\(\b_1,\dots,\b_k\)\in P_k$ and $a=\a_1+\dots+\a_k$, $b=\b_1+\dots,\b_k$, $p=2k-a-b+1$. Let $\X=\{x_1,\dots,x_p\}$ be such that $x_1+\dots+x_p=k$. Then 
    \begin{equation}
        \abs{\nplx\(\vaa,\vbb,\X\)} = k^2(p-1)!\cdot\dfrac{(a-1)!}{\a_1!\dots\a_k!}\cdot \dfrac{(b-1)!}{\b_1!\dots\b_k!}.
    \end{equation}
\end{cor}

\comm{
\begin{thm}\label{thm:coefficient-step-3}
    \begin{equation}
        C\(\vaa,\vbb\)=(-1)^{a+b-k-1}\cdot\dfrac{(a-1)!}{\a_1!\dots\a_k!}\cdot \dfrac{(b-1)!}{\b_1!\dots\b_k!}\cdot k\cdot \binom{a+b-2}{k-1}\,.
    \end{equation}
\end{thm}
}

Now we are ready to prove \Cref{thm:moments-ab-step2}.
\momentMiddleStep*

\begin{proof}
By \Cref{thm:num-of-partitions-2} and \Cref{cor:sum-of-prod-of-catalan}, 
    \begin{align*}
        & \sum_{\substack{\pi\in \np(\vaa), \sigma\in \np(\vbb): \\ \pi\oplus\sigma \in NC(k)}} \prod_{x\in S_{\pi\oplus\sigma}} C_{x-1} 
        = \sum_{(\pi, \sigma)\in \np\(\vaa,\vbb\)}\, \prod_{i_j\in S_{\pi\oplus \sigma}} C_{i_j-1}\\
        &= \sum_{i_j\geq 1:\, i_1+\dots+i_p=k}\, \dfrac{\abs{\np\(\vaa,\vbb,i_1,\dots,i_p\)}}{\perm\(i_1,\dots,i_p\)} \cdot C_{i_1-1} \cdot\dots\cdot C_{i_p-1} \\
        &= \sum_{i_j\geq 1:\, i_1+\dots+i_p=k}\, \dfrac{\abs{\nplx\(\vaa,\vbb,i_1,\dots,i_p\)}}{p\,!}\cdot C_{i_1-1} \cdot\dots\cdot C_{i_p-1} \\
        &= \sum_{i_j\geq 1:\, i_1+\dots+i_p=k}\, \(\dfrac{(a-1)!}{\a_1!\dots\a_k!}\cdot\dfrac{(b-1)!}{\b_1!\dots\b_k!} \cdot k^2\cdot (p-1)!/ p\,!\) \cdot C_{i_1-1} \cdot\dots\cdot C_{i_p-1} \\
        &= \dfrac{(a-1)!}{\a_1!\dots\a_k!}\cdot\dfrac{(b-1)!}{\b_1!\dots\b_k!} \cdot k^2\cdot \dfrac{\,1\,}{p} \cdot \(p\cdot \dfrac{(2k-p-1)!}{(k-p)!\,k!}\) \\
        &= \dfrac{(a-1)!}{\a_1!\dots\a_k!}\cdot\dfrac{(b-1)!}{\b_1!\dots\b_k!} \cdot k\cdot \dfrac{\(a+b-2\)!}{(a+b-k-1)!\,(k-1)!}\\
        &= \dfrac{(a-1)!}{\a_1!\dots\a_k!}\cdot\dfrac{(b-1)!}{\b_1!\dots\b_k!} \cdot k\cdot \binom{a+b-2}{k-1}\,.
    \end{align*}
\end{proof}



\section{Special Case: \texorpdfstring{$\ozm\opr \o$}{Ωz(m) ºR Ω} }
\label{section:ozm-1}

\setlength{\parskip}{1.5mm}
\setlength{\baselineskip}{1.3em}

In this section, we analyze the random matrix $M_{Z(m),\o} = D'RD$ where $\o'=\ozm$ and $\o$ is an arbitrary distribution. i.e. $\o_{2i}' = C(i,m)$ for all $i$.


Our main result is as follows.

\begin{thm}\label{thm:moments-mzshape-1-distr}
Let $\ozm$ be defined as in \Cref{defn:ozm-ozms}. Then
\begin{equation}
    \(\ozm\opr \o\)_{2k} = \sum_{\vaa_i\in P_k}\, \binom{mk}{a-1}\cdot \dfrac{(a-1)!}{\a_1!\dots\a_k!}\cdot \vv{\o}^{\,\vaa}\,.
\end{equation}
\end{thm}

After proving this result, we prove \Cref{thm:main-3-circ} by plugging in $\o=\ozmm$ to the above theorem.
\mainthreecirc*

\subsection{Proof of the Formula}\label{subsection:ozm-1-pf}

\setlength{\parskip}{1.5mm}
\setlength{\baselineskip}{1.3em}

\begin{defn}
    Let $\o' = \ozm$. Given a distribution $\o$, we define $A_m(k,0)$ to be 
    \begin{align}
        A_m(k,0) &= \sum_{\vaa,\vbb\in P_k}\, C\(\vaa,\vbb\)\cdot \voo^{\,\vaa} \cdot \vv{\o'}^{\,\vbb} = \sum_{\vaa,\vbb\in P_k}\, C\(\vaa,\vbb\)\cdot \voo^{\,\vaa} \cdot \(\prod_{i=1}^k C(i,m)^{\,\b_i}\)\,.
    \end{align}
    where
    \begin{equation}
        C(\vaa,\vbb) = (-1)^{a+b-k-1}\cdot k\cdot \binom{a+b-2}{k-1}\cdot \dfrac{(a-1)!}{\a_1!\dots\a_k!}\cdot \dfrac{(b-1)!}{\b_1!\dots\b_k!}\,.
    \end{equation}
\end{defn}

Note that by \Cref{thm:main-1-circ}, $A_m(k,0) = \(\ozm\opr \o\)_{2k}$. A restatement of \Cref{thm:moments-mzshape-1-distr} is the following.


\begin{thm}\label{thm:expression-for-Amk}
    \begin{equation}
        A_m(k,0) = \sum_{\vaa\in P_k} \binom{mk}{a-1}\cdot\dfrac{(a-1)!}{\a_1!\dots\a_k!}\cdot \voo^{\vaa}\,.
    \end{equation}
\end{thm}

To prove the above theorem, we need the following identities.

\begin{prop}\label{prop:alternating-sum-of-binomial-coeff}
\begin{equation}
    \sum_{i=0}^{m-k} (-1)^{i+(m-k)}\cdot\binom{m-i}{k}\cdot \binom{n}{i} = \binom{n-k-1}{m-k}\,.
\end{equation}
\end{prop}
\begin{proof}
    Note that \Cref{prop:alternating-sum-of-binomial-coeff} is always true when $k>m$ since it becomes $0=0$. W.O.L.G we will assume $k\leq m$. We prove this by induction on $m$.
    \begin{enumerate}
        \item $m=0:$ Since $k\leq m<n$, $k=0$ and $n\geq 1$. Then $\displaystyle LHS = \binom{0}{0}\cdot \binom{n}{0} = 1$ and $\displaystyle RHS = \binom{n-1}{0}=1$.
        \item $m\implies (m+1)$: Assume \Cref{prop:alternating-sum-of-binomial-coeff} is true for all triples $(n',m',k')$ where $m'\leq m$ and $k'\leq m'<n'$. For any $n,k$ such that $k\leq m+1 < n$,
        \begin{align*}
            & \sum_{i=0}^{m+1-k}\, (-1)^{i+m+1-k}\cdot  \binom{m+1-i}{k}\cdot \binom{n}{i}\\
            =& \sum_{i=0}^{m+1-k}\, (-1)^{i+m+1-k}\cdot \(\binom{m-i}{k}+ \binom{m-i}{k-1}\)\cdot \binom{n}{i}\\
            =& (-1)\cdot \(\sum_{i=0}^{m-k}\, (-1)^{i+m-k}\cdot \binom{m-i}{k}\cdot \binom{n}{i}\) + \binom{k-1}{k}\cdot \binom{n}{m+1-k} \\ 
            &\hspace{1cm} + \sum_{i=0}^{m-(k-1)}\, (-1)^{i+m-(k-1)}\cdot \binom{m-i}{k-1}\cdot \binom{n}{i}\\
            =& -\binom{n-k-1}{m-k} + 0 + \binom{n-(k-1)-1}{m-(k-1)}\\
            =& -\binom{n-k-1}{m-k} + \binom{n-k}{m-k+1}
            = \binom{n-k-1}{m-k+1} = \binom{n-k-1}{(m+1)-k}\,.
        \end{align*}
    \end{enumerate}
\end{proof}

\begin{prop}\label{prop:alternating-sum-of-binomial-coeff-2}
\begin{equation}
    \sum_{i=0}^n\, (-1)^{n-i}\cdot \binom{m+i}{k}\binom{n}{i} = \binom{m}{k-n}\,.
\end{equation}
\end{prop}
\begin{proof}
    We prove this by induction on $n$.
    \begin{enumerate}
        \item $n=0$: $\displaystyle LHS = \binom{m+0}{k}\binom{0}{0} = \binom{m}{k} = RHS$.
        \item $n\implies (n+1)$: Assume $\displaystyle\sum_{i=0}^n\, (-1)^{n-i}\cdot \binom{m+i}{k}\binom{n}{i} = \binom{m}{k-n}$. We want to prove that $\displaystyle\sum_{i=0}^{n+1}\, (-1)^{n+1-i}\cdot \binom{m+i}{k}\binom{n+1}{i} = \binom{m}{k-(n+1)}$:
        \begin{align*}
             &\sum_{i=0}^{n+1}\, (-1)^{n+1-i}\cdot \binom{m+i}{k} \binom{n+1}{i} \\
            =& \sum_{i=0}^{n+1}\, (-1)^{n+1-i}\cdot \binom{m+i}{k}\(\binom{n}{i} + \binom{n}{i-1}\)\\
            =& -\sum_{i=0}^{n}\, (-1)^{n-i}\cdot \binom{m+i}{k}\binom{n}{i} +  \sum_{i=1}^{n+1} (-1)^{(n+1)-i}\cdot \binom{(m+1)+(i-1)}{k} \binom{n}{i-1}\\
            =& -\binom{m}{k-n} + \binom{m+1}{k-n} 
            = \binom{m}{k-n-1} = \binom{m}{k-(n+1)}.
        \end{align*}
    \end{enumerate}
\end{proof}

\begin{lemma}\label{lem:sum-of-prod-of-m-catalan}
Assume $k\geq t\geq 0$. Then
    \begin{equation}
        \sum_{i_j\geq 1\,:\,i_1+\dots+i_t=k} C(i_1,m)\dots C(i_t,m) = \dfrac{\,t\,}{k}\cdot \binom{(m+1)k}{k-t}
    \end{equation}
\end{lemma}

\begin{figure}[hbt!]
    \centering
    \includegraphics[scale=0.35]{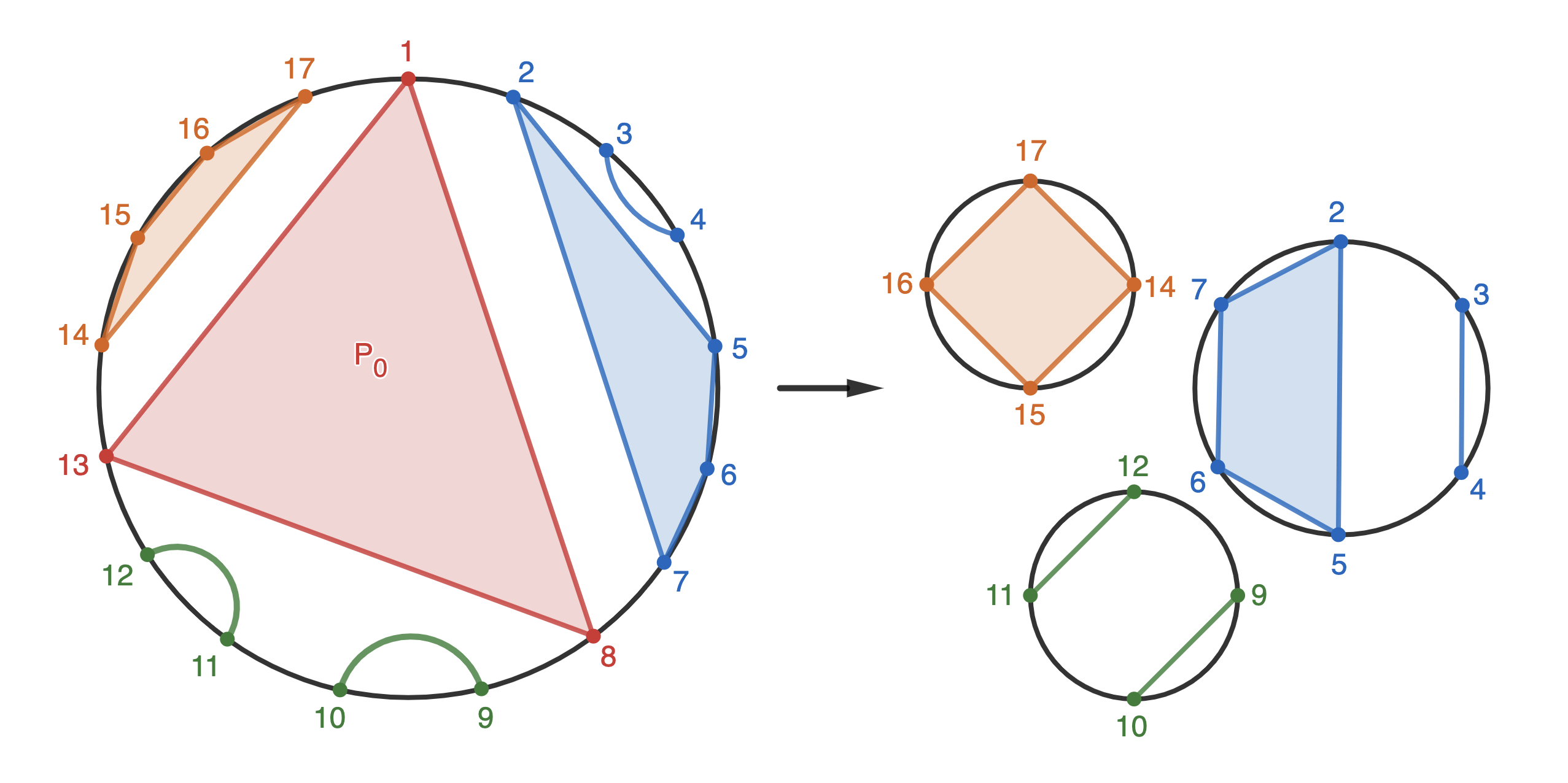}
    \caption{Illustration of \Cref{lem:sum-of-prod-of-m-catalan}: here $t=3$, $m=2$, $k=7$. For this configuration, $i_1=3$ (blue), $i_2=2$ (green), $i_3=2$ (orange).}
    \label{fig:sum-prod-catalan}
\end{figure}

\begin{proof}
    By \Cref{cor:m-catalan-nnp}, $C(k,m)$ is the number of non-crossing partitions of $[(m+1)k]$ with $k$ parts of sizes $(m+1)$. Therefore we can view LHS as the number of ways to put $k$ polygons of size $(m+1)$ and one polygon $P_0$ of size $t$ on the cycle of length $(m+1)k+t$ non-crossingly, with the following requirements:
    \begin{enumerate}
        \item $P_0$ must contain point $1$, and
        \item the $t$ parts on the cycle divided by $P_0$ are nonempty and contain some $(m+1)$-gon. i.e. if $P_0=\{p_1,\dots,p_t\}$ where $1=p_1<p_2<\dots<p_t$, then $p_{i+1}-p_i>1$ for $i=1,\dots,t-1$ and $p_t<k$.
    \end{enumerate}
    
    We can count the LHS in an alternative way: 
    \begin{enumerate}
        \item Count the number of non-crossing partitions of $[(m+1)k+t]$ using $k$ polygons of size $(m+1)$ and one polygon of size $t$. By \Cref{thm:num-of-np}, this is $\displaystyle \binom{(m+1)k+t}{(k+1)-1}\cdot \dfrac{k!}{k!1!} = \binom{(m+1)k+t}{k}$.
        \item We want to restrict that the $t$-gon has to contain $1$, this gives us a factor of $\dfrac{\,t\,}{(m+1)k+t}$.
        \item The above also counts the case when some of the parts from the $t$-gon are empty, so we will use inclusion-exclusion principle for the counting. If a fixed $(t-i)$ parts out of $t$ parts have to be empty, then it is equivalent to counting such non-crossing partitions of $[(m+1)k+i]$ with $k$ $(m+1)$-gons and one $i$-gon. This is $\displaystyle \dfrac{i}{(m+1)k+i}\cdot \binom{(m+1)k+i}{k}$. There are $\displaystyle \binom{t}{t-i}=\binom{t}{i}$ ways for choosing the $(t-i)$ empty parts. 
    \end{enumerate}
    
    Thus
    \begin{align*}
        LHS &= \sum_{i=1}^t\, (-1)^{t-i}\cdot \dfrac{i}{(m+1)k+i}\cdot \binom{(m+1)k+i}{k}\cdot \binom{t}{i}\\
        &= \sum_{i=1}^t\, (-1)^{t-i}\cdot \dfrac{\((m+1)k+i-1\)!}{k!(mk+i)!}\cdot \dfrac{t!}{(i-1)!(t-i)!}\\
        &= \sum_{i=1}^t\, (-1)^{t-i}\cdot \dfrac{\,t\,}{k}\cdot \binom{(m+1)k+i-1}{k-1} \binom{t-1}{i-1}\\
        &= \dfrac{\,t\,}{k}\cdot \(\sum_{j=0}^{t-1}\, (-1)^{t-1-j}\cdot \binom{(m+1)k+ j}{k-1} \binom{t-1}{j}\) \hspace{0.5cm}\(\text{Let } j=i-1\) \\
        &= \dfrac{\,t\,}{k}\cdot \binom{(m+1)k}{k-t} \hspace{1cm} \(\text{By } \Cref{prop:alternating-sum-of-binomial-coeff-2} \)
    \end{align*}
\end{proof}

\begin{proof}[Proof of \Cref{thm:expression-for-Amk}]
    \begin{align*}
         A_m(k,0) 
         &= \sum_{\vaa,\vbb\in P_k}\, C\(\vaa,\vbb\)\cdot \voo^{\,\vaa} \cdot \(\prod_{i=1}^k C(i,m)^{\,\b_i}\) \\
         &= \sum_{\vaa,\vbb\in P_k}\, (-1)^{a+b-k-1}\cdot k\cdot\binom{a+b-2}{k-1}\cdot \dfrac{(a-1)!}{\a_1!\dots\a_k!}\cdot \dfrac{(b-1)!}{\b_1!\dots\b_k!}\cdot \prod_{i=1}^k C(i,m)^{\,\b_i} \cdot \voo^{\vaa}\\
         &= \sum_{\vaa\in P_k}\,\(\sum_{\vbb\in P_k} (-1)^{a+b-k-1}\cdot k\cdot\binom{a+b-2}{k-1}\cdot \dfrac{(b-1)!}{\b_1!\dots\b_k!}\cdot \prod_{i=1}^k C(i,m)^{\,\b_i}\)\cdot \dfrac{(a-1)!}{\a_1!\dots\a_k!}\cdot \voo^{\vaa}\\
    \end{align*}
    
    Thus it is equivalent to prove that 
    \begin{align}
        \sum_{\vbb\in P_k}(-1)^{a+b-k-1}\cdot\dfrac{(b-1)!}{\b_1!\dots\b_k!}\cdot k\cdot\binom{a+b-2}{k-1}\cdot \prod_{i=1}^k C(i,m)^{\,\b_i} = \binom{mk}{a-1}\,.
    \end{align}
    
    By \Cref{lem:sum-of-prod-of-m-catalan},
    \begin{align*}
         & \sum_{\vbb\in P_k} (-1)^{a+b-k-1}\cdot\dfrac{(b-1)!}{\b_1!\dots\b_k!}\cdot k\cdot\binom{a+b-2}{k-1} \cdot\prod_{i=1}^k C(i,m)^{\,\b_i}\\
        =& \sum_{b=1}^k\,\sum_{\vbb\in P_k:\sum \b_i=b}\, (-1)^{a+b-k-1}\cdot k\cdot\binom{a+b-2}{k-1}\cdot\dfrac{\,1\,}{b}\cdot\(\dfrac{b!}{\b_1!\dots\b_k!}\cdot\prod_{i=1}^k C(i,m)^{\,\b_i}\)\\
        =& \sum_{b=1}^k\, (-1)^{a+b-k-1}\cdot \binom{a+b-2}{k-1}\cdot\dfrac{\,k\,}{b}\cdot \(\sum_{i_j\geq 1\,:\, i_1+\dots+i_b=k}\; \prod_{j=1}^b\,C\(i_j,m\)\)\\
        =& \sum_{b=1}^k\, (-1)^{a+b-k-1}\cdot \binom{a+b-2}{k-1}\cdot\dfrac{\,k\,}{b}\cdot \dfrac{\,b\,}{k} \cdot \binom{(m+1)k}{k-b}\\
        =& \sum_{b=k-(a-1)}^k\, (-1)^{a+b-k-1}\cdot \binom{a+b-2}{k-1}\cdot \binom{(m+1)k}{k-b}\\
        =& \sum_{i=0}^{a-1} (-1)^{(a-1)-i}\, \binom{(a-1)+(k-1)-i}{k-1}\cdot \binom{(m+1)k}{i} \hspace{1cm} (\text{Let } i=k-b)\\
    \end{align*}
    
    Plugging in $n'=(m+1)k$, $m'=(a-1)+(k-1)$, $k'=(k-1)$ for \Cref{prop:alternating-sum-of-binomial-coeff}, 
    \begin{align*}
         &\sum_{i=0}^{a-1} (-1)^{(a-1)-i}\, \binom{(a-1)+(k-1)-i}{k-1}\cdot \binom{(m+1)k}{i} \hspace{1cm} \\
        =& \sum_{i=0}^{m'-k'} (-1)^{(m'-k')-i}\, \binom{m'-i}{k'}\cdot \binom{n'}{i}\\
        =& \binom{n'-k'-1}{m'-k'} = \binom{(m+1)k-(k-1)-1}{a-1} = \binom{mk}{a-1}
    \end{align*}
    as needed.
\end{proof}

\subsection{Application when \texorpdfstring{$\o=\o_{Z(m')}$}{Ω = Ωz(m)} }\label{subsection:ozm-1-zn}

\setlength{\parskip}{1.5mm}
\setlength{\baselineskip}{1.3em}



As a direct application of the results from the previous section, in this section we prove \Cref{thm:main-3-circ}.
\mainthreecirc*

To prove \Cref{thm:main-3-circ}, we recall the main result from the paper \cite{CP20}.
\begin{thm}\label{thm:moments-zmshape}
    Let $M_{\mzshape}$ be the $Z(m)$-shape graph matrix. Let $\displaystyle M_{n,m} = \dfrac{1}{n^{m/2}} M_{\mzshape}$ and $\displaystyle r(n,m) = \dfrac{n!}{(n-m)!}$. Then 
    \begin{equation}
        \lim_{n\to\infty}\, \dfrac{1}{r(n,m)}\cdot \Ebb\[\trace\(\(M_{n,m}^TM_{n,m}\)^k\)\] = C(k,m)
    \end{equation}
    for all $k\in\Nbb$.
\end{thm}

Rewriting the above result in terms of $\ozm$, we get the following corollary.
\begin{cor}\label{cor:moments-zmshape}
    For any $k,m\geq 0$,
    \begin{equation}
        \(\ozm\)_{2k} = C(k,m).
    \end{equation}
\end{cor}


By the Trace Power Method and \Cref{cor:moments-zmshape}, \Cref{thm:main-3-circ} can be proved by proving the following statement.

\begin{restatable}{thm}{mainthree}
\label{thm:main-3}

\begin{equation}
    \(\ozm\opr \ozmm\)_{2k} = C\(k,m+m'\) = \(\o_{Z\(m+m'\)}\)_{2k}\,.
\end{equation}

\end{restatable}

By \Cref{thm:main-1-circ} and \Cref{cor:moments-zmshape}, it suffices to prove the following.

\begin{thm}\label{thm:moments-zmznshape}
    Let $\o_{2k}=C(k,m)$ and $\o'_{2k}=C(k,m')$  for all $k$. Then
    \begin{equation}
        \sum_{\vaa,\vbb\in P_k} C\(\vaa,\vbb\)\cdot \voo^{\vaa} \cdot\vv{\o'}^{\vbb} = C\(k,m+m'\) \,.
    \end{equation}
\end{thm}

An easy combinatorial identity is needed for the proof.
\begin{prop}\label{prop:combinatorial-id-1}
For any $m,n,k\in \Zbb_{\geq 0}$, 
    \begin{equation}
        \sum_{i=0}^k \binom{m}{i} \binom{n}{k-i} = \binom{m+n}{k}.
    \end{equation}
\end{prop}

\begin{proof}
    The RHS is the number of ways to choose $k$ elements from a collection of $m+n$ elements. Each summand in the LHS is the number of ways to choose $i$ elements from the first $m$ elements, and $k-i$ elements from the last $n$ elements. Summing up all the possible $i$'s gives the equality.
\end{proof}

Now we are ready to prove the theorem.
\begin{proof}[Proof of \Cref{thm:moments-zmznshape}]
    By \Cref{thm:expression-for-Amk}, when $\o'_{2k}=C\(k,m'\)$,
    \begin{align*}
        \sum_{\vaa,\vbb\in P_k} C\(\vaa,\vbb\)\cdot \voo^{\vaa} \cdot\vv{\o'}^{\vbb}
        = \sum_{\vaa\in P_k}\, \binom{m'k}{a-1}\cdot \dfrac{(a-1)!}{\a_1!\dots\a_k!}\cdot \voo^{\vaa}.
    \end{align*}
    
    Plugging in $\o_{2k}=C(k,m)$, we get
    \begin{align*}
         & \sum_{\vaa,\vbb\in P_k}\, C\(\vaa,\vbb\)\cdot \voo^{\vaa} \cdot\vv{\o'}^{\vbb} \\
        =& \sum_{\vaa\in P_k}\, \binom{m'k}{a-1}\cdot \dfrac{(a-1)!}{\a_1!\dots\a_k!}\cdot\prod_{i=1}^k C(i,m)^{\a_i} \\
        =& \sum_{a=1}^k\, \binom{m'k}{a-1}\cdot\dfrac{\,1\,}{a}\cdot \(\sum_{i_j\geq 1: i_1+\dots+i_a = k}\, \prod_{j=1}^a C(i_j,m)\) \\
        =& \sum_{a=1}^k\, \binom{m'k}{a-1}\cdot\dfrac{\,1\,}{a}\cdot \dfrac{\,a\,}{k}\cdot \binom{(m+1)k}{k-a} \hspace{1cm} \(\text{ By  \Cref{lem:sum-of-prod-of-m-catalan}}\) \\
        =& \dfrac{\,1\,}{k}\cdot\sum_{a=1}^k\, \binom{m'k}{a-1}\cdot \binom{(m+1)k}{k-a}\\
        =& \dfrac{\,1\,}{k}\cdot\sum_{i,j\geq 0: i+j=k-1}\, \binom{m'k}{i}\cdot \binom{(m+1)k}{j} \\
        =& \dfrac{\,1\,}{k}\cdot \binom{(m'+m+1)k}{k-1} \hspace{1cm} \(\text{ By  \Cref{prop:combinatorial-id-1}}\)\\
        =& \dfrac{\((m+m'+1)k\)!}{k!\cdot \((m+m')k+1\)!}
        = \dfrac{1}{(m+m')k+1}\cdot \binom{(m+m'+1)k}{k}
        = C\(k,m+m'\)\,.
    \end{align*}
\end{proof}

\section{ General Case: \texorpdfstring{$\ozm\opr \o^{(1)}\opr \dots \opr\o^{(s)}$}{Ωz(m) ºR Ω1 ºR …ºR Ωs}} \label{section:ozm-s}

\setlength{\parskip}{1.5mm}
\setlength{\baselineskip}{1.3em}



In this section, we generalize the result from the previous section. 

\begin{defn}\label{defn:Cm(vaa)}
    Given $\vaa_1,\dots,\vaa_s\in P_k$, we will denote $C_m\(\vaa_1,\dots,\vaa_s\)$ to be 
    \begin{equation}
        C_m\(\vaa_1,\dots,\vaa_s\) = \binom{(m-s+1)k}{mk - (a_1+\dots+a_s) + 1}\cdot k^{s-1}\cdot \prod_{i=1}^s \dfrac{(a_i-1)!}{\a_{i1}!\dots\a_{ik}!}
    \end{equation}
\end{defn}

\begin{thm}\label{thm:moments-mzshape-s-distr}
Given $\o^{(1)},\dots,\o^{(s)}$, 
\begin{equation}
    \(\ozm\opr \o^{(1)}\opr\dots \opr \o^{(s)}\)_{2k} = \sum_{\vaa_i\in P_k}\, C_m\(\vaa_1,\dots,\vaa_s\)\cdot \vv{\o^{(1)}}^{\,\vaa_1}\dots \vv{\o^{(s)}}^{\,\vaa_s}\,.
\end{equation}
\end{thm}

We will prove the theorem by induction on $s$. If we know the formula for \\ $\displaystyle \(\ozm\opr \o^{(1)}\opr\dots \opr \o^{(s-1)}\)_{2k}$, we can compute $\displaystyle\(\ozm\opr \o^{(1)}\opr\dots \opr \o^{(s)}\)_{2k}$ by using the formula from \Cref{thm:main-1-circ}:
\begin{align*}
    \(\ozm\opr \o^{(1)}\opr\dots \opr \o^{(s)}\)_{2k} 
    &= \(\(\ozm\opr \o^{(1)}\opr\dots \opr \o^{(s-1)}\)\opr \o^{(s)}\)_{2k} \\
    &= \sum_{\vaa,\vbb\in P_k}\, C\(\vaa,\vbb\)\cdot \(\vv{\ozm\opr \o^{(1)}\opr\dots \opr \o^{(s-1)}}\)^{\vaa}\cdot \vv{\o^{(s)}}^{\,\vbb}\,.
\end{align*}

\subsection{Extended Definition of Non-crossing Partitions}
\label{subsection:ozm-s-np-general}

\setlength{\parskip}{1.5mm}
\setlength{\baselineskip}{1.3em}

Recall that $\np\(\vaa\)$ is the set of noncrossing partitions of $[k]$ such that there are $\a_i$ many parts of size $i$ and $\np\(m\vaa\)$ is the set of noncrossing partitions of $[k]$ such that there are $\a_i$ many parts of size $mi$.

To prove our main results, we need to first extend our definitions and results on $\np\(\vaa\)$ and $\np\(m\vaa\)$.

\begin{defn}
    Let $V\subseteq[n]$ and $\vaa\in P_k$. We say that $\P\in\np\(V, m\vaa\)$ if $\abs{V}=mk$ and if $a_1< \dots < a_{mk}$ are the elements of $V$ in ascending order, and $\P' = \big\{\{i_1,\dots,i_p\}: \{a_{i_1},\dots,a_{i_p}\}\in \P\big\}$, then $\P'\in\np\(m\vaa\)$.
\end{defn}

\begin{eg}
    Assume $V=\{1,2,4,5,7,8\}\subseteq [9]$ and $m=2$, $k=3$, $\vaa=(1,1,0)$. Then $\P=\left\{ \{1,2\}, \{4,5,7,8\} \right\}\in\np\(V,m\vaa\)$ since $(a_1,\dots,a_6) = (1,2,4,5,7,8)$ and  $\P'=\left\{\{1,2\},\{3,4,5,6\}\right\}\in\np\(m\vaa\)$.
\end{eg}

\begin{defn}\label{defn:np-general}
    Given $\vaa_1,\dots,\vaa_s\in P_k$ and $m_1,\dots,m_s\in \Zbb^{+}$, let $m=m_1+\dots+m_s$. We define $\N\P\(m_1\vaa_1,\dots,m_s\vaa_s\)$ to be the set of non-crossing partitions of $[mk]$ corresponding to $m_1\vaa_1,\dots,m_s\vaa_s$, defined as follows:
    \begin{itemize}
        \item For each $t\in[s]$, let $V_t = \big\{ \(m_1+\dots+m_{t-1}\)+m(i-1)+j:i\in[k], j\in[m_t] \big\}$. Note that $\abs{V_t} = m_tk$. 
        \item $\P\in \np\(m_1\vaa_1,\dots,m_s\vaa_s\) \iff$ $\P = \P_1\bigsqcup\dots \bigsqcup \P_s$ where $\P_i\in\N\P\(V_i,m_i\vaa_i\)$ for each $i\in[s]$, and $\P$ is noncrossing (i.e. $\P_i$'s do not cross).
    \end{itemize} 
    
    See \Cref{fig:np-general} for an illustration.
    
    \begin{figure}[hbt!]
    \centering
    \includegraphics[scale=0.32]{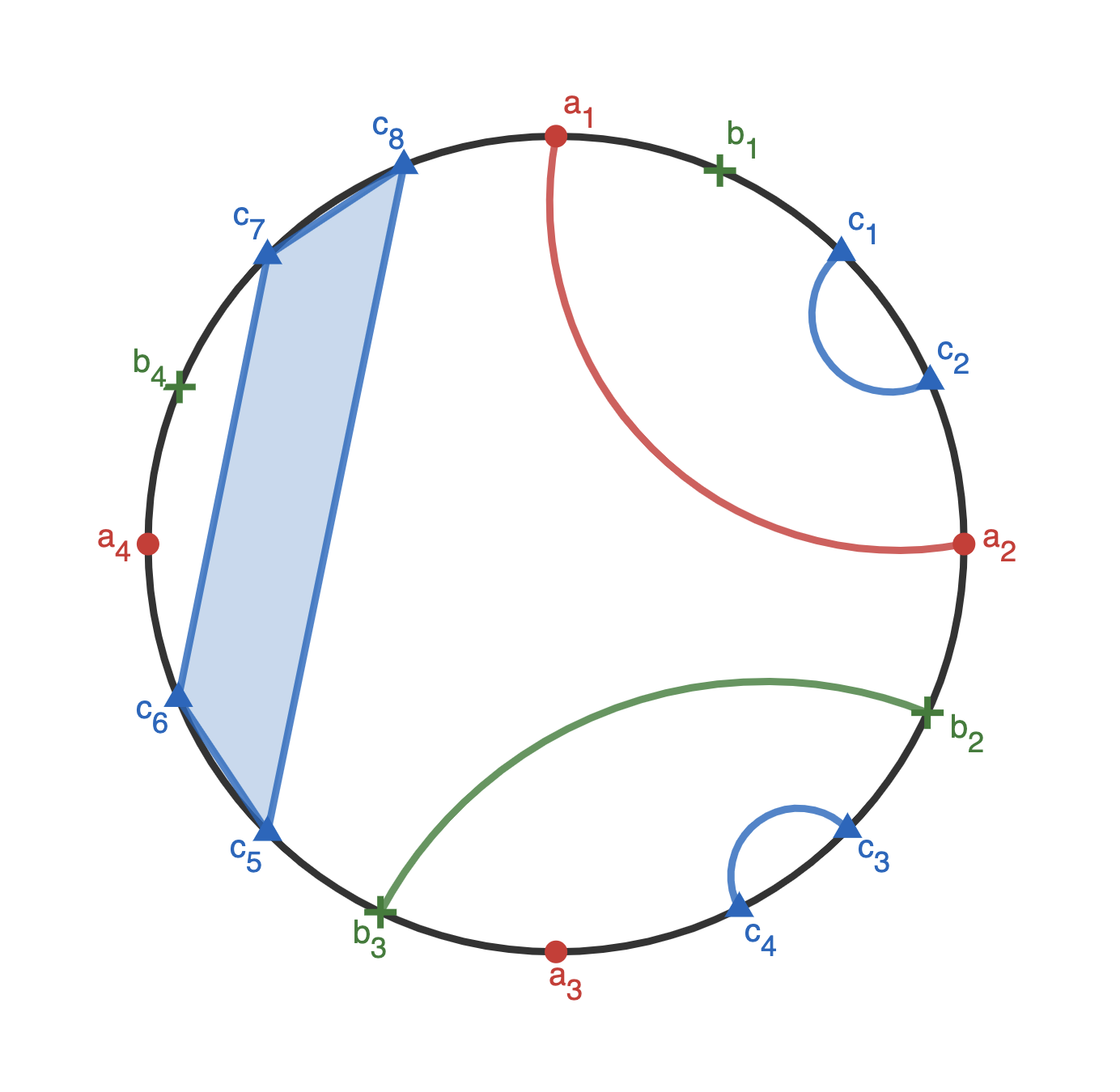}
    \caption{Illustration of \Cref{defn:np-general}: here $s=3, m=4,k=4$ and $m_1=m_2=1, m_3=2$. $V_1=\{a_1,a_2,a_3,a_4\}$, $V_2=\{b_1,b_2,b_3,b_4\}$ and $V_3=\{c_1,\dots,c_8\}$. }
    \label{fig:np-general}
\end{figure}
\end{defn}

\begin{rmk}
    Notice that when $s=1$, the definition of $\N\P\(m_1\vaa_1,\dots,m_s\vaa_s\)$ coincides with $\np\(m\vaa\)$ from Section 4.
\end{rmk}

\begin{defn}
   Given $\vaa_1,\dots,\vaa_s\in P_k$, we denote $\npm\(\vaa_1,\dots,\vaa_s\)$ to be $\np\(\(m-s+1\)\vaa_1,\vaa_2,\dots,\vaa_s\)$.
\end{defn}

\begin{defn}
     We denote 
     \begin{enumerate}
         \item $\np_k = \big\{\P\in \np(\vaa): \vaa\in P_k \big\} =$ the set of noncrossing partitions on $[k]$,
         \item $\np_{k,m} = \big\{\P\in \np(m\vaa): \vaa\in P_k \big\} =$ the set of noncrossing partitions on $[mk]$ where the size of each partition set is a multiple of $[m]$, and
         \item $\np_{k,m,s} = \big\{\P_1\cup \dots\cup \P_s \in \np_m(\vaa_1,\dots,\vaa_s): \vaa_i\in P_k \text{ for } i\in[s]\big\}$.
     \end{enumerate}
\end{defn}

\begin{rmk}
    We will add $\L_{A}$ in the end for all the above notations (i.e. $\pla\(\vaa\)$, $\npla\(\vaa\)$, $\npmla\(\vaa_1,\dots,\vaa_s\)$) to indicate that all polygons are distinct, including the ones of the same size. eg. $\abs{\npla\((k,0,\dots,0)\)} = k!$ while $\abs{\np\((k,0,\dots,0)\)} = 1$.
\end{rmk}

\begin{thm}\label{thm:np-general}
Let $\vaa_1,\dots,\vaa_s\in P_k$. 
    \begin{equation}
        \big|\np_m\(\vaa_1,\dots,\vaa_s\) \big| = C_m\(\vaa_1,\dots,\vaa_s\)\,.
    \end{equation}
    
    Equivalently,
    \begin{equation}
        \Big|\npmla\(\vaa_1,\dots,\vaa_s\) \Big| = \binom{(m-s+1)k}{mk - (a_1+\dots+a_s) + 1}\cdot k^{s-1}\cdot \prod_{i=1}^s (a_i-1)!
    \end{equation}
\end{thm}

To prove \Cref{thm:np-general}, we need the following middle step, which requires an equivalent definition of $\np_m\(\vaa_1,\dots,\vaa_s\)$, similar to the one we have seen in \Cref{subsection:moments-general}.

\begin{defn}\label{defn:np-general-touch}
    Let $\P_i=\lcurb P_{i1},\dots,P_{ik_i}\rcurb$ be partitions of $[n]$ for $i=1,\dots,s$. We say that $\P_1\cup\dots\cup \P_s$ is \emph{non-crossing} if 
    \begin{enumerate}
        \item Each $\P_i$ is non-crossing, and
        \item For any $P_i\in\P_i$ and $P_j\in\P_j$, they do not cross when placed on the cycle $\C_n$, and they touch at at most one point.
    \end{enumerate} 
\end{defn}

\begin{figure}[hbt!]
        \centering
        \includegraphics[scale=0.32]{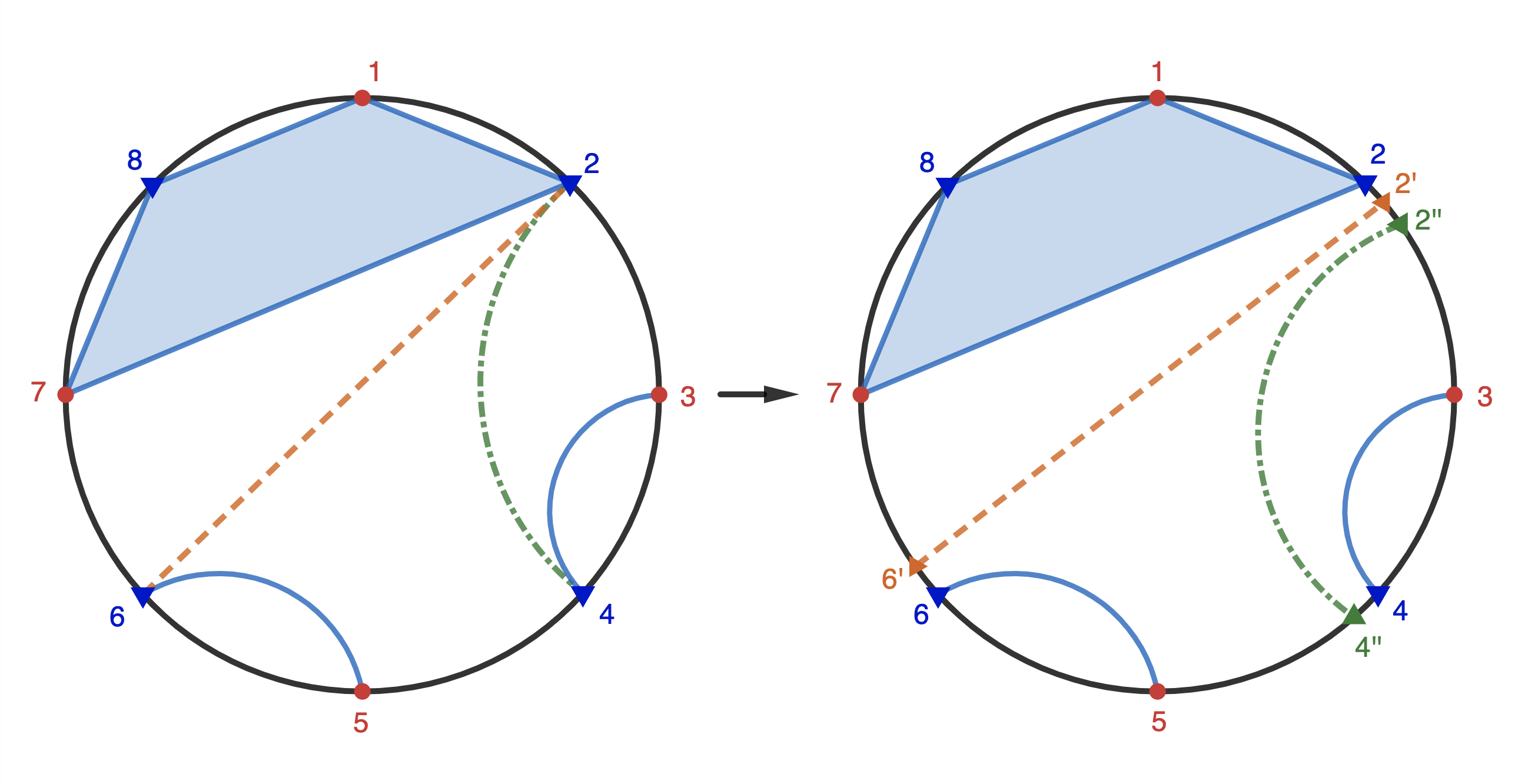}
        \caption{Illustration of \Cref{defn:np-general-touch} and \Cref{defn:np-general-vaa-touch}: here $m=4,s=3,k=4$, $V=\{2,4,6,8\}$. $\vaa_1= \vaa_2= \vaa_3=(2,1,0,0)$. $\P_1=\big\{\{1,2,7,8\}, \{3,4\},\{5,6\} \big\}\in\np\(2\vaa_1\)$ is marked blue, $\P_2=\big\{\{2,6\},\{4\},\{8\} \big\}\in\np\(V,\vaa_2\)$ is marked orange, $\P_3=\big\{\{2,4\},\{6\},\{8\} \big\}\in\np\(V,\vaa_3\)$ is marked green. }
        \label{fig:np-general-touch}
    \end{figure}

\begin{defn}\label{defn:np-general-vaa-touch}
    Let $\vaa_1,\dots,\vaa_s\in P_k$. Let $V=\big\{(m-s+1)i:i\in[k]\big\}$. We say $\P=\P_1\cup\dots \cup\P_s\in\npt_m\(\vaa_1,\dots,\vaa_s\)$ if
    \begin{enumerate}
        \item $\P_1\in\np\big((m-s+1)\vaa_1\big)$ and  $\P_i\in\np\(V, \vaa_i\)$ for $i=2,\dots,s$,
        
        \item $\P_1\cup\dots \cup\P_s$ is non-crossing,
        
        \item For any $P_i\in\P_i$ and $P_j\in\P_j$ where $i<j$ that touch at a point $t_0$, we can order $P_i$ and $Q_j$ as $P_j=\{t_0,p_1,\dots,p_x\}$ and $P_i=\{q_1,\dots,q_y,t_0\}$ such that $t_0,p_1,\dots,p_x,q_1,\dots,q_y$ are ordered in the clockwise direction on $\C_k$.
        
        
        Pictorially, for any point $t_0$ where some $P_j$'s touch together, if we perturb the ``$t_0$'' vertex of each $P_j$ in the clockwise direction $(j-1)\epsilon$ distance for some small enough $\epsilon$, the resulting partition consisting of $P_j'$s is noncrossing.
        
    \end{enumerate}
    
    For an example, see \Cref{fig:np-general-touch}. $P_1=\{1,2,7,8\}\in\P_1$ and $P_2=\{2,6\}\in\P_2$ touch at $2$. If we perturb vertex $2$ of $P_2$ in the clockwise direction to the $2'$ position on the right then $P_2'$ and $P_1$ do not cross.
    
    
\end{defn}

\begin{rmk}
    $\npt_2\(\vaa_1,\vaa_2\) = \np\(\vaa_1,\vaa_2\)$ from \Cref{subsection:moments-general}.
\end{rmk}

\begin{prop}\label{prop:np-general-touch-nontouch-equiv}
\begin{equation}
    \big|\np_m\(\vaa_1,\dots,\vaa_s\)\big| = \big|\npt_m\(\vaa_1,\dots,\vaa_s\)\big|\,.
\end{equation}
\end{prop}

\begin{figure}[hbt!]
    \centering
    \includegraphics[scale=0.32]{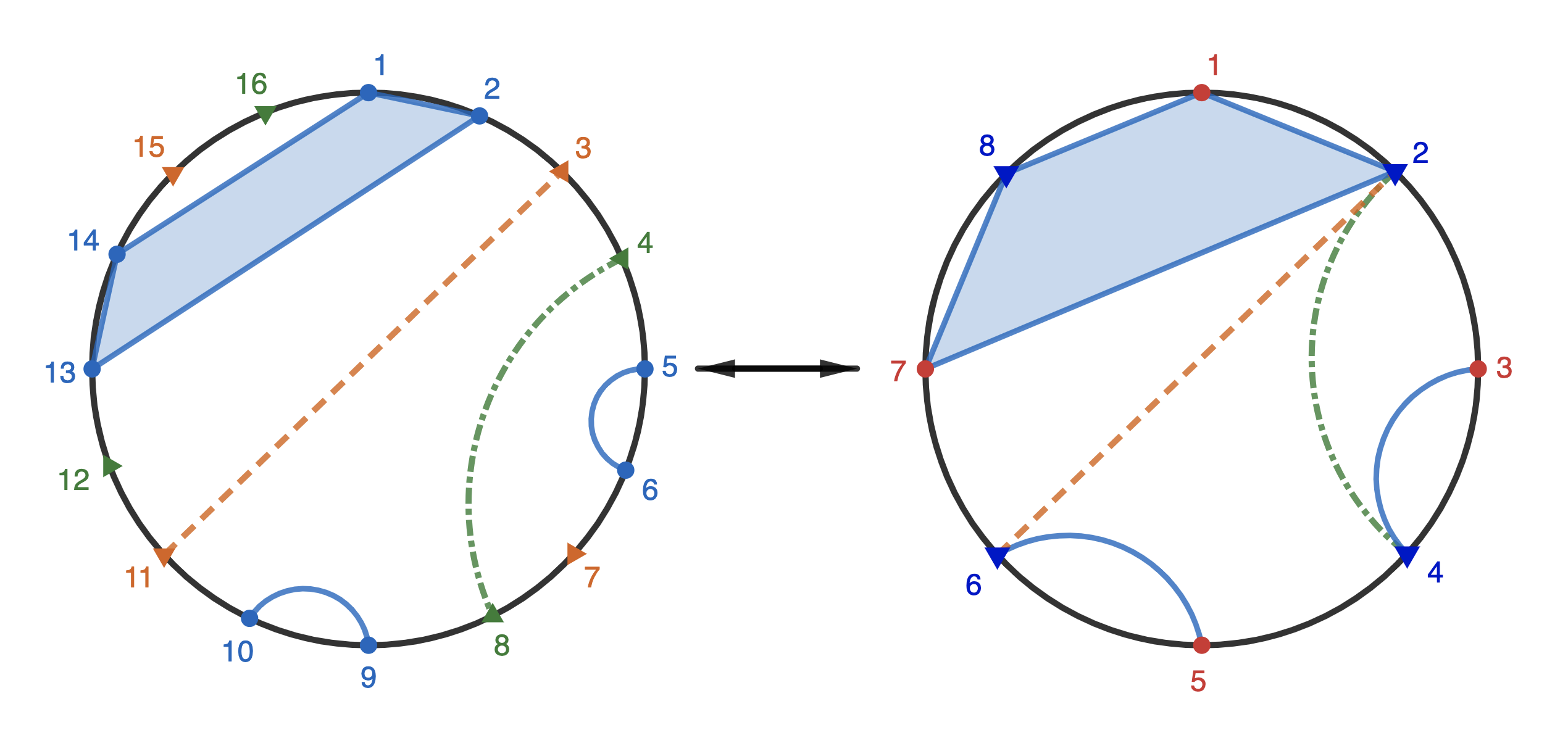}
    \caption{Illustration of \Cref{prop:np-general-touch-nontouch-equiv}: On the left is a partition in $\np_m\(\vaa_1,\dots,\vaa_s\)$ and on the right is a partition in $\npt_m\(\vaa_1,\dots,\vaa_s\)$. Here $m=4$, $s=3$, $k=4$.}
    \label{fig:np-general-equiv}
\end{figure}

\begin{proof}
We give a bijection between $\npm\(\vaa_1,\dots,\vaa_s\)$ and $\nptm\(\vaa_1,\dots,\vaa_s\)$. The idea is similar to the pictorial explanation for the third condition in \Cref{defn:np-general-vaa-touch}.
\begin{enumerate}
    \item $\npm\(\vaa_1,\dots,\vaa_s\)\to \nptm\(\vaa_1,\dots,\vaa_s\)$: Let $\P_1\bigsqcup \dots \bigsqcup\P_t\in \npm\(\vaa_1,\dots,\vaa_s\)$. By definition, they locate on the cycle of length $mk$. We will shrink the cycle length to $(m-s+1)k$ by grouping together the last $s$ points for every $m$ points. More precisely,
    \begin{enumerate}[i.]
        \item For each $P\in\P_1$, let $P' = \big\{(m-s+1)i+j: mi+j\in P \text{ for some } i\in\{0,1,\dots,k-1\}, j\in [m]\big\}$. Let $\displaystyle \P_1' = \underset{P\in\P_1}{\bigcup}\, P'$.
        \item For each $P_t\in\P_t$ where $t>1$, let $P_t' = \big\{(m-s+1)i: m(i-1)+j\in P \text{ for some } i\in[k], j\in [m]\big\}$. For each $t>1$, let $\displaystyle \P_t' = \underset{P_t\in\P_t}{\bigcup}\, P_t'$.
        \item We can see that $\P_1'\cup\dots\cup\P_t'\in \nptm\(\vaa_1,\dots,\vaa_s\)$.
    \end{enumerate}
    
    We can visualize the above procedure as follows: for each $i\in[k]$, take the consecutive $s$ points $m(i-1)+(m-s)+1, m(i-1)+(m-s)+2,\dots, mi$ and push all of them counterclockwise to the first point $m(i-1)+(m-s)+1$. This results in gluing together consecutive vertices of the polygons belonging to $\P_1, \P_2,\dots,\P_s$ all to one point, which is still allowed under $\nptm\(\vaa_1,\dots,\vaa_s\)$. See See \Cref{fig:np-general-equiv} for illustration.
    
    \item $\nptm\(\vaa_1,\dots,\vaa_s\)\to \npm\(\vaa_1,\dots,\vaa_s\)$: Let $\displaystyle \P_1\cup \dots \cup \P_s\in \npt_m\(\vaa_1,\dots,\vaa_s\)$. By definition, they locate on the cycle of length $(m-s+1)k$. We can expand the cycle length to $mk$ and separate the touching points in the expanded cycle to get a nontouching version of the noncrossing partition.
    \begin{enumerate}[i.]
        \item For each $P\in\P_1$, let $P' = \big\{mi+j: (m-s+1)i+j\in P \text{ for some } i\in\{0,1,\dots,k-1\}, j\in [m-s+1]\big\}$. Let $\displaystyle \P_1' = \underset{P\in\P_1}{\bigcup}\, P'$.
        \item For each $P_t\in\P_t$ where $t>1$, let $P_t' = \big\{mi-(s-t): (m-s+1)i\in P \text{ for some } i\in[k]\big\}$. For each $t>1$, let $\displaystyle \P_t' = \underset{P_t\in\P_t}{\bigcup}\, P_t'$.
        \item We can see that $\P_1'\bigsqcup\dots\bigsqcup\P_t'\in \npm\(\vaa_1,\dots,\vaa_s\)$.
    \end{enumerate}
    
    We can visualize the above procedure as follows: take each point $(m-s+1)i$ where the vertices of the polygons of different groups $\P_t$ can possibly touch, and expand this point into $s$ points in the clockwise direction, where polygon belonging to $\P_t$ takes the $t^{th}$ expanded point as the new vertex. See \Cref{fig:np-general-equiv} for an illustration.
    
    \item It is not hard to see that the above procedures are inverses of each other.
\end{enumerate}
\end{proof}

\begin{defn}
    Let $\vaa_1,\dots,\vaa_s\in P_k$ and $m\geq s$. Let $\X=\lcurb x_1,\dots,x_p\rcurb$ be such that $x_1+\dots+x_p=(m-s+1)k$. we define 
    \begin{enumerate}
        \item $\nptm\(\vaa_1,\dots,\vaa_s,\X\)$, the set of noncrossing partitions corresponding to $(m-s+1)\vaa_1,\dots,\vaa_s$ and part sizes $\X$, to be 
        \begin{equation*}
            \big\{ \P\in\np_{m}\(\vaa_1,\dots,\vaa_s\): S_{\P}=\{x_1,\dots,x_p\} \big\},
        \end{equation*}
        
        \item $\nptmlx\(\vaa_1,\dots,\vaa_s,\X\)$, the set of noncrossing partitions corresponding to $\vaa_1,\dots,\vaa_s$ and $\X$ as part labels, to be
        \begin{equation*}
            \big\{ (\P,L):\P\in\npm\(\vaa_1,\dots,\vaa_s\), S_{\P}=\X, L \text{ is a label of }\C_k/\P \text{ corresponding to } \X \big\}.
        \end{equation*}
        
        \item $\nptmlax\(\vaa_1,\dots,\vaa_s,\X\)$ to be the set of partitions corresponding to $\vaa_1,\dots,\vaa_s$ and labels $\X$, but where all the parts are distinct.
    \end{enumerate}
\end{defn}

\begin{lemma}\label{lem:np-general-X}
Let $\vaa_1,\dots,\vaa_s\in P_k$ and $m\geq s$. Let $p=mk-(a_1+\dots+a_s)+1$. Let $\X=\lcurb x_1,\dots,x_p\rcurb$ be such that $x_1+\dots+x_p=(m-s+1)k$. Then
\begin{equation}
    \nptmlax\(\vaa_1,\dots,\vaa_s, \X\) = \(m-s+1\)k^s\(p-1\)!(a_1-1)!\dots(a_s-1)!
\end{equation}
\end{lemma}

\begin{proof}
    We will prove this result using the method from \Cref{subsection:moments-general}.
    
    By the same argument from \Cref{subsection:moments-base-2}, we can reduce the number of polygons from $\vaa_2,\dots,\vaa_s$ to one for each $\vaa_j$, by combining polygons of size $x$ and $y$ into a single polygon of size $x+y-1$. i.e.
    \begin{align*}
        \Big|\nptmlax\(\vaa_1,\dots,\vaa_s, \X\)\Big| = \abs{\nptmlax\(\vaa_1, {\o_2^{(2)}}^{\b_2}\o_{2t_2}^{(2)}, \dots, {\o_2^{(s)}}^{\b_s}\o_{2t_s}^{(s)} \)}
    \end{align*}
    where $\displaystyle t_j = \(\sum_{i=2}^k i\cdot\a_{ji}\)-\(\a_{j2}+\dots+\a_{jk}-1\) = \(\sum_{i=1}^k i\cdot\a_{ji}\)-\(\a_{j1}+\a_{j2}+\dots+\a_{jk}\) + 1 =  k-a_j+1$ and $\b_j = k-t_j = a_j-1$ for each $j=2,\dots,s$.
    
    Now we can further combine the $(s-1)$ polygons $\o_{2t_j}^{(j)}$ of different types for $j=2,\dots,s$ (which can possibly touch each other) into one single polygon, using the method from \Cref{lem:num-of-partitions-polygon-2}. Since we do not combine the type-1 polygons, we keep regarding them as distinct ($\L_{\A}$), whereas we view polygons of types $j=2,\dots,s$ as not distinct. We give the special notation $\vaa_1\(\L_{\A,\Z}\)$ for this. i.e.  
    \begin{align*}
        \abs{\nptmlx\(\vaa_1 (\L_{\A,\X}), {\o_2^{(2)}}^{\b_2}\o_{2t_2}^{(2)}, \dots, {\o_2^{(s)}}^{\b_s}\o_{2t_s}^{(s)}, \X \)} = k^{s-2}\cdot\abs{\nptmlx\(\vaa_1 (\L_{\A,\X}), {\o_2'}^{\c}{\o'}_{2t}, \X\)} 
    \end{align*}
    where $t=t_2+\dots+t_s-(s-2)$ and $\c = k-t$. Here we obtain a factor of $k^{s-2}$ because each time we combine $\o_{2t'_j}^{(j)}$ and $\o_{2t_{j-1}}^{(j-1)}$, we get a factor of $k$ and we do the combination step $s-2$ times for $s-1$ polygons. Thus
    \begin{align*}
        &\abs{\nptmlax\(\vaa_1, {\o_2^{(2)}}^{\b_2}\o_{2t_2}^{(2)}, \dots, {\o_2^{(s)}}^{\b_s}\o_{2t_s}^{(s)}, \X \)} \\
        =& k^{s-2}\cdot\abs{\nptmlx\(\vaa_1 (\L_{\A,\X}), {\o_2'}^{\c}{\o'}_{2t}, \X\)}\cdot (a_2-1)!\dots(a_s-1)!
    \end{align*}
    since we can permute the dots ${\o_2^{(j)}}^{\a_j-1}$ for each $j=2,\dots,s$.
    
    Now we will count $\nptmlax\(\vaa_1, {\o_2'}^{\c}{\o'}_{2t}, \X\)$ by counting $\nplx\(\vbb, {\o_2'}^{\c'}{\o'}_{2t}, \X\)$ where if $\displaystyle \voo^{\vbb} = \prod_{j=1}^k \o_{2(m-s+1)j}^{\a_j}$ and $\c' = (m-s+1)k-t$. i.e. We are considering the number of noncrossing partitions corresponding to augmented $\vaa_1$ (by a factor of $(m-s+1)$) and a polygon of size $t$ of different type on the cycle of length $(m-s+1)k$. Notice that for $\nptmlx\(\vaa_1, {\o_2'}^{\c}{\o'}_{2t}, \X\)$, ${\o_2'}^{\c}{\o'}_{2t}$ is only allowed to use vertices from $\{(m-s+1)i:i\in[k]\}$, but in $\nplx\(\vbb, {\o_2'}^{\c'}{\o'}_{2t}, \X\)$, ${\o_2'}^{\c'}{\o'}_{2t}$ is allowed to use all vertices from $[(m-s+1)k]$. Thus for each configuration of polygons in $\nplx\(\vaa_1, {\o_2'}^{\c}{\o'}_{2t}, \X\)$, we can rotate it in the clockwise direction $1, 2,\dots m-s$ units to get an extra $(m-s)$ configurations in $\nplx\(\vbb, {\o_2'}^{\c'}{\o'}_{2t}, \X\)$. Thus
    \begin{equation*}
        \abs{\nptmlx\(\vaa_1 (\L_{\A,\X}), {\o_2'}^{\c}{\o'}_{2t}, \X\)} = \dfrac{1}{m-s+1}\cdot \abs{\nplx\(\vbb (\L_{\A,\X}), {\o_2'}^{\c'}{\o'}_{2t}, \X\)}.
    \end{equation*}
    
    Further reducing the polygons in $\vaa_1$, we get
    \begin{equation*}
        \abs{\nplx\(\vbb (\L_{\A,\X}), {\o_2'}^{\c'}{\o'}_{2t}, \X\)} = \abs{\nplx\(\o_2^{a_1-1}\o_{2r} (\L_{\A,\X}), {\o_2'}^{\c'}{\o'}_{2t}, \X\)} 
    \end{equation*}
    where $t=(m-s+1)k-a_1+1$.
    
    By the counting method in \Cref{subsection:moments-general}, we know that 
    \begin{equation*}
        \abs{\nplx\(\o_2^{a_1-1}\o_{2r}, {\o_2'}^{\c'}{\o'}_{2t}, \X\)} = \((m-s+1)k\)^2(p-1)!
    \end{equation*}
    
    Thus 
    \begin{align*}
        &\abs{\nplx\(\vbb (\L_{\A,\X}), {\o_2'}^{\c'}{\o'}_{2t}, \X\)} = \((m-s+1)k\)^2(p-1)!(a_1-1)!\\
        \implies & \abs{\nptmlx\(\vaa_1 (\L_{\A,\X}), {\o_2'}^{\c}{\o'}_{2t}, \X\)} = (m-s+1)k^2(p-1)!(a_1-1)! \\
    \end{align*}
    and 
    \begin{align*}
        \Big|\nptmlax\(\vaa_1,\dots,\vaa_s, \X\)\Big|
        &=\abs{\nptmlax\(\vaa_1, {\o_2^{(2)}}^{\b_2}\o_{2t_2}^{(2)}, \dots, {\o_2^{(s)}}^{\b_s}\o_{2t_s}^{(s)}, \X \)} \\
        &= k^{s-2}\cdot\abs{\nptmlx\(\vaa_1 (\L_{\A,\X}), {\o_2'}^{\c}{\o'}_{2t}, \X\)}\cdot (a_2-1)!\dots(a_s-1)!\\
        &= k^{s-2}\cdot\((m-s+1)k^2(p-1)!(a_1-1)!\)\cdot (a_2-1)!\dots(a_s-1)!\\
        &= (m-s+1)k^s(p-1)!(a_1-1)!\dots(a_s-1)!
    \end{align*}
    as needed.
\end{proof}

\begin{proof}[Proof of \Cref{thm:np-general}]

By \Cref{lem:np-general-X},
\begin{align*}
    \npmla\(\vaa_1,\dots,\vaa_s\)
    &= \sum_{\substack{x_i\geq 1:\\ x_1+\dots+x_p = (m-s+1)k}}\, \dfrac{\npmla\(\vaa_1,\dots,\vaa_s,\X\)}{\perm\(x_1,\dots,x_p\)}\\
    &= \sum_{\substack{x_i\geq 1:\\ x_1+\dots+x_p = (m-s+1)k}}\, \dfrac{\npmlax\(\vaa_1,\dots,\vaa_s,\X\)}{p!}\\
    &= \sum_{\substack{x_i\geq 1:\\ x_1+\dots+x_p = (m-s+1)k}}\, (m-s+1)k^s\cdot(p-1)!(a_1-1)!\dots(a_s-1)!\cdot \dfrac{\,1\,}{p!}\\
    &= \binom{(m-s+1)k-1}{p-1}\cdot (m-s+1)k^s\cdot (a_1-1)!\dots(a_s-1)!\cdot \dfrac{\,1\,}{p}\\
    &= \binom{(m-s+1)k-1}{p-1}\cdot \dfrac{(m-s+1)k}{p}\cdot k^{s-1}\cdot (a_1-1)!\dots(a_s-1)!\\
    &= \binom{(m-s+1)k}{p}\cdot k^{s-1}\cdot (a_1-1)!\dots(a_s-1)!\\
    &= \binom{(m-s+1)k}{mk-(a_1+\dots+a_s)+1}\cdot k^{s-1}\cdot (a_1-1)!\dots(a_s-1)!
\end{align*}

\end{proof}

We now extend the definition further to $\np_m\(\vaa_1,\dots,\vaa_s, \o_{2c}, \X\)$.

\begin{defn}\label{defn:np-general-c-gon}
     Let $\vaa_1,\dots,\vaa_s\in P_k$. Let $V=\big\{(m-s+1)i:i\in[k]\big\}$. We say $\P=\P_1\cup\dots \cup\P_s \cup\hat{\P}\in\nptm\(\vaa_1,\dots,\vaa_s, \o_{2c}\)$ if
    \begin{enumerate}
        \item $\P \in \npt_{m+1}\(\vaa_1,\dots,\vaa_s, \vaa_{s+1}\)$ where $\vaa_{s+1}$ corresponds to $\o_2^{k-c}\o_{2c}$.
        \item $\hat{\P}$ has to contain $(m-s+1)k$ (i.e. the last element of the cycle).
    \end{enumerate}
\end{defn}


\begin{defn}
    We define $\npt_{k,m,s}(\o_{2c}) = \big\{ \P=\P_1\cup\dots \cup\P_s \cup\hat{\P}\in\nptm\(\vaa_1,\dots,\vaa_s, \o_{2c}\): \vaa_i\in \P_k\big\}$.
\end{defn}

\begin{thm}\label{thm:np-general-c-gon}
    
    \begin{equation}
        \Big| \nptm\(\vaa_1,\dots,\vaa_s,\o_{2c}\) \Big| = c\cdot \binom{(m-s+1)k}{mk - (a_1+\dots+a_s) + c}\cdot k^{s-1}\cdot \prod_{i=1}^s \dfrac{(a_i-1)!}{\a_{i1}!\dots\a_{ik}!}.
    \end{equation}
\end{thm}

\begin{proof}

We first view the $c$-gon as adding an extra $\voo^{\vbb} = \o_2^{k-c}\o_{2c}$, and then factor out the double counting.

Let $m'=m+1$, $s'=s+1$ and $b=k-c+1$. By \Cref{thm:np-general}, 
\begin{align*}
    \Big| \npt_{m'}\(\vaa_1,\dots,\vaa_s,\vbb\) \Big|
    &= \binom{(m'-s'+1)k}{m'k - (a_1+\dots+a_s+b) + 1}\cdot k^{s'-1}\cdot \(\prod_{i=1}^s \dfrac{(a_i-1)!}{\a_{i1}!\dots\a_{ik}!}\)\cdot \dfrac{(b-1)!}{(k-c)!}\\
    &= \binom{(m-s+1)k}{(m+1)k - (a_1+\dots+a_s) -(k-c+1) + 1}\cdot k^{s}\cdot \prod_{i=1}^s \dfrac{(a_i-1)!}{\a_{i1}!\dots\a_{ik}!}\\
    &= \binom{(m-s+1)k}{mk - (a_1+\dots+a_s) +c}\cdot k^{s}\cdot \prod_{i=1}^s \dfrac{(a_i-1)!}{\a_{i1}!\dots\a_{ik}!}.
\end{align*}

Note that the original $c$-gon has to contain $(m-s+1)k$, while in the counting with the term $\voo^{\vbb}$, there is no such requirement on the $c$-gon. This gives a factor of $\dfrac{\,c\,}{k}$. Thus 
\begin{align*}
    \Big| \nptm\(\vaa_1,\dots,\vaa_s,\o_{2c}\) \Big| 
    &= \dfrac{\,c\,}{k}\cdot \Big| \npt_{m'}\(\vaa_1,\dots,\vaa_s,\vbb\) \Big| \\
    &= \dfrac{\,c\,}{k}\cdot \binom{(m-s+1)k}{mk - (a_1+\dots+a_s) +c}\cdot k^{s}\cdot \prod_{i=1}^s \dfrac{(a_i-1)!}{\a_{i1}!\dots\a_{ik}!} \\
    &= c\cdot \binom{(m-s+1)k}{mk - (a_1+\dots+a_s) + c}\cdot k^{s-1}\cdot \prod_{i=1}^s \dfrac{(a_i-1)!}{\a_{i1}!\dots\a_{ik}!}
\end{align*}
as needed
\end{proof}

\subsection{Formula for Moments of \texorpdfstring{$\ozm\opr \o^{(1)}\opr \dots, \opr\o^{(s)}$}{Ωz(m)ºR Ω1 ºR … ºR Ωs}}
\label{subsection:ozm-s-pf}

\setlength{\parskip}{1.5mm}
\setlength{\baselineskip}{1.3em}

First we are going to give another equivalent notation for the expression $\(\ozm\opr \o^{(1)}\opr\dots \opr \o^{(s)}\)_{2k}$ defined inductively based on \Cref{thm:main-1-circ}.

\begin{defn}\label{defn:Amsk}
We define $A_{m}^{(s)}\(k,\underbrace{0,\dots,0}_{s}\)$ associated with $\o^{(1)},\dots,\o^{(s)}$ inductively as the following:
\begin{enumerate}
    \item $\displaystyle A_{m}^{(0)}(k) = C(k,m)$.
    \item $\displaystyle A_{m}^{(1)}(k,0) = \sum_{\vaa,\vbb\in P_k} C\(\vaa,\vbb\)\cdot \(\prod_{i=1}^k C(i,m)^{\,\a_i}\) \cdot \vv{\o^{(1)}}^{\vbb}$.
    \item Let $\o'$ be the distribution where $\o'_{2i} = A_m^{(s-1)}\(i,\underbrace{0,\dots,0}_{s-1}\)$ for all $i$. Define 
    \begin{equation}
    \begin{aligned}
        A_{m}^{(s)}(k,0,\dots,0) 
        &= \sum_{\vaa,\vbb\in P_k} C\(\vaa,\vbb\)\cdot \vv{\o'}^{\vaa} \cdot \vv{\o^{(s)}}^{\vbb}\\
        &= \sum_{\vaa,\vbb\in P_k} C\(\vaa,\vbb\)\cdot \(\prod_{i=1}^k \(A_m^{(s-1)}(i,0,\dots,0)\)^{\,\a_i}\) \cdot \vv{\o^{(s)}}^{\vbb}
    \end{aligned}
    \end{equation}
    where 
    \begin{equation*}
        C\(\vaa,\vbb\) = (-1)^{a+b-k-1}\cdot k\cdot \binom{a+b-2}{k-1}\cdot \dfrac{(a-1)!}{\a_1!\dots\a_k!}\cdot \dfrac{(b-1)!}{\b_1!\dots\b_k!}
    \end{equation*}
\end{enumerate}
\end{defn}

\begin{rmk}
    By \Cref{thm:main-1-circ} and \Cref{cor:moments-zmshape}, $A_m^{(s)}\(k,0,\dots,0\) = \(\ozm\opr \o^{(1)}\opr\dots \opr \o^{(s)}\)_{2k} $ for all $k$.
\end{rmk}

Thus to prove \Cref{thm:moments-mzshape-s-distr}, it is equivalent to prove the following theorem.
\begin{restatable}{thm}{thmAmsFormula}
\label{thm:expression-for-mzshape-sdistr}
\begin{equation}
\begin{aligned}
    A_{m}^{(s)}\(k,\underbrace{0,\dots,0}_{s}\) 
    &= \sum_{\vaa_i\in P_k}\, C_m(\vaa_1,\dots,\vaa_s)\cdot \prod_{j=1}^s\, \vv{\o^{(i)}}^{\vaa_i} \\
    &= \sum_{\vaa_i\in P_k}\, \binom{(m-s+1)k}{mk - (a_1+\dots+a_s) + 1}\cdot k^{s-1}\cdot \(\prod_{j=1}^s \dfrac{(a_i-1)!}{\a_{i1}!\dots\a_{ik}!}\cdot \vv{\o^{(i)}}^{\vaa_i}\)
\end{aligned}
\end{equation}
where $\vaa_i=\(\a_{i1},\dots,\a_{ik}\)$, $\displaystyle a_i = \sum_{j=1} \a_{ij}$ for all $i=1,\dots,s$.
\end{restatable}

\begin{defn}
    Define $B_m^{(s)}(k)$ to be
    \begin{equation}
        B_m^{(s)}(k) = \sum_{\vaa_i\in P_k}\, C_m(\vaa_1,\dots,\vaa_s)\cdot \prod_{j=1}^s\, \vv{\o^{(i)}}^{\vaa_i}.
    \end{equation}
\end{defn}

\begin{defn}
    Given $\P=\{P_1,\dots, P_a\}\in \np_k$ and a distribution $\o$, we denote $\voo_{\P}$ to be
    \begin{equation}
        \voo_{\P} = \o_{2\,\abs{P_1}}\cdot \dots\cdot \o_{2\,\abs{P_a}} 
    \end{equation}
    and if the size of every $P_i$'s is a multiple of $m$, we denote $\voo_{\P/m}$ to be
    \begin{equation}
        \voo_{\P/m} = \o_{2\,\abs{P_1}/m}\cdot \dots\cdot \o_{2\,\abs{P_a}/m} 
    \end{equation}
\end{defn}

\begin{prop}\label{prop:Bms-alternative}
Assume $m\geq s$. Let $\mh = m-s+1$. Then 
    
    \begin{equation}
        \bms(k) = \sum_{\substack{\P_1\sqcup\dots \sqcup\P_s \\ \in\, \np_{k,m,s}}}\, \vv{\o^{(1)}}_{\P_1/\mh} \cdot \vv{\o^{(2)}}_{\P_2}\cdot \dots\cdot \vv{\o^{(s)}}_{\P_s} 
    \end{equation}
    where $\P_i=\big\{P_{i,1},\dots,P_{i,a_i}\big\}$ for each $i\in[s]$.
    
    Moreover, we can replace $\np_{k,m,s}$ in the summand with $\npt_{k,m,s}$.
\end{prop}

\begin{proof}
    By \Cref{thm:num-of-np}, $\displaystyle C_m\(\vaa_1,\dots,\vaa_s\) = \big|\np_m\(\vaa_1,\dots,\vaa_s\)\big|$. In particular, the size of any partition set corresponding to $\vaa_1$ is a multiple of $\mh=m-s+1$. Thus
    \begin{align*}
        \bms(k) 
        &= \sum_{\vaa_i\in P_k} C_m(\vaa_1,\dots,\vaa_s)\cdot \vv{\o^{(1)}}^{\vaa_1}\cdot\dots\cdot \vv{\o^{(s)}}^{\vaa_s} \\
        &= \sum_{\vaa_i\in P_k} \big|\np_m\(\vaa_1,\dots,\vaa_s\)\big|\cdot \vv{\o^{(1)}}^{\vaa_1}\cdot\dots\cdot \vv{\o^{(s)}}^{\vaa_s}\\
        &= \sum_{\vaa_i\in P_k}\; \sum_{\substack{\P_1\sqcup\dots \sqcup\P_s \in \\ \npm(\vaa_1,\dots,\vaa_s)}}\, \({\o^{(1)}_2}^{\a_{11}}\dots{\o^{(1)}}_{2k}^{\a_{1k}}\)\cdot\dots\cdot \({\o^{(s)}_2}^{\a_{s1}}\dots{\o^{(s)}}_{2k}^{\a_{sk}}\) \,\\
        &= \sum_{\vaa_i\in P_k}\; \sum_{\substack{\P_1\sqcup\dots \sqcup\P_s \in \\ \npm(\vaa_1,\dots,\vaa_s)}}\, \(\o^{(1)}_{2\,|P_{1,1}|/\mh}\cdot\dots\cdot \o^{(1)}_{2\,|P_{1,a_1}|/\mh}\)\cdot \prod_{j=2}^{s} \(\o^{(j)}_{2\,|P_{j,1}|}\cdot\dots\cdot \o^{(j)}_{2\,|P_{j,a_j}|}\) \\
        &= \sum_{\substack{\P_1\sqcup\dots \sqcup\P_s \\ \in\, \np_{k,m,s}}}\, \(\o^{(1)}_{2\,|P_{1,1}|/\mh}\cdot\dots\cdot \o^{(1)}_{2\,|P_{1,a_1}|/\mh}\)\cdot \prod_{j=2}^{s} \(\o^{(j)}_{2\,|P_{j,1}|}\cdot\dots\cdot \o^{(j)}_{2\,|P_{j,a_j}|}\)\, \\
        &= \sum_{\substack{\P_1\sqcup\dots \sqcup\P_s \\ \in\, \np_{k,m,s}}}\, \vv{\o^{(1)}}_{\P_1/\mh} \cdot \vv{\o^{(2)}}_{\P_2}\cdot \dots\cdot \vv{\o^{(s)}}_{\P_s} \,.
    \end{align*}
\end{proof}

\begin{lemma}\label{lem:sum-prod-general-np}
    \begin{equation}
    \begin{aligned}
         &\sum_{\substack{ k_i\geq 1: \\ k_1+\cdots+k_c = k}} B_{m}^{(s)}(k_1)\dots B_{m}^{(s)}(k_c) \\
        =& \,c\cdot \sum_{\a_i\in P_k} \binom{(m-s+1)k}{mk-(a_1+\dots+a_s)+c}\cdot k^{s-1}\cdot \(\prod_{i=1}^s \dfrac{(a_i-1)!}{\a_{i1}!\dots\a_{ik}!}\cdot \vv{\o^{(i)}}^{\vaa_i}\)
    \end{aligned}
    \end{equation}
\end{lemma}

\begin{proof}
    By \Cref{thm:np-general-c-gon}, $\displaystyle \Big| \nptm\(\vaa_1,\dots,\vaa_s,\o_{2c}\) \Big| = c\cdot \binom{(m-s+1)k}{mk - (a_1+\dots+a_s) + c}\cdot k^{s-1}\cdot \prod_{i=1}^s \dfrac{(a_i-1)!}{\a_{i1}!\dots\a_{ik}!}$\,. Thus we can rewrite the RHS as  $\displaystyle \sum_{\a_i\in P_k} \Big| \nptm\(\vaa_1,\dots,\vaa_s,\o_{2c}\) \Big| \cdot \prod_{i=1}^s\, \vv{\o^{(i)}}^{\vaa_i}$\,.
    
    By the same argument as \Cref{prop:Bms-alternative},
    \begin{align}
        RHS 
        &= \sum_{\a_i\in P_k} \Big| \nptm\(\vaa_1,\dots,\vaa_s,\o_{2c}\) \Big| \cdot \prod_{i=1}^s\, \vv{\o^{(i)}}^{\vaa_i} \\
        &= \sum_{\substack{\P_1\cup\dots\cup \P_s\cup \hat{\P} \\ \in\, \npt_{k,m,s}(\o_{2c})}}\, \vv{\o^{(1)}}_{\P_1/\mh}\cdot \vv{\o^{(2)}}_{\P_2}\cdot \dots\cdot \vv{\o^{(s)}}_{\P_s} \,.
    \end{align}
    where $\mh = m-s+1$.
    
    We are going to prove that there is a bijection between $\displaystyle \bigsqcup_{\substack{ k_i\geq 1: \\ k_1+\dots+k_c=k}}\, \npt_{k_1,m,s}\times \dots\times \npt_{k_c,m,s}$ and $\npt_{k,m,s}(\o_{2c})$.
    \begin{enumerate}
        \item $\displaystyle \Phi: \npt_{k,m,s}(\o_{2c}) \to  \bigsqcup_{\substack{ k_i\geq 1: \\ k_1+\dots+k_c=k}}\, \npt_{k_1,m,s}\times \dots\times \npt_{k_c,m,s}$: Let $\P_1\cup\dots \cup \P_s\cup \hat{\P} \in \npt_{k,m,s}(\o_{2c})$. $\hat{\P}$ contains a polygon $P_c$ of size $c$. Moreover, since $\mh k\in P_c$ by definition, we can assume $P_c = \{\mh i_1,\dots, \mh i_c\}$ for some $1\leq i_1< \dots< i_c=k$. Let $i_0=0$. Since $\P_1\cup\dots \cup \P_s\cup \hat{\P}$ is noncrossing, for each $j\in[c]$, everything in between $\mh i_{j-1}$ (exclusive) and $\mh i_{j}$ (inclusive) can be treated as a collection of noncrossing partitions in $\npt_{k_j,m,s}$ where $k_j = i_j-i_{j-1}$. More precisely, for each $j\in[c]$, let $V_j=\big\{\mh i_{j-1}+1, \mh i_{j-1}+2, \dots \mh i_{j} \big\}$. Then we can divide $\P_1,\dots,\P_s$ into $\P_{ij}$ for $i\in[s]$, $j\in[c]$ such that
        \begin{enumerate}[i.]
            \item $\P_{i1}\cup\dots\cup \P_{ic} = \P_{i}$, and
            \item for each $j\in[c]$, $\P_{1j}\cup\dots \cup \P_{sj}$ only uses vertices in $V_j$. i.e. $\P_{1j}\cup\dots \cup \P_{sj}\in \npt_{k_j,m,s}$ where $k_j = i_j-i_{j-1}$. See \Cref{fig:sum-prod-general-np} for an illustration.
        \end{enumerate}
        
        Thus the resulting $\big(\P_{11}\cup\dots \cup \P_{s1}, \dots, \P_{1c}\cup\dots \cup \P_{sc} \big)\in \npt_{k_1,m,s}\times \dots\times \npt_{k_c,m,s}$.
        
    \begin{figure}[hbt!]
        \centering
        \includegraphics[scale=0.32]{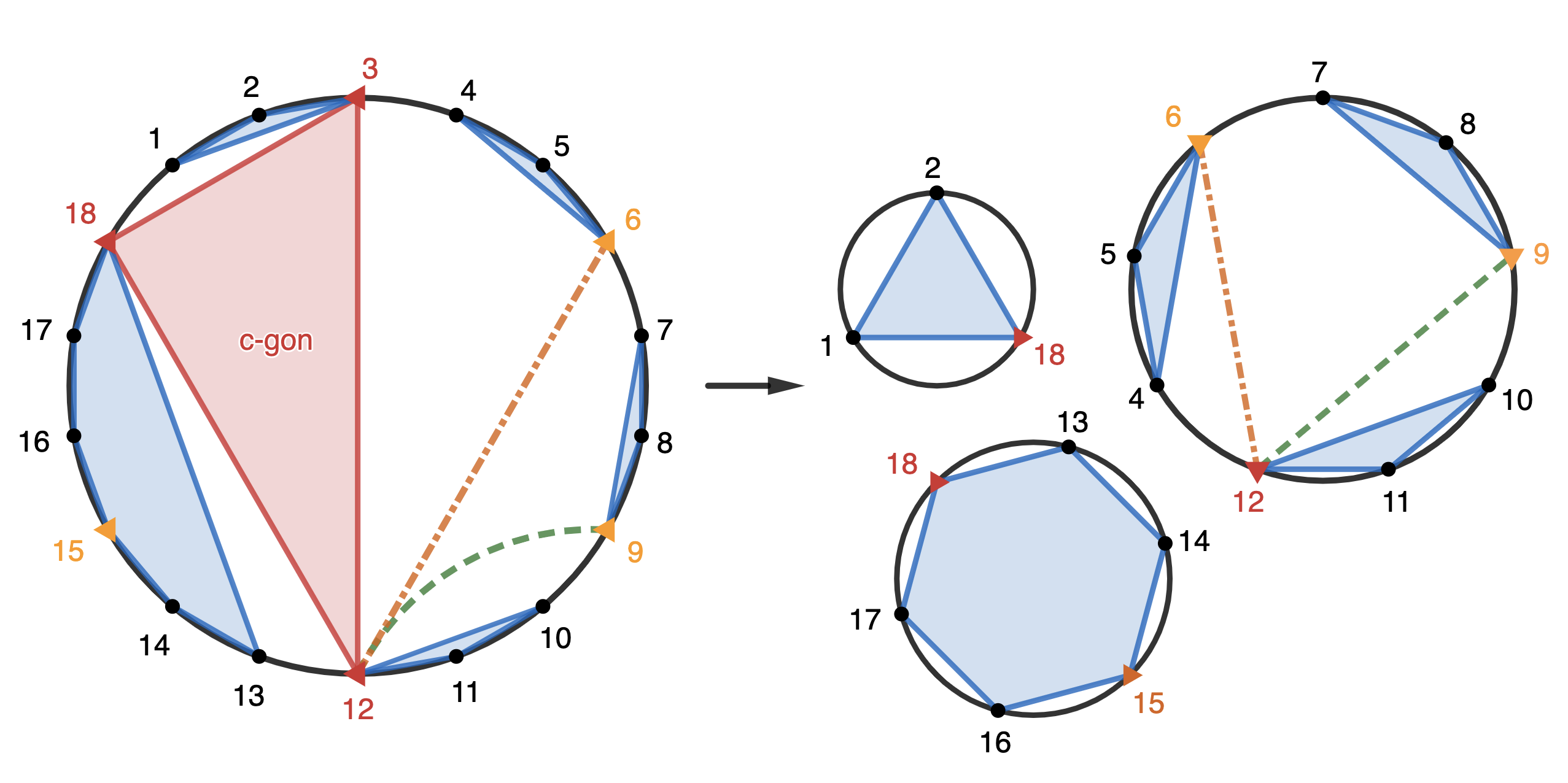}
        \caption{Illustration of \Cref{lem:sum-prod-general-np}: here $m=5, s=3, k=6$, $c=3$ and $\mh = m-s+1 = 3$. For this configuration $k_1=1$, $k_2=3$ and $k_3=2$.}
        \label{fig:sum-prod-general-np}
    \end{figure}
        
        \item $\displaystyle \Psi: \bigsqcup_{\substack{ k_i\geq 1: \\ k_1+\dots+k_c=k}}\, \npt_{k_1,m,s}\times \dots\times \npt_{k_c,m,s}\to \npt_{k,m,s}(\o_{2c})$: Let $\big(\P_{11}\cup\dots \cup \P_{s1}, \dots, \P_{1c}\cup\dots \cup \P_{sc} \big)\in \npt_{k_1,m,s}\times \dots\times \npt_{k_c,m,s}$ for some $k_i\geq 1$ where $k_1+\dots+k_c=k$. We will ``glue" them together in the following way:
        \begin{enumerate}[i.]
            \item For each $i\in[s], j\in[c]$, let $\P_{ij}' = \P_{ij}+\mh (k_1+\dots+k_{j-1}) := \big\{\{v+\mh (k_1+\dots+k_{j-1}): v\in P \}: P\in \P_{ij} \big\}$. i.e. we shift the indices of the polygons in each $\P_{ij}$ up by $\mh (k_1+\dots+k_{j-1})$.
            \item For each $i\in[s]$, let $\P_i = \P_{i1}\cup\dots \cup \P_{ic}$. 
            \item Let $P_c = \big\{\mh(k_1+\dots+k_j): j\in[c] \big\}$ and $\hat{P} = \big\{\{v\}, P_c: v\in [\mh k]\setminus P_c \big\}$. 
        \end{enumerate}
        
        We can see that $\P_1\cup\dots\cup \P_s\cup \hat{\P} \in \npt_{k,m,s}(\o_{2c})$.
        
        \item It is easy to see that $\Phi$ and $\Psi$ are inverses of each other.
    \end{enumerate}
    
    Since the bijection only shifts indices of the polygons, number of indices in the polygons stay the same. i.e. for each $\P_1\cup\dots\cup \P_s\cup \hat{\P}\in \npt_{k,m,s}(\o_{2c})$ where $\P_{i}=\{P_{i,1},\dots, P_{i,a_i}\}$, $\Phi\(\P_1\cup\dots\cup \P_s\cup \hat{\P}\) = \big(\P_{11}\cup\dots \cup \P_{s1}, \dots, \P_{1c}\cup\dots \cup \P_{sc} \big)$ for some $\P_{ij}$ where $\P_{ij}=\{P_{ij,1},\dots,P_{ij,a_{ij}}\}$. Then
    \begin{align*}
         &\(\o^{(1)}_{2\,|P_{1,1}|/\mh}\cdot\dots\cdot \o^{(1)}_{2\,|P_{1,a_1}|/\mh}\)\cdot \prod_{j=2}^{s} \(\o^{(j)}_{2\,|P_{j,1}|}\cdot\dots\cdot \o^{(j)}_{2\,|P_{j,a_j}|}\) \\
        =& \prod_{t=1}^c\, \( \(\o^{(1)}_{2\,|P_{1t,1}|/\mh}\cdot\dots\cdot \o^{(1)}_{2\,|P_{1t,a_{1t}}|/\mh}\)\cdot \prod_{j=2}^{s} \(\o^{(j)}_{2\,|P_{jt,1}|}\cdot\dots\cdot \o^{(j)}_{2\,|P_{jt,a_{jt}}|}\) \)
    \end{align*}
    
    Thus combining everything,
    \begin{align*}
         & \sum_{\substack{ k_i\geq 1: \\ k_1+\cdots+k_c = k}} B_{m}^{(s)}(k_1)\dots B_{m}^{(s)}(k_c) \\
        =& \sum_{\substack{ k_i\geq 1: \\ k_1+\cdots+k_c = k}} \prod_{t=1}^c\, \sum_{\substack{\P_{1t}\sqcup\dots \sqcup\P_{st} \\ \in\, \np_{k_t,m,s}}}\, \vv{\o^{(1)}}_{\P_{1t}/\mh}\cdot \vv{\o^{(2)}}_{\P_{2t}}\cdot \dots\cdot \vv{\o^{(s)}}_{\P_{st}} \\
        =& \sum_{\substack{ k_i\geq 1: \\ k_1+\cdots+k_c = k}} \sum_{\substack{\(\P_{11}\cup\dots \cup \P_{s1}, \dots, \P_{1c}\cup\dots \cup \P_{sc}\)\\ \in \npt_{k_1,m,s}\times \dots\times \npt_{k_c,m,s}}} \prod_{t=1}^c\, \vv{\o^{(1)}}_{\P_{1t}/\mh}\cdot \vv{\o^{(2)}}_{\P_{2t}}\cdot \dots\cdot \vv{\o^{(s)}}_{\P_{st}} \\
        =& \sum_{\substack{ k_i\geq 1: \\ k_1+\cdots+k_c = k}} \sum_{\substack{\(\P_{11}\cup\dots \cup \P_{s1}, \dots, \P_{1c}\cup\dots \cup \P_{sc}\)\\ \in \npt_{k_1,m,s}\times \dots\times \npt_{k_c,m,s}}}\, \vv{\o^{(1)}}_{\P_{1}/\mh}\cdot \vv{\o^{(2)}}_{\P_{2}}\cdot \dots\cdot \vv{\o^{(s)}}_{\P_{s}} \\ 
         &\hspace{1.5cm} \text{where } \P_1\cup\dots \cup\P_s\cup \hat{\P} = \Psi\big(\(\P_{11}\cup\dots \cup \P_{s1}, \dots, \P_{1c}\cup\dots \cup \P_{sc}\)\big)
         \text{ and } \P_i=\big\{P_{i,1},\dots, P_{i,a_i}\big\} \\
        =& \sum_{\substack{\P_1\cup\dots \cup\P_s\cup \hat{\P} \\ \in \npt_{k,m,s}(\o_{2c})}}\, \vv{\o^{(1)}}_{\P_{1}/\mh}\cdot \vv{\o^{(2)}}_{\P_{2}}\cdot \dots\cdot \vv{\o^{(s)}}_{\P_{s}} \\ 
        =& \,c\cdot \sum_{\a_i\in P_k} \binom{(m-s+1)k}{mk-(a_1+\dots+a_s)+c}\cdot k^{s-1}\cdot \(\prod_{i=1}^s \dfrac{(a_i-1)!}{\a_{i1}!\dots\a_{ik}!}\cdot \vv{\o^{(i)}}^{\vaa_i}\)
    \end{align*}
    as needed.
\end{proof}
Now we are ready to prove the main result for this section,  \Cref{thm:expression-for-mzshape-sdistr}.

\thmAmsFormula*

\begin{proof}[Proof of \Cref{thm:expression-for-mzshape-sdistr}]
We prove this by induction on $s$.
\begin{enumerate}
    \item Base case $s=0$: 
    \begin{align*}
        &\binom{(m-s+1)k}{mk - (a_1+\dots+a_s) + 1}\cdot k^{s-1}\cdot \prod_{i=1}^s \dfrac{(a_i-1)!}{\a_{i1}!\dots\a_{ik}!} \\
        =& \binom{(m+1)k}{mk+1}\cdot k^{-1}
        = \dfrac{\((m+1)k\)!}{(mk+1)!(k-1)!}\cdot\dfrac{\,1\,}{k}
        = \dfrac{1}{mk+1}\cdot\dfrac{\((m+1)k\)!}{(mk)!k!}\\
        =& \dfrac{1}{mk+1}\cdot \binom{(m+1)k}{k}
        = C(k,m)
    \end{align*}
    
    \item Assume $\displaystyle A_{m}^{(s)}\(k,\underbrace{0,\dots,0}_{s}\) = \sum_{\vaa_i\in P_k}\, C_m\(\vaa_1,\dots,\vaa_s\)\cdot \prod_{i=1}^s\, \vv{\o^{(i)}}^{\vaa_i}$, we will prove $\displaystyle A_m^{(s+1)}\(k,\underbrace{0,\dots,0}_{s}\) = \sum_{\vaa_i\in P_k}\, C_m\(\vaa_1,\dots,\vaa_{s+1}\)\cdot \prod_{i=1}^{s+1}\, \vv{\o^{(i)}}^{\vaa_i}$.
    
    For the LHS, by the definition of $A_{m}^{(s+1)}(k,0,\dots,0)$,
    \begin{align*}
        A_{m}^{(s+1)}(k,0,\dots,0) 
        &= \sum_{\vcc,\vbb\in P_k} C\(\vaa,\vbb\)\cdot \vv{A_m^{(s)}}^{\,\vcc} \cdot \vv{\o^{(s+1)}}^{\,\vbb} \\
        &= \sum_{\vcc,\vbb\in P_k} (-1)^{c+b-k-1}\cdot k \cdot \binom{b+c-2}{k-1}\cdot \dfrac{(c-1)!}{\c_1!\dots\c_k!}\cdot \vv{A_m^{(s)}}^{\,\vcc}\cdot \dfrac{(b-1)!}{\b_1!\dots\b_k!}\cdot \vv{\o^{(s+1)}}^{\vbb}\\
        &=\sum_{\vbb\in P_k} k\cdot \( \sum_{\vcc\in P_k} (-1)^{c+b-k-1}\cdot \binom{b+c-2}{k-1}\cdot \dfrac{(c-1)!}{\c_1!\dots\c_k!}\cdot \vv{A_m^{(s)}}^{\,\vcc}\) \cdot \dfrac{(b-1)!}{\b_1!\dots\b_k!}\cdot \vv{\o^{(s+1)}}^{\vbb}\\
    \end{align*}
    
    On the other hand, for the RHS,
    \begin{align*}
        & \sum_{\vaa_i\in P_k}\, C_m\(\vaa_1,\dots,\vaa_{s+1}\)\cdot \prod_{i=1}^{s+1}\, \vv{\o^{(i)}}^{\vaa_i} \\
        =& \sum_{\vaa_i\in P_k}\, \binom{(m-s)k}{mk - (a_1+\dots+a_{s+1}) + 1}\cdot k^{s}\cdot \(\prod_{i=1}^{s+1} \dfrac{(a_i-1)!}{\a_{i1}!\dots\a_{ik}!}\cdot \vv{\o^{(i)}}^{\vaa_i} \) \\
        =& \sum_{\vbb\in P_k} k\cdot \(\sum_{\vaa_i\in P_k} \binom{(m-s)k}{mk - \(a_1+\dots+a_{s}+b\) + 1}\cdot k^{s-1}\cdot \prod_{i=1}^{s} \dfrac{(a_i-1)!}{\a_{i1}!\dots\a_{ik}!}\cdot \vv{\o^{(i)}}^{\vaa_i} \) \cdot \dfrac{(b-1)!}{\b_1!\dots\b_k!}\cdot \vv{\o^{(s+1)}}^{\vbb}\\
    \end{align*}
    
    Thus it suffices to prove that
    \begin{align*}
         &\sum_{\vcc\in P_k} (-1)^{c+b-k-1}\cdot \binom{b+c-2}{k-1}\cdot \dfrac{(c-1)!}{\c_1!\dots\c_k!}\cdot \vv{A_m^{(s)}}^{\,\vcc}\\
        =& \sum_{\vaa_i\in P_k} \binom{(m-s)k}{mk - \(a_1+\dots+a_{s}+b\) + 1}\cdot k^{s-1}\cdot \(\prod_{i=1}^{s} \dfrac{(a_i-1)!}{\a_{i1}!\dots\a_{ik}!}\cdot \vv{\o^{(i)}}^{\vaa_i}\)
    \end{align*}
    
    Rewriting the LHS and using \Cref{lem:sum-prod-general-np}, we get
    \begin{align*}
        &\sum_{\vcc\in P_k} (-1)^{c+b-k-1}\cdot \binom{b+c-2}{k-1}\cdot \dfrac{(c-1)!}{\c_1!\dots\c_k!}\cdot \vv{A_m^{(s)}}^{\,\vcc}\\
        =& \sum_{\vcc\in P_k} (-1)^{c+b-k-1}\cdot \binom{b+c-2}{k-1}\cdot \dfrac{\,1\,}{c}\cdot  \dfrac{c\,!}{\c_1!\dots\c_k!}\cdot \(\prod_{i=1}^k A_m^{(s)}(i,0,\dots,0)^{\c_i}\) \\
        =& \sum_{c=1}^{k} (-1)^{c+b-k-1}\cdot \binom{b+c-2}{k-1}\cdot \dfrac{\,1\,}{c}\cdot \(\sum_{\substack{ k_i\geq 1: \\ k_1+\cdots+k_c = k}} A_{m}^{(s)}(k_1,0,\dots,0)\dots A_{m}^{(s)}(k_c,0,\dots,0)\) \\
        =& \sum_{c=1}^{k} (-1)^{c+b-k-1}\cdot \binom{b+c-2}{k-1}\cdot \(\sum_{\a_i\in P_k} \binom{(m-s+1)k}{a_1+\dots+a_s-(s-1)k-c}\cdot k^{s-1}\cdot \prod_{i=1}^s \dfrac{(a_i-1)!}{\a_{i1}!\dots\a_{ik}!}\cdot \vv{\o^{(i)}}^{\vaa_i} \) \\
        =& \sum_{\vaa_i\in P_k} \(\sum_{c=1}^{k} (-1)^{c+b-k-1}\cdot \binom{b+c-2}{k-1}\cdot \binom{(m-s+1)k}{a_1+\dots+a_s-(s-1)k-c}\)\cdot k^{s-1}\cdot \prod_{i=1}^{s} \dfrac{(a_i-1)!}{\a_{i1}!\dots\a_{ik}!}\cdot \vv{\o^{(i)}}^{\vaa_i}
    \end{align*}
    
    Therefore it is equivalent to prove that 
    \begin{align}
        \sum_{c=1}^{k} (-1)^{c+b-k-1}\cdot \binom{b+c-2}{k-1}\cdot \binom{(m-s+1)k}{(a_1+\dots+a_s)-(s-1)k-c}
        = \binom{(m-s)k}{mk-\(a_1+\dots+a_{s}+b\)+1} \label{eqn:key identity}
    \end{align}
    
    Recall \Cref{prop:alternating-sum-of-binomial-coeff}
    \begin{equation}
        \sum_{i=0}^{m-k} (-1)^{i+(m-k)}\cdot \binom{m-i}{k}\cdot \binom{n}{i} = \binom{n-k-1}{m-k}\,.
    \end{equation}
    
    Let $i=(a_1+\dots+a_s)-(s-1)k-c$ in \eqref{eqn:key identity}. Notice that for the summands to be nonzero, we require
    \begin{enumerate}[i.]
        \item $b+c-2\geq k-1$ which implies that $c\geq k-b+1$
        \item Since $\displaystyle\binom{(m-s+1)k}{(a_1+\dots+a_s)-(s-1)k-c} = \binom{(m-s+1)k}{mk-(a_1+\dots+a_s)+c}$, $(m-s+1)k \geq mk-(a_1+\dots+a_s)+c$ which implies that $c\leq (a_1+\dots+a_s)-(s-1)k$.
    \end{enumerate}
    
    Thus 
    \begin{align*}
         & \sum_{c=1}^{k}\, (-1)^{c+b-k-1}\cdot \binom{b+c-2}{k-1}\cdot \binom{(m-s+1)k}{(a_1+\dots+a_s)-(s-1)k-c}\\
        =& \sum_{c=k-b+1}^{(a_1+\dots+a_s)-(s-1)k}\, (-1)^{c+b-k-1}\cdot \binom{b+c-2}{k-1}\cdot \binom{(m-s+1)k}{(a_1+\dots+a_s)-(s-1)k-c}\\
        =& \sum_{i=0}^{(a_1+\dots+a_s)+b-sk-1}\, \binom{b+(a_1+\dots+a_s)+(s-1)k-2-i}{k-1} \binom{(m-s+1)k}{i}\\
        =& \binom{(m-s+1)k-(k-1)-1}{\(b+(a_1+\dots+a_s)+(s-1)k-2\)-(k-1)}\\
        =& \binom{(m-s)k}{(a_1+\dots+a_s)+b-sk-1}
        =  \binom{(m-s)k}{mk-(a_1+\dots+a_s+b)+1}
    \end{align*}
    as needed.
\end{enumerate}
    
\end{proof}


\section{Connection with Graph Matrices}
\label{section:graph-matrices}

\setlength{\parskip}{1.5mm}
\setlength{\baselineskip}{1.3em}


\comm{
\begin{thm}
    $\mzshapesdistr$ and $\mzshapesdistrg$ have the same moments. i.e.  Let $\displaystyle r(n,m)=\dfrac{n!}{(n-m)!}$, then
    \begin{equation}
        \lim_{n\to\infty}\, \dfrac{\,1\,}{n}\cdot\Ebb\[\trace\(\(\mzshapesdistr{\mzshapesdistr}^T\)^{k}\)\] = \lim_{n\to\infty}\, \dfrac{\,1\,}{r(n,m)}\cdot\Ebb\[\trace\(\(M^{(G)}_{Z(m),s}{M^{(G)}_{Z(m),s}}^T\)^{k}\)\] 
    \end{equation}
    for all $k\geq 0$.
\end{thm}
}

In this section, we will prove \Cref{thm:main-2-circ}. Let $\mzshapesdistrg$ be defined as in \Cref{defn:mzshape-sdistr-graph-matrix}.
\maintwocirc*

To prove this theorem, we use the trace power method. The following theorem follows from Lemma B.1 and B.3 in \cite{bai2010spectral}.

\begin{thm}\label{thm:trace-power-nonsym}
    Let $\left\{M_n: n \in \mathbb{N}\right\}$ be a family of random $a(n)\times b(n)$ matrices. Let $r(n)= \min\{a(n), b(n)\}$. If $F(x)$ is some distribution function such that the following are true:
    \begin{enumerate}
        \item For each $k\geq 0$, $\displaystyle \beta_k = \lim_{n\to\infty}\, \dfrac{\,1\,}{r(n)}\, \Ebb\[ \trace\(\(M_n{M_n}^T\)^{k}\) \]$  exists and is finite. Moreover, $\displaystyle  \int x^{2k} \,dF(x) = \beta_k$ for all $k\geq 0$.
        \item $\displaystyle \dfrac{\,1\,}{r(n)}\, \Ebb\[ \trace\(\({M_n{M_n}^T}\)^{k}\) \] < \infty$ for all $n\geq 1, k\geq 0$. 
        \item $\displaystyle \sum_{n}\, \var\(\dfrac{\,1\,}{r(n)}\, \trace\(\(M_n{M_n}^T\)^{k}\)\) <\infty$.
        \item (Carleman condition) $\displaystyle \sum {\beta_{2k}}^{-1/2k} = \infty$.
    \end{enumerate}
    Then $F_n(x)$, the distribution of singular values of $M_n$, converges to $F(x)$ almost surely.
\end{thm}

In this version of the paper, we focus on the first condition and only briefly sketch how to verify the remaining conditions. The second condition is trivial as all of the distributions have bounded moments. For the third condition, it is sufficient to show that when we expand out $\dfrac{\,1\,}{r(n)^2}\, \Ebb\[ \trace\(\({M_n{M_n}^T}\)^{k}\)^2 \]$, the nonzero terms which do not appear in $\dfrac{\,1\,}{r(n)^2}\, \Ebb\[ \trace\(\({M_n{M_n}^T}\)^{k}\) \]^2$ have at least two fewer distinct indices than the dominant terms in $\dfrac{\,1\,}{r(n)^2}\, \Ebb\[ \trace\(\({M_n{M_n}^T}\)^{k}\) \]^2$. This implies that for each $k$, $\var\(\dfrac{\,1\,}{r(n)}\, \trace\(\(M_n{M_n}^T\)^{k}\)\)$ is $O\(\frac{1}{n^2}\)$. For the fourth condition, it is sufficient to show that if $\o$ and $\o'$ satisfy Carleman's condition then so does $\o \opr \o'$ and show that $\ozm$ satisfies Carleman's condition.

To show that the first condition holds, it suffices to prove the following statement.

\begin{restatable}{thm}{maintwo}
\label{thm:main-2}
    $\ozm\opr \o^{(1)}\opr \dots \opr \o^{(s)}$ and the limiting distributions of singular values of $\mzshapesdistrg$ have the same moments. In other words, if we let $\displaystyle r(n,m)=\dfrac{n!}{(n-m)!}$ then
    \begin{equation}
    \begin{aligned}
        \(\ozm\opr \o^{(1)}\opr \dots \opr \o^{(s)}\)_{2k} 
        &= \lim_{n\to\infty}\, \dfrac{\,1\,}{r(n,m)}\cdot\Ebb\[\trace\(\(M^{(G)}_{Z(m),s}{M^{(G)}_{Z(m),s}}^T\)^{k}\)\] \\
        &= \(\ozms\)_{2k}\,.
    \end{aligned}
    \end{equation}
\end{restatable}

We will prove the above result by proving their moments have the same recurrence relations. First, we will introduce the techniques for analyzing the trace power moments of the graph matrices. We will then find the recurrence relation for the trace power moments of $\mzshapesdistrg$ and show that it matches the recurrence relation for $\(\ozm\opr \o^{(1)}\opr\dots\opr \o^{(s)}\)_{2k}$.

\subsection{Preliminary Analysis for Graph Matrices}

\setlength{\parskip}{1.5mm}
\setlength{\baselineskip}{1.3em}

\begin{figure}[hbt!]
    \centering
    \includegraphics[scale=0.32]{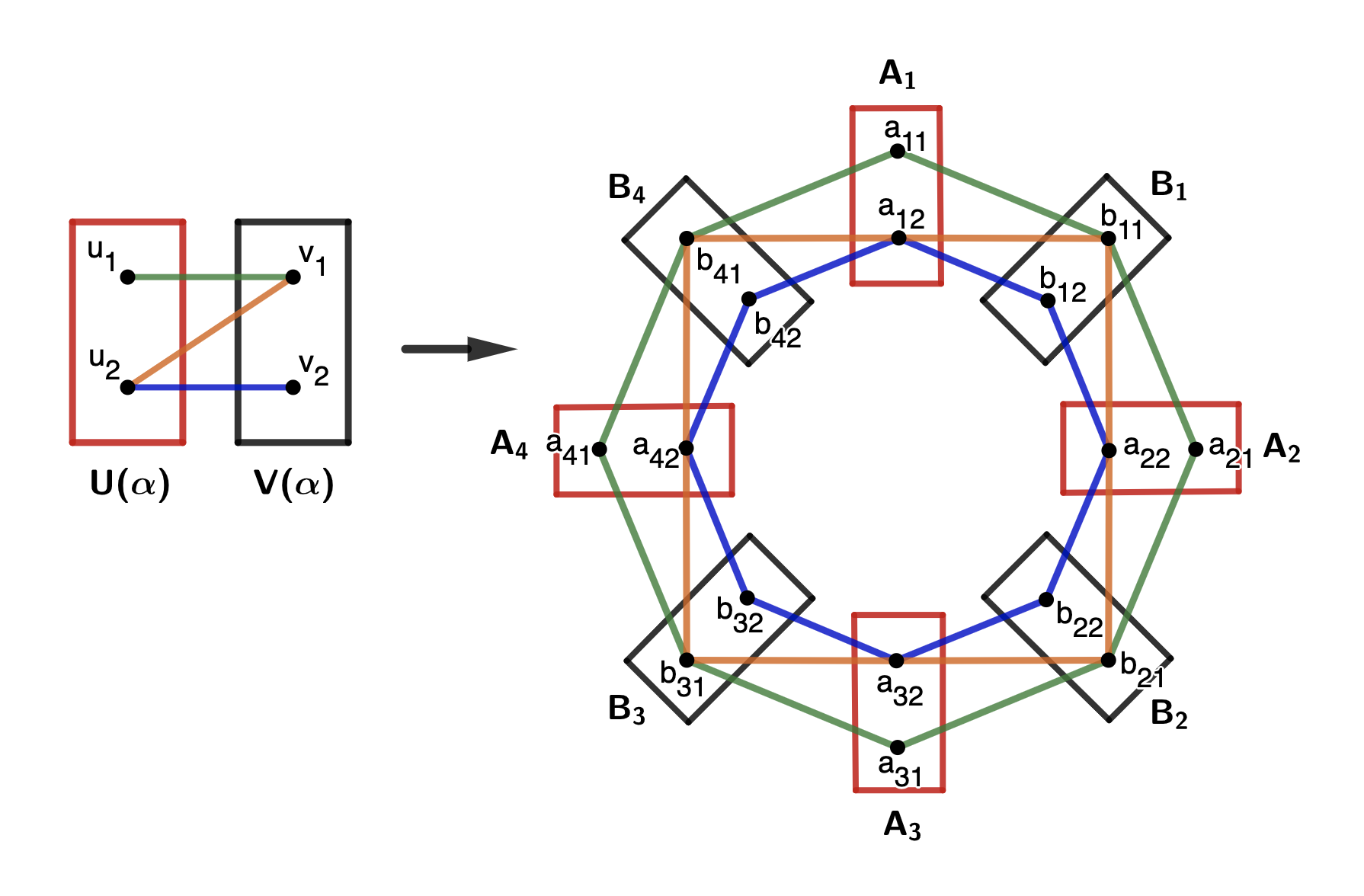}
    \caption{On the left is a shape $\a$, on the right is $H(\a,2q)$ where $q=4$.}
    \label{fig:H-a-2q)}
\end{figure}

\begin{defn}[Definition 3.2 of \cite{AMP20}] \label{defn:copies}
	Given a shape $\a$ and a $q \in \Nbb$, we define $H(\a,2q)$ to be the multi-graph which is formed as follows:
	\begin{enumerate}
		\item Take $q$ copies $\a_1,\ldots,\a_{q}$ of $\a$ and take $q$ copies $\a^{T}_1,\ldots,\a^{T}_q$ of $\a^{T}$, where $\a^{T}$ is the shape obtained from $\a$ by switching the role of $U_{\a}$ and $V_{\a}$.
		\item For all $i \in [q]$, we glue them together by setting $V_{\a_i} = U_{\a^{T}_i}$ and $V_{\a^{T}_i} = U_{\a_{i+1}}$ (where $\a_{q+1} = \a_1$).
	\end{enumerate}
	We define $V(\a,2q) = V(H(\a,2q))$ and we define $E(\a,2q) = E(H(\a,2q))$. See Figure \ref{fig:H-a-2q)} for an illustration.
\end{defn}

\begin{rmk}
	$H(\a,2q)$ is defined as a multi-graph because edges will be duplicated if $U_{\a}$ or $V_{\a}$ contains one or more edges. That said, in this paper we only consider $\a$ such that $U_{\a}$ and $V_{\a}$ do not contain any edges, so here $H(\a,2q)$ will always be a graph.
\end{rmk}
	
\begin{defn}[Definition 3.8 of \cite{AMP20}: Constraint graphs on $H(\alpha,2q)$]
	We define \\$\mathcal{C}_{\a,\oa,2q} = \left\{C(\phi): \phi: V(\alpha,2q) \to [n] \text{ is piecewise injective}\right\}$ to be the set of all possible constraint graphs on $V(\alpha,2q)$ which come from a piecewise injective map $\phi: V(\alpha,2q) \to [n]$.
		
	Given a constraint graph $C \in \mathcal{C}_{\a,\oa,2q}$, we make the following definitions:
	\begin{enumerate}
		\item We define $N(C) = \abs{\left\{\phi: V(\alpha,2q) \to [n]: \phi \text{ is piecewise injective}, C(\phi) = C\right\}}$.
		\item We define $\val(C) = \Ebb\[\chi_{\phi\(E\(\alpha,2q\)\)}(G)\]$ where $\phi: V(\alpha,2q) \to [n]$ is any piecewise injective map such that $C(\phi)=C$.
	\end{enumerate}
		
	We say that a constraint graph $C$ on $H(\a,2q)$ is \textit{nonzero-valued} if $\val(C)\neq 0$.
\end{defn}

As observed in \cite{AMP20}, with these definitions $\displaystyle \Ebb\[\trace\(\(M_{\alpha}M_{\alpha}^T\)^q\)\]$ can be re-expressed as follows.

\begin{prop}[Proposition 3.9 of \cite{AMP20}] \label{prop:exp-value-diffrep}
    For all shapes $\alpha$ and all $q \in \Nbb$,
    \begin{equation*}
		\Ebb\[\trace\(\(M_{\a}M_{\a}^T\)^q\)\] = \sum_{C \in \mathcal{C}_{\a,\oa, 2q}}{N(C)\val(C)}.
    \end{equation*}
\end{prop}

\begin{prop}\label{prop:num-phi-maps-under-constr-graph}
	For every constraint graph $ C \in \mathcal{C}_{\a,\oa, 2q}$, $\displaystyle N(C) = \frac{n!}{\(n - \abs{V(\a,2q)} + \abs{E(C)}\)!}$. 
\end{prop}
\begin{proof}
    Observe that choosing a piecewise injective map $\phi$ such that $C(\phi) = C$ is equivalent to choosing a distinct element of $[n]$ for each of the $n - \abs{V(\alpha,2q)} + \abs{E(C)}$ connected components of $C$.
\end{proof}

Since the number of constraint graphs in $\mathcal{C}_{\a,\oa, 2q}$ depends on $q$ but not on $n$, as $n \to \infty$ we only care about the nonzero-valued constraint graphs in $\mathcal{C}_{\a,\oa, 2q}$ which have the minimum possible number of edges. We call such constraint graphs dominant.

\begin{defn}[Dominant Constraint Graphs]\label{defn:dominant-constraint-graph}
We say a constraint graph $C \in \C_{\a,\oa,2q}$ is a \textit{dominant constraint graph} if $\val(C) \neq 0$ and $
\abs{E(C)}=\min\left\{\abs{E(C')}: C' \in \C_{\a,\oa,2q}, \val(C') \neq 0\right\}$. We denote the set of dominant constraint graphs in $\C_{\a,\oa,2q}$ to be $\D_{\a,\oa,2q}$.
\end{defn}

\begin{defn}\label{defn:weighted-sum-of-dominant-constraint-graph}
    We define the \emph{weighted sum of the dominant constraint graphs $C\in\C_{\a,\oa,2q}$} to be 
    \begin{equation*}
        W\(\D_{\a,\oa,2q}\) = \sum_{C\in\D_{\a,\oa,2q}} \val(C).
    \end{equation*}
\end{defn}

\begin{rmk}
    When $\oa=\o_{\pm1}$, $W\(\D_{\a,\oa,2q}\)=
    \abs{\D_{\a,\oa,2q}} =$ the number of dominant constraint graphs in $\C_{\a,\oa,2q}$.
\end{rmk}

We now state the number of edges in dominant constraint graphs.
\begin{defn}[Vertex Separators]
    We say that $S \subseteq V(\alpha)$ is a \emph{vertex separator} of $\alpha$ if every path from a vertex $u \in U_{\alpha}$ to a vertex $v \in V_{\alpha}$ contains at least one vertex in $S$.
\end{defn}

\begin{defn}
	Given a shape $\alpha$, define $s_{\alpha}$ to be the minimum size of a vertex separator of $\alpha$.
\end{defn}

\begin{lemma}[Follows from Lemma 6.4 of \cite{AMP20}]\label{lem:min-constr-edges-for-bipartite-shape}
    For any bipartite shape $\a$, for any nonzero-valued $C \in \mathcal{C}_{\a,\oa, 2q}$, $\abs{E(C)} \geq (q-1)s_{\a}$. Moreover, the bound is tight. i.e. There exists a nonzero-valued $C \in \C_{\a,\oa, 2q}$ such that $\abs{E(C)} = (q-1)s_{\a}$.
\end{lemma}

\begin{cor}\label{cor:dominant-number-of-edges}
For all bipartite shapes $\a$, for all dominant constraint graphs $C \in \mathcal{C}_{(\alpha,2q)}$, $\abs{E(C)} = (q-1)s_{\alpha}$.
\end{cor}

The following Corollary follows from \Cref{prop:exp-value-diffrep}, \Cref{prop:num-phi-maps-under-constr-graph} and \Cref{cor:dominant-number-of-edges}.

\begin{cor}\label{cor:expected-trace-for-bipartite-shapes}
    For all bipartite shapes $\alpha$ associated with distributions $\oa$, taking $\displaystyle r_{approx}(n) = \frac{n!}{(n - s_{\alpha})!}$,
    \begin{equation}
        \lim_{n \to \infty}{\frac{1}{r_{approx}(n)}\Ebb\[\trace\(\(\frac{M_{\a,\oa}M_{\a,\oa}^T}{n^{|V(\a)| - s_{\a}}}\)^q\)\]} = W\(\D_{\a,\oa,2q}\)
    \end{equation}
\end{cor}

Thus, to determine the spectrum of the singular values of $M_{\alpha}$ for a bipartite shape $\alpha$, we need to count the number of constraint graphs $C \in \mathcal{C}_{(\alpha,2q)}$ such that $C$ is dominant.

\begin{rmk}
    We write $r_{approx}$ rather than $r$ here because if $s_{\alpha} \leq \min{\left\{|U_{\alpha}|,|V_{\alpha}|\right\}}$ then the rank of $M_{\alpha}$ will generally be $\displaystyle\frac{n!}{\(n - \min{\left\{|U_{\alpha}|,|V_{\alpha}|\right\}}\)!}$ rather than $\displaystyle\frac{n!}{\(n - s_{\alpha}\)!}$.
\end{rmk}

\begin{rmk}
    The same statement is true for general $\alpha$ except that the number of edges in a dominant constraint graph $C \in \mathcal{C}_{(\alpha,2q)}$ is $q\abs{V\(\alpha\) \setminus \(U_{\alpha} \cup V_{\alpha}\)} + (q-1)(s_{\alpha} - \abs{U_{\alpha} \cap V_{\alpha}})$ rather than $(q-1)s_{\alpha}$.
\end{rmk}

Applying \Cref{cor:expected-trace-for-bipartite-shapes} to $\mzshapesdistrg$, we get the following result.

\begin{cor}\label{cor:trace-power-weight-dominated-constraint}
Let $\oa=\o_{\mzshape,s}$ be as in \Cref{defn:mzshape-sdistr-graph-matrix}. Let $\displaystyle r(n,m)=\dfrac{n!}{(n-m)!}$. Then
    \begin{equation}
        \lim_{n \to \infty} {\dfrac{1}{r(n,m)}\cdot \Ebb\[\trace\(\(\mzshapesdistrg{\mzshapesdistrg}^T\)^q\)\]} = W\(\D_{\mzshape,\oa,2q}\)
    \end{equation}
\end{cor}

\begin{proof}
Let $\a = \mzshape$. Then $s_{\a} = m$, $\displaystyle r_{approx}(n) = \dfrac{n!}{(n-m)!}$, and $n^{\abs{V(\a)-s_{\a}}} = n^{m}$. Thus $\displaystyle \dfrac{M_{\a,\oa}{M_{\a,\oa}}^T}{n^{|V(\a)| - s_{\a}}} = \dfrac{M_{\a,\oa}{M_{\a,\oa}}^T}{n^{m}} = \(\dfrac{M_{\mzshape,s}}{n^{m/2}}\)\(\dfrac{M_{\mzshape,s}}{n^{m/2}}\)^T = \mzshapesdistrg{\mzshapesdistrg}^T$. Plugging these into \Cref{cor:expected-trace-for-bipartite-shapes} we get the above statement.
\end{proof}

\begin{rmk}
    By \Cref{cor:trace-power-weight-dominated-constraint}, to prove \Cref{thm:main-2}, it suffices to prove
    \begin{equation}
        \(\ozm\opr \o^{(1)}\opr \dots \opr \o^{(s)}\)_{2k} 
        = W\(\D_{\mzshape,\oams,2k}\)
    \end{equation}
    for all $k\geq 0$.
\end{rmk}

In the following two subsections, we prove this by proving they have the same recurrence relations,




\subsection{Recurrence Relations for \texorpdfstring{$\(\ozm\opr \o^{(1)}\opr \dots \opr\o^{(s)}\)_{2k}$}{(Ωz(m) ºR Ω1 ºR … ºR Ωs)2k}}

\setlength{\parskip}{1.5mm}
\setlength{\baselineskip}{1.3em}

We will use the notation $A_m^{(s)}\(k,0,\dots,0\)$ from \Cref{section:ozm-s}.
\begin{equation}
\begin{aligned}
    A_m^{(s)}\(k,0,\dots,0\) 
    &= \(\ozm\opr \o^{(1)}\opr\dots \opr \o^{(s)}\)_{2k} \\
    &= \sum_{\vaa_i\in P_k}\, C_m\(\vaa_1,\dots,\vaa_s\)\cdot \vv{\o^{(1)}}^{\vaa_1}\dots \vv{\o^{(s)}}^{\vaa_s}
\end{aligned}
\end{equation}
where 
\begin{equation}
    C_m\(\vaa_1,\dots,\vaa_s\) = \binom{(m-s+1)k}{mk-(a_1+\dots+a_s)+1}\cdot k^{s-1} \cdot \prod_{i=1}^s \dfrac{(a_i-1)!}{\a_{i1}!\dots\a_{ik}!}\,.
\end{equation}

\begin{rmk}
    Note by definition of $ C_m\(\vaa_1,\dots,\vaa_s\)$, $ C_m\(\vaa_1,\dots,\vaa_s\)\neq 0 \implies m-s+1 > 0 \implies m\geq s$. Thus we can always assume that $m\geq s$.
\end{rmk}

\begin{prop}\label{prop:Ams-alternative}
Assume $m\geq s$. Let $\mh = m-s+1$. Then 
    \begin{equation}
        A_m^{(s)}(k,0,\dots,0) = \sum_{\vaa_i\in P_k}\, \sum_{\substack{\P_1\sqcup\dots \sqcup\P_s \in \\ \npm(\vaa_1,\dots,\vaa_s)}}\, \(\o^{(1)}_{2\,|P_{1,1}|/\mh}\cdot\dots\cdot \o^{(1)}_{2\,|P_{1,a_1}|/\mh}\)\cdot \prod_{j=2}^{s} \(\o^{(j)}_{2\,|P_{j,1}|}\cdot\dots\cdot \o^{(j)}_{2\,|P_{j,a_j}|}\)
    \end{equation}
    where $\P_i=\big\{P_{i,1},\dots,P_{i,a_i}\big\}$ for each $i\in[s]$.
\end{prop}

\begin{proof}
    By \Cref{thm:num-of-np}, $\displaystyle C_m\(\vaa_1,\dots,\vaa_s\) = \big|\np_m\(\vaa_1,\dots,\vaa_s\)\big|$. In particular, the size of any partition set corresponding to $\vaa_1$ is a multiple of $\mh=m-s+1$. Thus
    \begin{align*}
        A(k,0) 
        &= \sum_{\vaa_i\in P_k} C_m(\vaa_1,\dots,\vaa_s)\cdot \vv{\o^{(1)}}^{\vaa_1}\cdot\dots\cdot \vv{\o^{(s)}}^{\vaa_s} \\
        &= \sum_{\vaa_i\in P_k} \big|\np_m\(\vaa_1,\dots,\vaa_s\)\big|\cdot \vv{\o^{(1)}}^{\vaa_1}\cdot\dots\cdot \vv{\o^{(s)}}^{\vaa_s}\\
        &= \sum_{\vaa_i\in P_k}\; \sum_{\substack{\P_1\sqcup\dots \sqcup\P_s \in \\ \npm(\vaa_1,\dots,\vaa_s)}}\, \({\o^{(1)}_2}^{\a_{11}}\dots,{\o^{(1)}}_{2k}^{\a_{1k}}\)\cdot\dots\cdot \({\o^{(s)}_2}^{\a_{s1}}\dots,{\o^{(s)}}_{2k}^{\a_{sk}}\) \,\\
        &= \sum_{\vaa_i\in P_k}\; \sum_{\substack{\P_1\sqcup\dots \sqcup\P_s \in \\ \npm(\vaa_1,\dots,\vaa_s)}}\, \(\o^{(1)}_{2\,|P_{1,1}|/\mh}\cdot\dots\cdot \o^{(1)}_{2\,|P_{1,a_1}|/\mh}\)\cdot \prod_{j=2}^{s} \(\o^{(j)}_{2\,|P_{j,1}|}\cdot\dots\cdot \o^{(j)}_{2\,|P_{j,a_j}|}\) \,.
    \end{align*}
\end{proof}

Next we will define $A_m^{(s)}\(k,r_1,\dots,r_s\)$. For simplicity, we will only give the definition for the case $s=1$ here and defer the full definition to \Cref{subsec:recur-grid-general} (see \Cref{defn:Amskr}). 

When $s=1$, we denote $A_m^{(s)}\(k,r_1,\dots,r_s\)$ as $A_m\(k,r\)$.

\begin{defn}\label{defn:Amkr}
    Let $r\geq 0$. We define $A_m(k,r)$ to be
    \begin{equation}
        A_m\(k,r\) = \sum_{\vaa\in P_k}\, \sum_{\substack{\P=\{P_1,\dots,P_a\}\in \\ \np\(m\vaa\):\, 1\in P_1 }}\,  \(\o_{2\,\(|P_1|/m+r\)}\o_{2\,|P_2|/m} \cdot\dots\cdot \o_{2\,|P_a|/m}\)
    \end{equation}
\end{defn}

\begin{rmk}
    When $m=1$, we further denote $A_m\(k,r\)$ as $A\(k,r\)$. Then
    \begin{equation}
        A\(k,r\) = \sum_{\vaa\in P_k}\, \sum_{\substack{\P=\{P_1,\dots,P_a\}\in \\ \np(\vaa):\, 1\in P_1 }}\,  \(\o_{2\,\(|P_1|+r\)}\o_{2\,|P_2|} \cdot\dots\cdot \o_{2\,|P_a|}\)
    \end{equation}
\end{rmk}

\begin{rmk}
    When $r=0$, $A_m\(k,r\)$ coincides with $A_m\(k,0\)$.
\end{rmk}

We will prove the following main result.
\begin{restatable}{thm}{recurGrid}
\label{thm:recur-grid-general}
\begin{equation}
\begin{aligned}
    A_{m}^{(s)}\(k,r_1,\dots,r_s\) 
    = \sum_{\substack{i_1,\dots,i_{m+1}\geq 0: \\ i_1+\dots+i_{m+1} = k-1}}\, & A_{m}^{(s)}\(i_1,r_1+1,0,\dots,0\)\dots A_{m}^{(s)}\(i_s,0,\dots,0,r_s+1\)\cdot\\
    & A_{m}^{(s)}\(i_{s+1},0,\dots,0\)\dots A_{m}^{(s)}\(i_{m+1},0,\dots,0\)
\end{aligned}
\end{equation}
\end{restatable}

\subsubsection{Base Case s=1, m=1}\label{subsec:recur-grid-s1m1}

\setlength{\parskip}{1.5mm}
\setlength{\baselineskip}{1.3em}

\begin{defn}\label{defn:line shape}
    Let $\a_{0}$ be the bipartite shape with vertices $V(\a_{0})=\{u,v\}$ and a single edge $\{u,v\}$ with distinguished tuples of vertices $U_{\a_{0}}=(u)$ and $V_{\a_{0}}=(v)$. We call $\a_0$ the \textit{line shape}.
\end{defn}

\begin{rmk}
    By definition, $\a_{Z(1)} = \a_0$. We will let $\o_{\a_0}$ denote $\o_{Z(1)}$.
\end{rmk}

When $m=1, s=1$, applying \Cref{prop:Ams-alternative}, we have that
\begin{equation}
\begin{aligned}
    A(k,0) 
    &= \(\o_{Z(1)}\opr \o\)_{2k} = \(\o_{\a_0}\opr \o\)_{2k} = \sum_{\vaa\in P_k} \binom{k}{a-1}\cdot \dfrac{(a-1)!}{\a_1!\dots\a_k!}\cdot \voo^{\vaa} \\
    &= \sum_{\vaa\in P_k}\, \sum_{\P=\{P_1,\dots,P_a\}\in  \np(\vaa)}\, \(\o_{2\,|P_1|}\cdot\dots\cdot \o_{2\,|P_a|}\)
\end{aligned}
\end{equation}
and
\begin{equation}
    A\(k,r\) = \sum_{\vaa\in P_k}\, \sum_{\substack{\P=\{P_1,\dots,P_a\}\in \\ \np(\vaa):\, 1\in P_1 }}\,  \(\o_{2\,\(|P_1|+r\)}\o_{2\,|P_2|} \cdot\dots\cdot \o_{2\,|P_a|}\)
\end{equation}

We will prove the following recurrence relation.
\begin{thm}\label{thm:recur-grid-base-line}
    \begin{equation}
    A\(k,r\) 
    = \sum_{\substack{i, j\geq 0: \\ i+j = k-1}}\, A\(i,0\)\cdot A\(j,r+1\)\,
    \end{equation}
and 
\begin{equation}
    A(0,r) = \o_{2r}\,.
\end{equation}
\end{thm}

The following result is not needed for the proof of \Cref{thm:recur-grid-base-line}, but gives us an alternative way of looking at $A_m(k,r)$ when $s=1$. More precisely, we can deduce an explicit formula for $A_m(k,r)$.
\begin{defn}
    Given $\vaa\in P_k$, we denote $\voo^{\vaa}(i,+r)$ to be $\o_2^{\a_1}\cdot\dots\cdot\(\o_{2i}^{\a_i-1}\o_{2(i+r)}\)\cdot\dots\cdot\o_{2k}^{\a_k}$.
\end{defn}

\begin{prop}
    When $s=1$, $\displaystyle A_m(k,0) = \sum_{\vaa\in P_{k}}\,  C_m\(\vaa\)\cdot\voo^{\vaa}$ and
    \begin{align}
        A_m(k,r) = \sum_{\vaa\in P_{k}}\,  C_m\(\vaa\)\cdot \( \sum_{i}^k\,  \dfrac{i\cdot\a_{i}}{k} \cdot \voo^{\vaa}\(i,+r\)\)\,.
    \end{align}
\end{prop}


\begin{eg}
$ $

\begin{enumerate}
    \item $A(2,0) = \o_2^2 + \o_4$, $A(2,1) = \o_2\o_4 + \o_6$, $A(2,2) = \o_2\o_6 + \o_8$.
    \item $A(3,0) = \o_2^3 + 3\,\o_2\o_4 + \o_6$, $A(3,1) = \o_2^2\o_4 + 3\(\dfrac{\,1\,}{3}\,\o_4\o_4 + \dfrac{\,2\,}{3}\,\o_2\o_6\) + \o_8 = \o_2^2\o_4 +\o_4^2 + 2\,\o_2\o_6 + \o_8$.
\end{enumerate}
    
\end{eg}

\begin{proof}[Proof of \Cref{thm:recur-grid-base-line}]
\begin{figure}[hbt!]
    \centering
    \includegraphics[scale=0.32]{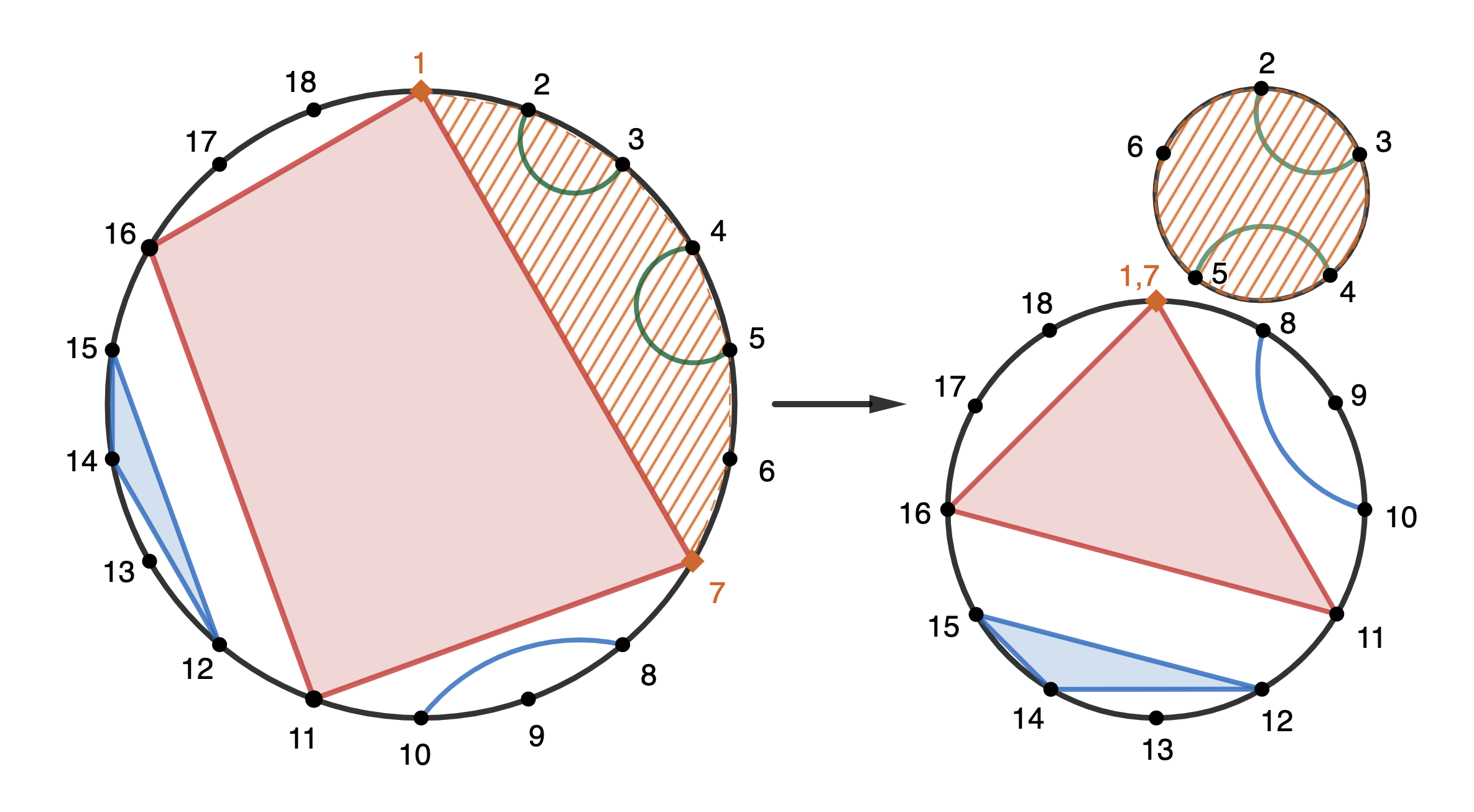}
    \caption{Illustration of \Cref{thm:recur-grid-base-line}: Here $i=6$. The orange shaded part on the right represents $A(i-1,0)$, and the other part represents $A(k-i,r+1)$.}
    \label{fig:recurrence-grid-base-1}
\end{figure}

We obtain this recurrence relation with the following steps:
\begin{enumerate}[i.]
    \item Given a noncrossing partition $\P$ of $[k]$, we will look at the smallest $i\geq 1$ such that $1$ and $i+1$ are in the same partition set (i.e. $1$ and $i+1$ are vertices of the same polygon $P_1$, and $i+1$ is immediately after $1$ in the clockwise direction). In the case when $1$ is isolated, $i$ is chosen to be $k$.
    \item Since $\P$ is noncrossing, everything in between $1$ and $i+1$ can be viewed a noncrossing partition on the cycle with elements $\{2,3,\dots,i\}$ (i.e. cycle of length $i-1$). This gives rise to $A(i-1,0)$.
    \item Now we identify $1$ with $i+1$ and shrink the size of the polygon $P_1$ by $1$. The remaining part from $i+1$ to $k$ is a noncrossing partition on $k-i$ elements, and this gives rise to $A(k-i,r+1)$.
\end{enumerate}

See \Cref{fig:recurrence-grid-base-1} for an illustration.

More precisely,
\begin{align*}
    A(k,r)
    &= \sum_{\substack{\P = \{P_1,\dots,P_a\}\in\np_k:\\ 1\in P_1}}\, \o_{2\,\(\abs{P_1}+r\)} \o_{2\,\abs{P_2}}\cdot\dots \cdot \o_{2\,\abs{P_a}} \\
    &= \sum_{i=1}^k \sum_{\substack{\P = \{P_1,\dots,P_a\}\in\np_k:\\ 1,i+1\in P_1,\, j\not\in P_1 \text{ for any } 2\leq j\leq i}}\, \o_{2\,\(\abs{P_1}+r\)} \o_{2\,\abs{P_2}}\cdot\dots \cdot \o_{2\,\abs{P_a}} \\
    &= \sum_{i=1}^k \sum_{\substack{\P = \{P_1,\dots,P_a\}\in\np_k:\\ 1,i+1\in P_1,\, j\not\in P_1 \text{ for any } 2\leq j\leq i}}\, \(\o_{2\,\abs{P_{i_1}}}\cdot\dots \cdot \o_{2\,\abs{P_{i_b}}}\)\cdot \(\o_{2\,\(\abs{P_1}+r\)}\o_{2\,\abs{P_{j_1}}}\cdot\dots\cdot \o_{2\,\abs{P_{j_{c-1}}}} \) \\
    &\hspace{2cm} \text{where } \{i_1,\dots, i_b\}\sqcup \{1,j_1,\dots,j_{c-1}\} = [a], \\
    &\hspace{3cm} P_{i_1}\sqcup\dots\sqcup P_{i_b} = \{2,\dots,i\}, P_1\sqcup P_{j_1}\sqcup\dots\sqcup P_{j_{c-1}} = \{1,i+1,\dots,k\} \\
    &= \sum_{i=1}^k \(\sum_{\substack{\Q=\{Q_1,\dots,Q_b\}\\ \in\np(i-1)}}\, \o_{2\,\abs{Q_1}}\cdot\dots\cdot \o_{2\,\abs{Q_b}}\)\cdot \(\sum_{\substack{ \R=\{R_1,\dots, R_c\}\in \\ \np(k-i):\, 1\in R_1 }} \o_{2\,\(\abs{R_1}+r+1\)}\cdot\dots\cdot \o_{2\,\abs{R_c}}\)\\
    &= \sum_{i=1}^k \, A(i-1,0)\cdot A(k-i,r+1) 
    = \sum_{i,j\geq 0:\, i+j=k-1} \, A(i,0)\cdot A(j,r+1)
\end{align*}

\end{proof}

\begin{defn}\label{defn:gamma-shape}
    Let $\c$ be the shape with vertices $V(\c)=\{u,w,v\}$ and edges \\ $E(\c)=\left\{\{u,v\},\{v,w\}\right\}$ with distinguished tuples of vertices $U_{\a_{Z}}=(u,w)$ and $V_{\a_{Z}}=(v,w)$. Let $\oc$ be the distributions associated with $\c$ where $\o_{\{u,v\}} = \o_{\pm1}$ and $\o_{\{v,w\}} = \o$.
    
    \begin{figure}[hbt!]
        \centering
        \begin{subfigure}[t]{.3\textwidth}
            \includegraphics[width=1\linewidth]{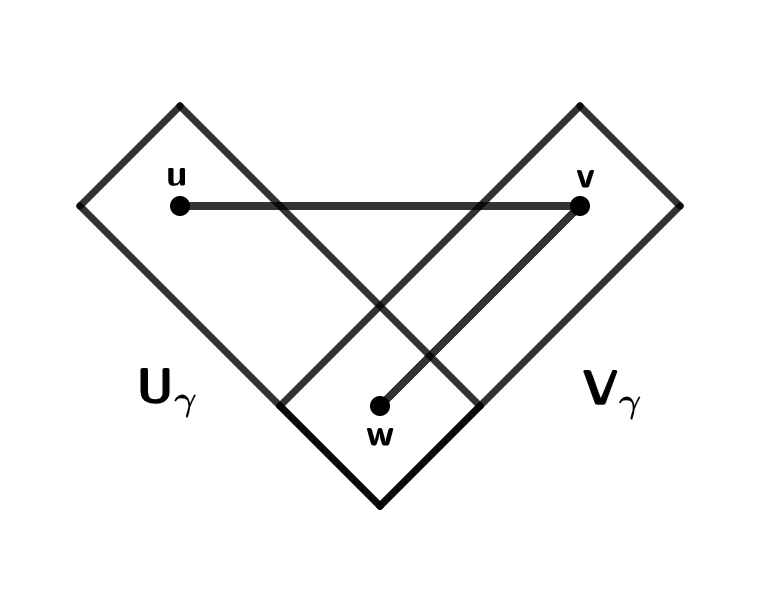}
            \caption{Shape $\c$.}
            \label{fig:gamma-shape}
        \end{subfigure}\hspace{0.7cm}
        \begin{subfigure}[t]{.3\textwidth}
            \includegraphics[width=1\linewidth]{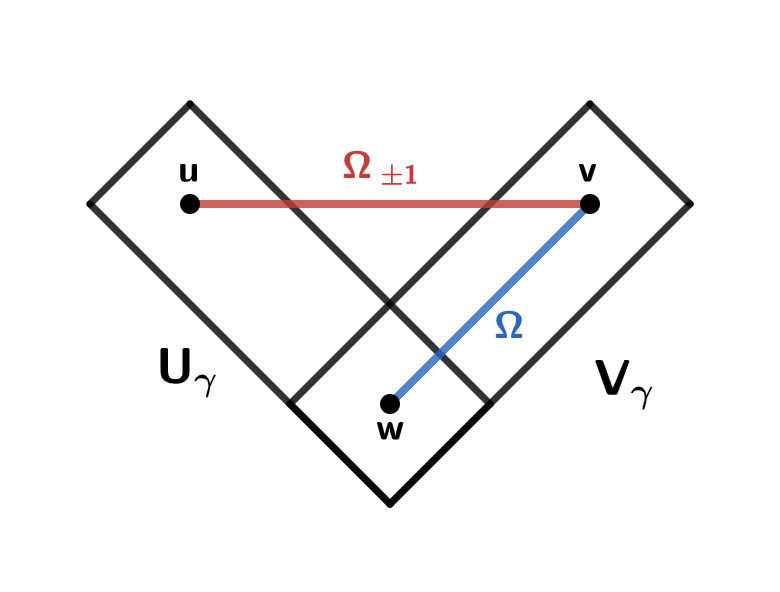}
            \caption{Associated distributions $\oc$ for shape $\c$.}
            \label{fig:gamma-distribution-2}
        \end{subfigure}
        \caption{Shape $\c$ and its associated distributions $\oc$.}
        \label{fig:gamma-shape-distribution}
    \end{figure}
\end{defn}

\begin{rmk}
    When $m=1$, $\o_{Z(1)} = \o_{\a_0}$, the line shape distribution. $\o_{\a_0}\opr \o$ will have the same limiting distribution of singular value moments as $M_{\c,\oc}$, the graph matrix of the shape $\c$ associated with distributions $\oc$. This fact can also be proved by matching their recurrence relations of the moments.
\end{rmk}

\subsubsection{Base Case s=1}\label{subsec:recur-grid-s1}

\setlength{\parskip}{1.3mm}
\setlength{\baselineskip}{1.3em}

When $s=1$, applying \Cref{prop:Ams-alternative}, we have that
\begin{equation}
\begin{aligned}
    A_m(k,0) 
    &= \(\ozm\opr \o\)_{2k} = \sum_{\vaa\in P_k} \binom{mk}{a-1}\cdot \dfrac{(a-1)!}{\a_1!\dots\a_k!}\cdot \voo^{\vaa} \\
    &= \sum_{\vaa\in P_k}\; \sum_{\substack{\P = \{P_1,\dots,P_a\}\\ \in\np(m\vaa)}}\, \o_{2\,\abs{P_1}/m}\cdot\dots\cdot \o_{2\,\abs{P_a}/m} 
\end{aligned}
\end{equation}
and 
\begin{equation}
    A_m(k,r) =\sum_{\vaa\in P_k}\; \sum_{\substack{\P = \{P_1,\dots,P_a\}\\ \in\np(m\vaa):\, 1\in P_1}}\, \o_{2\,\(\abs{P_1}/m+r\)} \o_{2\,\abs{P_2}/m}\cdot\dots \cdot \o_{2\,\abs{P_a}/m}\,.
\end{equation}

We want to prove the following.
\begin{thm}\label{thm:recur-grid-base-mzshape}
\begin{equation}
    A_m\(k,r\) 
    = \sum_{\substack{i_1,\dots,i_{m+1}\geq 0: \\ i_1+\dots+i_{m+1} = k-1}}\, A_m\(i_1,0\)\dots A_m\(i_m,0\)\cdot A_m\(i_{m+1},r+1\)\,
\end{equation}
and 
\begin{equation}
    A_m(0,r) = \o_{2r}\,.
\end{equation}
\end{thm}

\begin{proof}
\begin{figure}[hbt!]
    \centering
    \includegraphics[scale=0.3]{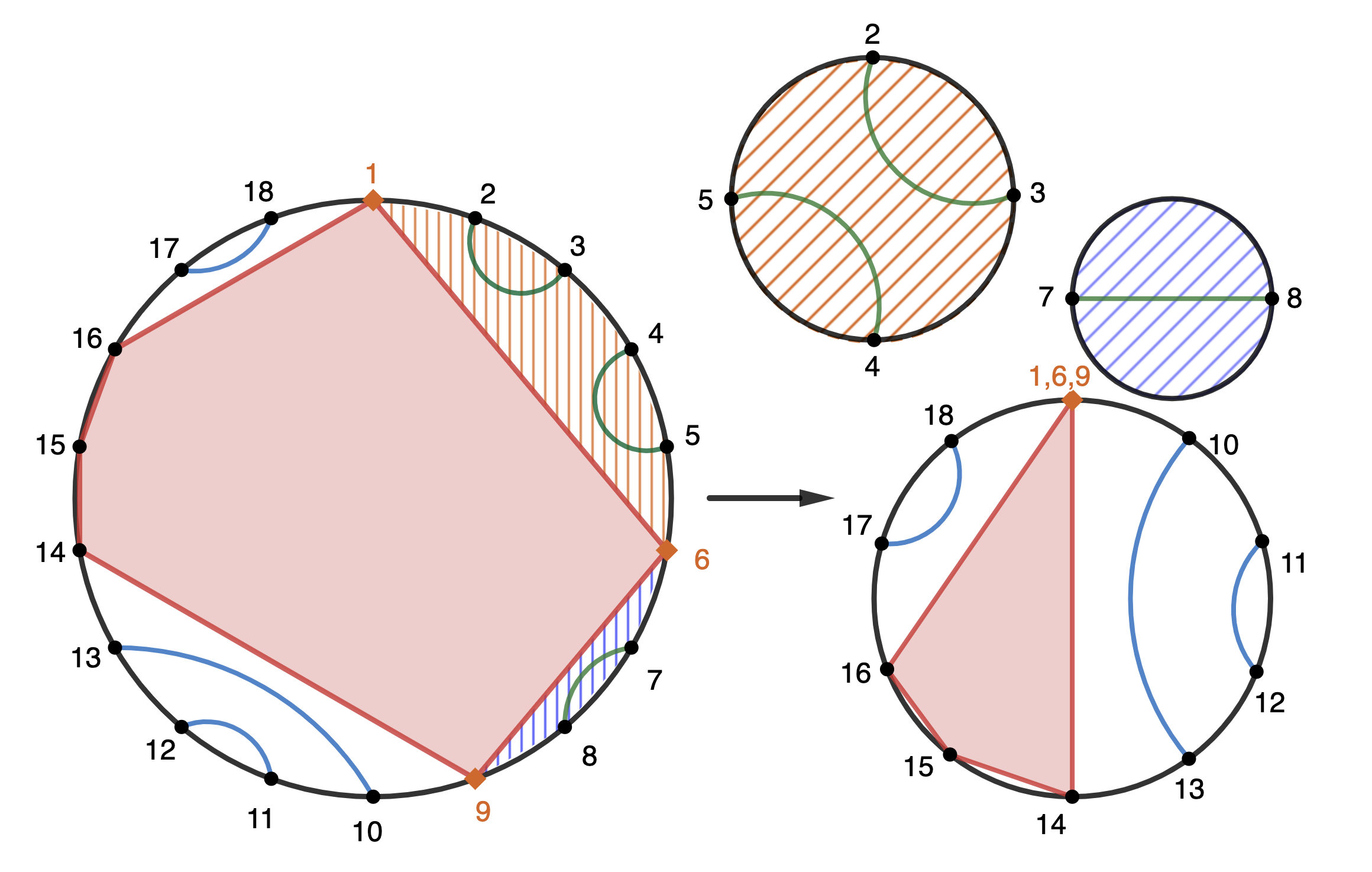}
    \caption{Illustration of \Cref{thm:recur-grid-base-mzshape}: Here $m=2, k=9$, $i_1=2$, $i_2=1$. The orange and blue shaded parts on the right represent $A_2(i_1,0)$ and $A_2(i_2,0)$, respectively. The remaining part represents $A_2(k-i_1-i_2,r+1)$.}
    \label{fig:recurrence-grid-base-2}
\end{figure}

\begin{figure}[hbt!]
    \centering
    \includegraphics[scale=0.3]{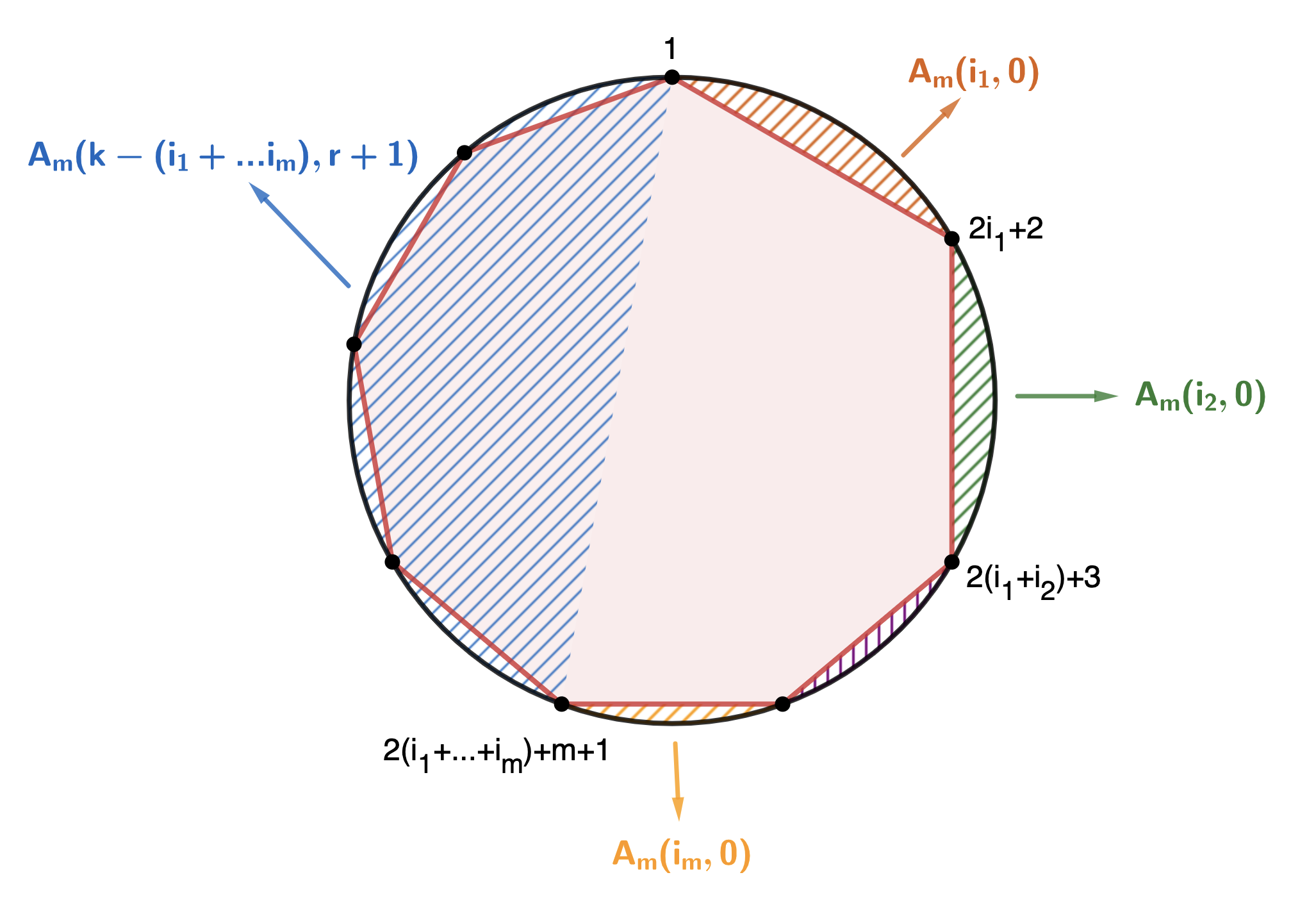}
    \caption{Illustration of \Cref{thm:recur-grid-base-mzshape}: }
    \label{fig:recurrence-grid-base-m}
\end{figure}

The proof is similar to the proof for \Cref{thm:recur-grid-base-line}, with the modification that now the size of every polygon in a partition is augmented by a factor of $m$. 

Consider a $\P=\{P_1,\dots,P_a\}\in\np(m\vaa)$ where $\vaa\in P_k$. Then $\abs{P_i}$ is a multiple of $m$ for all $i$. Since $\P$ is noncrossing, we can deduce that if $x<y$ are adjacent in the same partition set, then $y-x-1$ must be a multiple of $m$.

We will obtain this recurrence relation with the following steps:
\begin{enumerate}[i.]
    \item Let $\P\in \np(m\vaa)$ for some $\vaa\in P_k$. Since size of each partition set in $\P$ is a multiple of $m$, there are at least $m$ elements in $P_1$. W.L.O.G., we can assume $1\in P_1$. Denote $i_0=0$. Let $i_1,\dots,i_m\geq 0$ be such that $V=\{m(i_0+\dots+i_s)+s+1: s=0,1,\dots,m\}$ are the first $m+1$ elements in the clockwise direction that are in $P_1$. In the case when there are exactly $m$ elements in $P_1$, $m(i_1+\dots+i_m)+m+1 = mk+1 \iff i_m = k-i_1-\dots-i_{m-1}-1$. 
    \item Since $\P$ is noncrossing, everything in between $m(i_1+\dots+i_{s-1})+s$ and $m(i_1+\dots+i_s)+s+1$ can be viewed as an augmented noncrossing partition on the cycle with elements $\{m(i_1+\dots+i_{s-1})+s+1, m(i_1+\dots+i_{s-1})+s+2,\dots, m(i_1+\dots+i_{s})+s\}$ (i.e. cycle of length $mi_s$). This gives rise to $A(i_s,0)$.
    \item Now we identify all the vertices in $V$ and shrink the size of the polygon $P_1$ by $m$. The remaining part from $m(i_1+\dots+i_m)+m+1$ to $mk$ is an augmented noncrossing partition on $m(k-i_1-\dots-i_m-1)$ elements, and this gives rise to $A(k-i_1-\dots-i_m-1,r+1)$.
\end{enumerate}

See \Cref{fig:recurrence-grid-base-1} for an illustration.

More precisely, given $i_1,\dots,i_m\geq 0$, we denote $\vi = \(i_1,\dots,i_m\)$ and
\begin{itemize}
    \item $\(V_{\;\vi}\)_j = \big\{m(i_1+\dots+i_{j-1})+j+1, \dots, m(i_0+\dots+i_j)+j \big\}$ for $j\in [m]$, and 
    \item $\(V_{\;\vi}\)_0 = \big\{m(i_1+\dots+i_j)+j+1: j\in[m] \big\} \bigcup \big\{1,m(i_0+\dots+i_m)+m+2,\dots, mk \big\}$. i.e. $\displaystyle \(V_{\;\vi}\)_0 = [mk]\setminus \bigcup_{j=1}^m \(V_{\;\vec{i}}\)_j$.
\end{itemize}

Note that $\abs{\(V_{\;\vi}\)_j} = mi_j$ for $j\in[m]$, and $\abs{\(V_{\;\vi}\)_0} = m(k-i_1-\dots-i_m)$.

Then
\begin{align*}
     & A_m(k,r)
    = \sum_{\vaa\in P_k}\; \sum_{\substack{\P = \{P_1,\dots,P_a\}\\ \in\np(m\vaa):\, 1\in P_1}}\, \o_{2\,\(\abs{P_1}/m+r\)} \o_{2\,\abs{P_2}/m}\cdot\dots \cdot \o_{2\,\abs{P_a}/m} \\
    =& \sum_{i_1=0}^{k-1}\dots \sum_{i_m=0}^{k-1}\; \sum_{\vaa\in P_k}\; \sum_{\substack{\P = \{P_1,\dots,P_a\}\\ \in\np(m\vaa): 1\in P_1}}\, \(\prod_{j=1}^m \o_{2\,|P_{i_{j,1}}|/m} \dots \o_{2\,|P_{i_{j,a_j}}|/m}\)\cdot \(\o_{2\,\(|P_1|/m+r\)}\o_{2\,|P_{j_1}|/m}\dots \o_{2\,|P_{j_{c-1}}|/m} \) \\
    &\hspace{1cm} \text{where } \(\bigsqcup_{j=1}^m \big\{i_{j,1},\dots, i_{j, a_j}\big\}\) \bigsqcup\, \big\{1,j_1,\dots,j_{c-1}\big\} = [a],\, \vi = \(i_1,\dots,i_m\), \text{ and } \\
    &\hspace{2.2cm} P_{i_{j,1}}\sqcup\dots \sqcup P_{i_{j,a_j}} = \(V_{\;\vi}\)_{j}  \text{ for } j\in[m],\, P_1\sqcup P_{j_1}\sqcup\dots \sqcup P_{j_{c-1}} =\(V_{\;\vi}\)_{0} \\
    =& \sum_{i_1=0}^{k-1}\dots \sum_{i_m=0}^{k-1}\; \prod_{j=1}^m\,  \(\sum_{\substack{\Q_j=\{Q_{j,1},\dots,Q_{j,a_j}\}\\ \in\np(mi_j)}}\, \o_{2\,|Q_{j,1}|/m} \dots  \o_{2\,|Q_{j,a_j}|/m}\)\cdot \\
    &\hspace{4cm} \(\sum_{\substack{ \R=\{R_1,\dots, R_b\}\in \\ \np\(m(k-i_1-\dots-i_m-1)\):\, 1\in R_1 }} \o_{2\,\(\abs{R_1}/m+r+1\)}\cdot\dots\cdot \o_{2\,\abs{R_c}/m}\)\\
    =& \sum_{i_1=0}^{k-1}\dots \sum_{i_m=0}^{k-1}\; A(i_1,0)\cdot\dots\cdot A(i_m,0)\cdot A(k-i_1-\dots-i_m-1,r+1)\\ 
    =& \sum_{\substack{i_1,\dots, i_{m+1}\geq 0:\, \\ i_1+\dots+i_{m+1}=k-1}} \, A(i_1,0)\cdot\dots \cdot A(i_m,0)\cdot A(i_{m+1},r+1)
\end{align*}
as needed.

\end{proof}

\subsubsection{General Case}\label{subsec:recur-grid-general}

\setlength{\parskip}{1.5mm}
\setlength{\baselineskip}{1.3em}

Recall the definition of $A_m^(s)(k,0,\dots,0)$ and \Cref{prop:Ams-alternative}.
\begin{equation}
\begin{aligned}
    A_m^{(s)}&\(k,0,\dots,0\) 
    = \sum_{\vaa_i\in P_k}\, C_m\(\vaa_1,\dots,\vaa_s\)\cdot \vv{\o^{(1)}}^{\vaa_1}\dots \vv{\o^{(s)}}^{\vaa_s} \\
    &= \sum_{\vaa_i\in P_k}\, \sum_{\substack{\P_1\sqcup\dots \sqcup\P_s \in \\ \npm(\vaa_1,\dots,\vaa_s)}}\, \(\o^{(1)}_{2\,|P_{1,1}|/\mh}\cdot\dots\cdot \o^{(1)}_{2\,|P_{1,a_1}|/\mh}\)\cdot \prod_{j=2}^{s} \(\o^{(j)}_{2\,|P_{j,1}|}\cdot\dots\cdot \o^{(j)}_{2\,|P_{j,a_j}|}\)
\end{aligned}
\end{equation}
where $\mh = m-s+1$ and $\P_i=\big\{P_{i,1},\dots,P_{i,a_i}\big\}$ for each $i\in[s]$.

\begin{rmk}
    Let $\mh=m-s+1$. We will label the vertices in the cycle of length $mk$ by $v_{1,1},\dots,v_{1,\mh}, v_{2,1},\dots, v_{s,1}$, $v_{1,\mh+1},\dots,v_{1,2\mh}, v_{2,2},\dots, v_{s,2},\dots\dots v_{1, (k-1)\mh+1},\dots, v_{1,k\mh}, v_{2,k},\dots, v_{s,k}$ in the clockwise direction, starting from $1$. i.e. The vertices on the cycle are $\big\{v_{1,(t-1)\mh+1},\dots,v_{1,t\mh}, v_{2,t},\dots, v_{s,t}: t\in[k]\big\}$. More precisely, the set of vertices $\{v_{j,1},\dots,v_{j,k}\}$ represents the vertices used by the $j^{th}$ type partition. See \Cref{fig:np-ms-label} for illustration.
    
    \begin{figure}[hbt!]
    \centering
    \begin{subfigure}[t]{.45\textwidth}
        \includegraphics[width=1\linewidth]{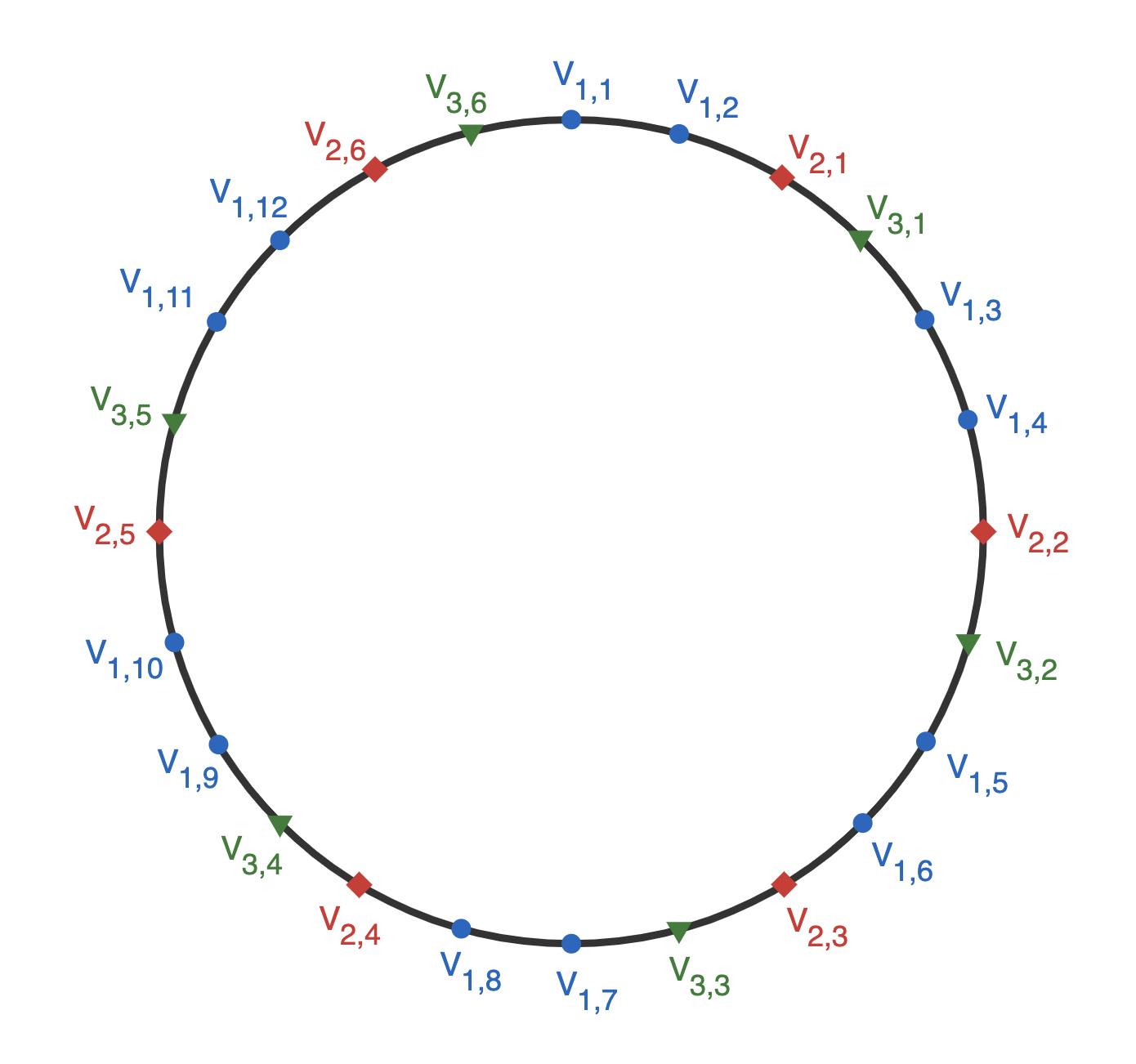}
        \caption{Label.}
        \label{fig:np-ms-label-1}
    \end{subfigure}\hspace{0.5cm}
    \begin{subfigure}[t]{.45\textwidth}
        \includegraphics[width=1\linewidth]{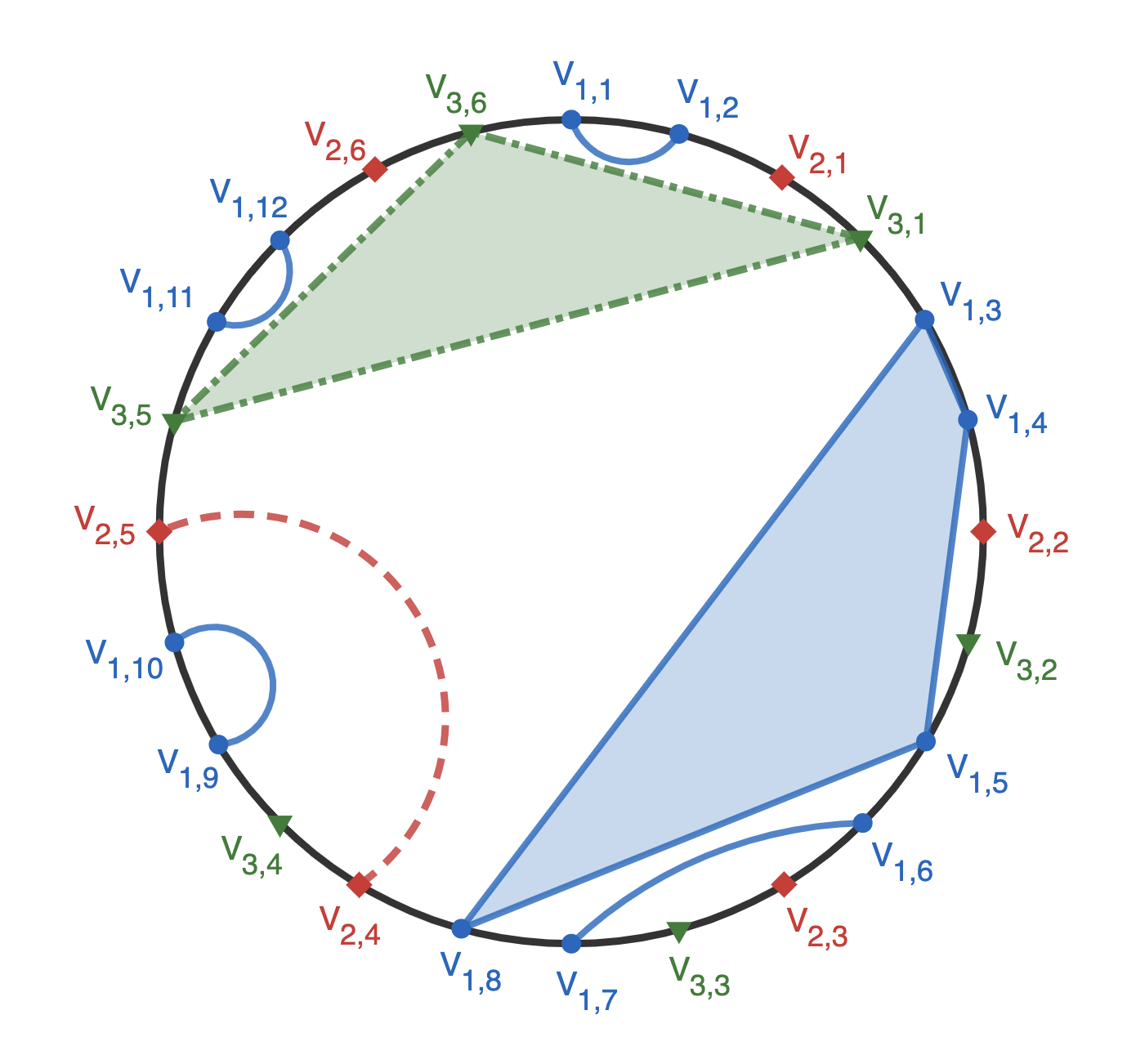}
        \caption{An example partition.}
        \label{fig:np-ms-label-2}
    \end{subfigure}
    \caption{Labelling of the cycle corresponding to $\npm\(\vaa_1,\dots,\vaa_s\)$. Here $s=3$, $m=4$, $\mh=m-s+1=2$ and $k=6$.}
    \label{fig:np-ms-label}
    \end{figure}
\end{rmk}

Now we will give the full definition of $A_m^{(s)}\(k,r_1,\dots,r_s\)$. 
\begin{defn}\label{defn:Amskr}
    Let $m\geq s$ and $\mh=m-s+1$. We define $A_m^{(s)}\(k,r_1,\dots,r_s\)$ to be
    \begin{equation}
    \begin{aligned}
        \sum_{\vaa_i\in P_k}\, \sum_{\substack{\P_1\sqcup\dots \sqcup\P_s \in \\ \npm(\vaa_1,\dots,\vaa_s)}}\, & \(\o^{(1)}_{2\,\(|P_{1,1}|/\mh+r_1\)}\o^{(1)}_{2\,|P_{1,2}|/\mh} \cdot\dots\cdot \o^{(1)}_{2\,|P_{1,a_1}|/m}\)\cdot \\ &\hspace{2cm} \prod_{j=2}^{s} \(\o^{(j)}_{2\,\(|P_{j,1}|+r_j\)}\cdot \o^{(j)}_{2\,|P_{j,2}|}\cdot \dots\cdot \o^{(j)}_{2\,|P_{j,a_j}|}\)
    \end{aligned}
    \end{equation}
    where $\P_i=\big\{P_{i,1},\dots,P_{i,\, a_i}\big\}$ for each $i\in[s]$, and $P_{i,1}$ are such that 
    \begin{enumerate}
        \item $v_{s,1}\in P_{s,1}$.
        \item Let $v_{s,\, j_s}$ be the first vertex in the clockwise direction from $v_{s,1}$ that is in $P_{s,1}$. In the case when $v_{s,1}$ is isolated, let $v_{s,\, j_s} = v_{s,1}$.
        \item Inductively, assume that we have picked $P_{s,1}, P_{s-1,1}, \dots, P_{i,1}$ and $v_{s,\, j_s}, v_{s-1,\, j_{s-1}}, \dots, v_{i,\, j_i}$, now we are going to pick $P_{i-1,1}$ and $v_{i-1,\, j_{i-1}}$. Let $P_{i-1,1}$ be the set containing $v_{i-1,j_i}$ (which is the vertex immediately next to $v_{i,j_i}$ in the counterclockwise direction). Let $v_{i-1,j_{i-1}}$ be the first vertex in the clockwise direction from $v_{i-1,j_i}$. In the case when $v_{i-1,j_i}$ is isolated, let $v_{i-1,\, j_{i-1}} = v_{i-1,j_{i}}$.
        \item Do the above process for $i=3,\dots,s$, resulting in $P_{2,1},\dots, P_{s,1}$ and $v_{2,j_2},\dots,v_{s,j_s}$. $v_{1,\mh j_2}$ is the vertex immediately next to $v_{2,j_2}$ in the counterclockwise direction. Now we let $P_{1,1}$ be the set containing $v_{1,\mh j_2}$.
    \end{enumerate}
    
    \begin{figure}[hbt!]
        \centering
        \includegraphics[scale=0.32]{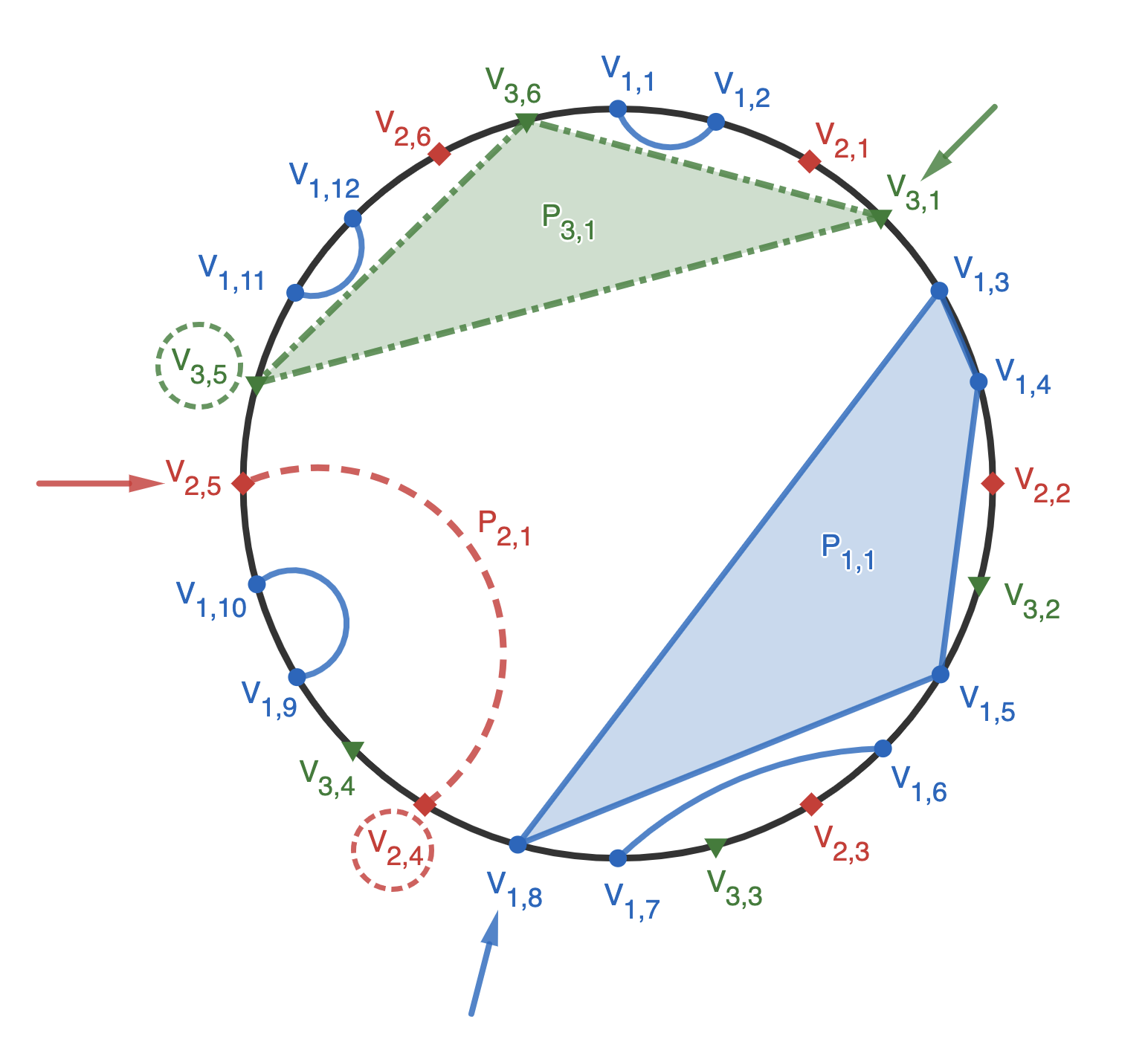}
        \caption{Illustration of \Cref{defn:Amkr}: Here $s=3$, $m=4$, $\mh=m-s+1=2$, and $k=6$. We start with $v_{s,1} = v_{3,1}$. Then $v_{3,j_3} = v_{3,5}$, $v_{2,j_3}=v_{2,5}$,  $v_{2,j_2} = v_{2,4}$ and $v_{1,\mh j_2} = v_{1,8}$. $P_{1,1}, P_{2,1}, P_{3,1}$ are as labeled. }
        \label{fig:np-Amsr}
    \end{figure}
    
    See \Cref{fig:np-Amsr} for illustration.
\end{defn}

\begin{rmk}
    When $r_1=\dots=r_s=0$, $A_m^{(s)}\(k,r_1,\dots,r_s\)$ coincides with $A_m^{(s)}\(k,0,\dots,0\)$.
\end{rmk}

We are now ready to prove the general case, \Cref{thm:recur-grid-general}.
\recurGrid*

\begin{proof}

We obtain the recurrence relation through the following steps.
\begin{enumerate}
    \item Let $v_{i,j_i}$ and $P_{i,1}$ for $i\in[s]$ be as in \Cref{defn:Amskr}. 
    
    \item For each $i\in\{2,\dots,s\}$, Identify $v_{i,j_{i+1}}$ with $v_{i,j_i}$. $j_{s+1}$ is defined to be $1$. Since for each $i$, $v_{i,j_i}$ is the first vertex in the clockwise direction from $v_{i,j_{i+1}}$ in $P_{i,1}$, everything in between $v_{i,j_i}$ and $v_{i,j_{i+1}}$ in the clockwise direction can be viewed as $\ams\(k_i,0,\dots,r_i+1,0\dots, 0\)$, where $k_i= j_{i+1}-j_i \mod k$.
    

\begin{figure}[hbt!]
    \centering
    \includegraphics[scale=0.35]{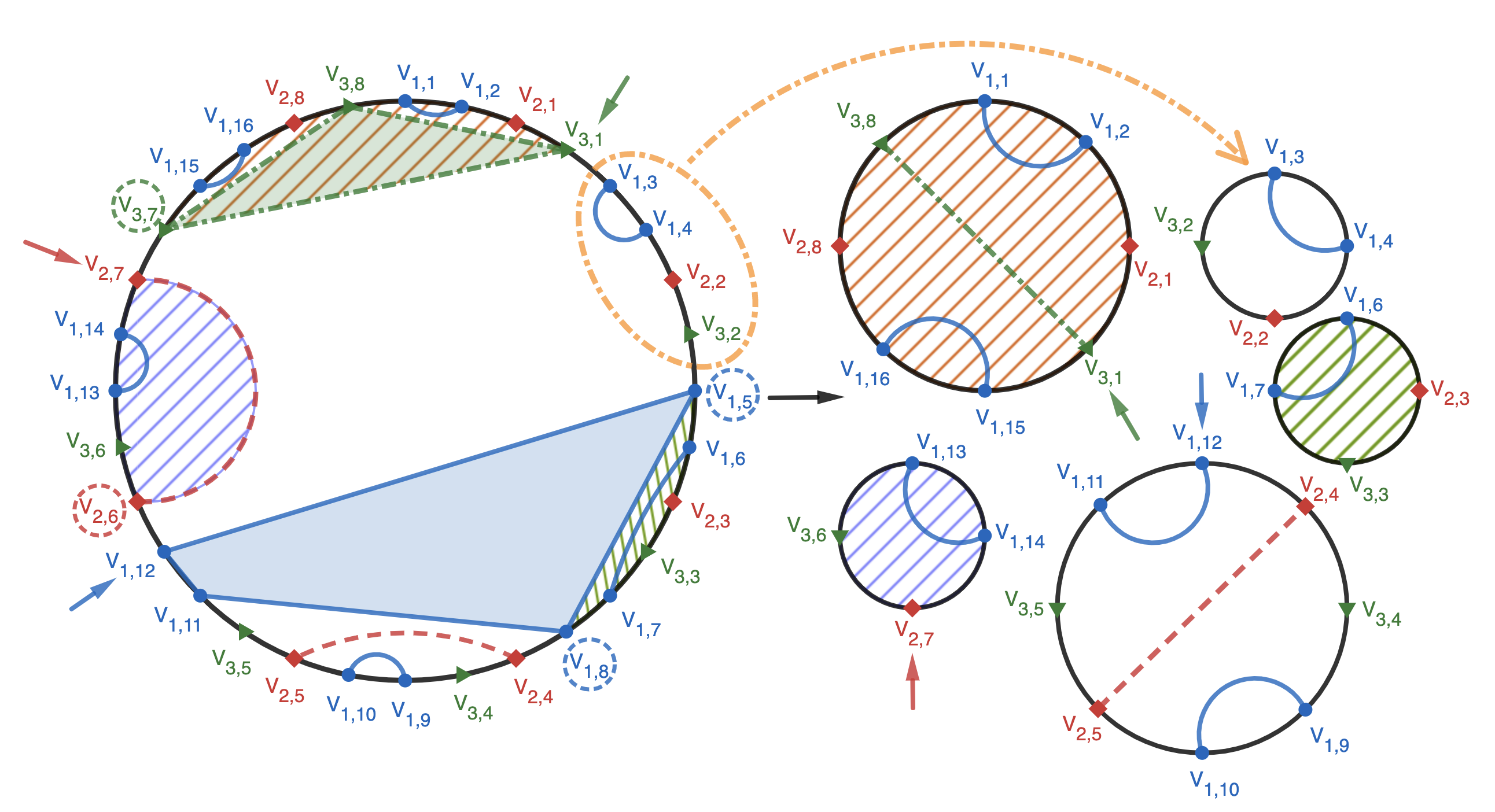}
    \caption{Illustration of the proof of \Cref{thm:recur-grid-general}: Here $s=3$, $m=4$, $\mh=m-s+1=2$ and $k=8$. In the clockwise direction, $v_{3,7}$ is the first vertex in the same partition set as $v_{3,1}$, and $v_{2,6}$ is the first vertex in the same partition set as $v_{2,7}$. Thus $j_4=1, j_3=7$ and $j_2=6$. The resulting orange and blue shaded parts correspond to $A_4^{(3)}(k_3,0,0,r_3+1)$ and $A_4^{(3)}(k_2,0,r_2+1,0)$, respectively, where $k_3=2$, $k_2=1$. $v_{1,5}$ and $v_{1,8}$ are the first two vertices in the same partition set as $v_{1,12}$ in the clockwise direction, thus $t_1=2, t_2=3$. The resulting yellow-circled part and the green shaded part correspond to $A_4^{(3)}(q_1,0,0,0)$ and $A_4^{(3)}(q_2,0,0,0)$, respectively, where $q_1=q_2=1$. The last cycle corresponds to $A_4^{(3)}(j_2-t_2-1,r_1+1,0,0)$ where $j_2-t_2-1 = 6-3-1 = 2$.}
    \label{fig:recur-grid-general}
\end{figure}

\comm{ 
\begin{figure}[hbt!]
    \centering
    \includegraphics[scale = 0.32]{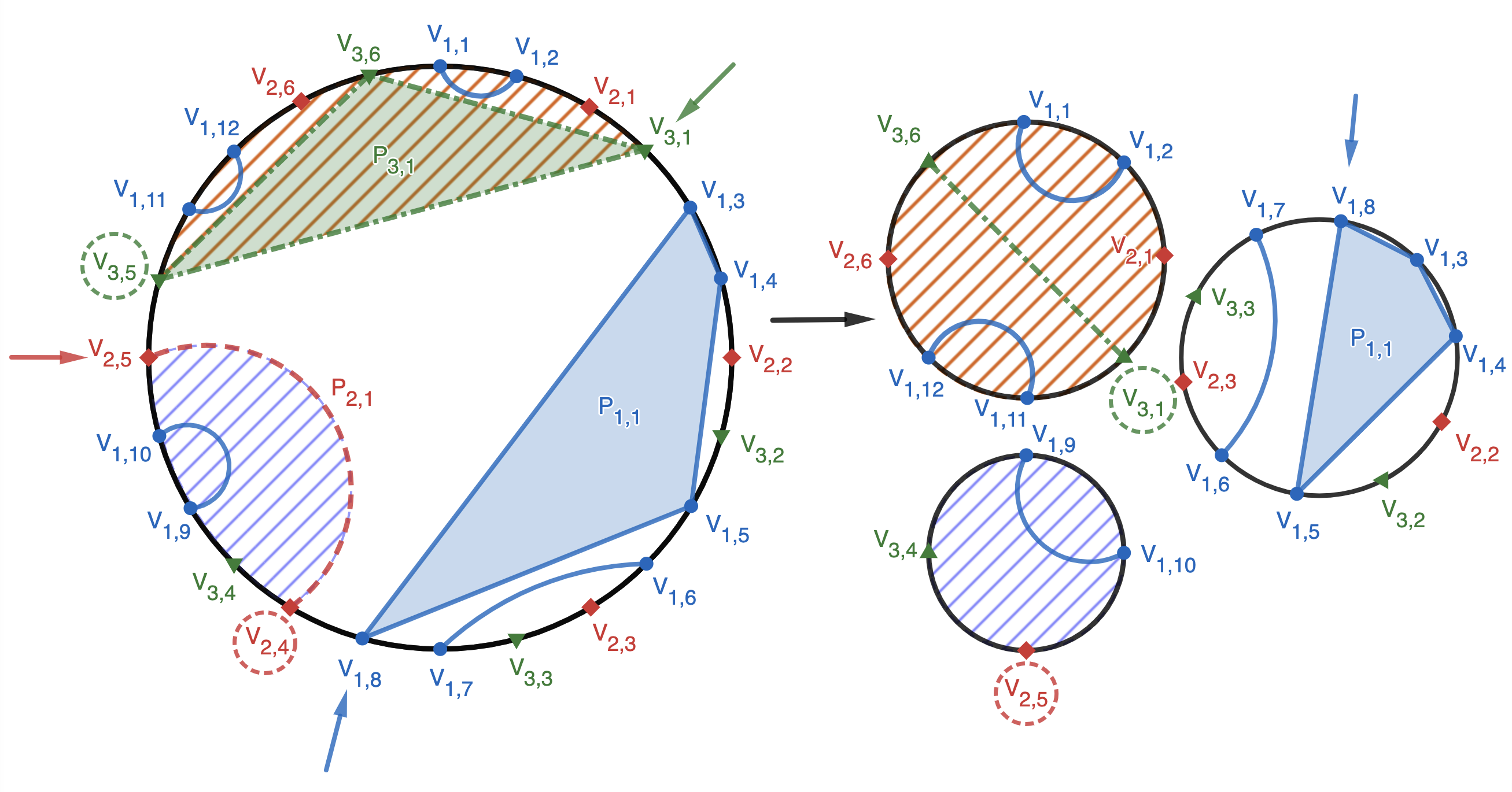}
    \caption{Illustration of Step 1-3 in the proof of \Cref{thm:recur-grid-general}: Here $s=3$, $m=4$, $\mh=m-s+1=2$ and $k=6$. In the clockwise direction, $v_{3,5}$ is the first vertex in the same partition set as $v_{3,1}$, and $v_{2,4}$ is the first vertex in the same partition set as $v_{2,5}$. The resulting orange and blue shaded parts correspond to $A_4^{(3)}(k_3,0,0,r_3+1)$ and $A_4^{(3)}(k_2,0,r_2+1,0)$, respectively, where $k_3=2$, $k_2=1$. The last cycle on the right is the remaining part described in Step 3, where $j_2=4$.}
    \label{fig:np-Amsr-split-1}
\end{figure}}
    
    \item Let $v_{1,t_1\mh+1}, v_{1,t_2\mh+2},\dots, v_{1,t_{\mh}\mh+\mh}$ be the first $\mh$ vertices in the clockwise direction that are in the same partition set as $v_{1,j_2\mh}$, for some $1\leq t_1\leq \dots \leq t_{\mh}$. Note that it takes the form $v_{1,t_{j}\mh+j}$ because the size of any partition set of the type $v_1$ is a multiple of $\mh$, and  we are only considering noncrossing partitions, thus the number of type $v_1$ vertices strictly in between $v_{1,t_{j}\mh+j}$ and $v_{1,t_{j+1}\mh+j+1}$ is a multiple of $\mh$.
    
    \item Let $t_0=j_2$. For each $i\in[\mh]$, identify $v_{1,t_{i-1}\mh+i-1}$ and $v_{1,t_{i}\mh+i}$. Everything in between $v_{1,t_{i-1}\mh+i-1}$ and $v_{1,t_{i}\mh+i}$ in the clockwise direction can be viewed as $\ams\(q_i,0,\dots, 0\)$, where $q_i= t_i-t_{i-1}$ for $i>1$ and $q_1 = t_1-1$. The last remaining part can be viewed as $\ams\(j_2-t_{\mh}-1, r_1+1,0,\dots,0\)$.

\comm{ 
\begin{figure}[hbt!]
    \centering
    \includegraphics[scale = 0.32]{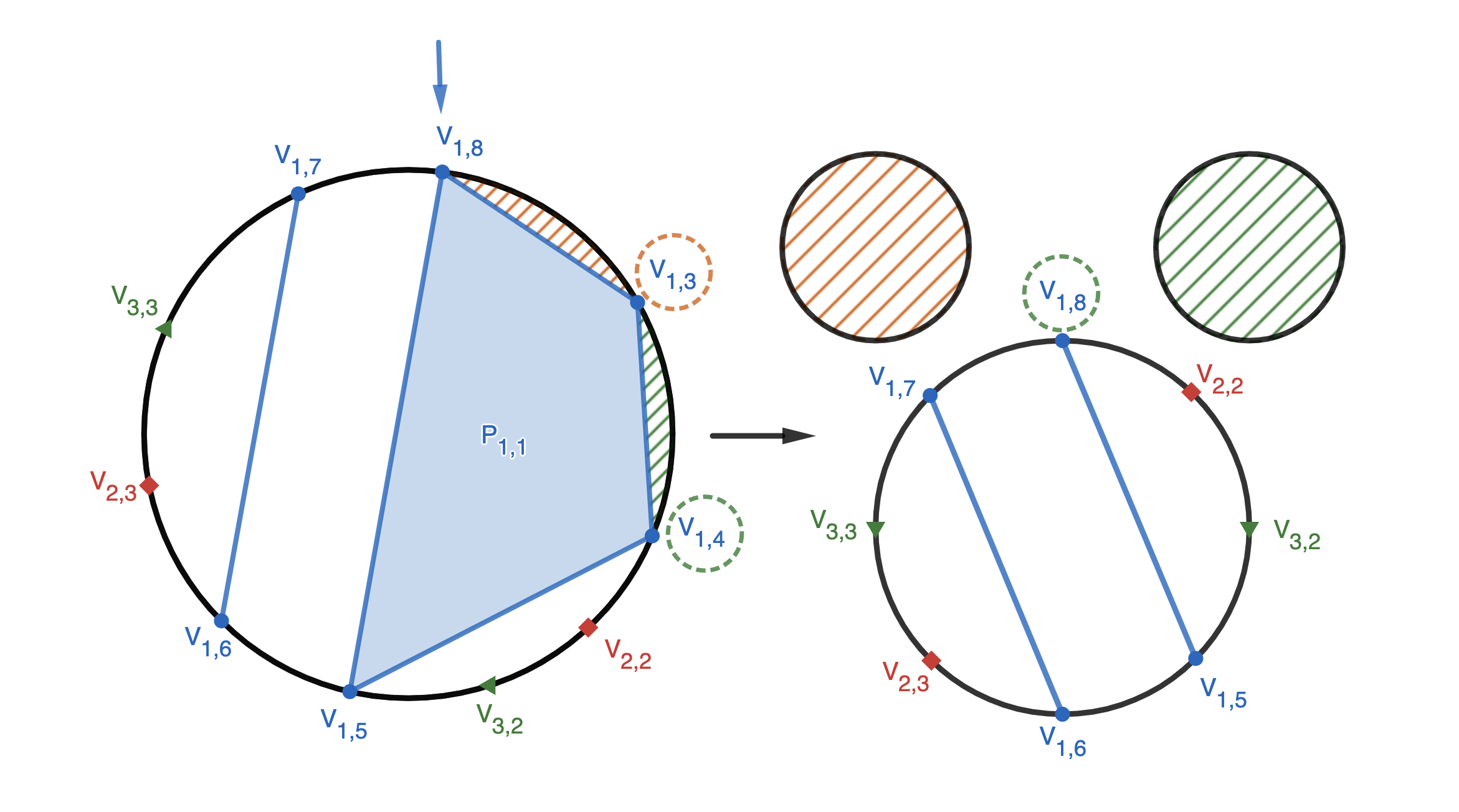}
    \caption{Illustration of Step 4-5 in the proof of \Cref{thm:recur-grid-general}: Here $v_{1,3}$ and $v_{1,4}$ are the first two vertices in the same partition set as $v_{1,8}$, thus $t_1=t_2=1$. The orange and green shaded parts correspond to $A_4^{(3)}(q_1,0,0,0)$ and $A_4^{(3)}(q_2,0,0,0)$, respectively, where $q_1=q_2=0$. The last cycle corresponds to $A_4^{(3)}(j_2-t_2-1,r_1+1,0,0)$ where $j_2-t_2-1 = 4-1-1 = 2$.}
    \label{fig:np-Amsr-split-2}
\end{figure}}
\end{enumerate}

Note that 
\begin{align*}
     & (k_2+\dots+k_s) + (q_1+\dots+q_{\mh}) + (j_2-t_{\mh} - 1) \\
    =& \(\(k+1-j_s\) + \(j_s-j_{s-1}\) +\dots + \(j_3-j_2\)\) \\ &\hspace{2cm} + \((t_{\mh}-t_{\mh-1}) + (t_{\mh-1} - t_{\mh-2}) + \dots + (t_2-t_1) + (t_1-1)\) + (j_2-t_{\mh} - 1) \\
    =& \(k+1 - j_2\) + \(t_{\mh}-1\) + (j_2-t_{\mh} - 1) = k-1.
\end{align*}

Combining everything, we have

\begin{align*}
     & A_m^{(s)}\(k,r_1,\dots,r_s\) \\
    =& \sum_{0\leq j_2\leq \dots\leq j_s}\,\sum_{0\leq t_1\leq \dots\leq t_{\mh}}\, \(\prod_{i=2}^s A_m^{(s)}\(k_i,0,\dots,r_i+1,\dots, 0\)\)\cdot \(\prod_{i=1}^{\mh} A_m^{(s)}\(q_i,0,\dots, 0\)\)\cdot \\
     &\hspace{1.5 cm} A_m^{(s)}\(j_2-t_{\mh}-1, r_1+1,0,\dots,0\) \text{ where } k_i= j_{i+1}-j_i \mod k \text{ and } q_i= t_i-t_{i-1}, q_1=t_1-1 \\
    =& \sum_{\substack{i_1,\dots,i_{m+1}\geq 0: \\ i_1+\dots+i_{m+1} = k-1}}\, \(\prod_{x=1}^s A_m^{(s)}\(i_x,0,\dots,r_x+1,\dots, 0\)\)\cdot
    \(\prod_{y=s+1}^{s+\mh} A_m^{(s)}\(i_y,0,\dots, 0\)\) \\
    =& \sum_{\substack{i_1,\dots,i_{m+1}\geq 0: \\ i_1+\dots+i_{m+1} = k-1}}\, \(\prod_{x=1}^s A_m^{(s)}\(i_x,0,\dots,r_x+1,\dots, 0\)\)\cdot
    \(\prod_{y=s+1}^{m+1} A_m^{(s)}\(i_y,0,\dots, 0\)\)
\end{align*}
as needed.

\end{proof}

\subsection{Recurrence Relations for the Trace Power Moments of  \texorpdfstring{$\mzshapesdistrg$}{Mz(m),s}}

\setlength{\parskip}{1.5mm}
\setlength{\baselineskip}{1.3em}

\begin{defn}\label{defn:Bms}
    Given distributions $\o^{(1)},\dots, \o^{(s)}$, we define $B_m^{(s)}\(k,0,\dots,0\)$ to be 
    \begin{equation}
        \bms\(k,0,\dots,0\) = W\(\D_{\mzshape,\oams,2k}\).
    \end{equation}
\end{defn}

\begin{defn}
    We define $\bms\(k,r_1,\dots,r_s\)$ to be the weight of dominant constraint graphs of length $k$ associated with $\mzshape$ and $\oams$ where the first spoke on the $i^{th}$ layer has an multi-edge of multiplicity $(2r_i+1)$ for each $i\in[s]$.
\end{defn}

\begin{rmk}
    When $r_1=\dots=r_s=0$,  $\bms\(k,r_1,\dots,r_s\)$ coincides with  $\bms\(k,0,\dots,0\)$.
\end{rmk}

We will prove the following main result.
\begin{restatable}{thm}{recurGraph}
\label{thm:recur-graph-general}
    \begin{equation}
    \begin{aligned}
        \bms\(k,r_1,\dots,r_s\) 
        = \sum_{\substack{i_1,\dots,i_{m+1}\geq 0: \\ i_1+\dots+i_{m+1} = k-1}}\, & \bms\(i_1,r_1+1,0,\dots,0\)\dots \bms\(i_s,0,\dots,0,r_s+1\)\cdot\\
        & \bms\(i_{s+1},0,\dots,0\)\dots \bms\(i_{m+1},0,\dots,0\)
    \end{aligned}
    \end{equation}
\end{restatable}


We list below some of the definitions and properties of dominant constraint graphs from \cite{CP20}, which we will need to prove the main result. 

\begin{defn}[Well-behaved]\label{defn:well-behaved}
    Given a shape $\a$, we say that a constraint graph $C \in \mathcal{C}_{(\a,2q)}$ is \emph{well-behaved} if whenever $u \= v$ in $C$, $u$ and $v$ are copies of the same vertex in $\a$ or $\a^T$.
\end{defn}

\begin{thm}\label{thm:mzshape-split}
    If $C$ is a dominant constraint graph in $\C_{\mzshape,\oams, 2q}$, then the following statements are true:
    \begin{enumerate}
        \item If $a_{i,1}\=a_{j,1}$ for some $1\leq i<j\leq q$, then $a_{i,k}\=a_{j,k}$ for all $k\in[m]$. Similarly, if $b_{i,m}\=a_{j,m}$ for some $1\leq i<j\leq q$, then $b_{i,k}\=b_{j,k}$ for all $k\in[m]$. More generally,
        \begin{enumerate}[i.]
            \item if $a_{i,t}\=a_{j,t}$, for some $t\in[m]$, then $a_{i,l}\=a_{j,l}$ for all $t\leq l\leq m$.
            \item if $b_{i,t}\=b_{j,t}$ for some $t\in[m]$, then $b_{i,l}\=b_{j,l}$ for all $1\leq l\leq t$.
        \end{enumerate}
        \item For a given wheel $V_t$, we relabel its vertices as $v_1,\dots,v_{2q}$, starting from $a_{1,t}$. If $v_j$ is the first vertex $v_1$ is constrained to (i.e. if $j$ is the smallest index such that $v_1\=v_j$), then $v_2\=v_{j-1}$. More generally,
        \begin{itemize}
            \item for any given $i\in[2q]$, if $v_j$ is the first vertex in the \emph{clockwise} direction such that $v_j\= v_i$, then $v_{i+1}\= v_{j-1}$, where $i+1$ and $j-1$ are taking $\mod 2q$. e.g. if $v_6$ is the first vertex constrained to $v_{2q}$ in the clockwise direction, then $v_{2q+1} = v_1\= v_5$.
        \end{itemize}
    \end{enumerate}
\end{thm}

The proof idea of the first part of this theorem is that in dominant constraint graphs, constrained spokes on the same layer should not cross each other. More precisely, if spokes $e_1,e_2,e_3,e_4$ are ordered in the clockwise direction on the same layer, and $e_1 \= e_3$, $e_2\= e_4$, then all four of them should be constrained together. 

\begin{defn}[Line Shape]\label{defn:line-shape}
    Let $\a_{0}$ be the bipartite shape with vertices $V(\a_{0})=\{u,v\}$ and a single edge $\{u,v\}$ with distinguished tuples of vertices $U_{\a_{0}}=(u)$ and $V_{\a_{0}}=(v)$. We call $\a_0$ the \textit{line shape}.
\end{defn}

\begin{defn}\label{defn:non-crossing-line-shape}
    Let $\a_0$ be the line shape as in \Cref{defn:line-shape}. Let $H(\a_0,2q)$ be the multi-graph as in \Cref{defn:copies}. We label the vertices of $H(\a_0,2q)$ as $\left\{i_j:j\in[2q]\right\}$.
    
    We say a constraint graph $C \in \C_{(\a_0,2q)}$ is \textit{non-crossing} if it is possible to draw the constraint edges (within the circle containing the vertices) so that they do not cross. Equivalently, $C$ is non-crossing if whenever $i_{x}\=i_y$ and $i_{v}\=i_w$ for some $x\leq v\leq y\leq w$, we have that $i_{x}\=i_y\=i_{v}\=i_w$.
    
    
\end{defn}

\begin{defn}\label{defn:mzshape-induced-constraint-graph}
   Let $\mzshape$ be the multi-layer Z-shape. Let $C$ be a constraint graph on $H(\mzshape,2q)$. For $i\in[m]$, we denote $C_i$ the induced subgraph of $C$ on vertices $V_i$. We call $C_i$ the \textit{induced constraint graph} of $C$ on $V_i$. 
\end{defn}

\begin{defn}\label{defn:well-defined-mzshape-constraint-graph}
    Let $C$ be a constraint graph in $\C_{\mzshape,\oams, 2q}$. We say $C$ is \emph{non-crossing} if the induced constraint graphs $C_i$'s are non-crossing.
\end{defn}

\begin{thm}\label{thm:mzshape-dominant-properties}
    If $C$ is a dominant constraint graph in $\C_{\mzshape,\oams, 2q}$, then
    \begin{enumerate}
        \item $C$ is well-behaved and non-crossing.
        \item The induced constraint graphs $C_i$ on wheels $W_i$ are dominant constraint graphs in $\C_{\(\a_0,2q\)}$.
    \end{enumerate}
\end{thm}

\begin{thm}\label{thm:line-shape-cross}
    Let $C \in \C_{(\a_0,2q)}$ be a dominant constraint graph. Then $C$ is well-behaved and non-crossing.
\end{thm}

\subsubsection{Base Case s=1, m=2}

\setlength{\parskip}{1.5mm}
\setlength{\baselineskip}{1.3em}

\begin{defn}
    Given $\o$, we will denote $\B_{k,r}$ to be the set of dominant constraint graphs of length $k$ associated with $\zshape$ and $\o_{\zshape,1}$ where the first spoke on the first layer has an multi-edge of multiplicity $(2r+1)$. 
    Note that $\B_{k,0} = \D_{\zshape,\oz,2k}$.
\end{defn}

When $s=1$, for simplicity denote $\o^{(1)}$ as $\o$ and  $B_m^{(1)}(k,r)$ as $B_m(k,r)$. When $s=1$ and $m=2$, 
\begin{equation}
    B_2(k,r) = W\(\B_{k,r}\)\,.
\end{equation}


\begin{thm}\label{thm:recur-graph-zshape}
\begin{equation}
    B_2\(k,r\) 
    = \sum_{\substack{i_1,i_2,i_3\geq 0: \\ i_1+i_2+i_3 = k-1}}\, B_2\(i_1, 0\)\cdot B_2\(i_2, 0\)\cdot B_2\(i_3, r+1\)
\end{equation}
and 
\begin{equation}
    B_2(0,r) = \o_{2r}\,.
\end{equation}
\end{thm}

\begin{proof}

\begin{figure}[hbt!]
    \centering
    \begin{subfigure}[t]{.48\textwidth}
        \includegraphics[width=1\linewidth]{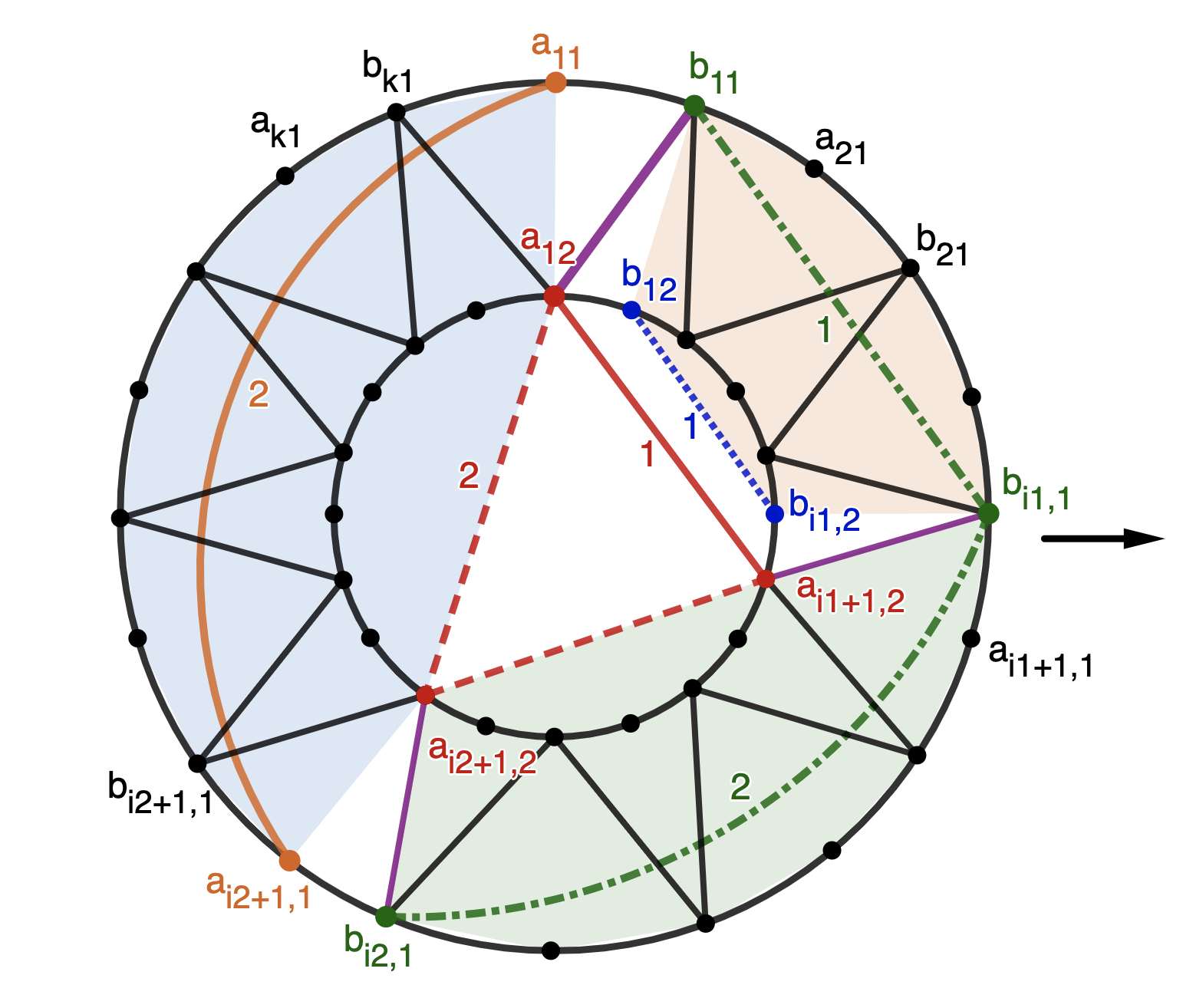}
    \end{subfigure}
    \begin{subfigure}[t]{.48\textwidth}
        \includegraphics[width=1\linewidth]{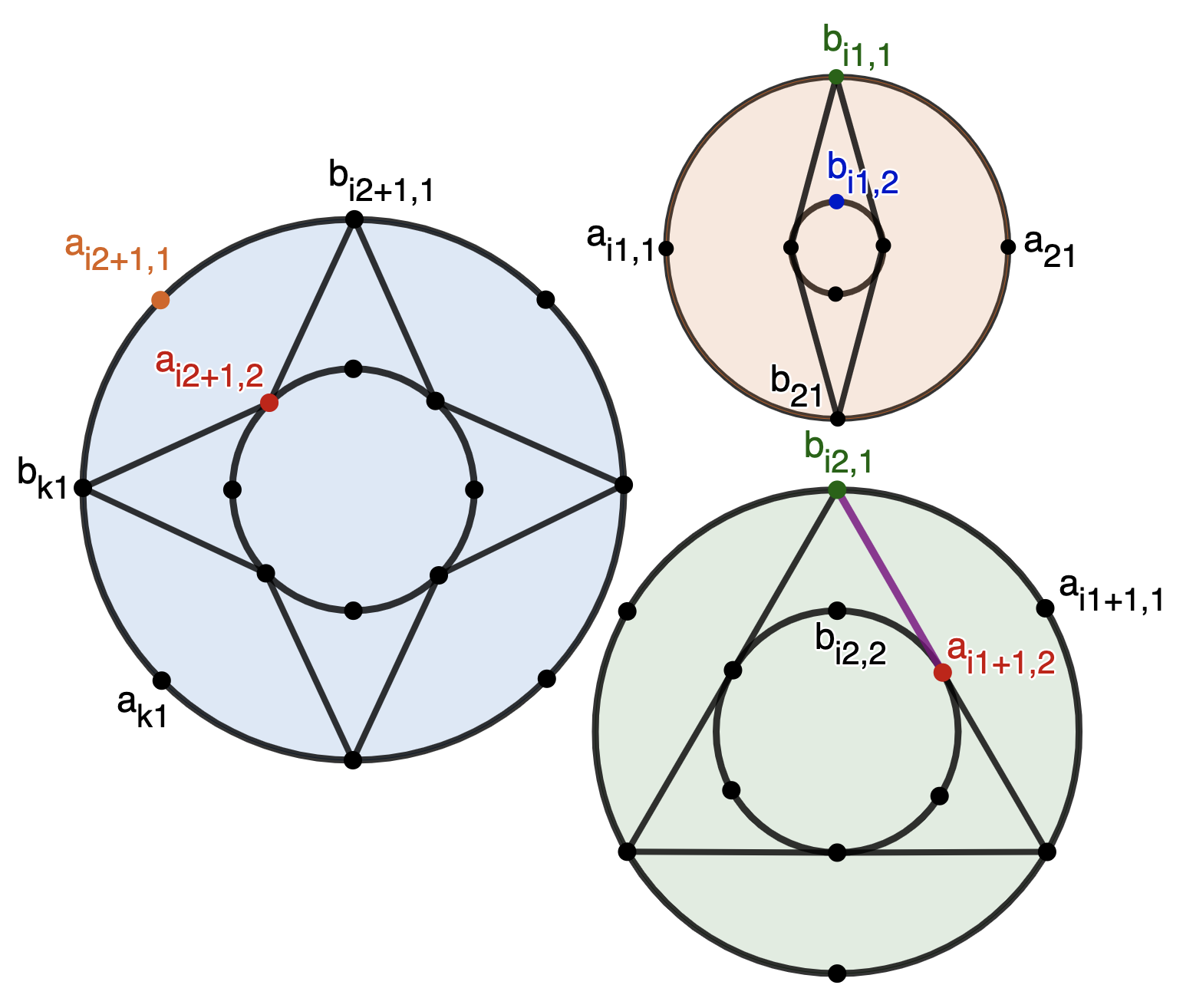}
    \end{subfigure}
    \caption{Illustration of the proof of \Cref{thm:recur-graph-zshape}: The spoke $\{a_{12},b_{11}\}$, colored purple, has multiplicity $2r+1$. For constraint edges, solid constraint edges imply dashed ones with the same numbering. For example, for the constraint edges numbered 1, $a_{1,2}\longleftrightarrow a_{i_1+1,2} \implies b_{1,2}\longleftrightarrow b_{i_1,2} \implies b_{1,1}\longleftrightarrow b_{i_1,1}$.}
    \label{fig:recur-zshape-split}
\end{figure}

We will first prove that there is a bijection between $\B_{k,r}$ and $\displaystyle \bigcup_{\substack{ k_i\geq 0: \\ k_1+k_2+k_3= k-1}}\, \B_{k_1,0}\times \B_{k_2,0} \times \B_{k_3,r+1}$.

\begin{enumerate}
    \item $\B_{k,r} \to \displaystyle \bigcup_{\substack{ k_i\geq 0: \\ k_1+k_2+k_3= k-1}}\, \B_{k_1,0}\times \B_{k_2,0} \times \B_{k_3,r+1}$: Let $C\in \B_{k,r}$, with the spoke $\{a_{1,2}, b_{1,1}\}$ having multiplicity $2r+1$.  We split $C$ into three parts with the following steps. 
    
    \begin{enumerate}[i.]
        \item Let $i_1 \in [k]$ be the first index such that $a_{1,2} \= a_{i_1+1,2}$. In the case of $a_{1,2}$ being isolated, $i_1$ is chosen to be $k$. This implies that $b_{1,2} \= b_{i_1,2}$, which further implies that $b_{1,1}\=b_{i_1,1}$.
    
        \item Let $i_2\in [k]$ be the first index such that $a_{1,1}\= a_{i_2+1,1}$. In the case of $a_{1,1}$ being isolated, $i_2$ is chosen to be $k$. This implies that $b_{1,1}\= b_{i_2,1}$ and $a_{1,2}\= a_{i_2+1,2}$. 
    
        \item Observe that $i_2\geq i_1$: otherwise if $i_2<i_1$, $a_{1,1}\= a_{i_2+1,1}$ and $b_{1,1}\= b_{i_1,1}$ would imply $a_{1,1}\= a_{i_2+1,1} \= b_{1,1} \= b_{i_1,1}$ since $C$ is dominant and they cross. But $a_{1,1}\= b_{1,1}$ contradicts that $C$ is parity preserving.
    
        \item Since $a_{1,2} \= a_{i_1+1,2}$ and $a_{1,2}\= a_{i_2+1,2}$, $a_{i_1+1,2}\= a_{i_2+1,2}$. Since $b_{1,1}\=b_{i_1,1}$ and $b_{1,1}\= b_{i_2,1}$, $b_{i_1,1}\= b_{i_2,1}$. 
    
        \item Contracting $b_{1,1}$ with $b_{i_1,1}$ and $b_{1,2}$ with $b_{i_1,2}$ give $H\(\a_Z,2(i_1-1)\)$. The induced constraint graph $C_1 \in \C_{\(\a_Z,2(i_1-1)\)}$ is dominant, thus $C_1\in \B_{i_1-1,0}$. 
    
        \item Contracting $a_{i_1+1,2}$ with $a_{i_2+1,2}$ and $b_{i_1,1}$ with $b_{i_2,1}$ give $H\(\a_Z,2(i_2-i_1)\)$. Moreover, the spokes $\{a_{1,2}, b_{1,1}\}$, $\{a_{i_1+1,2}, b_{i_1,1}\}$ and $\{a_{i_2+1,2}, b_{i_2,1}\}$ are identified together, resulting in multiplicity $(2r+1)+2 = 2(r+1)+1$. thus the induced constraint graph $C_2$ is in $\B_{k,r+1}$.
    
        \item Contracting $a_{1,1}$ with $a_{i_2+1,1}$ and $a_{1,2}$ with $a_{i_2+1,2}$ give $H\(\a_Z,2(k-i_2)\)$. The induced constraint graph $C_3 \in \C_{\(\a_Z,2(k-i_2)\)}$ is dominant, thus $C_3 \in \B_{k-i_2,0}$.
\end{enumerate}
    
    Since $(i_1-1)+(i_2-i_1)+(k-i_2) = k-1$, we conclude that $(C_1,C_2,C_3)\in \displaystyle \bigcup_{\substack{ k_i\geq 0: \\ k_1+k_2+k_3= k-1}}\, \B_{k_1,0}\times \B_{k_2,0} \times \B_{k_3,r+1}$.
    
    \item $\displaystyle \bigcup_{\substack{ k_i\geq 0: \\ k_1+k_2+k_3= k-1}}\, \B_{k_1,0}\times \B_{k_2,0} \times \B_{k_3,r+1} \to \B_{k,r}$: Conversely, given $C_1\in \B_{k_1,0}, C_2\in \B_{k_2,0}$ and $C_3\in \B_{k_3,r+1}$ for some $k_1+k_2+k_3$, we can reverse the above steps and glue them together to get $C\in \B_{k,r}$.
\end{enumerate}

    Thus, 
    \begin{align*}
    B(k,r) 
    &= W\(\B_{k,r}\)
     = W\(\bigcup_{\substack{ k_i\geq 0:\, k_1+k_2+k_3 = k-1}} \B_{k_1,0}\times \B_{k_2,0} \times \B_{k_3,r+1} \) \\
    &= \sum_{k_i\geq 0:\, k_1+k_2+k_3 = k-1}\, W\(\B_{k_1,0}\times \B_{k_2,0} \times \B_{k_3,r+1}\) \\
    &= \sum_{k_i\geq 0:\, k_1+k_2+k_3 = k-1}\, W\(\B_{k_1,0}\)\cdot W\(\B_{k_2,0}\)\cdot W\(\B_{k_3,r+1}\) \\
    &= \sum_{k_i\geq 0:\, k_1+k_2+k_3 = k-1}\, B(k_1,0)B(k_2,0)B(k_3,r+1).
    \end{align*}
    as needed.
\end{proof}


\subsubsection{General Case}

\setlength{\parskip}{1.5mm}
\setlength{\baselineskip}{1.3em}

\begin{defn}
    Given a distribution $\o^{(1)},\dots, \o^{(s)}$ and $t\leq s$, we will denote $\B_{m,k,r_1,\dots,r_t}$ to be the set of dominant constraint graphs of length $k$ associated with $\mzshape$ and $\oams$ where the first spoke on the $j^{th}$ layer has an multi-edge of multiplicity $(2r_j+1)$ for each $j\in[t]$. Note that $\B_{m,k,r_1,\dots,r_s} = \D_{\mzshape,\oams,2k}$ and $\B_{m,k} = \D_{\mzshape, 2k}$.
        
    
\end{defn}

\begin{defn}
    For $t\in[s]$, we will use $\B_{m,k,r\ve_t}$ to denote $\B_{m,k,r_1,\dots,r_s}$ where $r_t=t$ and $r_j=0$ for all $j\neq t$.
\end{defn}

\begin{obs}
    $\bms(k,r_1,\dots,r_s) = W\(\B_{m,k,r_1,\dots,r_s}\)$, $\bms\(k,0,\dots,r_j,\dots,0\) = W\(\B_{m,k,r_j\ve_j}\)$ and $\bms(k,0,\dots,0) = W\(\B_{m,k}\)$.
\end{obs}

Now we are ready to prove \Cref{thm:recur-graph-general}.
\recurGraph*

\begin{proof}

\begin{figure}[hbt!]
    \centering
    \begin{subfigure}[t]{.48\textwidth}
        \includegraphics[width=1\linewidth]{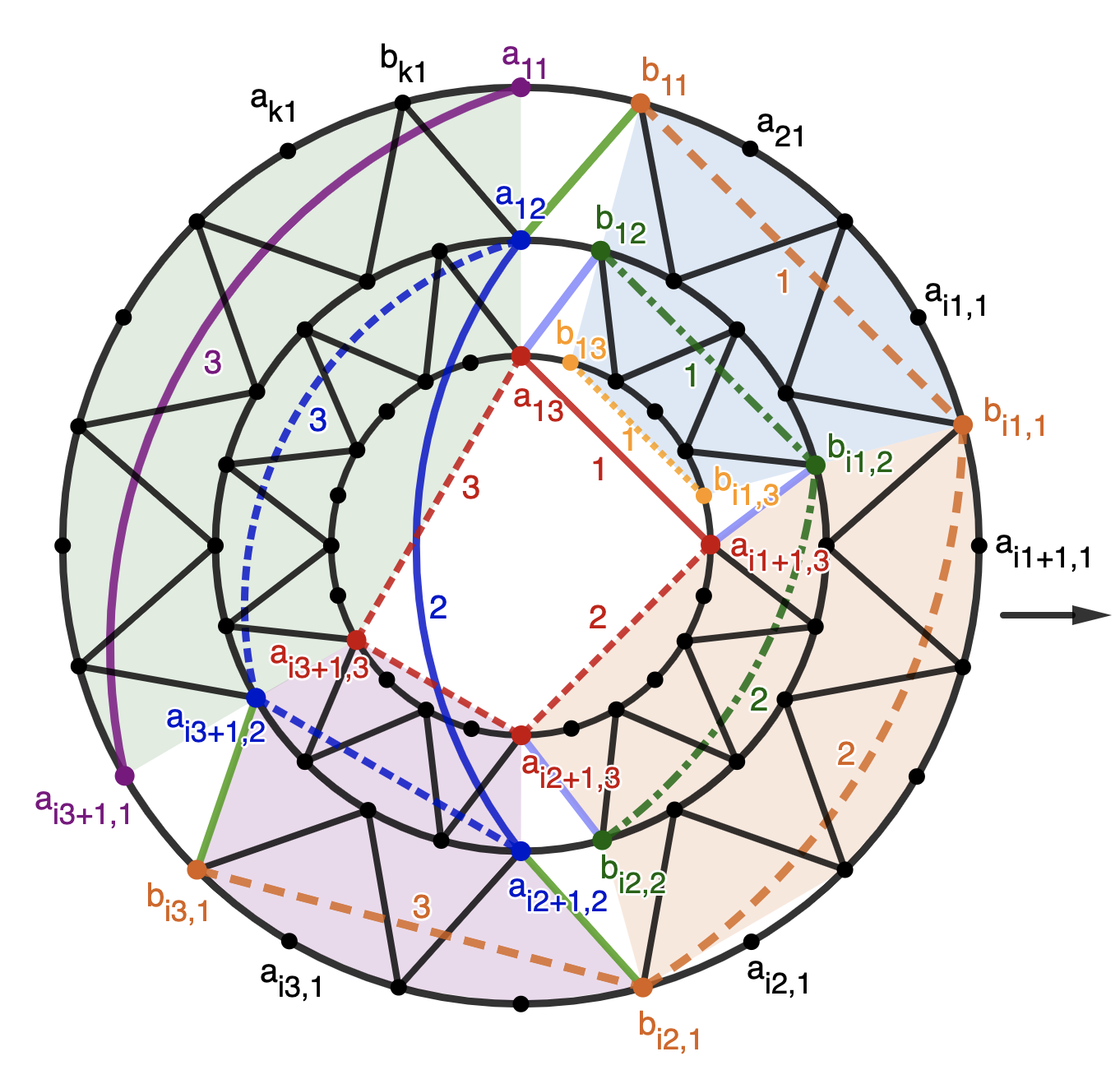}
    \end{subfigure}
    \begin{subfigure}[t]{.45\textwidth}
        \includegraphics[width=1\linewidth]{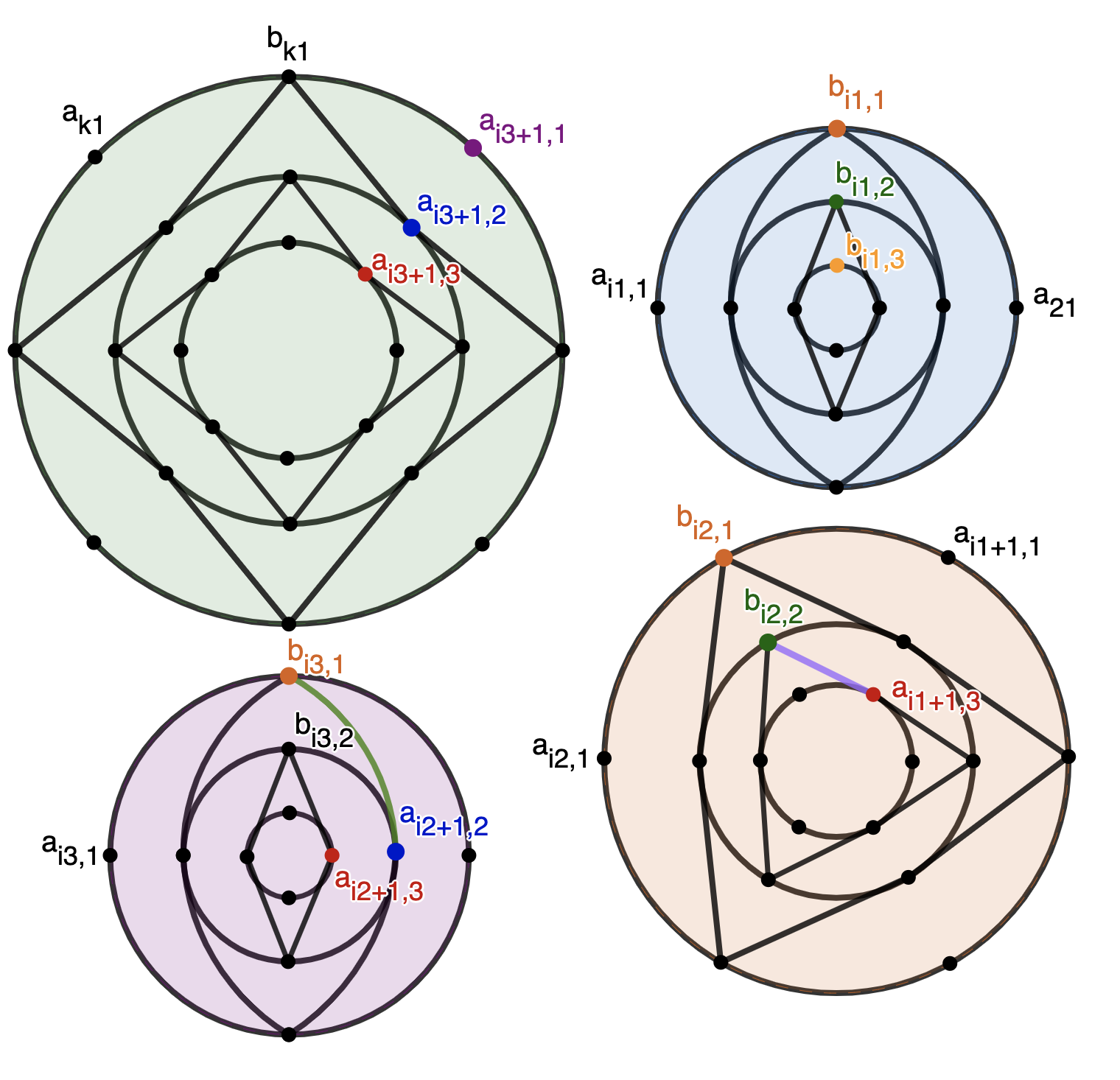}
    \end{subfigure}
    \caption{Illustration of \Cref{thm:recur-graph-general}. For constraint edges, solid constraint edges imply dashed ones with the same numbering. }
    \label{fig:recur-mzhape-split}
\end{figure}

Similar to the proof of \Cref{thm:recur-graph-zshape}, we now will prove that there is a bijection between $\B_{m,k,r_1,\dots,r_s}$ and $\displaystyle \bigcup_{\substack{ k_i\geq 0: \\ k_1+\dots+k_{m+1}= k-1}}\, \B_{m,k_1,(r_1+1)\ve_1}\times \dots\times \B_{m,k_s,(r_s+1)\ve_s} \times \B_{m,k_{s+1}}\times\dots \B_{m,k_{m+1}}$.

\begin{itemize}
    \item $\displaystyle \B_{m,k,r_1,\dots,r_s} \to  \bigcup_{\substack{ k_i\geq 0: \\ k_1+\dots+k_{m+1}= k-1}}\, \B_{m,k_1,(r_1+1)\ve_1}\times \dots\times \B_{m,k_s,(r_s+1)\ve_s} \times \B_{m,k_{s+1}}\times\dots \B_{m,k_{m+1}}$: 
    
    Let $C\in \B_{m,k,r_1,\dots,r_s}$, with the spokes $\{a_{1,j+1}, b_{1,j}\}$ having multiplicity $2r_j+1$ for $j\in[s]$.  We split $C$ into $m+1$ parts with the following steps. 
    
    \begin{enumerate}
        \item Let $i_1 \in [k]$ be the first index such that $a_{1,m} \= a_{i_1+1,m}$. In the case of $a_{1,m}$ being isolated, $i_1$ is chosen to be $k$. This implies that $b_{1,m} \= b_{i_1,m}$, which further implies that $b_{1,j}\=b_{i_1,j}$ for all $j\in[m]$.
    
        \item Let $i_2\in [k]$ be the first index such that $a_{1,m-1}\= a_{i_2+1,m-1}$. In the case of $a_{1,m-1}$ being isolated, $i_2$ is chosen to be $k$. This implies that 
        \begin{enumerate}[i.]
            \item $b_{1,m-1}\= b_{i_2,m-1}$, which further implies that $b_{1,j}\= b_{i_2+1,j}$ for all $j\in[m-1]$. 
            \item $a_{1,m} \= a_{i_2+1, m}$.
        \end{enumerate}
        
        \item Continuing this way, at the $t^{th}$ step where $t\in[m]$, let $i_{t}\in [k]$ be the first index such that $\a_{1,m-t+1} \= a_{i_t+1, m-t+1}$. In the case of $a_{1,m-t+1}$ being isolated, $i_t$ is chosen to be $k$. This implies that 
        \begin{enumerate}[i.]
            \item $b_{1,m-t+1}\= b_{i_t,m-t+1}$, which further implies that $b_{1,j}\= b_{i_t+1,j}$ for all $j\in[m-t+1]$. 
            \item $a_{1,j} \= a_{i_t+1,j}$ for all $j\in \{m-t+1, \dots, m\}$.
        \end{enumerate}
    
        \item Observe that $i_{j+1}\geq i_j$ for all $j\in[m-1]$: otherwise if $i_{j+1}<i_j$ for some $j$, $a_{1,m-j}\= a_{i_{j+1}+1,m-j}$ and $b_{1,m-j}\= b_{i_j,m-j}$ would imply $a_{1,m-j}\= a_{i_{j+1}+1,m-j} \= b_{1,m-j} \= b_{i_j,m-j}$ since $C$ is dominant and they cross. But $a_{1,m-j}\= b_{1,m-j}$ contradicts that $C$ is parity preserving.
    
        \item Let $t\in[m]$. Since $a_{1,j} \= a_{i_t+1,j}$ for all $j\in \{m-t+1,\dots, m\}$, for each $j\in [m]$, $a_{1,j}\= a_{i_{m-j+1}+1,j}\= \dots \= a_{i_m+1,j}$. Since $b_{1,j}\=b_{i_t,j}$ for all $j\in[m-t+1]$, for each $j\in[m]$, $b_{1,j}\= b_{i_1,j}\= \dots\= b_{i_{m-j+1},j}$. 
    
        \item Contracting $b_{1,j}$ with $b_{i_1,j}$ for all $j\in[m]$ gives $H\(\mzshape,2(i_1-1)\)$. The induced constraint graph $C_0 \in \C_{\(\a_Z,2(i_1-1)\)}$ is dominant, thus $C_0\in \B_{m,i_1-1}$. 
        
        \item Let $t\in[m-1]$. Contracting $b_{i_t,j}$ with $b_{i_{t+1},j}$ for all $j\in[m-t]$ and contracting $a_{i_t+1,j}$ with $a_{i_{t+1}+1,j}$ for all $j\in \{m-t+1,\dots, m\}$ give $H\(\a_Z,2(i_{t+1}-i_t)\)$. Moreover, the spokes $\{a_{1,m-t+1}, b_{1,m-t}\}$, $\{a_{i_t+1,m-t+1}, b_{i_t,m-t}\}$ and $\{a_{i_{t+1}+1,m-t+1}, b_{i_{t+1},m-t}\}$ are identified together, resulting in multiplicity $(2r_{m-t+1}+1)+2 = 2(r_{m-t+1}+1)+1$ if $m-s\leq t\leq m-1$. Thus the induced constraint graph $C_t$ is in $\B_{m,i_{t+1}-i_t}$ if $t\in [m-s-1]$ and in $\B_{m,i_{t+1}-i_t,(r_{m-t}+1)\ve_{m-t}}$ if $m-s\leq t\leq m-1$.
    
        \item Finally, contracting $a_{1,j}$ with $a_{i_m+1,j}$ for all $j\in[m]$ gives $H\(\a_Z,2(k-i_m)\)$. The induced constraint graph $C_m \in \C_{\(\a_Z,2(k-i_m)\)}$ is dominant, thus $C_m \in \B_{m,k-i_m}$.
\end{enumerate}
    
    Since $\displaystyle (i_1-1)+ \sum_{j=1}^{m-1} \(i_{j+1}-i_j\) +(k-i_m) = k-1$, and there are $s$ many $C_j$'s (i.e. $C_{m-s}, \dots, C_{m-1}$) belonging to $\B_{m,k_1, (r_1+1)\ve_1},\dots, \B_{m,k_s, (r_s+1)\ve_s}$ for $k_j= i_{m-j+1}-i_{m-j}$, we conclude that
    \begin{equation*}
        (C_0,C_1,\dots, C_m)\in \bigcup_{\substack{ k_i\geq 0: \\ k_1+\dots +k_{m+1}= k-1}}\, \B_{m, k_1, (r_1+1)\ve_1}\times\dots \times \B_{m, k_s, (r_s+1)\ve_s} \times \B_{m,k_{s+1}}\times \dots\times \B_{m,k_{m+1}}\,.
    \end{equation*}
    
    \item $\displaystyle \bigcup_{\substack{ k_i\geq 0: \\ k_1+\dots+k_{m+1}= k-1}}\, \B_{m,k_1,(r_1+1)\ve_1}\times \dots\times \B_{m,k_s,(r_s+1)\ve_s} \times \B_{m,k_{s+1}}\times\dots \B_{m,k_{m+1}} \to \B_{m,k,r_1,\dots,r_s}$: 
    
    Conversely, given $C_i\in \B_{m,k_i,(r_i+1)\ve_j}$ for $i\in[s]$ and $C_j\in \B_{m,k_j}$ for some $j\in\{s+1,\dots, m+1\}$, we can reverse the above steps and glue them together to get $C\in \B_{m,k,r_1,\dots,r_s}$.
\end{itemize}

    Thus 
    \begin{align*}
     & B_m(k,r_1,\dots, r_s) 
    =  W\(\B_{m,k,r_1,\dots, r_s}\) \\
    =& W\(\bigcup_{\substack{ k_i\geq 0: \\ k_1+\dots+k_{m+1}= k-1}}\, \B_{m,k_1,(r_1+1)\ve_1}\times \dots\times \B_{m,k_s,(r_s+1)\ve_s} \times \B_{m,k_{s+1}}\times\dots \B_{m,k_{m+1}} \) \\
    =& \sum_{\substack{ k_i\geq 0:\\ k_1+\dots+k_{m+1} = k-1}}\, W\(\B_{m,k_1,(r_1+1)\ve_1}\times \dots\times \B_{m,k_s,(r_s+1)\ve_s} \times \B_{m,k_{s+1}}\times\dots \B_{m,k_{m+1}}\) \\
    =& \sum_{\substack{k_i\geq 0:\,\\ k_1+\dots+k_{m+1} = k-1}}\, \prod_{i=1}^{s}\, W\(\B_{m,k_i,\(r_i+1\)\ve_i}\)\cdot \prod_{j=s+1}^{m+1} W\(\B_{m,k_j}\) \\
    =& \sum_{\substack{k_i\geq 0:\,\\ k_1+\dots+k_{m+1} = k-1}}\, \(\prod_{i=1}^{s}\, B_m(m,k_i,0,\dots, r_i+1,\dots, 0)\)\cdot \(\prod_{j=s+1}^{m+1} B_m(m, k_j,0,\dots,0)\)
    \end{align*}
    as needed.

\end{proof}

\begin{proof}[Proof of \Cref{thm:main-2}]
    This follows directly from \Cref{thm:recur-grid-general} and \Cref{thm:recur-graph-general}.
\end{proof}


\bibliographystyle{plainnat}
\bibliography{main}

\appendix

\section{Analyzing Inner Products with Random Vectors}

\setlength{\parskip}{1.5mm}
\setlength{\baselineskip}{1.3em}

\begin{defn}
    A \emph{perfect matching of $[n]$} is $\displaystyle \M = \left\{\{a_i,b_i\}: i\in[n/2], a_i, b_i\in [n], \bigsqcup_{i=1}^{n/2} \{a_i\}\sqcup \{b_i\} = [n] \right\}$\,.
\end{defn}

We use the following calculation from Jones and Potechin \cite{JP21} which is a spherical analogue of Isserlis' Theorem/Wick's Theorem.
\begin{thm}
    For any vectors $\vd_1,\dots,\vd_k \in \Rbb^n$, 
    \begin{equation}
        \Ebb_{\vvv\in S^{n-1}}\[\;\prod_{i=1}^{k}\<\vvv,\vd_i\>\;\] = \(\prod_{j=1}^{k/2}\; n+2j-2\)^{-1}\cdot \sum_{\substack{\M\text{ perfect } \\ \text{matchings of } [k] }}  \prod_{e=(i,j)\in \M} \<\vd_i,\vd_j\>.
    \end{equation}
\end{thm}

\begin{prop}
    The number of perfect matchings between $2k$ elements is $\displaystyle\(\binom{2k}{2}\binom{2k-2}{2}\dots \binom{2}{2}\)/k! = \dfrac{(2k)!}{2^{k}\cdot k!}=(2k-1)\cdot(2k-3)\dots 3\cdot 1$.
\end{prop}

\begin{eg}
Assume $k$ is even, then
\begin{enumerate}
    \item $\displaystyle\Ebb_{\vvv\in S^{n-1}} \<\vvv,\vd_1\>\<\vvv,\vd_2\> = \dfrac{\,1\,}{n}\cdot\<\vd_1,\vd_2\>$.
    \item $\displaystyle\Ebb_{\vvv\in S^{n-1}} \<\vvv,\vd_1\>^2\<\vvv,\vd_2\>^2 = \dfrac{1}{n(n+2)}\(2\<\vd_1,\vd_2\>^2+\<\vd_1,\vd_1\>\<\vd_2,\vd_2\>\)$.
    \item $\displaystyle\Ebb_{\vvv\in S^{n-1}}\<\vvv,\vd_1\>^k = \dfrac{(k-1)\cdot(k-3)\dots 3\cdot 1}{n(n+2)\dots(n+k-2)}\cdot\<\vd_1,\vd_1\>^{k/2}$.
    \item Let $\vd_1=\ve_1=(1,0,\dots,0)$ from 3, then $\displaystyle\Ebb\[\, v_1^k \,\] = \dfrac{(k-1)\cdot(k-3)\dots 3\cdot 1}{n(n+2)\dots(n+k-2)}$ where $v_1$ is the first entry of a random unit vector $\vvv\in\Rbb^n$.
    \item if $x,y$ are two entries of different rows and columns from a random orthogonal matrix, then
    \begin{equation*}
        \Ebb\[x^2y^2\]=\underset{\vvv\in S^{n-1}}{\Ebb}\[\<\vvv,\ve_1\>^2\cdot\(\underset{\vw\in S^{n-2}}{\Ebb}\<\vw,\ve_2^{\;\perp}\>^2\)\]
    \end{equation*}
    where $\ve_2^{\;\perp}=\ve_2-\<\ve_2,\vvv\>\vvv$.
\end{enumerate}
\end{eg}

\section{Alternative Proof without Free Probability}

\setlength{\parskip}{1.5mm}
\setlength{\baselineskip}{1.3em}

In this section we give an alternative proof for \Cref{thm:main-1-circ} without using the results from the free probability theory.
\mainonecirc*

Using the Trace Power Method, $\displaystyle \(\o\circ_R\o'\)_{2k} =  \lim_{n\to\infty}\dfrac{\,1\,}{n}\cdot\Ebb\[\trace\(\(MM^T\)^k\)\]$ where $M=DRD'$. So proving \Cref{thm:main-1-circ} is equivalent to proving the following.

\begin{restatable}{thm}{mainone}
\label{thm:main-1}
Let $M=DRD'$ where $R$ is an $n\times n$ random orthogonal matrix, $D$ and $D'$ are $n\times n$ random diagonal matrices where the diagonal elements are drawn independently from distributions $\Omega$ and $\Omega'$, respectively. Denote $\vec{\a}=(a_1,\dots,\a_k)$ and $\vec{\b}=(\b_1,\dots,\b_k)$. Then
    \begin{align*}
        \lim_{n\to\infty}\dfrac{\,1\,}{n}\cdot\Ebb\[\trace\(\(MM^T\)^k\)\]
        = \sum_{\vec{\a}, \vec{\b}\in P_k} C\(\vaa,\vbb\)\cdot \voo^{\,\vaa}\cdot \vv{\o'}^{\,\vbb}
    \end{align*}
    where if we let $\a_1+\dots+\a_k=a$ and $\b_1+\dots+\b_k=b$, then
    \begin{equation}
        C\(\vaa,\vbb\) = (-1)^{k+a+b-1}\cdot k\cdot \binom{a+b-2}{k-1}\cdot \dfrac{(a-1)!}{\a_1!\dots\a_k!}\cdot \dfrac{(b-1)!}{\b_1!\dots\b_k!}\;.
    \end{equation}
    
    In particular, note that $C\(\vaa,\vbb\)=0$ if $a+b\leq k$.
\end{restatable}

More specifically, we are going to prove the equivalent statement of \Cref{prop:moments-ab-step1} in the language of random matrices alternatively, and the rest of the proofs remain the same.

\begin{thm}\label{thm:app-step1}
 Let $M$ be defined as before. Then
\begin{equation}
    \lim_{n\to\infty}\dfrac{\,1\,}{n}\cdot\Ebb\[\trace\(\(MM^T\)^k\)\] = \sum_{\vaa, \vbb\in P_k}\, (-1)^{\,a+b-k-1} \cdot\( \sum_{(\pi, \sigma)\in \np(\vaa,\vbb)}\,  \prod_{x\in S_{\pi\oplus \sigma}}  C_{x-1}\) \cdot \voo^{\,\vaa}\cdot \vv{\o'}^{\,\vbb}
\end{equation}
where $a=\a_1+\dots+\a_k$ and $b=\b_1+\dots+\b_k$.
\end{thm}

\begin{rmk}
    The above theorem is equivalent to the statement
    \begin{equation}
        \varphi\((ab)^k\) = \sum_{\vaa, \vbb\in P_k}\, (-1)^{\,A+B-k-1} \cdot\( \sum_{(\pi, \sigma)\in \np(\vaa,\vbb)}\,  \prod_{x\in S_{\pi\oplus \sigma}}  C_{x-1}\) \cdot \vec{\varphi_a}^{\vaa} \cdot \vec{\varphi_b}^{\vbb}
\end{equation}
\end{rmk}

\begin{proof}[Proof of \Cref{thm:main-1}]
    With the same steps in \Cref{section:moments}, we can prove that
\begin{equation}
    \sum_{(\pi, \sigma)\in \np(\vaa,\vbb)}\,  \prod_{x\in S_{\pi\oplus \sigma}}  C_{x-1} = k\cdot \binom{a+b-2}{k-1}\cdot\dfrac{(a-1)!}{\a_1!\dots\a_k!}\cdot \dfrac{(b-1)!}{\b_1!\dots\b_k!}
\end{equation}
and thus \Cref{thm:app-step1} follows.
\end{proof}

\begin{defn}
    Let $\vi=\(i_1,\dots,i_k\)\subseteq [n]^k$ and $\vaa = \(\a_1,\dots,\a_k\)\in P_k$. We define $n_{j} = \abs{\lcurb s\in[k]:i_s=i_j\rcurb}$ to be the number of elements in $\vi$ that are equal to $i_j$, including $i_j$ itself. We say that $P\(\,\vi\,\) = \vaa$ if $\abs{\lcurb j: n_{j}=i \rcurb } = i\cdot\a_i$ for all $i\in[k]$.
\end{defn}

\begin{eg}
    $P(1,1,1)=(0,0,1)$, $P(1,2,1,4,4,3,1) = (2,1,1,0,0,0,0)$.
\end{eg}

\begin{prop}
    Let $M=DRD'$ where $R$ is an $n\times n$ random orthogonal matrix, $D$ and $D'$ are $n\times n$ random $\o$ and $\o'$-distribution diagonal matrices, respectively. Then
    \begin{align*}
        \Ebb\[\trace\(MM^T\)^k\]
        &=\sum_{i_m\in[n], j_m\in[n]} \underset{d_{i_m}\sim\o}{\Ebb}\[\, d_{i_1}^2\dots d_{i_{k}}^2\,\]\cdot
        \Ebb\[\,\prod_{m=1}^k R(i_{m},j_{m})R(i_{m+1},j_{m})\]\cdot
        \underset{d_{j_m}'\sim\o'}{\Ebb}\[\,d_{j_1}'^2\dots d_{j_{k}}'^2\,\]
    \end{align*}
\end{prop}

\begin{prop}
    Let $M=DRD'$ and $D\(\vaa,\vbb\)$ be defined as
    \begin{equation}
        D\(\vaa, \vbb\) 
        = \lim_{n\to\infty}\dfrac{\,1\,}{n}\cdot\(\sum_{\substack{ \vi,\,\vj\subseteq [n]^k: \\ P(\,\vi\,)=\vaa,\, P(\,\vj\,)=\vbb}} \Ebb\[\,\prod_{m=1}^k R(i_{m},j_{m})R(i_{m+1},j_{m})\]\)
    \end{equation}
    Then 
    \begin{equation}
        \lim_{n\to\infty}\, \dfrac{\,1\,}{n}\cdot \Ebb\[\trace\(MM^T\)^k\] = \sum_{\vaa, \vbb\in P_k} D\(\vaa, \vbb\) \cdot \voo^{\,\vaa}\cdot \vv{\o'}^{\,\vbb}
    \end{equation}
\end{prop}

\begin{proof}

\begin{align*}
    \Ebb\[\trace\(\(MM^T\)^k\)\] 
    & =\sum_{i_m\in[n], j_m\in[n]} \underset{d_{i_m}\sim\o}{\Ebb}\[\, d_{i_1}^2\dots d_{i_{k}}^2\,\]\cdot
    \Ebb\[\,\prod_{m=1}^k R(i_{m},j_{m})R(i_{m+1},j_{m})\]\cdot
    \underset{d_{j_m}'\sim\o'}{\Ebb}\[\,d_{j_1}'^2\dots d_{j_{k}}'^2\,\] \\
    & =\sum_{\vec{\a}, \vec{\b}\in P_k} \sum_{\substack{P(\,\vec{i}\,)=\vec{\a}, \\ P(\,\vec{j}\,)=\vec{\b}}} 
    \underset{d_{i_m}\sim\o}{\Ebb}\[\, \prod_{m=1}^k d_{i_m}^2\,\]\cdot
    \Ebb\[\,\prod_{m=1}^k R(i_{m},j_{m})R(i_{m+1},j_{m})\]\cdot
    \underset{d_{j_m}'\sim\o'}{\Ebb}\[\,\prod_{m=1}^k d_{j_m}'^2\,\] \\
    & =\sum_{\vec{\a}, \vec{\b}\in P_k} 
    \(\prod_{i=1}^k\o_{2i}^{\;\a_i}\)
    \cdot \(\sum_{\substack{P(\,\vec{i}\,)=\vec{\a}, \\ P(\,\vec{j}\,)=\vec{\b}}} 
    \Ebb\[\,\prod_{m=1}^k R(i_{m},j_{m})R(i_{m+1},j_{m})\]\)
    \cdot \(\prod_{j=1}^k\o_{2j}'^{\;\b_j}\) \\
\end{align*}

Dividing both sides by $n$ and evaluating the limit as $n\to\infty$, we get the conclusion.
\end{proof}

By the above proposition, to prove \Cref{thm:app-step1}, it suffices to prove the following.
\begin{thm}\label{thm:app-step2}
For all $\vaa, \vbb\in P_k$,
    \begin{equation}
        D\(\vaa,\vbb\) = (-1)^{A+B-k-1}\cdot \(\sum_{(\pi, \sigma)\in \np(\vaa,\vbb)}\, \prod_{x\in S_{\pi\oplus \sigma}}  C_{x-1} \)
    \end{equation}
\end{thm}

\begin{rmk}
    For given $\vaa,\vbb$, we will use notations $D\(\vaa,\vbb\)$ and $D\(\voo^{\vaa},\vv{\o'}^{\vbb}\)$ interchangeably.
\end{rmk}

We will prove \Cref{thm:app-step2} by first proving the base cases first, and then prove it generally.

\subsection{ Case when \texorpdfstring{$D\(\vaa,\vbb\) = D\(\o_2^{\,k},\o_2'^{\,k}\)$}{D(a,ß)= D(Ω2ˆk,Ω2'ˆk)}} 

\setlength{\parskip}{1.5mm}
\setlength{\baselineskip}{1.3em}

We want to show that 
\begin{equation}
    C\(\o_2^{k}, \o_2'^{k}\) 
    = (-1)^{k-1}\cdot \dfrac{\,1\,}{k}\cdot \binom{2k-2}{k-1} = (-1)^{k-1}\cdot C_{k-1}
\end{equation}
where $C_i$ is the $i^{th}$ Catalan Number.

\begin{defn}
Given $s_1,\dots,s_k\in[k]$ such that $s_1+\dots+s_m=k$, we define $\val\(G_{s_1},\dots, G_{s_m}\)$ to be 
\begin{equation*}
    \val\(G_{s_1},\dots,G_{s_m}\)
    =\Ebb\[\prod_{m=1}^k R(i_{m},j_{m})R(i_{m+1},j_{m})\]
\end{equation*}
where $i_1,\dots, i_k\in[n]$ and $j_1,\dots, j_k\in[n]$ are chosen such that $\displaystyle i_{l_1}=i_{l_2} \iff l_r=1+\sum_{j=1}^{x_r} s_j$ for some $x_r\in\{0,1,\dots,m\}$ for each $r=1,2$ and $j_1\neq j_2\neq\dots\neq j_k$\,. See \Cref{fig:grid-matrix-combine} for illustration.

\begin{figure}[hbt!]
    \centering
    \includegraphics[scale=0.32]{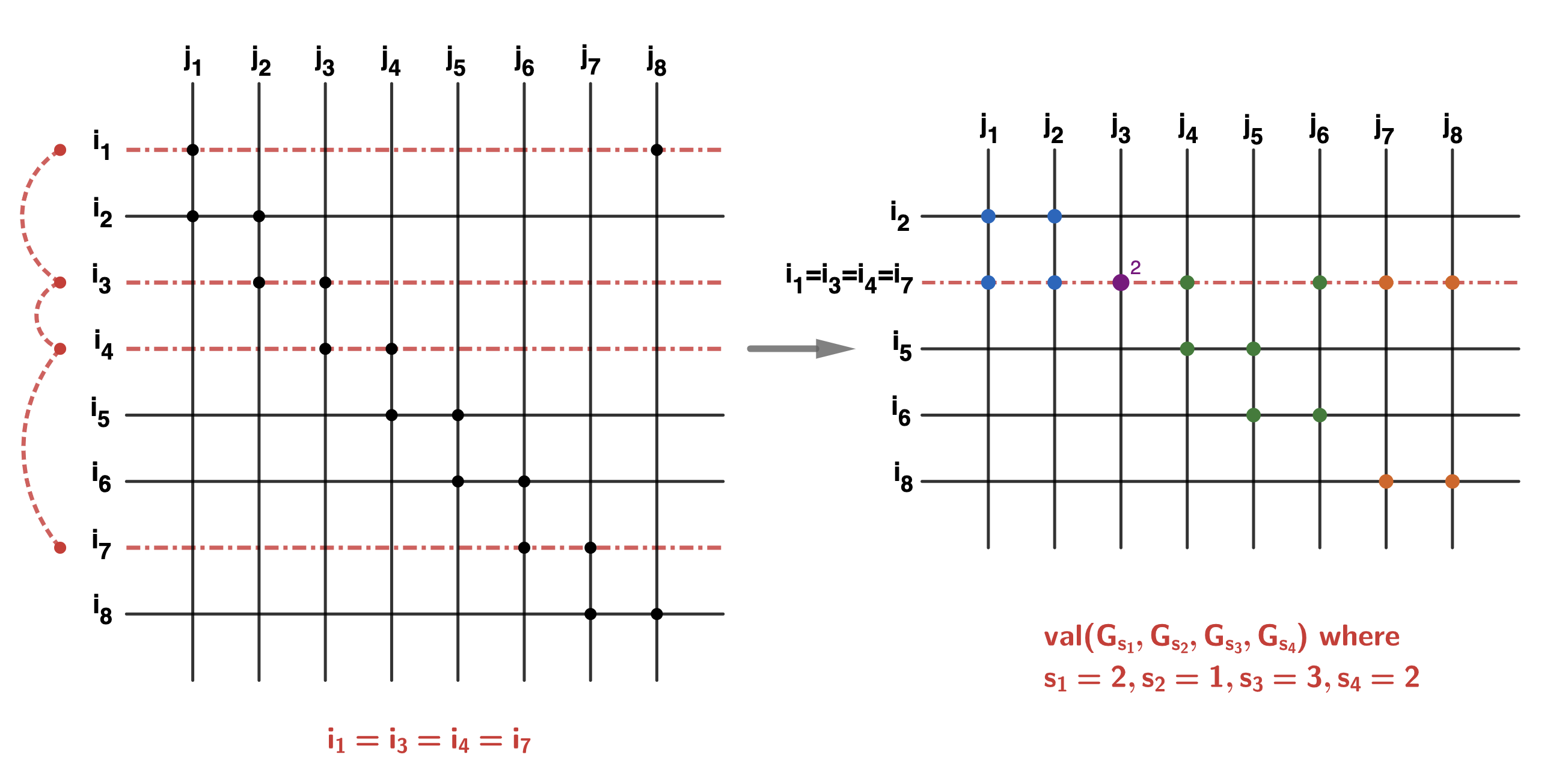}
    \caption{Example of $\val\(G_{s_1},\dots, G_{s_m}\)$ when $n=8$, $m=4$, $s_1=2$, $s_2=1$, $s_3=3$ and $s_4=2$. Note that the purple dot is doubled.}
    \label{fig:grid-matrix-combine}
\end{figure}

\medskip
In particular, $\displaystyle\val(G_k) = \Ebb\[\prod_{m=1}^k R(i_{m},j_{m})R(i_{m+1},j_{m})\]$ where $i_1\neq\dots\neq i_k$ and  $j_1\neq\dots\neq j_k$.
\end{defn}
We can compute $\val\(G_{s_1},\dots, G_{s_m}\)$ using the following theorem which we show in Section \ref{subsection:checkingnoncrossingoptimality}.
\begin{thm}\label{thm:disjointstaircases}
    \begin{equation*}
        \val\(G_{s_1},\dots,G_{s_m}\) \sim \prod_{j=1}^m \val(G_{s_j}) \hspace{1cm} \(\text{ as } n\to\infty \)
    \end{equation*}
\end{thm}
\begin{thm}\label{thm:staircaseCatalan}
\begin{equation*}
    \val(G_k) \sim (-1)^{k-1}\cdot  n^{-2k+1} \cdot C_{k-1} \hspace{1cm} \(\text{ as } n\to\infty \)
\end{equation*}
where $C_{k-1}$ is the $(k-1)^{th}$ Catalan number.
\end{thm}

\begin{proof}
We prove this by induction on $k$.
\begin{enumerate}
    \item Base case $k=1$: $\displaystyle\val(G_k) = \Ebb\[R(i_1,j_1)^2\] = \underset{\vvv\in S^{n-1}}{\Ebb}\[\<\vvv, \ve_1\>^2\] = \dfrac{\,1\,}{n}\cdot \<\ve_1,\ve_1\> = \dfrac{\,1\,}{n} = n^{-1}\cdot C_0$.
        
    \item Inductive case: assume $\val(G_i) \sim (-1)^{i-1}\cdot  n^{-2i+1} \cdot C_{i-1}$ for all $i<k$, we want to prove that $\val(G_k) \sim (-1)^{k-1}\cdot  n^{-2k+1} \cdot C_{k-1}$.
    
    We will define $\E_1,\dots,\E_k$ inductively starting from $\E_1$. Note that $\E_j$ is a function of $\vvv_1,\ldots,\vvv_{j-1}$.
    
    \begin{enumerate}[i.]
        \item We define $\displaystyle \E_1=\underset{\vvv_1\in S^{n-k}}{\Ebb}\[\<\vvv_1,\ve_k^{\;\perp\times (k-1)}\>\<\vvv_1,\ve_{1}^{\;\perp\times (k-1)}\>\]$ where $\ve^{\;\perp\times i}$ is the orthogonal projection of $\ve$ onto the orthogonal complement of $\spann\(\vvv_{k-i+1},\dots,\vvv_k\)$, and
            
        \item for each $2\leq i \leq k$, $\displaystyle \E_i=\underset{\vvv_i\in S^{n-1-(k-i)}}{\Ebb}\[\<\vvv_i,\ve_{i-1}^{\;\perp\times (k-i)}\>\<\vvv_i,\ve_{i}^{\;\perp\times (k-i)}\>\cdot  \E_{i-1}\]$.
    \end{enumerate} 
        
    Defined this way, $\E_k = \val(G_k)$.
    
    Since 
    \begin{align*}
        \ve^{\;\perp\times i} 
        = \ve^{\;\perp\times (i-1)}-\<\vvv_{k-i+1},\ve^{\;\perp\times (i-1)}\>\cdot \vvv_{k-i+1}\,,
    \end{align*}
    
    we have that 
    \begin{align*}
        \E_1
        &= \underset{\vvv_1\in S^{n-k}}{\Ebb}\[\<\vvv_1,\ve_k^{\;\perp\times (k-1)}\>\<\vvv_1,\ve_{1}^{\;\perp\times (k-1)}\>\]\\
        &= \dfrac{1}{n-k+1} \cdot \<\,\ve_k^{\;\perp\times (k-1)}, \ve_{1}^{\;\perp\times (k-1)}\,\> \\
        &=\dfrac{1}{n-k+1} \cdot\(\<\,\ve_k^{\;\perp\times (k-2)}, \ve_{1}^{\;\perp\times (k-2)}\,\> - \<\vvv_{2}, \ve_k^{\;\perp\times(k-2)}\>\<\vvv_{2},\ve_1^{\;\perp\times(k-2)}\> \)\\
        &=\dfrac{1}{n-k+1} \cdot \bigg(\<\,\ve_k^{\;\perp\times (k-3)}, \ve_{1}^{\;\perp\times (k-3)}\,\> \\
        &\hspace{2cm} - \<\vvv_{3}, \ve_k^{\;\perp\times(k-3)}\> \<\vvv_{3},\ve_1^{\;\perp\times(k-3)}\> - \<\vvv_{2}, \ve_k^{\;\perp\times(k-2)}\> \<\vvv_{2},\ve_1^{\;\perp\times(k-2)}\> \bigg) \\
        &= \dots\dots \\
        &= -\dfrac{1}{n-k+1}\cdot \(\sum_{i=1}^{k-1} \<\vvv_{k-i+1}, \ve_k^{\;\perp\times(i-1)}\> \<\vvv_{k-i+1},\ve_{1}^{\;\perp\times(i-1)}\>\)\\
        &= -\dfrac{1}{n-k+1}\cdot \(\sum_{i=1}^{k-1} \<\vvv_{k-i+1}, \ve_k\> \<\vvv_{k-i+1},\ve_{1}\>\)
    \end{align*}
    
    Thus
    \begin{align*}
        \E_k
        &=\dfrac{-1}{n-k+1}\big(\;\Ebb\[\vvv_1\sim \vvv_2\]+ \Ebb\[\vvv_1\sim \vvv_3\]+ \dots+ \Ebb\[\vvv_1\sim \vvv_{k}\]\;\big)\\
        &\\
        &\sim -n^{-1}\cdot \big(\val\(G_1,G_{k-1}\)+\val\(G_2,G_{k-2}\) + \dots + \val\(G_{k-1},G_1\)\big) \\
        &\sim -n^{-1}\cdot \(\sum_{i=1}^{k-1} \((-1)^{i-1}\cdot  n^{-2i+1} \cdot C_{i-1}\) \cdot \( (-1)^{k-i-1}\cdot  n^{-2(k-i)+1} \cdot C_{k-i-1}\)\)\\
        &= (-1)^{k-1}\cdot n^{-2k+1}\cdot\(\sum_{i=1}^{k-1} C_{i-1}C_{k-i-1}\) \\
        &=(-1)^{k-1}\cdot n^{-2k+1}\cdot C_{k-1} \,.
    \end{align*}
    where $\vvv_1\sim \vvv_i$ means folding the first row of the matrix to the $i^{th}$ row. See \Cref{fig:grid-combine-pf} for illustration.
\end{enumerate}

    \begin{figure}[hbt!]
        \centering
        \includegraphics[scale = 0.32]{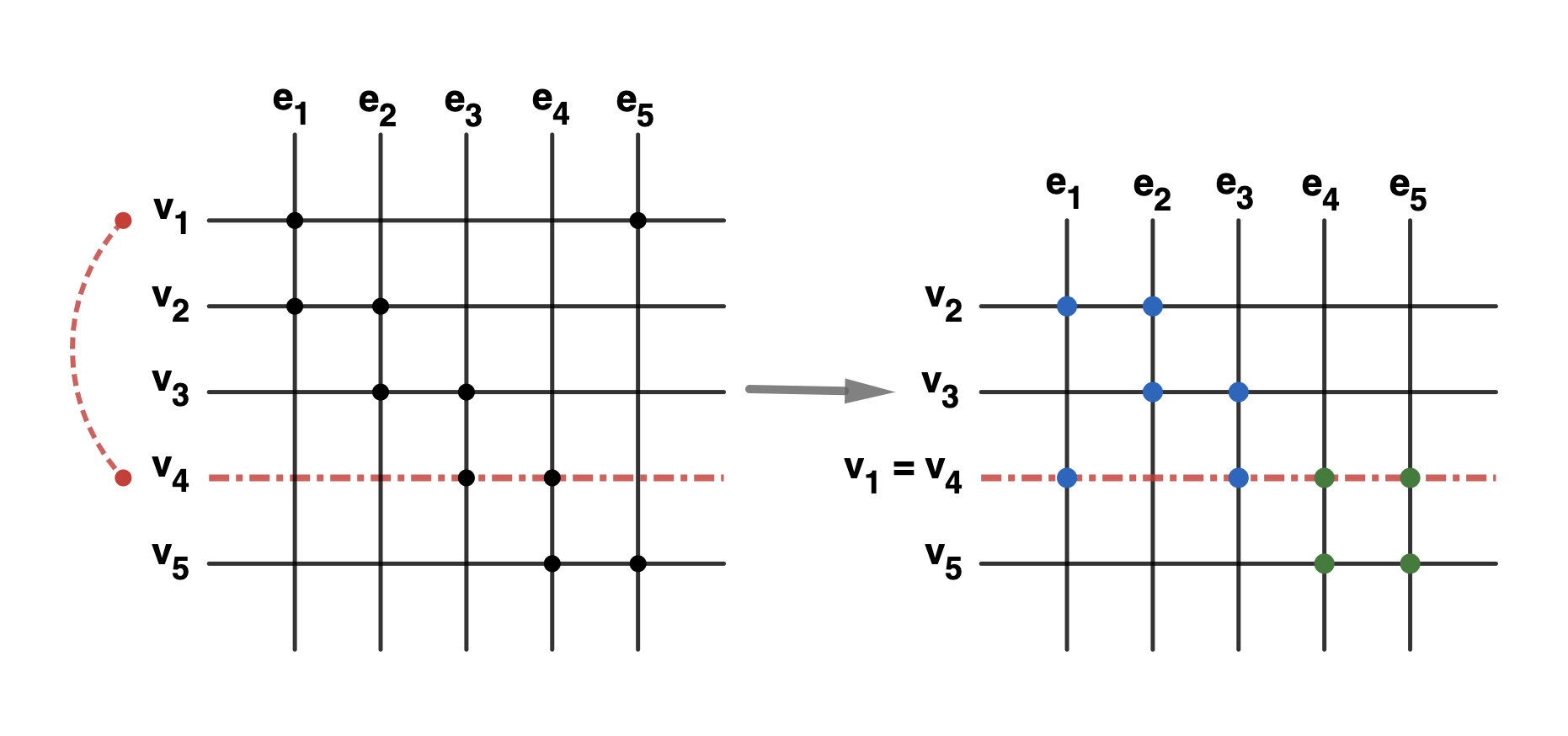}
        \caption{Illustration of $\Ebb\[\vvv_1\sim \vvv_4\]$.}
        \label{fig:grid-combine-pf}
    \end{figure}
    
\end{proof}

\begin{thm}
    $\displaystyle D\(\Omega_2^k, \,{\Omega'_2}^k\)=(-1)^{k-1}\cdot C_{k-1}$, where $C_k$ is the $\(k-1\)^{th}$ Catalan number.
\end{thm}
\begin{proof}
    \begin{equation*}
        \begin{aligned}
        D\(\Omega_2^k,\,{\Omega'_2}^k\)
        &=\lim_{n\to\infty}\dfrac{\,1\,}{n}\cdot\( \sum_{\substack{i_1\neq\dots\neq i_k,\\ j_1\neq\dots\neq j_k}} \Ebb\[\prod_{m=1}^k R(i_{m},j_{m})R(i_{m+1},j_{m})\]\) \\
        &=\lim_{n\to\infty}\dfrac{\,1\,}{n}\cdot\( \sum_{\substack{i_1\neq\dots\neq i_k,\\ j_1\neq\dots\neq j_k}} \val\(G_k\) \) \\
        &=\lim_{n\to\infty}\dfrac{\,1\,}{n}\cdot n^{2k}\cdot \val\(G_k\)\\
        &=\lim_{n\to\infty} n^{2k-1}\cdot n^{-2k+1}\cdot (-1)^{k-1}\cdot C_{k-1}
         = (-1)^{k-1}\cdot C_{k-1}\,.
        \end{aligned}
    \end{equation*}
    
    \smallskip
    Thus $D\(\Omega_2^k,\,{\Omega'_2}^k\) = (-1)^{k-1}\cdot C_{k-1}$.
\end{proof}

\subsection{ Case when \texorpdfstring{$D\(\vaa,\vbb\) = D\(\vaa, \o_2'^{\,k}\)$}{D(a,ß)= D(a,Ω2'ˆk)}} 

\setlength{\parskip}{1.5mm}
\setlength{\baselineskip}{1.3em}

We want to show that

\begin{restatable}{thm}{CoefficientStepTwo}
\label{thm:coefficient-step-2}
    Let $\vaa = \(\a_1,\dots,\a_k\)\in P_k$ and $a=\a_1+\dots+\a_k$. Then
    \begin{equation}
        D\(\voo^{\vaa},\, \Omega_2'^{\;k}\) = (-1)^{a-1}\cdot \dfrac{(a-1)!}{\a_1!\dots\a_k!}\cdot \binom{a+k-2}{k-1} \,.
    \end{equation}
\end{restatable}

\begin{defn}
    We say that an assignment $i:[k]\to[n]$ (conventionally denote $i(m)$ as $i_m$) 
    \emph{respects} a partition $\P$ of $[k]$ if for all $j_1,j_2\in[k]$, $i_{j_1}=i_{j_2} \iff j_1$ and $j_2$ are in the same parts under $\P$.
\end{defn}

\begin{defn}
    Given $\vaa = \(\a_1,\dots,\a_k\), \vbb = \(\b_1,\dots,\b_k\)\in P_k$ and $\P_r\in \P\(\vaa\)$, $\P_c\in\P\(\vbb\)$, we define
    \begin{enumerate}
        \item $\displaystyle \val\(\P_r,\P_c\) = \underset{}{\Ebb}\[\prod_{m=1}^k R\(i_m,j_m\)R\(i_{m+1},j_m\)\]$, where $i:[k]\to[n]$ and $j:[k]\to[n]$ are some assignments that respect the partitions $\P_r$ and $\P_c$, respectively.
        \item $\displaystyle N\(\P_r,\P_c\) = \abs{\lcurb\(i,j\): i:[k]\to[n] \text{ respects } \P_r \text{ and } j:[k]\to[n] \text{ respects } \P_c \rcurb}$.
    \end{enumerate}
    
    In the special cases when $\P_c$ or $\P_r$ is $\P_0\in \P(k,0,\dots,0)$, we denote $\val\(\P_r\)$ to be $\val\(\P_r,\P_0\)$ and $\val\(\P_c\)$ to be $\val\(\P_0,\P_c\)$.
\end{defn}

\begin{prop}
Let $\vaa=\(\a_1,\dots,\a_k\), \vbb=\(\b_1,\dots,\b_k\)\in P_k$ and $\P_r\in \P\(\vaa\)$, $\P_c\in\P\(\vbb\)$. Let $a=\a_1+\dots+\a_k$ and $b=\b_1+\dots+\b_k$. Then
    \begin{equation*}
        N\(\P_r,\P_c\) = \dfrac{n!}{(n-a)!}\cdot \dfrac{n!}{(n-b)!} \sim n^{a+b}\,.
    \end{equation*}
\end{prop}

By definitions of $\val\(\P_r,\P_c\)$ and $N\(\P_r,\P_c\)$, we can express $D\(\vaa,\vbb\)$ as the following. 
\begin{prop}
Let $\vaa=\(\a_1,\dots,\a_k\)\in P_k$ and $\vbb=\(\b_1,\dots,\b_k\)\in P_k$. Then
\begin{equation}
    D\(\vaa,\vbb\) = \lim_{n\to\infty}\, \dfrac{\,1\,}{n}\cdot \(\sum_{\substack{\P_r\in\P\(\vaa\), \P_c\in\P\(\vbb\)}} N\(\P_r,\P_c\)\cdot\val\(\P_r,\P_c\)\)\,.
\end{equation}
\end{prop}
As we show in Section \ref{subsection:checkingnoncrossingoptimality}, we only need to consider partitions $\P_r$ which are non-crossing.
\begin{thm}\label{thm:noncrossing-dominant-row}
    Given $\vaa = \(\a_1,\dots,\a_k\)\in P_k$, let $a=\a_1+\dots+\a_k$. Let $\P_r\in \P\(\vaa\)$. Then
    \begin{equation}
        n^{a+k}\cdot\val\(\P_r\)=\Theta(n)
        \iff \P_r\in\np\(\vaa\).
    \end{equation}
    Moreover, if $\P_r$ is not non-crossing, then $n^{a+k}\cdot\val\(\P_r\)=O(1)$.
\end{thm}

\begin{cor}
    Let $\vaa=\(\a_1,\dots,\a_k\)\in P_k$. Then
    \begin{equation}
        D\(\o_2^{\;\a_1}\dots\o_{2k}^{\;\a_k},\, \Omega_2'^{\;k}\) = \lim_{n\to\infty}\, \dfrac{\,1\,}{n}\cdot \(\sum_{\P_r\in\np\(\vaa\)} N\(\P_r\)\cdot\val\(\P_r\)\)\,.
    \end{equation}
\end{cor}

\begin{defn}
    Let $\P=\lcurb P_1,\dots,P_m \rcurb$ be a partition of $[k]$. Let $\C_k$ be the cycle graph with vertices $\{1,\dots,k\}$. Define $\C_k/\P$ to be the graph obtained by identifying vertices together under $\P$. 
\end{defn}

\begin{prop}\label{prop:grid-cycle-equiv-1}
    Let $\vaa=\(\a_1,\dots,\a_k\)\in P_k$,  $\P_r\in\np(\vaa)$ and $a=\a_1+\dots+\a_k$. Let $S_{\P}=\lcurb i_1,\dots,i_p \rcurb$. Then
    \begin{equation}
        \lim_{n\to\infty}\, n^{a+k-1}\cdot\val\(\P_r\) = (-1)^{a-1}\cdot\prod_{j=1}^p C_{i_j-1}
    \end{equation}
    where $C_{i_j-1}$ is the $\(i_j-1\)^{th}$ Catalan number.
\end{prop}

\begin{figure}[hbt!]
    \centering
    \includegraphics[scale=0.3]{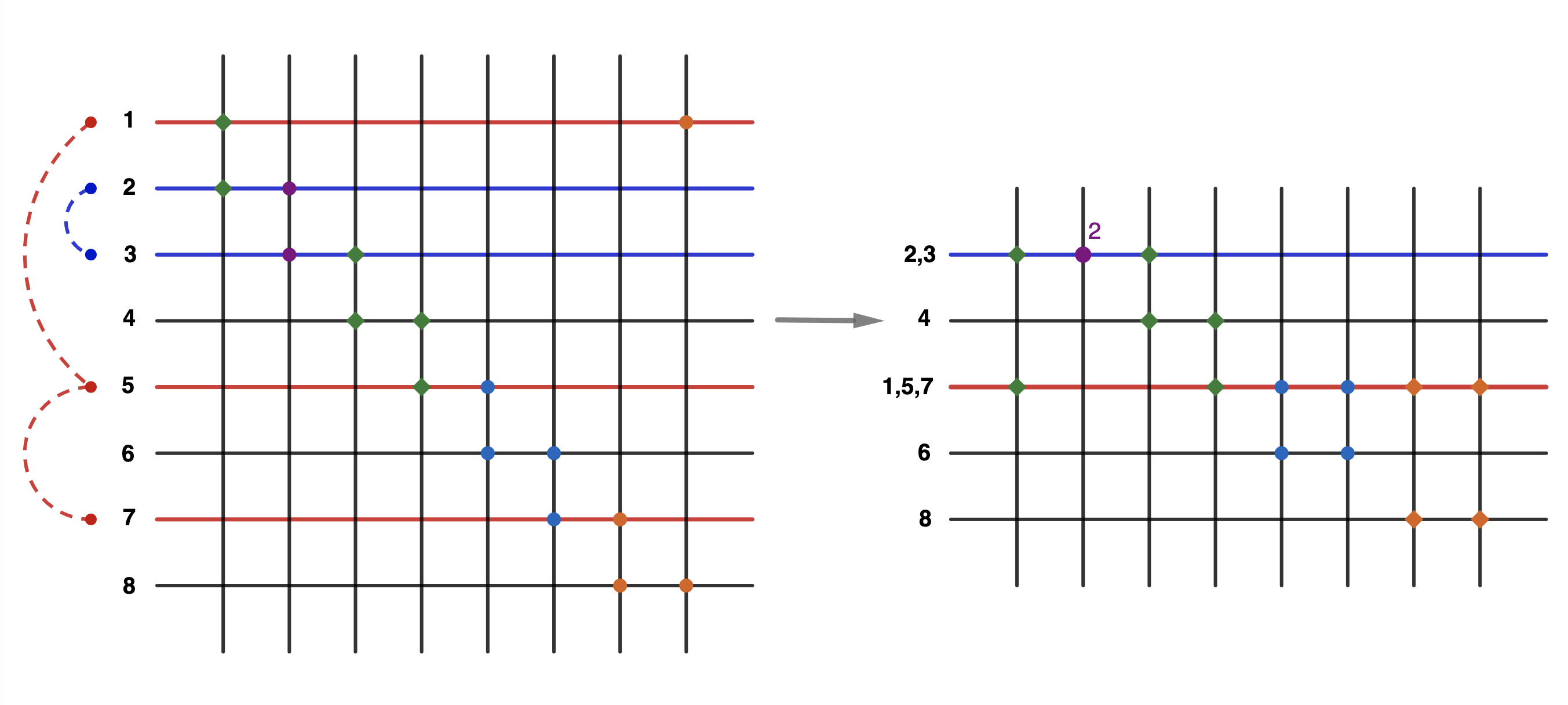}
    \caption{Illustration of \Cref{prop:grid-cycle-equiv-1}: $\P_r=\lcurb\{4\},\{6\},\{8\},\{2,3\},\{1,5,7\}\rcurb$, $\val\(\P_r\) = \val\(G_{i_1},\dots,G_{i_4}\)$ where $i_1=1$, $i_2=3$, $i_3=i_4=2$.}
    \label{fig:np-grid}
\end{figure}

\begin{figure}[hbt!]
    \centering
    \includegraphics[scale=0.3]{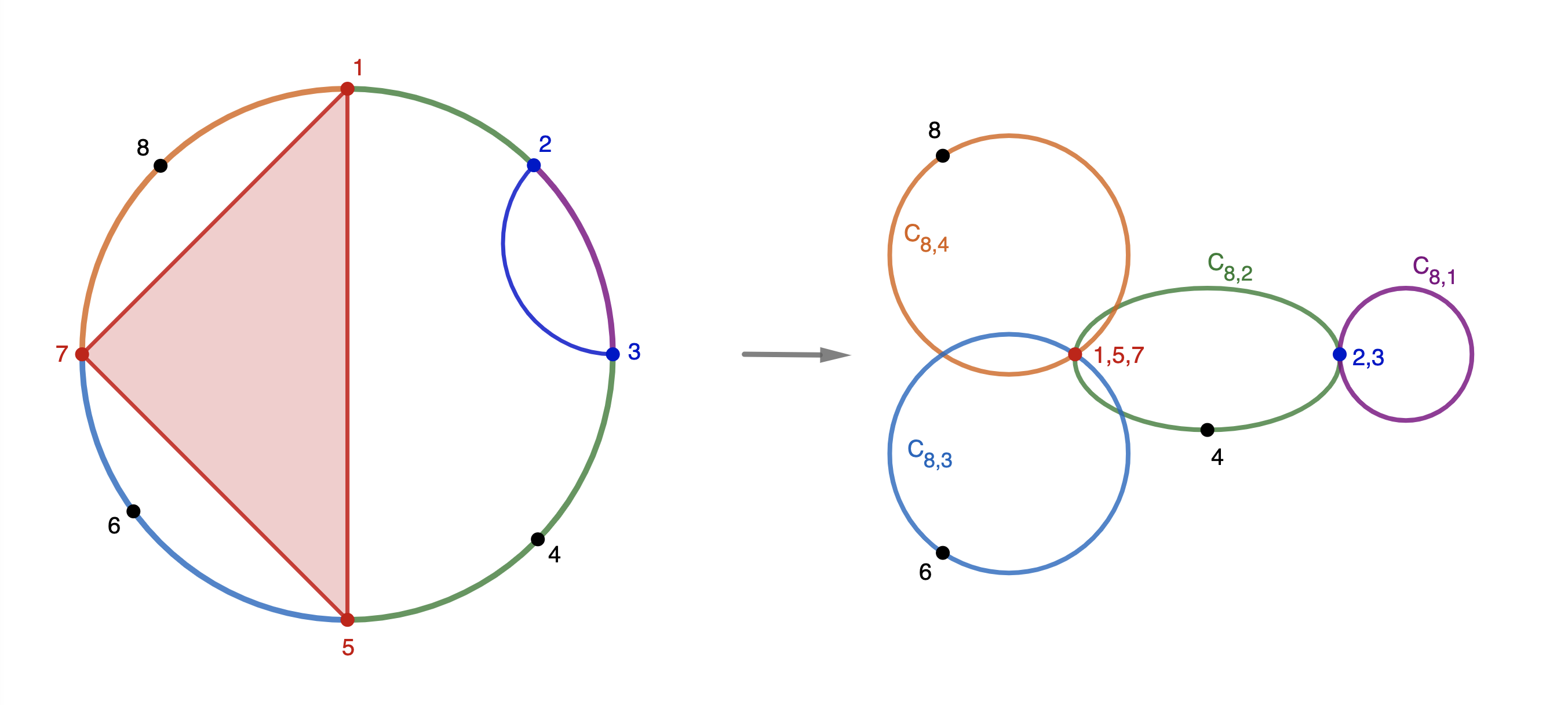}
    \caption{Illustration of \Cref{prop:grid-cycle-equiv-1}: $\P=\lcurb\{4\},\{6\},\{8\},\{2,3\},\{1,5,7\}\rcurb$, $S_{\P} = \{1,3,2,2\}$}
    \label{fig:np-cycle}
\end{figure}

\begin{proof}
By \Cref{prop:noncross-partition-cycles}, $\C_{k}/\P_r = \C_{k,1}\cup\dots\cup\C_{k,p}$ where $p=k-a+1$ and $i_j=\abs{\C_{k,j}}$ for each $j\in[p]$. Then
    \begin{align*}
        \lim_{n\to\infty}\, n^{a+k-1}\cdot \val\(\P_r\)
        &= \lim_{n\to\infty}\, n^{a+k-1}\cdot \val\(G_{i_1},\dots,G_{i_p}\)\\
        &= \lim_{n\to\infty}\, n^{a+k-1}\cdot\prod_{j=1}^p\, \val\(G_{i_j}\)\\
        &= \lim_{n\to\infty}\, n^{a+k-1}\cdot\prod_{j=1}^p\,(-1)^{i_j-1}\cdot n^{-2i_j+1}\cdot C_{i_j-1}\\
        &= \lim_{n\to\infty}\, n^{a+k-1}\cdot (-1)^{k-p}\cdot n^{-2k+p}\cdot\prod_{j=1}^p\,C_{i_j-1}\\
        &= (-1)^{a-1}\cdot\lim_{n\to\infty}\, \(n^{a-k-1+p}\cdot \prod_{j=1}^p\,C_{i_j-1}\)
         = (-1)^{a-1}\cdot\prod_{j=1}^p\,C_{i_j-1}\,.
    \end{align*}
\end{proof}

\subsection{The General Case}

\setlength{\parskip}{1.5mm}
\setlength{\baselineskip}{1.3em}

Now we are ready to show the general case.

\begin{restatable}{thm}{CoefficientStepThree}
\label{thm:coefficient-step-3}
    Let $\vaa, \vbb\in P_k$ and $a=\a_1+\dots+\a_k$, $b=\b_1+\dots+\b_k$. Then
    \begin{equation}
        D\(\vaa, \vbb\) = (-1)^{a+b-k-1}\cdot\dfrac{(a-1)!}{\a_1!\dots\a_k!}\cdot\dfrac{(b-1)!}{\b_1!\dots\b_k!}\cdot k\cdot \binom{a+b-2}{k-1}.
    \end{equation}
\end{restatable}

\begin{defn}\label{defn:non-crossing-two-partitions}
    Let $\P=\lcurb P_1,\dots,P_k\rcurb$ and $\Q=\lcurb Q_1,\dots, Q_l\rcurb$ be two partitions of $[n]$. We say that $\P\cup \Q$ is \emph{non-crossing} if 
    \begin{enumerate}
        \item $\P$ and $\Q$ are non-crossing, and
        \item For any $P\in\P$ and $Q\in\Q$, they do not cross when placing on the cycle $\C_n$, and they touch at at most one point.
    \end{enumerate} 
    
    See \Cref{fig:np-two-partitions} for an illustration.
    
    \begin{figure}[hbt!]
        \centering
        \includegraphics[scale=0.35]{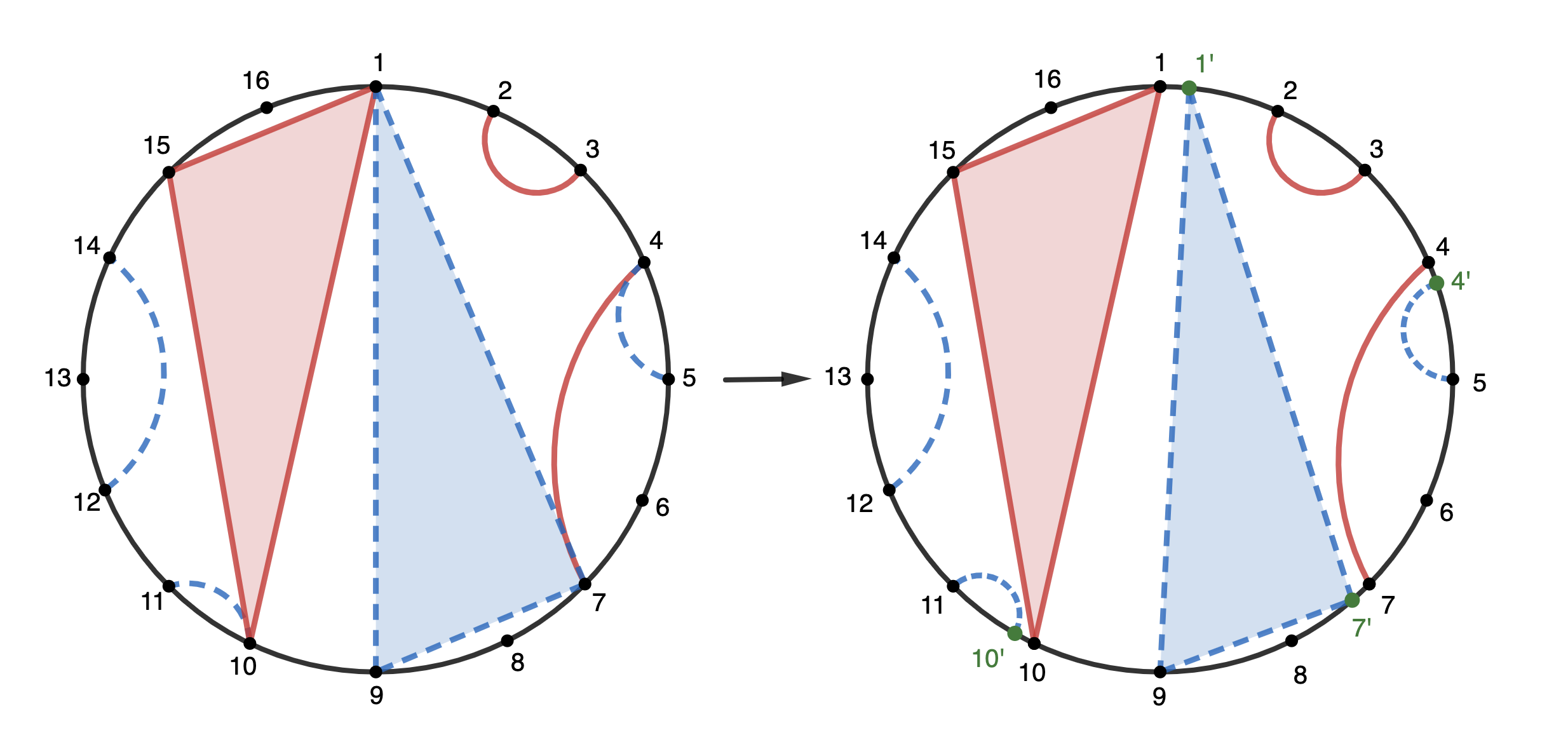}
        \caption{Illustration of \Cref{defn:non-crossing-two-partitions} and \Cref{defn:np-alpha-beta}: $\P$ is colored red and labelled with solid lines, $\Q$ is colored blue and labelled with dash lines. $\P\cup\Q$ is non-crossing and moreover, $\in\np\(\vaa,\vbb\)$ where $\vaa$ and $\vbb$ correspond to $\o_2^{9}\o_4^2\o_6$ and $\o_2^7\o_4^3\o_6$, respectively.}
        \label{fig:np-two-partitions}
    \end{figure}
\end{defn}

\begin{defn}\label{defn:np-alpha-beta}
    Let $\vaa,\vbb\in P_k$. We say $\P\cup\Q\in\np\(\vaa,\vbb\)$ if
    \begin{enumerate}
        \item $\P\in\np\(\vaa\)$, $\Q\in\np\(\vbb\)$,
        \item $\P\cup\Q$ is non-crossing,
        \item For any $P\in\P$ and $Q\in\Q$ that touches at a point $t_0$, we can order $P$ and $Q$ as $P=\{p_1,\dots,p_x,t_0\}$ and $Q=\{t_0, q_1,\dots,q_y,t_0\}$ such that $t_0,q_1,\dots,q_y,p_1,\dots,p_x$ are ordered in the clockwise direction on $\C_k$.
        
        Pictorially, for any $P\in\P$ and $Q\in\Q$ that touch at a point $t_0$, we can perturb the $t_0$ vertex of $Q$ in the clockwise direction a little so that the perturbed $Q$ does not cross $P$. See \Cref{fig:np-two-partitions} for an illustration.
        
        
    \end{enumerate}
    
\end{defn}

As we show in \Cref{subsection:checkingnoncrossingoptimality}, we only need to consider non-crossing partitions.
\begin{thm}\label{thm:dominant-noncrossing}
    Given $\vaa, \vbb\in P_k$, let $a=\a_1+\dots+\a_k$ and $b=\b_1+\dots+\b_k$. Let $\P_r\in \P\(\a_1,\dots,\a_k\)$ and $\P_c\in\P\(\b_1,\dots,\b_k\)$. Then
    \begin{equation*}
        n^{a+b}\cdot\val\(P_r,P_c\) = \Theta(n)
        \iff \P_r\cup \P_c\in\np\(\vaa,\vbb\).
    \end{equation*}
    Moreover, if $\P_r\cup\P_c$ is not non-crossing, then $n^{a+b}\cdot \val\(\P_r\cup\P_c\)=O(1)$.
\end{thm}

\begin{cor}
    Let $\vaa=\(\a_1,\dots,\a_k\), \vbb=\(\b_1,\dots,\b_k\)\in P_k$. Then
    \begin{align*}
        D\(\vaa,\vbb\)
        = \lim_{n\to\infty} \dfrac{\,1\,}{n}\cdot \(\sum_{\P_r\cup\P_c\in\np\(\vaa,\vbb\)}\, N\(\P_r,\P_c\)\cdot\val\(\P_r,\P_c\)\)\\
    \end{align*}
\end{cor}

\begin{defn}\label{defn:np-cycle-size-Spq}
    Let $\P\cup\Q$ be a non-crossing partition of $[k]$. Assume $\C_k/(\P\cup\Q) = \C_{k,1}\cup\dots\cup \C_{k,p}$. We define $S_{\P,\Q}$ to be the unordered sequence $\lcurb i_1,\dots,i_p\rcurb$ where $i_j$ is the size of the cycle $\C_{k,j}$ for each $j\in[p]$.
\end{defn}

\begin{eg}

\begin{figure}[hbt!]
    \centering
    \includegraphics[scale=0.32]{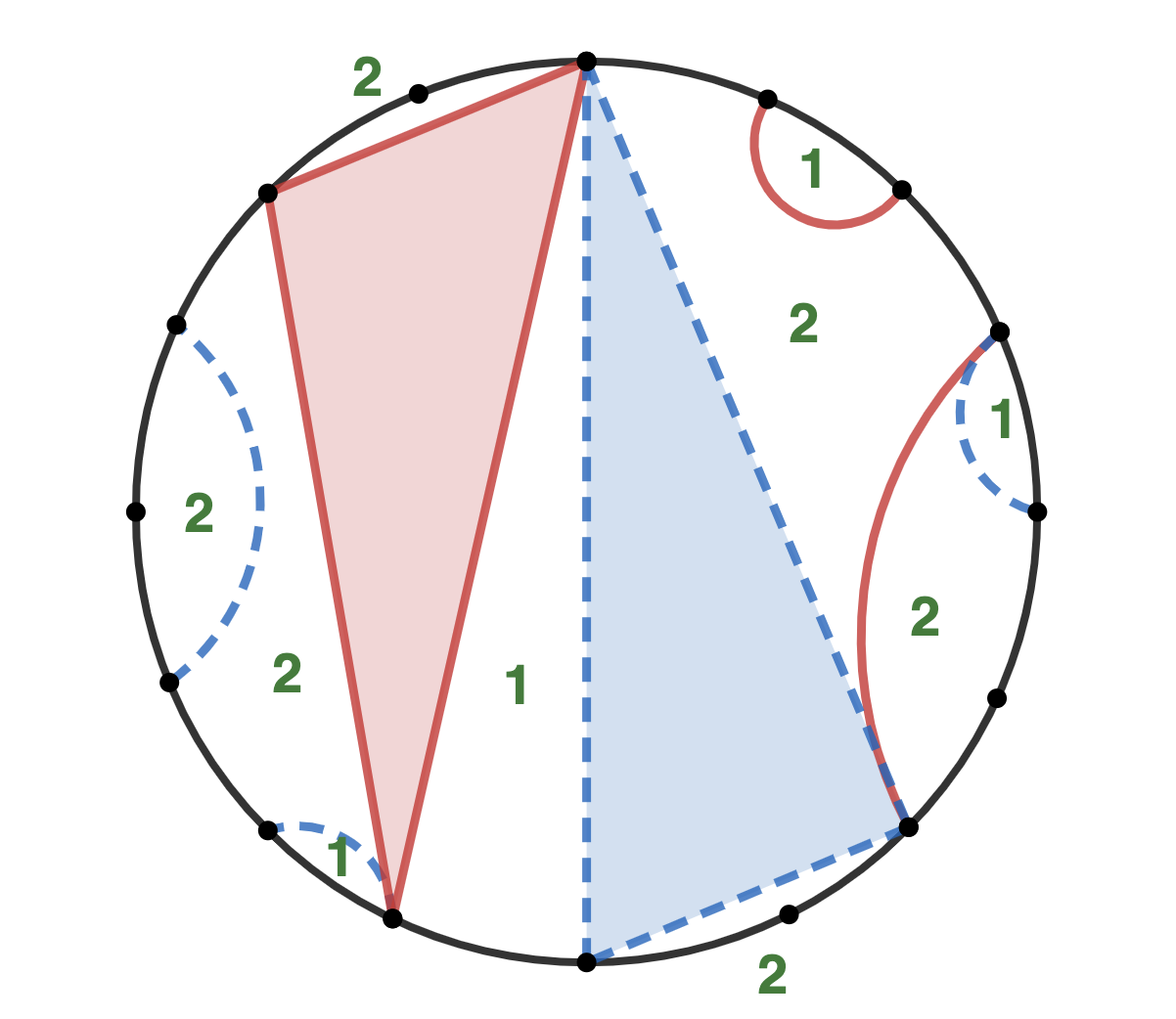}
    \caption{Illustration of \Cref{defn:np-cycle-size-Spq}.}
    \label{fig:np-two-partitions-size}
\end{figure}

    We can label the sizes for \Cref{fig:np-two-partitions}, then $S_{\P,\Q} = \{1,1,1,1,2,2,2,2,2,2\}$ as shown in \Cref{fig:np-two-partitions-size}.
\end{eg}

\begin{prop}\label{prop:grid-cycle-equiv-2}
    Let $\vaa=\(\a_1,\dots,\a_k\),\vbb=\(\b_1,\dots,\b_k\)\in P_k$, $\P_r\cup\P_c\in\np\(\vaa,\vbb\)$, and $a=\a_1+\dots+\a_k$, $b=\b_1+\dots+\b_k$. Let $S_{\P,\Q}=\lcurb i_1,\dots,i_p \rcurb$. Then
    \begin{equation}
        \lim_{n\to\infty}\, n^{a+b-1}\cdot\val\(\P_r,\P_c\) = (-1)^{a+b-k-1}\cdot\prod_{j=1}^p C_{i_j-1}
    \end{equation}
    where $C_{i_j-1}$ is the $\(i_j-1\)^{th}$ Catalan number.
\end{prop}
\begin{proof}
    By \Cref{prop:noncross-partition-cycles}, $\C_{k}/\(\P_r\cup\P_c\) = \C_{k,1}\cup\dots\cup\C_{k,p}$ where $p=2k-a-b+1$ and $i_j=\abs{\C_{k,j}}$ for each $j\in[p]$. Then

    \begin{align*}
        \lim_{n\to\infty}\, n^{a+b-1}\cdot \val\(\P_r,\P_c\)
        &= \lim_{n\to\infty}\, n^{a+b-1}\cdot \val\(G_{i_1},\dots,G_{i_p}\)\\
        &= \lim_{n\to\infty}\, n^{a+b-1}\cdot\prod_{j=1}^p\, \val\(G_{i_j}\)\\
        &= \lim_{n\to\infty}\, n^{a+b-1}\cdot\prod_{j=1}^p\,(-1)^{i_j-1}\cdot n^{-2i_j+1}\cdot C_{i_j-1}\\
        &= \lim_{n\to\infty}\, n^{a+b-1}\cdot (-1)^{k-p}\cdot n^{-2k+p}\cdot\prod_{j=1}^p\,C_{i_j-1}\\
        &= (-1)^{a+b-k-1}\cdot\lim_{n\to\infty}\, \(n^{a+b-1-2k+p}\cdot\prod_{j=1}^p\,C_{i_j-1}\)\\
        &= (-1)^{a+b-k-1}\cdot\prod_{j=1}^p\,C_{i_j-1}\,.
    \end{align*}
\end{proof}

\subsection{Optimality of Non-Crossing Partitions} \label{subsection:checkingnoncrossingoptimality}

\setlength{\parskip}{1.5mm}
\setlength{\baselineskip}{1.3em}

We now show Theorems \ref{thm:disjointstaircases}, \ref{thm:noncrossing-dominant-row}, and \ref{thm:dominant-noncrossing} which say that when we take the limit as $n \to \infty$, we only need to consider non-crossing partitions and we can ignore the interaction between parts of our grid which are not part of the same term.

\begin{figure}[hbt!]
    \centering
    \includegraphics[scale=0.32]{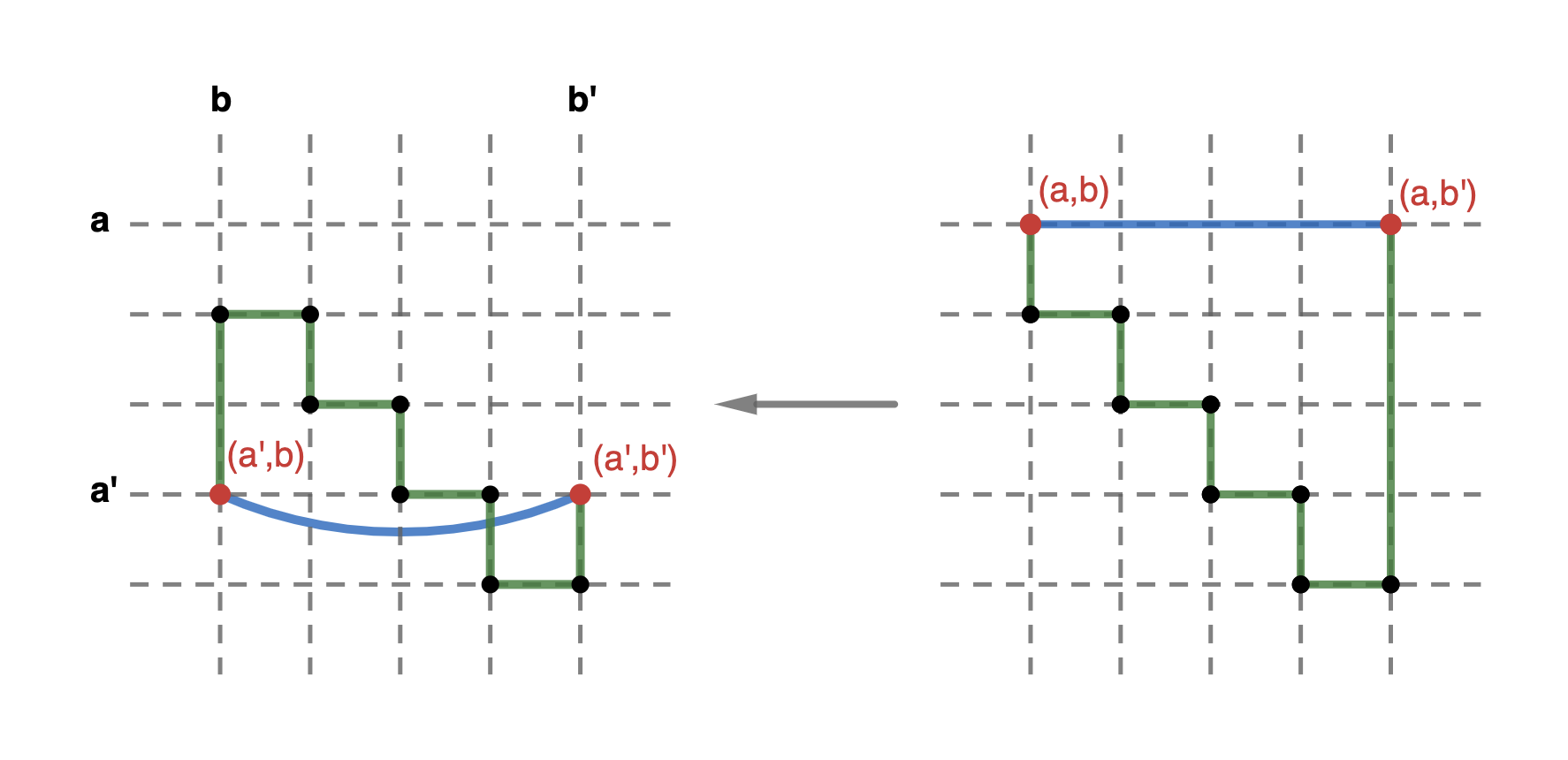}
    \caption{Illustration of \Cref{lem:staircase-matching}, 2: After we shift row $a$ to row $a'$, we have the row matching $\(M_r \setminus \left\{\{(a,b),(a,b')\}\right\}\)  \cup \left\{\{(a',b),(a',b')\}\right\}$. Note that while this is allowed in \Cref{lem:staircase-matching}, it will never happen in any efficient sequence of operations because the staircase is still one connected component.}
    \label{fig:shift-matching-1}
\end{figure}

\begin{figure}[hbt!]
    \centering
    \includegraphics[scale=0.32]{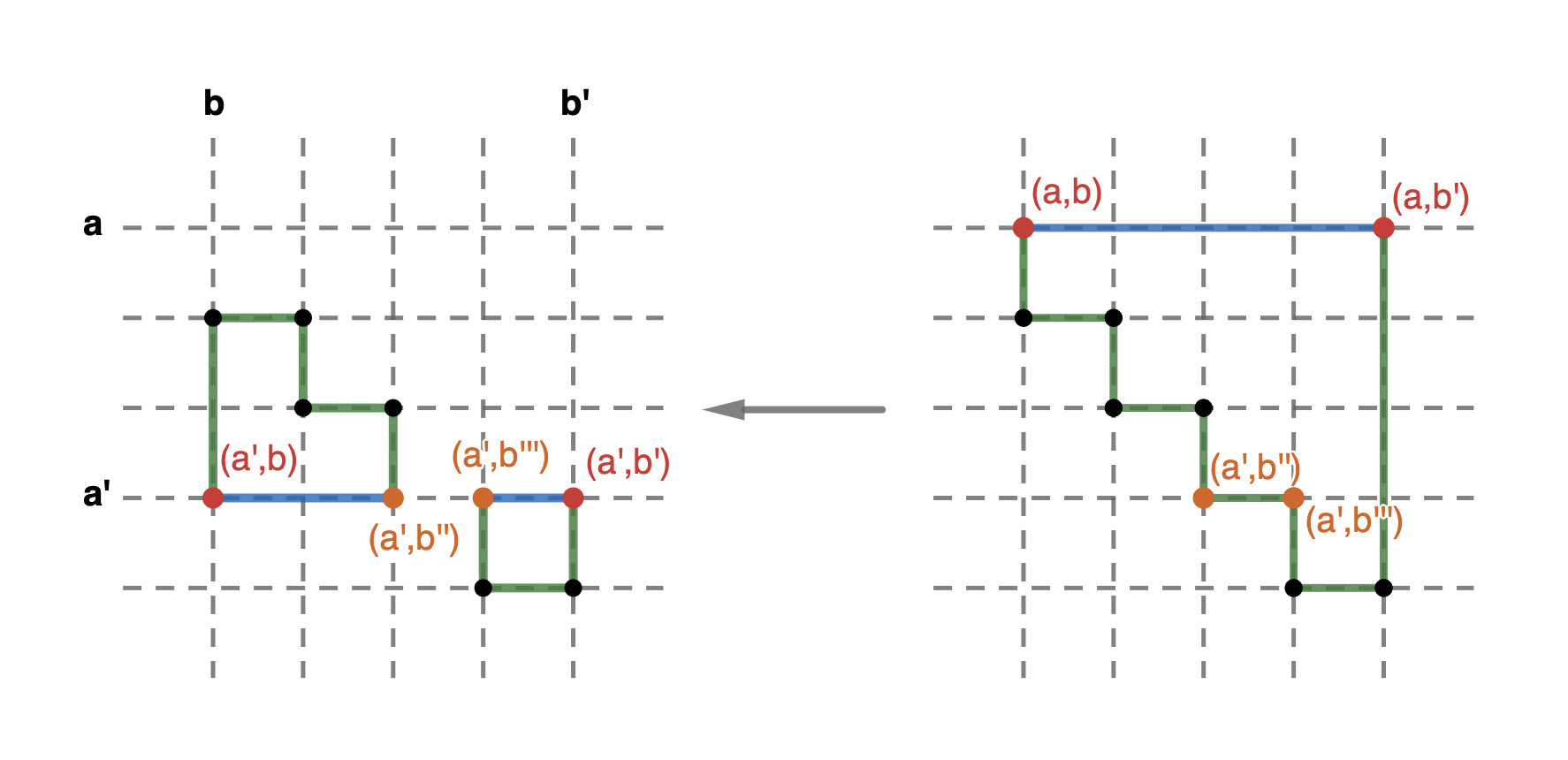}
    \caption{Illustration of \Cref{lem:staircase-matching}, 2: After we shift row $a$ to row $a'$, we have the row matching $\(M_r \setminus \left\{\{(a,b),(a,b')\},\{(a',b''),(a',b''')\}\right\}\)  \cup \left\{\{(a',b),(a',b'')\}, \{(a',b'),(a',b''')\}\right\}$.}
    \label{fig:shift-matching-2}
\end{figure}

Our setup is as follows. We start with the set of vertices $$V_0 = \left\{(i,i): i \in [k]\right\} \cup \left\{(i+1,i): i \in [k-1]\right\} \cup \{(1,k)\}$$
which we think of a staircase of length $k$ (see Definition \ref{def:staircase} below). We then perform the following kinds of operations on our current multi-set of vertices $V$.
\begin{defn}\label{def:mergeoperation}
In a \emph{row merge} operation, we merge a set of rows $A \subseteq [k]$ together. More precisely, we choose a representative $a \in A$ and then for each $a' \in A \setminus \{a\}$, we replace each vertex $(a',b) \in V$ with $(a,b)$. We define the weight of a row merge operation to be $|A| - 1$. The intuition for this is that the row merge operation reduces the number of distinct indices by $|A| - 1$.

\begin{figure}[hbt!]
    \centering
    \includegraphics[scale=0.3]{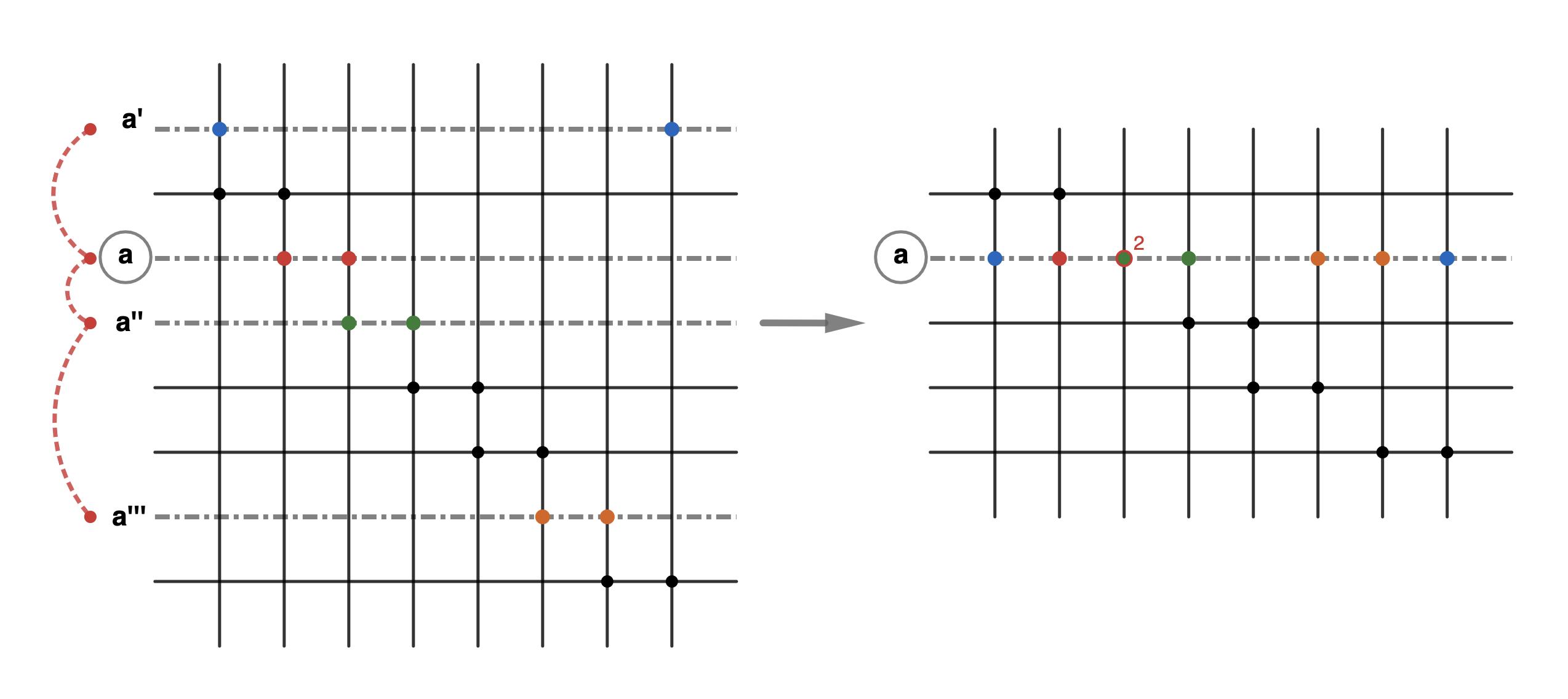}
    \caption{Illustration of \Cref{def:mergeoperation}: here $A=\{a,a',a'',a'''\}$, and we merge rows $A$ into row $a$.}
    \label{fig:grid-merge-operation}
\end{figure}

We define a \emph{column merge} operation and its weight in the same way.
\end{defn}
\begin{defn}\label{def:shiftoperation}
We define a \emph{shift} operation as follows. We take two vertices $(a,b)$ and $(a,b')$ in the same row $a$ and we shift them to a different row $a'$ by replacing them with the vertices $(a',b)$ and $(a',b')$. We write this shift operation as $\{(a,b),(a,b')\} \to \{(a',b),(a',b')\}$. We take all shift operations to have weight $1$.

\begin{figure}[hbt!]
    \centering
    \includegraphics[scale=0.32]{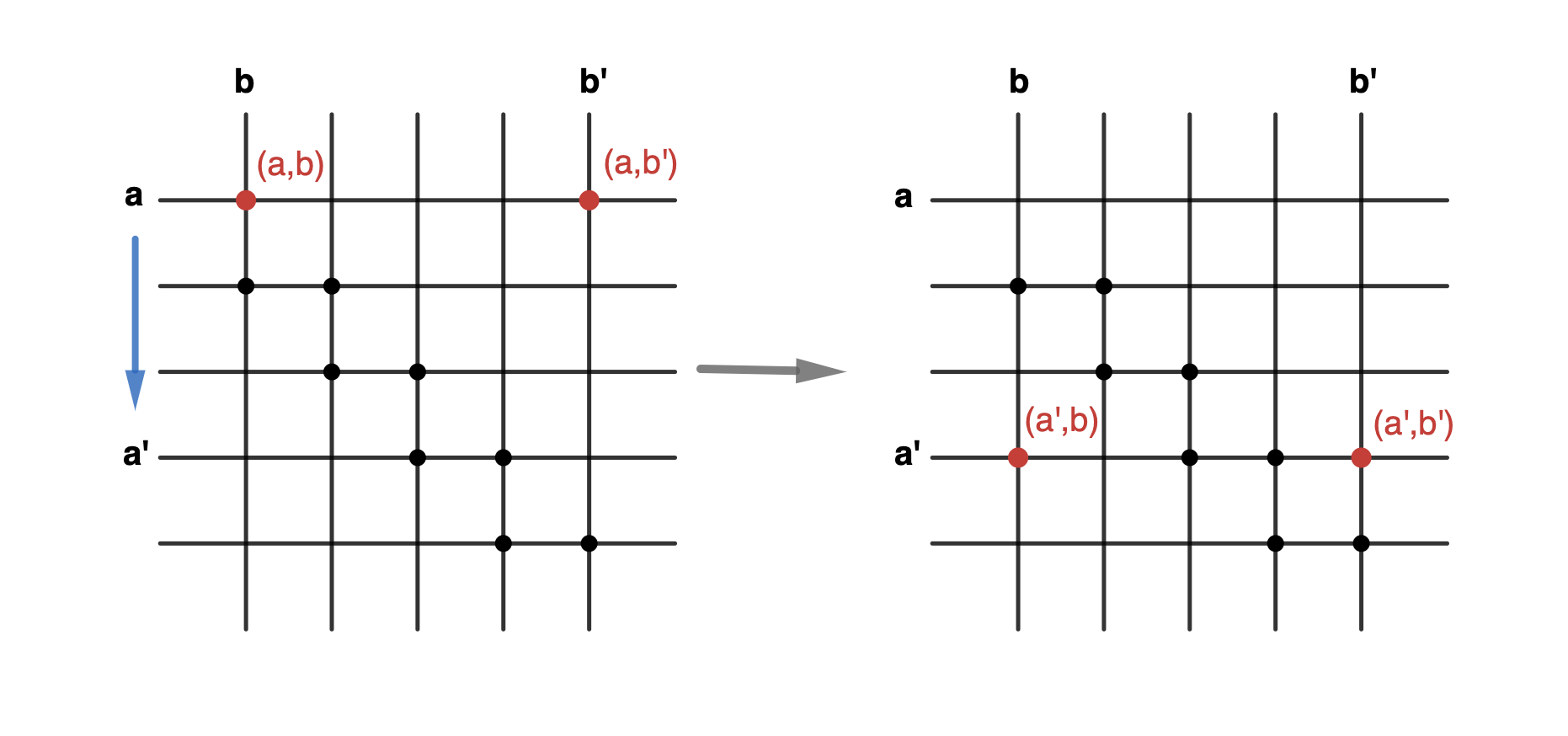}
    \caption{Illustration of \Cref{def:shiftoperation}: the shift operation $\{(a,b),(a,b')\} \to \{(a',b),(a',b')\}$.}
    \label{fig:grid-shift-operation}
\end{figure}
\end{defn}
In the analysis, we perform these operations in the following way.
\begin{enumerate}
    \item We first apply a sequence of row and column merge operations. For these operations, we never merge a row or column which was part of a previous merge operation. The order in which we apply these operations does not matter, so we can assume that we first apply all of the column merge operations and then apply the row merge operations.
    
    These merge operations correspond to creating the parts in the partitions $\P_r$ and $\P_c$.
    \item We then iteratively choose a row $a$, pair up the vertices in that row, and shift some or all of the pairs of vertices to different rows which have not yet been chosen.
    
    A shift operation $\{(a,b),(a,b')\} \to \{(a',b),(a',b')\}$ corresponds to the folding operation $\vvv_a \sim \vvv_{a'}$ in the proof of Theorem \ref{thm:staircaseCatalan}.
\end{enumerate}
Since we are taking the limit as $n \to \infty$, we focus on sequences of operations with minimum total weight.
\begin{defn}
We say that a sequence of row merge, column merge, and shift operations on a starting set of vertices $V_0$ is \emph{efficient} if it results in a multi-set of vertices where each vertex has multiplicity at least $2$ and the total weight of the operations is minimized.
\end{defn}
We show that in efficient sequences of operations, each operation is applied to a single staircase and breaks this staircase into smaller staircases.
\begin{defn}\label{def:staircase}
We define a \emph{staircase of length $j$} to be a set of vertices of the form 
$$\left\{(a_i,b_i): i \in [j]\right\} \cup \left\{(a_{i+1},b_i): i \in [j-1]\right\} \cup \{(a_1,b_j)\}$$
together with the edges 
$$
\left\{\{(a_i,b_i),(a_{i+1},b_i)\}: i \in [j-1]\right\} \cup \left\{\{(a_{i+1},b_i),(a_{i+1},b_{i+1})\}: i \in [j-1]\right\} \cup \left\{\{(a_j,b_j),(a_1,b_j)\}, \{(a_1,b_j),(a_1,b_1)\}\right\}
$$
where $a_1,\dots,a_j$ are distinct indices in $[k]$ and $b_1,\dots,b_j$ are distinct indices in $[k]$.
\begin{figure}[hbt!]
    \centering
    \includegraphics[scale=0.32]{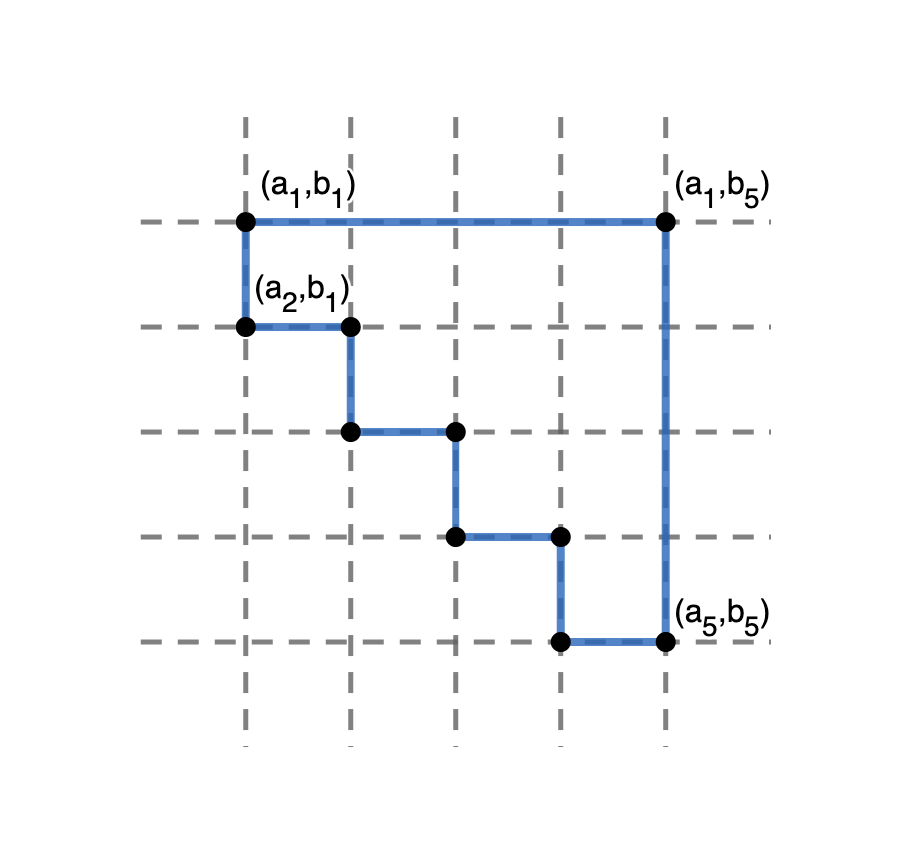}
    \caption{Illustration of \Cref{def:staircase}.}
    \label{fig:grid-staircase}
\end{figure}
\end{defn}
Row merges, column merges, and shifts on a single staircase split this staircase into $w+1$ new staircases where $w$ is the weight of the operation.
\begin{enumerate}
    \item If we have a staircase on the vertices $$\left\{(a_i,b_i): i \in [j]\right\} \cup \left\{(a_{i+1},b_i): i \in [j-1]\right\} \cup \{(a_1,b_j)\}$$ and apply a row merge operation on the indices $a_{i_1},\dots,a_{i_t}$ (where $a_{i_1} < \dots < a_{i_t}$) then we have the $t$ staircases 
    $$\left\{\left\{(a_i,b_i): i \in [i_x,i_{x+1}-1]\right\} \cup \left\{(a_{i+1},b_i): i \in [i_x,i_{x+1}-1]\right\} \cup \{(a_{i_x},b_{i_{x+1}-1})\}: x \in [t]\right\}$$
    where we set $a_{i_1} = \dots = a_{i_t}$, $a_{i+j} = a_i$, and $b_{i+j} = b_i$. Column merges have a similar effect.
    \item If we choose a row of the staircase and shift its two vertices to a different row, this has the same effect as merging the two rows.
\end{enumerate}
In order to show Theorems \ref{thm:disjointstaircases}, \ref{thm:noncrossing-dominant-row}, and \ref{thm:dominant-noncrossing} (i.e. that non-crossing partitions are optimal), we need to show the following. If we start with a staircase of length $k$, for all efficient sequences of operations,
\begin{enumerate}
    \item The total weight of the operations is $k-1$.
    \item After all of the operations, we are left with a multi-set of vertices $V$ where each vertex has multiplicity exactly $2$.
    \item Each operation only affects vertices in a single staircase and breaks this staircase into $w+1$ new staircases where $w$ is the weight of the operation. Thus, at each step we have a set of staircases. Moreover, these staircases are always disjoint, i.e. we never have a vertex which appears multiple times where one copy is in one staircase and another copy is in a different staircase.
\end{enumerate}
We now prove these statements. We start with the second statement.
\begin{lemma}
After any efficient sequence of operations, we are left with a multi-set of vertices $V$ where each vertex has multiplicity exactly $2$.
\end{lemma}
\begin{proof}
We observe if we have a sequence of operations where we are left with a multi-set of vertices where each vertex has multiplicity at least $2$ and some vertex $v$ has multiplicity at least $3$ then we can reduce the total weight of the operations as follows.

Let $b$ be the column $v$ is in and modify the sequence of operations so that all of the rows which have a vertex in column $b$ at the end are merged together. Since we have that at the end, the total number of vertices in column $b$ is even, each vertex has multiplicity at least $2$, and $v$ has multiplicity at least $3$,  if we merged $j$ rows together than we must have a vertex of multiplicity at least $2j+2$. Now observe that in order to have a vertex of multiplicity at least $2j+2$, at least $j+1$ columns must have been merged together. By removing this column merge, we can end with vertices with multiplicity $2$ and reduce the total weight of the operations by at least $j$. The additional row merge only increased the total weight of the operations by $j-1$, so the total weight of the operations decreased by at least $1$, as needed.
\end{proof}
To prove the first and third statements, we add horizontal and vertical edges between the vertices. We show that for any efficient sequence of operations, we can choose these edges so that they give us our disjoint staircases.
\begin{defn}
Given a multi-set of vertices $V \subseteq \{(a,b): a,b \in [m]\}$ such that every row and column has an even number of vertices, we define a \emph{row matching} $M_r$ to be a matching between the vertices of $V$ such that every vertex $(a,b) \in V$ is matched with a vertex $(a,b') \in V$ in the same row. Similarly, we define a \emph{column matching} $M_c$ to be a matching between the vertices of $V$ such that every vertex $(a,b) \in V$ is matched with a vertex $(a',b) \in V$ in the same column. If there are multiple copies of the same vertex $(a,b)$ in $V$ then these copies can be matched to each other or to other vertices.

We define $(V,M_r,M_c)$ to be the multi-graph with vertices $V$ and edges $M_r \cup M_c$.
\end{defn}
\begin{lemma}\label{lem:staircase-matching}
Given any sequence of operations where we start with
$$V_0 = \left\{(i,i): i \in [k]\right\} \cup \left\{(i+1,i): i \in [k-1]\right\} \cup \{(1,k)\}$$
and end with $k$ vertices of multiplicity $2$, we can choose a row matching $M_r$ and a column matching $M_c$ for each step such that 
\begin{enumerate}
    \item If we have a row merge operation where we merge rows $a_{i_1},\dots,a_{i_t}$ then before the merge, $M_r$ has the two vertices in each row $a_{i_x}$ paired up and after the merge, $M_r$ has a matching between these $2t$ vertices. No other edges of the row matching $M_r$ are changed and the column matching $M_c$ is unchanged.
    
    The analogous statement holds for column merges.
    \item If we have a shift operation $\{(a,b),(a,b')\} \to \{(a',b),(a',b')\}$ then before the shift operation, $\{(a,b),(a,b')\} \in M_r$. After the shift operation, we either have the row matching 
    $$\(M_r \setminus \left\{\{(a,b),(a,b')\}\right\}\)  \cup \left\{\{(a',b),(a',b')\}\right\}$$
    or the row matching 
    $$\(M_r \setminus \left\{\{(a,b),(a,b')\},\{(a',b''),(a',b''')\}\right\}\)  \cup \left\{\{(a',b),(a',b'')\}, \{(a',b'),(a',b''')\}\right\}$$
    where $\{(a',b''),(a',b''')\}$ is another edge which was in $M_r$. The column matching $M_c$ is unchanged. See \Cref{fig:shift-matching-1} and \Cref{fig:shift-matching-2} for illustration (we drew these figures from right to left because in the proof of \Cref{lem:staircase-matching} we will construct the row matchings by starting at the end and working backwards).

    \item At the beginning, $(V_0,M_r,M_c)$ is a staircase of length $k$. At the end, $M_r$ and $M_c$ consist of loops on the $k$ vertices of multiplicity $2$ (i.e. each vertex is paired with the other copy of that vertex).
\end{enumerate}
\end{lemma}

\begin{proof}
To prove this lemma, we start from the end and work backwards. More precisely, at the end we take $M_r$ and $M_c$ to be the row and column matchings which consist of loops on the $k$ vertices of multiplicity $2$. For each operation, we show that given the row and column matchings after the operation, we can find row and column matchings before the operation which satisfy the needed conditions.
\begin{enumerate}
    \item If we perform a shift operation $\{(a,b),(a,b')\} \to \{(a',b),(a',b')\}$ and after the shift we have that $\{(a',b),(a',b')\} \in M_r$ then we can take the row matching
    $$\left(M_r \setminus \left\{\{(a',b), (a',b')\}\right\}\right)  \cup \left\{\{(a,b), (a,b')\}\right\}$$ before the shift. Otherwise, we have that after the shift, $\{(a',b),(a',b'')\},\{(a',b'),(a',b''')\} \in M_r$ for some $b'',b'''$ so we can take the row matching
    $$\left(M_r \setminus \left\{\{(a',b),(a',b'')\}, \{(a',b'),(a',b''')\}\right\}\right) \cup \left\{\{(a,b),(a,b')\}, \{(a',b''),(a',b''')\}\right\}$$
    before the shift.
    \item If we have a row merge operation where we merge rows $a_{i_1},\dots,a_{i_t}$ then the only option for how to pair up the vertices in these rows in $M_r$ before the merge is to pair up the two vertices in each row $a_{i_x}$. We can leave the column matching $M_c$ and the edges in $M_r$ not involving the vertices in rows $a_{i_1},\dots,a_{i_t}$ unchanged.
    
    We handle column merges in a similar way. 
\end{enumerate}
Finally, we observe that at the start, each row and column have exactly two vertices so there is only one choice for $M_r$ and $M_c$ and this gives the staircase of length $k$.
\end{proof}
With these row and column matchings, we can now prove the first and third statements.
\begin{cor}
Efficient sequences of operations on a staircase of length $k$ have total weight $k-1$.
\end{cor}
\begin{proof}
We make the following observations about how the multi-graph $(V,M_r,M_c)$ changes at each step. 
\begin{enumerate}
    \item There is one connected component at the start and $k$ connected components at the end.
    \item Each operation increases the number of connected components by at most $w$ where $w$ is the weight of the operation.
\end{enumerate}
These observations immediately imply that the sequence of operations must have total weight at least $k-1$, as needed.
\end{proof}
\begin{cor}
For any efficient sequence of operations on a staircase of length $k$, each operation only affects vertices in a single staircase and breaks this staircase into $w+1$ new staircases where $w$ is the weight of the operation. Thus, at each step we have a set of staircases. Moreover, these staircases are always disjoint, i.e. we never have a vertex which appears multiple times where one copy is in one staircase and another copy is in a different staircase.
\end{cor}
\begin{proof}
To prove this statement, we consider how the multi-graph $(V,M_r,M_c)$ changes at each step. We show that for each operation of weight $w$, the only way for this operation to increase the number of connected components by $w$ is if it acts on a single staircase and breaks this staircase into $w+1$ disjoint staircases. We have the following cases:
\begin{enumerate}
    \item If we have a row merge operation where we merge rows $a_{i_1},\dots,a_{i_t}$ then consider the connected components involving the vertices in these rows. In order for the number of connected components to increase by $t-1$, there must be exactly one connected component before the merge and exactly $t$ connected components after the merge. Before the merge, we have a collection of disjoint staircases, so all vertices in these rows must be in the same staircase. Since we are acting on a single staircase, it is not hard to show that there is a unique choice for how to pair these vertices up after the merge which results in $t$ connected components and this choice splits the original staircase into $t$ disjoint staircases.
    
    Similar logic applies to column merges.
    \item If we have a shift operation $\{(a,b),(a,b')\} \to \{(a',b),(a',b')\}$, then before the shift we have that $\{(a,b),(a,b')\} \in M_r$ because of the way that we chose $M_r$ and $M_c$. After the shift, we will either have the row matching
    $$\(M_r \setminus \left\{\{(a,b),(a,b')\}\right\}\) \cup \left\{\{(a',b),(a',b')\}\right\}$$
    or the row matching 
    $$\(M_r \setminus \left\{\{(a,b),(a,b')\},\{(a',b''),(a',b''')\}\right\}\)  \cup \left\{\{(a',b),(a',b'')\},\{(a',b'),(a',b''')\}\right\}$$
    for some edge $\{(a',b''),(a',b''')\} \in M_r$. However, the only way that the shift operation can increase the number of connected components is if two edges in the same staircase are removed. This uniquely determines $\{(a',b''),(a',b''')\}$ (assuming that the staircase has row $a'$, otherwise the shift operation cannot increase the number of connected components). There are now two choices for $b''$ and $b'''$ as these indices can be swapped. One of these choices keeps the staircase as a single connected component while the other splits the staircase into two staircases. Thus, there is a unique choice for $M_r$ after the shift which increases the number of connected components and this choice breaks the staircase containing $\{(a,b),(a,b')\}$ into two staircases. 
\end{enumerate}
Finally, we show that we never have two different staircases with the same vertex. To see this, assume that there is an efficient sequence of operations such that at some point, there are two different staircases which each have a copy of some vertex $(a,b)$ and consider the first point where this occurs. Observe that we can swap these copies for the remainder of the sequence, both in terms of which vertices are chosen for shift operations and the edges in $M_r$ and $M_c$. If we make this swap then looking at the previous step, we will have a different column matching $M_c$. However, we showed above that for efficient sequences of operations, $M_r$ and $M_c$ are uniquely determined at each step by the operation we perform. This is a contradiction, so we cannot have two different staircases with the same vertex. 
\end{proof}

\end{document}